# OPTIMIZATION WITH DISCRETE SIMULTANEOUS PERTURBATION STOCHASTIC APPROXIMATION USING NOISY LOSS FUNCTION MEASUREMENTS

by
Qi Wang

A dissertation submitted to Johns Hopkins University in conformity with the requirements for the degree of Doctor of Philosophy.

Baltimore, Maryland
October, 2013

# Abstract

Discrete stochastic optimization considers the problem of minimizing (or maximizing) loss functions defined on discrete sets, where only noisy measurements of the loss functions are available. The discrete stochastic optimization problem is widely applicable in practice, and many algorithms have been considered to solve this kind of optimization problem. Motivated by the efficient algorithm of simultaneous perturbation stochastic approximation (SPSA) for continuous stochastic optimization problems, we introduce the middle point discrete simultaneous perturbation stochastic approximation (DSPSA) algorithm for the stochastic optimization of a loss function defined on a $p$-dimensional grid of points in Euclidean space.

We show that the sequence generated by DSPSA converges to the optimal point under some conditions. Consistent with other stochastic approximation methods, DSPSA formally accommodates noisy measurements of the loss function. We also show the rate of convergence analysis of DSPSA by solving an upper bound of the mean squared error of the generated sequence. In order to compare the performance of DSPSA with the other algorithms such as the stochastic ruler algorithm (SR) and the stochastic comparison algorithm (SC), we set up a bridge between DSPSA and the other two algorithms by comparing the probability in a big-$O$ sense of not achieving the optimal solution. We show the theoretical and numerical comparison results of DSPSA, SR, and SC. In addition, we consider an application of DSPSA towards developing optimal public health strategies for containing the spread of influenza given limited societal resources.

This dissertation also contains three appendices. The first appendix considers the analysis of practical step size selection in stochastic approximation algorithms for continuous problems. The second appendix discusses the rate of convergence analysis of SPSA for time-varying loss functions. The third appendix focuses on the numerical experiments on the properties of the upper bound of the mean squared errors for the sequence generated by DSPSA.

First Reader and Advisor: James C. Spall

Second Reader: Daniel P. Robinson

# Acknowledgements

My dissertation research had the usual challenges, and with the help of many people, I finally completed my thesis. I would like to acknowledge the contribution of some of them here.

Foremost, Professor Jim Spall, who is an excellent teacher, a famous scholar and a smart thinker, has made a huge contribution to the thesis. He introduced me to the field of stochastic approximation. He provided many good ideas and helped me to find the correct way. He also sacrificed much of his leisure time to read my work and write useful comments on my thesis. He gave a lot of effort to try to encourage me when I met difficulties, and I fully appreciate his kindness. There are many others things he has done for me. Even though I have no space to list all of them here, but I definitely remember all of them in my heart. I also want to thank Kathy Spall for her kindness to review the wording and grammar issues in some of the thesis for me.

I also want to greatly thank the second reader Professor Daniel Robinson for reading my thesis. He took large amount of time and effort in making this thesis better.

I appreciate the support I got from the Department of Applied Mathematics and Statistics of the Johns Hopkins University. All the professors and staffs helped me a lot in many areas. The fellowship and the teaching assistant opportunities helped me to pay the tuition and to keep independent from the financial support of my parents. The Duncan Fund provided financial support to me to attend several conferences. I want to thank Professor Daniel Naiman especially, who helped to apply to the Duncan Fund for me.

I also want to thank the Applied Physics Lab (APL) of the Johns Hopkins University, since I got research funds from them for some semesters. These funds help me to put more time on research.

Finally, I am grateful for the strong support from my husband, Ming Ye. He is very smart and told me many efficient ways to think of mathematical problems. He also helped to review some of the thesis and provided many useful comments.

# Contents

# List of Tables

# List of Figures

# Chapter 1

# Introduction

The purpose of this chapter is to introduce the research problem, provide a literature review, and summarize main research results. We discuss the motivation in Section 1.1, and introduce the problem and previous work in Section 1.2. These sections are followed by a summary of the main results of the research.

## 1.1 Motivation

The optimization of real-world stochastic systems typically involves the use of a mathematical algorithm that iteratively seeks out the solution. It is often the case that the domain of optimization is discrete and that only noisy objective function measurements are available to carry out the optimization process.

Some problems of interest within this framework include transmission problem in networks (Wieselthier et al., 1992, Cassandras and Julka, 1995, Mishra et al., 2007), spreading code design problem (Krishnamurthy et al., 2004), assignment problem of people evacuation (Francis, 1981), facility locating problem (Ermoliev and Leonardi, 1982), weapons assignment problem (Cullenbine et al., 2003), antenna selection problem (Liu et al, 2012), and, more generally, distributing a discrete amount of resources to a finite number of users in the face of uncertainty (Castanon and Wohletz, 2009, Wang and Gao, 2010).

In this paragraph, we briefly discuss some of the problems above. Wieselthier et al. (1992) discuss the model of a multiple service, multiple resource (MSMR) problem. In networks that support voice traffic, the acceptance of a call is based on the commitment of resources at all nodes along the path. The call may be blocked due to the lack of resource at any node along the path. The objective is to optimize the total performance of the networks, which is related to throughput and blocking probability. Cassandras and Julka (1994) discuss the problem where a single resource must provide service to a set of customer classes and the single resource can serve multiple customers simultaneously from the same class. The objective is to determine the assignment probabilities and use them to optimize the performance by solving the scheduling problem. Mishra et al. (2007) consider the problem of admission control of packets in communication networks under dependent service times. The parameters to be optimized are the thresholds of acceptance and rejection of the incoming packets. Krishnamurthy et al. (2004) consider the problem of optimization of the spreading codes of users in the CDMA system. The objective is to maximize the signal-to-interference-plus-noise ratio by choosing the optimal spreading codes. Francis (1981) discusses the problem of evacuating a building in the minimal time by assigning people to different evacuation routes. Liu et al. (2012) discuss the problem of sum rate maximization antenna selection in MIMO (multiple input, multiple output) two-way AF relay with imperfect CSI (channel state information). They want to determine the optimal antenna subset under the imperfect information, which is a discrete problem with noise. Castanon and Wohletz (2009) consider one class of stochastic resource allocation problems, and in their class, resources assigned may fail to complete the task with certain probability. The objective is to minimize the incomplete task value and the cost of using resources.

Many algorithms, such as ranking and selection, multiple comparisons, stochastic ruler and stochastic comparison (more discussions will be in Section 1.2.2), have been considered for such discrete stochastic optimization problems. By the “No Free Lunch Theorem” (Spall 2003, Subsection 1.2.2), we know that no algorithm can outperform all other algorithms for all kinds of functions. There is an inherent trade-off between the robustness and the efficiency of algorithms. In view of the trade-off, some algorithms may be more robust, which indicates that they may be applicable for many kinds of functions, but at the same time, these algorithms may not be efficient. On the other hand, some other algorithms may be very efficient for some specific kinds of functions, but they may only be applicable for fewer kinds of loss functions. In view of objective functions, for some loss functions we have to compare the function values in all points to decide the optimal solution, while for others we can use the special structure of the functions to find out the optimal solution more efficiently.

As will be discussed in Section 1.2.2, many algorithms for discrete optimization often focus on the comparisons of loss function values at different points or the comparisons of the ranks of candidate points, and these algorithms do not really take advantage of function structure. Furthermore, due to the noisy measurements of loss function, point-wise comparisons in each iteration may not lead to the right results. Hill et al. (2004) introduce a new algorithm, which uses the similar idea of the simultaneous perturbation stochastic approximation (SPSA) algorithm (Spall, 1992, 1998a). SPSA is one kind of stochastic approximation algorithm used for continuous stochastic optimization problems. The gradient estimate in SPSA is based on the idea of simultaneous perturbation instead of the general finite difference estimate in the finite difference stochastic approximation algorithm (FDSA). Compared with FDSA, SPSA can achieve the same level of accuracy by using $p$ times fewer noisy measurements of the loss function in each iteration, where $p$ is the dimension of the problem. Overall, SPSA has many advantages such as: 1) it is relatively easy to implement SPSA into real problems; 2) there are only two noisy measurements of the loss function in each iteration; 3) SPSA implicitly makes use of loss function structure; 4) SPSA can handle noise properly; 5) SPSA is efficient for high-dimensional problems. In Hill et al. (2004), they use the similar idea of SPSA and introduce a stochastic approximation type algorithm for discrete stochastic optimization problems. They also develop preliminary results associated with convergence for a separable discrete loss function under special conditions. However, their algorithm can be shown not to converge to the optimal solution in simple examples.

For some one-dimensional problems, this algorithm converges to a point next to the optimal solution, but not the optimal one.

We introduce a different form of discrete version of SPSA (different from Hill et al., 2004). This new version of discrete SPSA (DSPSA) applies to a broad range of problems, while potentially retaining the essential efficiency advantage of standard SPSA. We show almost sure convergence of this algorithm under some conditions, discuss the rate of convergence property, and compare it with other popular algorithms.

# 1.2 Discrete Stochastic Optimization

In this section, we introduce the discrete stochastic optimization problem that we will consider, and do the literature review for a broad range of existing algorithms.

## 1.2.1 Problem Statement

We consider a real-valued function $L(\boldsymbol{\theta}): \mathbb{Z}^p \to \mathbb{R}$, where $\boldsymbol{\theta}$ is a $p$-dimensional vector, $\mathbb{Z}^p$ is the space of all $p$-dimensional multivariate integer points, $\mathbb{R}$ is the space of reals, and $L$ is a loss function to be minimized. Suppose $\Theta \subseteq \mathbb{Z}^p$ is the domain of the allowable values of $\boldsymbol{\theta}$. We consider the optimization problem of

$$\min_{\boldsymbol{\theta} \in \Theta} L(\boldsymbol{\theta}).$$

Suppose the solution set of the problem above is $\Theta^*$. Then we have

$$\Theta^* = \arg\min_{\boldsymbol{\theta}\in\Theta} L(\boldsymbol{\theta}) = \left\{\boldsymbol{\theta}^* \in \Theta : L(\boldsymbol{\theta}^*) \le L(\boldsymbol{\theta})\ \forall \boldsymbol{\theta} \in \Theta\right\}.$$

But the measurements of the objective function $L(\boldsymbol{\theta})$ involve noise. We assume the noisy measurement of $L(\boldsymbol{\theta})$ is $y(\boldsymbol{\theta})$, and

$$y(\boldsymbol{\theta}) = L(\boldsymbol{\theta}) + \varepsilon(\boldsymbol{\theta}),$$

where $\varepsilon$ is the measurement of noise. It follows that the stochastic optimization problem can be regarded as

$$\min_{\boldsymbol{\theta}\in\Theta} E\left(y(\boldsymbol{\theta})\right).$$

Generally, we do not know the exact form of the objective function $L(\boldsymbol{\theta})$, and we need to do the optimization only based on the noisy measurements of the loss function through simulations.

Here we consider a general domain within $\mathbb{Z}^p$. However, there are often discrete problems with a countable number of candidates that are not at the integer-based grid points. For such non-grid problems, we may map these non-integer points onto the domain of $\mathbb{Z}^p$, and there are many possible mappings. We will further discuss this problem in Section 8.2 as a future research direction.

## 1.2.2 Previous Work on Algorithms

There are many papers discussing the discrete stochastic optimization problems. In this section, we review a broad range of these methods. Basically, we can divide these algorithms into several classes: statistical approach, random search method, stochastic approximation, and other algorithms. Hill (2013), Fu (2002) and Swisher et al. (2004) provide some literature reviews on the discrete stochastic optimization with noisy loss function measurements.

### 1.2.2.1 Statistical Approaches

The class of statistical approaches includes the methods of ranking and selection (R&S), multiple comparisons, ordinal optimization and others. The basic idea of these methods is to use statistical analysis to make decisions under a confidence level.

R&S and multiple comparisons procedures (Swisher et al., 2003, Goldsman and Nelson, 1994, Kim and Nelson, 2007, Bechhofer et al., 1995, Bechhofer, 1954) evaluate all candidates from a given fixed and finite set of feasible domain. Due to the noisy measurements, these methods need to collect a certain number of observations for each candidate. These algorithms mainly focus on the optimization problem with small number of feasible solutions (20 is a frequently recommended number [Spall, 2003, p. 302]).

For the R&S algorithm, based on observations, the goal is to select either the best candidate with a specific probability or a subset of candidates that contains the best candidate with a specific probability. These two types of selection procedures correspond to the two general classes of R&S: indifference-zone ranking and subset selection.

For multiple comparisons procedures, the goal is to quantify the difference between candidates and construct simultaneous confidence intervals of the differences between loss function values at different points. Swisher et al. (2003) divide this approach into three classes: all-pairwise comparisons approaches, multiple comparisons with a control, and multiple comparisons with the best.

The ordinal optimization algorithm is designed for the problem where it is easier to determine the relative order than the precise function values. Ho et al. (2000) and Ho et al. (2007) discuss the ideas of ordinal optimization, and these ideas are similar to the R&S algorithm and multiple comparisons algorithm. Ordinal optimization can work for the domains with a large number of candidates. In the ordinal optimization algorithm, a budget allocation problem is solved to decide the number of measurements needed for each candidate point to achieve some level of correct selection probability. Due to the noise, a large number of observations are needed for each point to achieve a certain confidence level of correct selection. The variance reduction techniques can greatly help to reduce the cost and improve the performance.

The algorithm discussed in Kleywegt et al. (2001) is not a typical R&S method, but we think it belongs to the class of statistical approach. In their paper, the loss function value is estimated by the sample average. They solve the sample average optimization problem for several replicates until a stopping criterion is satisfied. Compared with the R&S method, the algorithm in Kleywegt et al. (2001) does not do the screening and selection in each iteration, and it statistically picks the best solution from different replicates based on a confidence interval.

### 1.2.2.2 Random Search Type Algorithms

The basic ideas of the random search method are discussed in Chapter 2 of Spall (2003). These algorithms are applicable for the noise-free case, but with some modifications, such as adding a threshold, these algorithms are also applicable for the noisy case. Basically there are two classes of random search algorithms for discrete stochastic optimization. For the first class of algorithms, in each iteration the loss function is measured only at the candidate point and the current point, and for the second class of algorithms, the loss function is measured at points in a big set (promising area). Now we discuss these two classes of random search algorithms for discrete stochastic optimization in detail.

For the first class of random search type algorithms, suppose the generated sequence is $\left\{\hat{\boldsymbol{\theta}}_k\right\}$. Assume that $N(\hat{\boldsymbol{\theta}}_k)$ is a set of neighboring points of point $\hat{\boldsymbol{\theta}}_k$ excluding $\hat{\boldsymbol{\theta}}_k$ itself (the new candidate point should be generated from $N(\hat{\boldsymbol{\theta}}_k)$). For any point $\tilde{\boldsymbol{\theta}}_k$ in $N(\hat{\boldsymbol{\theta}}_k)$, suppose $R(\hat{\boldsymbol{\theta}}_k, \tilde{\boldsymbol{\theta}}_k)$ is the probability to generate $\tilde{\boldsymbol{\theta}}_k$ as the new candidate, where $R(\hat{\boldsymbol{\theta}}_k, \tilde{\boldsymbol{\theta}}_k) > 0$ and $\sum_{\tilde{\boldsymbol{\theta}}_k \in N(\hat{\boldsymbol{\theta}}_k)} R(\hat{\boldsymbol{\theta}}_k, \tilde{\boldsymbol{\theta}}_k) = 1$. The general random search algorithm in the first class can be described as:

Step 1: Pick the initial guess $\hat{\boldsymbol{\theta}}_0$, $k = 0$.

Step 2: Given $\hat{\boldsymbol{\theta}}_k$, generate a new candidate $\tilde{\boldsymbol{\theta}}_k \in N(\hat{\boldsymbol{\theta}}_k)$ with probability $R(\hat{\boldsymbol{\theta}}_k, \tilde{\boldsymbol{\theta}}_k)$.

Step 3: Do the comparisons based on some criterion to decide whether to accept $\tilde{\boldsymbol{\theta}}_k$. If accepted, let $\hat{\boldsymbol{\theta}}_{k+1} = \tilde{\boldsymbol{\theta}}_k$; otherwise $\hat{\boldsymbol{\theta}}_{k+1} = \hat{\boldsymbol{\theta}}_k$.

Step 4: Set $k = k + 1$, go to step 2.

In each iteration, a candidate point is generated from the neighborhood of the current point, and based on some comparisons, the current point is either replaced by the new candidate point or retained. The definition of the neighborhood $N(\hat{\boldsymbol{\theta}}_k)$ depends on the choice of algorithm. For example, the neighborhood of $\hat{\boldsymbol{\theta}}_k$ can be $\Theta \setminus \{\hat{\boldsymbol{\theta}}_k\}$ or a set containing points that are closed to $\hat{\boldsymbol{\theta}}_k$ under some norm.

The first class of random search type algorithms for discrete stochastic optimization includes the stochastic ruler algorithm (Yan and Mukai, 1992), the stochastic comparison algorithm (Gong et al., 1999), and the modified simulated annealing algorithm (Gelfand and Mitter, 1989, Fox and Heine, 1995, Gutjahr and Pflug, 1996, Alrefaei and Andradottir, 1999). Andradottir (1999) introduces an idea to accelerate the convergence speed of the first class of random search type algorithms. By recording all previous data and adding extra optimization step in each iteration, Andradottir (1999) improves the performance of the stochastic ruler algorithm, the stochastic comparison algorithm and the simulated annealing algorithm. In the following, we expand the discussion of each method in the first class of random search type algorithms for discrete stochastic optimization.

For the stochastic ruler algorithm, the idea is to change the minimization problem to a maximization problem. A stochastic ruler $U_{u,v}$ is defined and it is uniformly distributed over the interval $[u, v]$. The equivalent maximization problem is

$$\max_{\boldsymbol{\theta} \in \Theta} P\left( y(\boldsymbol{\theta}) \leq U_{u,v} \right).$$

The implementation of the algorithm is to do multiple comparisons in each iteration between the noisy measurements of loss function at the current point and the stochastic ruler. If any of the comparisons indicate that $U_{u,v}$ is smaller than $y(\boldsymbol{\theta})$, the current point will be kept; otherwise the current point will be replaced. More discussions on the stochastic ruler algorithm can be seen in Section 5.2.

For the stochastic comparison algorithm, the idea is similar to the stochastic ruler algorithm, where the minimization problem is also replaced by a maximization problem. In each iteration, the noisy measurements of loss function values at the current point and the candidate point are compared to decide whether the current point is kept or replaced. Gong et al. (1999) consider the case of $N(\boldsymbol{\theta}) = \Theta \setminus \boldsymbol{\theta}$ (set of $\Theta$ excluding $\boldsymbol{\theta}$), and the authors show that under some reasonable conditions, the generated sequence converges to the global optimal point. More discussions on the stochastic comparison algorithm can be seen in Section 5.3.

For the modified simulated annealing algorithm, when the variance of the noise goes to 0 as $k \to \infty$, Gelfand and Mitter(1989) and Gutjahr and Pflug (1996) show the result of convergence for the noisy case. Without the restrictive variance assumption, the most obvious way to deal with the noise is to use sample average to do the comparisons in each iteration. However, the cost may be dramatically increased. Fox and Heine (1995) consider the algorithm to calculate the average based on all previous iterations. Another way to handle the noise is to add a threshold (e.g. accepting the candidate point $\tilde{\boldsymbol{\theta}}_k$, if for the current point $\hat{\boldsymbol{\theta}}_k$, $y(\tilde{\boldsymbol{\theta}}_k) \leq y(\hat{\boldsymbol{\theta}}_k) + \iota$, where $\iota < 0$ is a threshold). In the Proposition 8.1 of Spall (2003), the author shows a result for the method of threshold, but this proposition does not provide clear convergence property for the full process. Alrefaei and Andradottir (1999) consider the modification of a simulated annealing algorithm by using constant temperature. But they record the noisy measurements in all previous iterations and set the current optimal point to be the point with the smallest average function value based on all historical information.

Andradottir (1999) uses the similar idea of Alrefaei and Andradottir (1999) and expands the idea to other random search type algorithms, such as the stochastic ruler algorithm and the stochastic comparison algorithm. Generally, Andradottir (1999) adds one more step into the original random search type algorithms. In this step, all old information is stored and the point with the smallest average noisy loss function value (based on all old information) is set as the current optimal point. The convergence properties can be proved through the strong law of large numbers and the central limit theorem. Even though lots of information needs to be stored, the extra optimization step increases the speed of convergence.

For the second class of random search type algorithms in discrete stochastic optimization, the loss function values are measured on a set of points in each iteration. This class of algorithms is sometimes called "adaptive random search," where a promising area is chosen in each iteration. Generally, for this class, the algorithms contain two critical steps: the sampling scheme that determines the points to be sampled and the estimation scheme that decides the points to be measured and the number of measurements for each point. This kind of class includes convergent optimization via most-promising-area stochastic search (Hong and Nelson, 2006, Li et al., 2009), adaptive hyperbox algorithm (Xu et al., 2013), locally convergent random-search algorithm (Hong and Nelson, 2007), and general nested partitions (Pichitlamken and Nelson, 2003). In the following, we expand the discussion of each method in the second class of random search type algorithms.

The basic idea of the convergent optimization via most-promising-area stochastic search (COMPASS) is to find the most promising area and to do more sampling within that area to increase the efficiency of the algorithm. All points sampled are to be evaluated for the noisy loss function values, and all observations up to the current iteration are used to decide the optimal solution in each iteration. The adaptive hyperbox algorithm (AHA) is a modification of COMPASS done by modifying the way to generate the most promising area. The method of locally convergent random-search algorithm generalizes the idea in COMPASS and provides a revised form of COMPASS, such that in the estimation scheme fewer points need to be evaluated. The nested partitions method for the noise-free case is discussed in Shi and Olafsson (2000), and they show that the most promising area shrinks to the optimal solution with probability one. The application of nested partitions in noisy loss function is considered in Pichitlamken and Nelson (2003) and they modify the algorithm in several parts: 1) the way that they choose the most promising area is based on all the observations up to current iteration, while Shi and Olafsson (2000) only use measurements in current iteration; 2) at search termination, they pick the optimal solution with the smallest sample mean up to that time. Based on these two key modifications, they can show the almost sure convergence under less restrictive conditions.

### 1.2.2.3 Stochastic Approximation Type Algorithms

The discrete version of the stochastic approximation algorithms mimics the idea of stochastic approximation algorithms for the continuous stochastic optimization problem. Stochastic approximation (SA) has been well developed for continuous parameter problems, and Kushner (2010) provides a survey of the SA methods. The general formula of SA is

$$\hat{\boldsymbol{\theta}}_{k+1} = \hat{\boldsymbol{\theta}}_k - a_k \hat{\boldsymbol{g}}_k(\hat{\boldsymbol{\theta}}_k),$$

where $\{a_k\}$ is the gain sequence and $\hat{\boldsymbol{g}}_k(\hat{\boldsymbol{\theta}}_k)$ is the estimate of the gradient or subgradient at point $\hat{\boldsymbol{\theta}}_k$. Few papers have considered using the idea of stochastic approximation in the discrete stochastic optimization problems.

Dupac and Herkenrath (1982) solve the problem of finding a root of a function defined on $\mathbb{Z}^p$. Their idea is to extend the discrete function to a continuous one by multilinear interpolation. Then the discrete root finding problem is transferred to the continuous root finding problem, which can be solved by standard SA algorithms. Bhatnagar et al. (2011) set up an algorithm by making use of the similar idea of Dupac and Herkenrath (1982) on multilinear interpolation and the idea of SPSA, but they mainly focus on the problem of minimizing the long-run average cost

$$\min_{\boldsymbol{\theta}\in\Theta} \frac{1}{n} E\left( \sum_{j=1}^{n} h(X_j) \right),$$

where $h(\cdot)$ is real-valued function and $\{X_j\}$ is the underlying Markov process that depends on $\boldsymbol{\theta}$. Gokbayrak and Cassandras (1999) also extend the discrete objective function into a continuous differentiable function, which can be solved by the general SA algorithm. The function value of each non-multivariate-integer point is defined as the convex combination of multivariate integer points. In order to estimate the gradient for each point, they need many loss function measurements in each iteration. Therefore, the expense is costly for this algorithm. Lim (2012) constructs a piecewise linear continuous extension based on the original discrete function. But in order to get the approximated gradient, $p+1$ loss function measurements are needed in each iteration. Thus, for high-dimension problems, the cost of the algorithm is really high. Hill et al. (2004) and Hill (2005) use a different way from Gokbayrak and Cassandras (1999) and Lim (2012) to calculate the gradient-like approximation, and they also provide some preliminary results associated with the convergence for a separable discrete loss function under some conditions. Compared with Gokbayrak and Cassandras (1999) and Lim (2012), Hill et al. (2004) only need two measurements in each iteration, which involves less cost. Wang and Spall (2011) suggest a new discrete version of SPSA, which generalizes the type of loss functions discussed in Hill et al. (2004).

### 1.2.2.4 Other Algorithms

Besides the types of algorithms discussed in Sections 1.2.2.1, 1.2.2.2, and 1.2.2.3, there are some other algorithms that solve the discrete stochastic optimization problem including the evolutionary policy selection – Monte Carlo (EPI-MC) (Hannah and Powell, 2010) and the stochastic branch and bound approach (Norkin et al., 1998a and Norkin et al., 1998b).

Hannah and Powell (2010) consider the algorithm of EPI-MC, which is designed for one-stage stochastic combinatorial optimization with a finite action space ($\Theta$) and a noisy cost function. Policy switching and mutation are two features used in the algorithm of EPI-MC. In each iteration, a set of policies (points) is chosen to evaluate the cost values by Monte Carlo simulation. Norkin et al. (1998a) consider the stochastic branch and bound approach, where three operations are executed iteratively: partition of set, estimation of loss function in each subset, and removal of some subset based on the statistical upper and lower bounds.

## 1.3 Summary of Main Results

The primary contribution of this dissertation to the field is that we introduce a new discrete version of simultaneous perturbation stochastic approximation algorithm (DSPSA) for the discrete stochastic optimization problem. We show that DSPSA is a convergent algorithm under some conditions that are similar to those for some other discrete algorithms. We also provide a rate of convergence analysis of DSPSA for both finite sample performance and asymptotical performance, and compare DSPSA with the stochastic ruler algorithm and the stochastic comparison algorithm. Besides these theoretical results, we apply DSPSA on an epidemic problem in public health to solve for the optimal intervention method, including vaccination priority, antiviral agent policy, and time of school closure.

Currently, most discrete stochastic optimization algorithms do not make use of the loss function's structure in a manner analogous to the gradient information for continuous problems, and these algorithms may be robust even when the feasible domain is very strange. By the "No Free Lunch Theorem", no algorithm can have both perfectly efficient performance and perfectly robust performance, which means there is a trade-off between efficiency and robustness for any algorithm. The discrete stochastic approximation algorithms implicitly make use of the structure information of the loss function through a "sort-of" gradient (subgradient), which can lead to a sequence with fast convergence performance under some assumptions. However, the robustness of the discrete stochastic approximation algorithms may be sacrificed because of the assumptions for convergence property.

There are few algorithms that are of the discrete stochastic approximation type, and the theoretical analysis of these algorithms is not well developed. DSPSA is one kind of discrete stochastic approximation algorithm, and we provide detailed theoretical analysis on its convergence properties. Moreover, the discrete stochastic approximation type algorithms are discrete analogues of continuous stochastic approximation algorithms, while most other discrete stochastic optimization algorithms are based on point-wise comparisons. The characteristics of discrete stochastic approximation type algorithms are very different from other discrete stochastic optimization algorithms. Thus, to our knowledge, there is no work on discussing the comparisons between them. We set up a bridge to do the comparisons of DSPSA and random search type algorithms (such as the stochastic ruler algorithm and the stochastic comparison algorithm).

Overall, there are many good properties of DSPSA, including: 1) DSPSA is a simple algorithm and it is easy to implement DSPSA into computer code; 2) The number of coefficients to be picked is small, and we provide practical guidelines on the choice of coefficients; 3) DSPSA implicitly makes use of loss function structure and it can lead to very efficient performance for some kinds of loss functions; 4) There are only two noisy measurements of the loss function in each iteration; 5) Theory on the convergence and rate of convergence may be developed using powerful methods in stochastic analysis, as carried out by this dissertation. However, no algorithm can be perfect, so DSPSA also has some difficulties: 1) DSPSA is only guaranteed to converge from within a local area around the optimum; 2) DSPSA may be restricted to some kinds of loss functions that satisfy some conditions, and not all conditions are easy to check; 3) For non-integer-grid domain, we need to do reformulation before using DSPSA. These difficulties are also true for many other discrete stochastic algorithms.

The dissertation is organized into eight chapters and three appendices. In Chapter 2, we introduce DSPSA and show the convergence of DSPSA. In Chapter 3, we show the convergence rate of DSPSA by providing an upper bound for the mean square error of the generated sequence. We also discuss the properties of the upper bound theoretically. Furthermore, we consider the guidelines on the choice of coefficients of the gain sequence. In Chapter 4, some numerical experiments are done to show the performances and properties of DSPSA. In Chapter 5, we show the rate of convergence of two random search type algorithms: the stochastic ruler algorithm (SR) and the stochastic comparison algorithm (SC). We then theoretically compare these two algorithms (SR and SC) with DSPSA. In Chapter 6, we compare these three algorithms (SR, SC, and DSPSA) in numerical tests. In Chapter 7, we discuss the application of DSPSA towards developing optimal public health strategies for containing the spread of influenza given limited societal resources. Chapter 8 concludes the work in this dissertation and discusses some remaining problems for further research.

This thesis also contains three appendices. In Appendix A, we consider the analysis of practical step size selection in stochastic approximation algorithms for continuous problems. The practical gain sequence selection is different from the optimal selection (theoretically derived from asymptotical performance). We provide a formal way to justify the reasons why we choose this gain sequence in practice. In Appendix B, we consider the rate of convergence of simultaneous perturbation stochastic approximation algorithm (SPSA) for time-varying loss functions. One important application of time-varying loss function is in the model-free adaptive control with nonlinear stochastic systems, and model-free adaptive control is useful in many practical areas. Therefore, the results in Appendix B show the reasonable performance of SPSA in model-free control in the big-$O$ sense. In Appendix C, we do the numerical experiments on the properties of the upper bound of $E\left\|\hat{\boldsymbol{\theta}}_k - \boldsymbol{\theta}^*\right\|^2$ in DSPSA. We show that the numerical results in Appendix C are consistent with the theoretical analysis in Section 3.2.

# Chapter 2

# Discrete Simultaneous Perturbation Stochastic Approximation Algorithm

In this chapter, we introduce the discrete simultaneous perturbation stochastic approximation algorithm, and discuss the almost sure convergence of DSPSA. We consider both the unconstrained and constrained domains. It is shown that the sequence generated by DSPSA converges to the optimal solution under some general conditions.

## 2.1 Algorithm Description

In this section, we introduce the discrete simultaneous perturbation stochastic approximation (DSPSA) algorithm. First we discuss the motivation of the algorithm by considering the one-dimensional case. Suppose we have a discrete function $L : \mathbb{Z} \to \mathbb{R}$, where $\mathbb{Z}$ denotes the set of integers {…, −2, −1, 0, 1, 2, …}. We want to find the minimal value of the loss function $L$. Let the noisy measurement of the loss function $L$ be $y$, where $y = L + \varepsilon$ and $\varepsilon$ indicates the noise. Figure 2.1 shows an example of a discrete function in one dimensional case with a line connecting the neighbor integer points. The piecewise linear function $\bar{L}$ is a continuous extension of $L$, but $\bar{L}$ is nondifferentiable at the integer points. For a point $\theta \in \mathbb{R} \setminus \mathbb{Z}$, the gradient of $\bar{L}(\theta)$ is

$$
\begin{aligned}
g(\theta) &= L(\lceil \theta \rceil) - L(\lfloor \theta \rfloor) \\
&= L\left(\pi(\theta) + \frac{1}{2}\right) - L\left(\pi(\theta) - \frac{1}{2}\right),
\end{aligned}
$$

where $\lfloor \cdot \rfloor$ is the floor function, $\lceil \cdot \rceil$ is the ceiling function, and $\pi(\theta) = \lfloor \theta \rfloor + 1/2$ is the middle point between $\lfloor \theta \rfloor$ and $\lceil \theta \rceil$. Here $\pi(\theta)$ is the middle point between two neighbor integer points, so $\pi(\theta) \pm 1/2$ must be an integer. If $\theta$ is an integer point, then $\pi(\theta)$ is the middle point between $\theta$ and $\theta + 1$. We see that $g(\theta)$ is also well defined at any integer point $\theta$, where $g(\theta)$ is a subgradient (a vector $\boldsymbol{g}$ is a subgradient of $L$ at $\boldsymbol{\theta}$ if $L(\boldsymbol{\theta}') - L(\boldsymbol{\theta}) \geq \boldsymbol{g}^T(\boldsymbol{\theta}' - \boldsymbol{\theta})$ for all $\boldsymbol{\theta}' \in \mathbb{R}^p$ ). For $\theta \in \mathbb{R}$, the estimated gradient (subgradient) based on noisy function measurements is

$$
\hat{g}(\theta) = \frac{y\left(\pi(\theta) + \frac{1}{2}\Delta\right) - y\left(\pi(\theta) - \frac{1}{2}\Delta\right)}{\Delta},
$$

where $\Delta$ is the random perturbation, and here we consider the special case when the $\Delta$ is a Bernoulli random variable taking values $\pm 1$. We see that the form of the estimated gradient (subgradient) for this piecewise linear continuous extension $\bar{L}$ is the same as the gradient (subgradient) estimate for $\bar{L}$ in the simultaneous perturbation stochastic approximation (SPSA) algorithm. The continuous extension $\bar{L}$ is nondifferentiable at the integer points, and He, Fu and Marcus (2003) has shown that SPSA method converges for nondifferentiable convex continuous functions on a compact and convex domain. Therefore, for the one-dimensional case, by using the similar idea of SPSA, we have a convergent sequence.

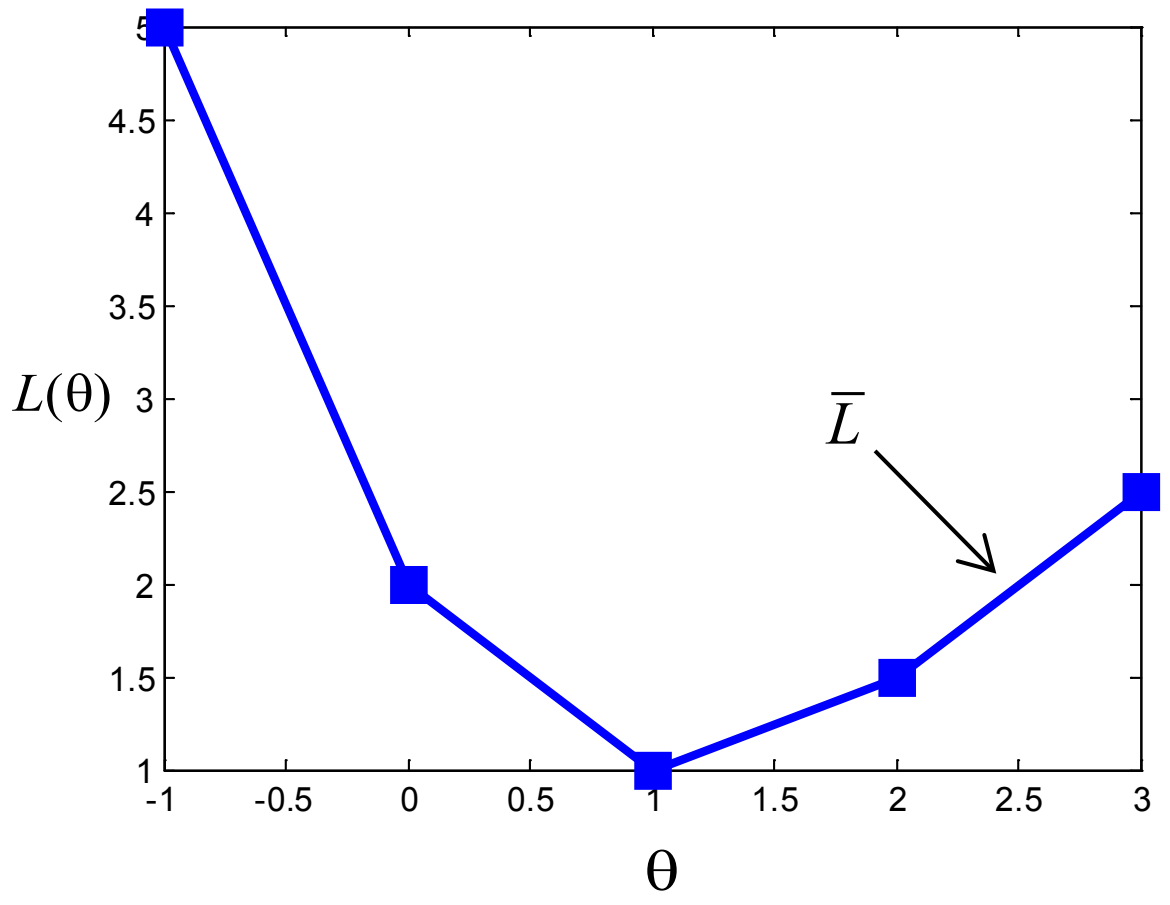

**Figure 2.1** Example of discrete function with strictly convex continuous extension $\bar{L}$.

Motivated by the one-dimensional case, we consider the case when $\boldsymbol{\theta}$ is a $p$-dimensional vector, where $p$ = 1, 2, 3, … . We have the general basic algorithm as below for the unconstrained problem, where the noisy loss function $y = L + \varepsilon$, $L$: $\mathbb{Z}^p \rightarrow \mathbb{R}$ and $\varepsilon$ is the noise. The basic algorithm is:

Step 1: Pick an initial guess $\hat{\boldsymbol{\theta}}_0 \in \mathbb{R}^p$, and set $k = 0$.

Step 2: Generate the random perturbation vector $\boldsymbol{\Delta}_k = \left[\Delta_{k1}, \Delta_{k2}, ..., \Delta_{kp}\right]^T$, where $\boldsymbol{\Delta}_k$ has a user-specified distribution satisfying conditions discussed in Section 2.2. A special case that we will focus on is when the $\Delta_{ki}$ are independent Bernoulli random variables taking the values $\pm 1$ with probability $1/2$.

Step 3: $\boldsymbol{\pi}(\hat{\boldsymbol{\theta}}_k) = \lfloor \hat{\boldsymbol{\theta}}_k \rfloor + \mathbf{1}_p/2$, where $\mathbf{1}_p$ is a $p$-dimensional vector with all components being unity, $\hat{\boldsymbol{\theta}}_k = \left[\hat{\theta}_{k1}, ..., \hat{\theta}_{kp}\right]^T$, and $\lfloor \hat{\boldsymbol{\theta}}_k \rfloor = \left[\lfloor \hat{\theta}_{k1} \rfloor, ..., \lfloor \hat{\theta}_{kp} \rfloor\right]^T$.

Step 4: Evaluate $y$ at $\hat{\boldsymbol{\theta}}_k^+ = \left[\boldsymbol{\pi}(\hat{\boldsymbol{\theta}}_k) + \boldsymbol{\Delta}_k/2\right]$ and $\hat{\boldsymbol{\theta}}_k^- = \left[\boldsymbol{\pi}(\hat{\boldsymbol{\theta}}_k) - \boldsymbol{\Delta}_k/2\right]$, where the $[\cdot]$ here is the round operator (In the following, if "$[\cdot]$" is the round operator, not simple brackets, we will note it.). Form the estimate of $\hat{\boldsymbol{g}}_k(\hat{\boldsymbol{\theta}}_k)$,

$$\hat{\boldsymbol{g}}_k(\hat{\boldsymbol{\theta}}_k) = \left(y(\hat{\boldsymbol{\theta}}_k^+) - y(\hat{\boldsymbol{\theta}}_k^-)\right)\boldsymbol{\Delta}_k^{-1},$$

where $\boldsymbol{\Delta}_k^{-1} = \left[\Delta_{k1}^{-1}, ..., \Delta_{kp}^{-1}\right]^T$. Note: since $\boldsymbol{\pi}(\hat{\boldsymbol{\theta}}_k)$ is the middle point of unit hypercube, then we have $\hat{\boldsymbol{\theta}}_k^+ \neq \hat{\boldsymbol{\theta}}_k^-$.

Step 5: Update the estimate according to the recursion

$$\hat{\boldsymbol{\theta}}_{k+1} = \hat{\boldsymbol{\theta}}_k - a_k \hat{\boldsymbol{g}}_k(\hat{\boldsymbol{\theta}}_k),$$

where $a_k = a/(1 + A + k)^{\alpha}$ is the step size, and $a > 0$, $A \geq 0$, $0.5 < \alpha \leq 1$ are the coefficients of the step size. Set $k = k + 1$. If $k < M$, go to Step 2.

Step 6: After $M$ iterations, set the approximated optimal solution to be $\left[\hat{\boldsymbol{\theta}}_M\right]$, where $M$ is the maximum number of allowed iterations based on the cost limit, and $[\cdot]$ is the round operator (each component of $\hat{\boldsymbol{\theta}}_M$ is rounded to its nearest integer value).

In the theoretical analysis below, we define $\bar{\boldsymbol{g}}(\boldsymbol{\pi}(\boldsymbol{\theta}))$ as the mean gradient-like expression centered at $\boldsymbol{\pi}(\boldsymbol{\theta})$:

$$\bar{\boldsymbol{g}}(\boldsymbol{\pi}(\boldsymbol{\theta})) \equiv E\left(\left(L(\boldsymbol{\theta}^+) - L(\boldsymbol{\theta}^-)\right)\boldsymbol{\Delta}^{-1}\middle|\boldsymbol{\theta}\right),$$

where $\boldsymbol{\Delta}$ is a $p$-dimensional vector that has the same definition as $\boldsymbol{\Delta}_k$ mentioned above. If all components of $\boldsymbol{\Delta}$ are independent Bernoulli random variables taking the values $\pm 1$ with probability $1/2$, then

$$\boldsymbol{\theta}^+ = \left[\boldsymbol{\pi}(\boldsymbol{\theta}) + \frac{\boldsymbol{\Delta}}{2}\right] = \boldsymbol{\pi}(\boldsymbol{\theta}) + \frac{\boldsymbol{\Delta}}{2},$$

$$\boldsymbol{\theta}^- = \left[\boldsymbol{\pi}(\boldsymbol{\theta}) - \frac{\boldsymbol{\Delta}}{2}\right] = \boldsymbol{\pi}(\boldsymbol{\theta}) - \frac{\boldsymbol{\Delta}}{2},$$

where $[\cdot]$ is the round operator, and $\bar{\boldsymbol{g}}(\boldsymbol{\pi}(\boldsymbol{\theta}))$ can be further written as

$$\bar{\boldsymbol{g}}(\boldsymbol{\pi}(\boldsymbol{\theta})) = \frac{1}{2^p}\sum\nolimits_{\boldsymbol{\Delta}}\left(L\left(\boldsymbol{\pi}(\boldsymbol{\theta}) + \frac{\boldsymbol{\Delta}}{2}\right) - L\left(\boldsymbol{\pi}(\boldsymbol{\theta}) - \frac{\boldsymbol{\Delta}}{2}\right)\right)\boldsymbol{\Delta},$$

where $\Sigma_{\boldsymbol{\Delta}}$ indicates the summation over all possible directions $\boldsymbol{\Delta}$. Note that $\boldsymbol{\Delta}^{-1} = \boldsymbol{\Delta}$ in the special case of Bernoulli $\pm 1$ perturbations.

Although most discussions in this thesis pertain to the unconstrained problem, let us offer some comments on the constrained problem with the lower and upper bounds in the $i$th coordinate of the feasible domain being $l_i \in \mathbb{Z}$ and $u_i \in \mathbb{Z}$, respectively. For this constrained case, we only have a finite number of feasible unit hypercubes (denote the feasible unit hypercube to be the unit hypercube with all its corner points being in the feasible domain), and the sequence $\{\hat{\boldsymbol{\theta}}_k\}$ generated by the general algorithm of DSPSA may be out of these feasible unit hypercubes. Thus, we need to modify the general algorithm to handle the bounded domain case. Suppose $\hat{\boldsymbol{\theta}}_k = \left[\hat{\theta}_{k1},...,\hat{\theta}_{kp}\right]^T$. Let $\boldsymbol{\Psi}(\hat{\boldsymbol{\theta}}_k)$ $= \left[\psi_1(\hat{\theta}_{k1}),...,\psi_p(\hat{\theta}_{kp})\right]^T$ be the projection to map the $\hat{\boldsymbol{\theta}}_k$ back to the feasible unit hypercubes, where we set

$$\psi_i(\hat{\theta}_{ki}) = \begin{cases} l_i & \hat{\theta}_{ki} < l_i \\ \hat{\theta}_{ki} & l_i \le \hat{\theta}_{ki} \le u_i \\ u_i - \tau & \hat{\theta}_{ki} > u_i, \end{cases}$$

where $\tau$ is a very small positive number (e.g. $\tau = 10^{-10}$, a simple example will be given in the following to show why $\tau$ is introduced). Due to the modified definition of $\boldsymbol{\pi}(\hat{\boldsymbol{\theta}}_k)$ = $\left\lfloor \boldsymbol{\Psi}(\hat{\boldsymbol{\theta}}_k) \right\rfloor + \mathbf{1}_p/2$, which is given in the modified step 3 below, we find that in order to make $\boldsymbol{\pi}(\hat{\boldsymbol{\theta}}_k)$ to be the middle point of one feasible unit hypercube, we have to add $\tau$ here. For example, for one-dimensional case, when the feasible domain is {0, 1} ($l_1$ = 0, $u_1$ = 1), then there is only one feasible unit hypercube with the middle point being 0.5.

Suppose $\hat{\boldsymbol{\theta}}_0 = 1.5$, then by the modified definition of $\boldsymbol{\pi}(\hat{\boldsymbol{\theta}}_0)$ we have $\boldsymbol{\pi}(\hat{\boldsymbol{\theta}}_0) = 0.5$. However, if we get rid of $\tau$, then $\boldsymbol{\pi}(\hat{\boldsymbol{\theta}}_0) = 1.5$, which is not within the only feasible unit hypercube. Besides making $\boldsymbol{\pi}(\hat{\boldsymbol{\theta}}_k)$ to be the middle point of one feasible unit hypercube, we also need to make $\hat{\boldsymbol{\theta}}_k^{\pm}$ to be feasible points. Therefore, for the bounded domain case, we modify $\hat{\boldsymbol{\theta}}_k^{\pm}$ as $\hat{\boldsymbol{\theta}}_k^{\pm} = \left[\boldsymbol{\Psi}(\boldsymbol{\pi}(\hat{\boldsymbol{\theta}}_k) \pm \boldsymbol{\Delta}_k/2)\right]$, where $[\cdot]$ is the round operator. We see that after adding the projection $\boldsymbol{\Psi}$, the modified values of $\hat{\boldsymbol{\theta}}_k^{\pm}$ are within the feasible domain. Then step 3, step 4 and step 6 of the general DSPSA can be modified as:

Step 3 (modified): Let $\boldsymbol{\pi}(\hat{\boldsymbol{\theta}}_k) = \left\lfloor \boldsymbol{\Psi}(\hat{\boldsymbol{\theta}}_k) \right\rfloor + \mathbf{1}_p/2$, where $\mathbf{1}_p$ is a $p$-dimensional vector with all components being unity, and $\left\lfloor \boldsymbol{\Psi}(\hat{\boldsymbol{\theta}}_k) \right\rfloor$ is equal to $\left[\left\lfloor \psi_1(\hat{\theta}_{k1}) \right\rfloor, \ldots, \left\lfloor \psi_p(\hat{\theta}_{kp}) \right\rfloor\right]^T$.

Step 4 (modified): Evaluate $y$ at the points $\hat{\boldsymbol{\theta}}_k^{+} = \left[\boldsymbol{\Psi}(\boldsymbol{\pi}(\hat{\boldsymbol{\theta}}_k) + \boldsymbol{\Delta}_k/2)\right]$ and $\hat{\boldsymbol{\theta}}_k^{-} = \left[\boldsymbol{\Psi}(\boldsymbol{\pi}(\hat{\boldsymbol{\theta}}_k) - \boldsymbol{\Delta}_k/2)\right]$, where $[\cdot]$ is the round operator. Form the estimate of $\hat{\boldsymbol{g}}_k(\hat{\boldsymbol{\theta}}_k)$,

$$\hat{\boldsymbol{g}}_k(\hat{\boldsymbol{\theta}}_k) = \left(y(\hat{\boldsymbol{\theta}}_k^{+}) - y(\hat{\boldsymbol{\theta}}_k^{-})\right)\boldsymbol{\Delta}_k^{-1},$$

where $\boldsymbol{\Delta}_k^{-1} = \left[\Delta_{k1}^{-1}, \ldots, \Delta_{kp}^{-1}\right]^T$. Note: since $\boldsymbol{\pi}(\hat{\boldsymbol{\theta}}_k)$ is the middle point of one feasible unit hypercube, then we have $\hat{\boldsymbol{\theta}}_k^{+} \neq \hat{\boldsymbol{\theta}}_k^{-}$.

Step 6 (modified): After $M$ iterations, set the approximated optimal solution to be $\left[\boldsymbol{\Psi}(\hat{\boldsymbol{\theta}}_M)\right]$, where $M$ is the maximum number of allowed iterations based on the cost limit and $[\cdot]$ is the round operator (as in the unconstrained algorithm).

We see that the value of $\boldsymbol{\pi}(\hat{\boldsymbol{\theta}}_k)$ is modified by the mapping, and the value of $\hat{\boldsymbol{\theta}}_k$ is not modified by the mapping, which means that the value of $\hat{\boldsymbol{\theta}}_k$ is allowed to be out of the feasible unit hypercubes. The reason to allow $\hat{\boldsymbol{\theta}}_k$ to be outside the feasible unit hypercubes is that we do not want to lose the information in $\hat{\boldsymbol{g}}_k(\hat{\boldsymbol{\theta}}_k)$ through the mapping.

## 2.2 Almost Sure Convergence

We now present an almost sure (a.s.) convergence result for $\hat{\boldsymbol{\theta}}_k$. First we introduce some definitions that will be used in the proof. For any point $\boldsymbol{\theta}$, we denote the set of middle points of all unit hypercubes containing $\boldsymbol{\theta}$ as $\mathcal{M}_{\boldsymbol{\theta}}$. If $\boldsymbol{\theta}$ lies strictly inside one unit hypercube, then $\mathcal{M}_{\boldsymbol{\theta}}$ contains only one point. However, if $\boldsymbol{\theta}$ lies on the boundary of unit hypercube, then $\mathcal{M}_{\boldsymbol{\theta}}$ contains more than one point (at most $2^p$). For any point $\boldsymbol{m}_{\boldsymbol{\theta}}$ in $\mathcal{M}_{\boldsymbol{\theta}}$, we have $|m_{\theta i} - t_i| \leq 1/2$ for $i = 1,\ldots,p$, where $\boldsymbol{\theta} = [t_1,\ldots,t_p]^T$ and $m_{\theta i}$ is the $i$th component of $\boldsymbol{m}_{\boldsymbol{\theta}}$. Furthermore, let $\Im_k = \{\hat{\boldsymbol{\theta}}_0, \hat{\boldsymbol{\theta}}_1, ..., \hat{\boldsymbol{\theta}}_k\}$, $\xi_k = \{\boldsymbol{\Delta}_0, \boldsymbol{\Delta}_1, ..., \boldsymbol{\Delta}_k\}$, $\boldsymbol{\Delta}_k^{-T} = (\boldsymbol{\Delta}_k^{-1})^T$, and $\Omega$ is the set of all possible events.

**Theorem 2.1.** Suppose $L$ is defined on $\mathbb{Z}^p$ and $L$ has unique minimal point $\boldsymbol{\theta}^*$. Assume (i) $a_k > 0$, $\lim_{k\to\infty} a_k = 0$, $\sum_{k=0}^{\infty} a_k = \infty$ and $\sum_{k=0}^{\infty} a_k^2 < \infty$; (ii) the $\{\boldsymbol{\Delta}_k\}$ are independent vectors, the components of $\boldsymbol{\Delta}_k$ are independently distributed random variables, and $\boldsymbol{\Delta}_k^{-T}\boldsymbol{\Delta}_k^{-1}$ is uniformly bounded in $k$; (iii) For all $k$, $E\left((\varepsilon_k^+ - \varepsilon_k^-) \mid \Im_k, \xi_k\right) = 0$ a.s., and $\text{var}(\varepsilon_k^{\pm})$ is uniformly bounded in $k$; (iv) $E\left(L(\hat{\boldsymbol{\theta}}_k^+) - L(\hat{\boldsymbol{\theta}}_k^-)\right)^2$ is uniformly bounded in $k$; (v) $\bar{\boldsymbol{g}}(\boldsymbol{m}_{\boldsymbol{\theta}})^T(\boldsymbol{\theta} - \boldsymbol{\theta}^*) > 0$ for all $\boldsymbol{m}_{\boldsymbol{\theta}} \in \mathcal{M}_{\boldsymbol{\theta}}$ and all $\boldsymbol{\theta} \in \mathbb{R}^p \setminus \{\boldsymbol{\theta}^*\}$; (vi) $\sup_{k\geq 0} \| \hat{\boldsymbol{\theta}}_k \| < \infty$ a.s. Then $\hat{\boldsymbol{\theta}}_k \to \boldsymbol{\theta}^*$ a.s.

*Remarks:*

1. The inner product condition (v) in Theorem 2.1 is a natural extension of the standard inner product condition for continuous problems (e.g. Spall, 2003, p. 106). The inner product condition (v) here is dependent on the distribution of the perturbation direction $\boldsymbol{\Delta}$, because the definition of $\bar{\boldsymbol{g}}(\boldsymbol{m}_{\boldsymbol{\theta}})$ is dependent on the distribution of $\boldsymbol{\Delta}$ (the value of $\bar{\boldsymbol{g}}(\boldsymbol{m}_{\boldsymbol{\theta}})$ is constructed based on both the loss function and the distribution of $\boldsymbol{\Delta}$). Moreover, the inner product condition (v) here is only a sufficient condition, which means that DSPSA may also be effective for some loss functions where condition (v) is not satisfied.

2. The combination of conditions (ii), (iii) and (iv) expresses the similar idea as the condition (B.4′′) on p. 183 of Spall (2003). Here the condition of uniformly bounded $\boldsymbol{\Delta}_k^{-T}\boldsymbol{\Delta}_k^{-1}$ is equivalent to the condition of uniformly bounded $\Delta_{ki}^{-2}$, so condition (ii) in Theorem 2.1 is an analogue of the finite inverse moments condition of standard SPSA.

3. Conditions (i) and (vi) are similar to the conditions (B.1′′) and (B.3′′) on p. 183 of Spall (2003).

*Proof.* By the algorithm, we have

$$\begin{aligned}\hat{\boldsymbol{\theta}}_{k+1} &= \hat{\boldsymbol{\theta}}_k - a_k \hat{\boldsymbol{g}}_k(\hat{\boldsymbol{\theta}}_k) \\ &= \hat{\boldsymbol{\theta}}_k - a_k \left( y(\hat{\boldsymbol{\theta}}_k^+) - y(\hat{\boldsymbol{\theta}}_k^-) \right) \boldsymbol{\Delta}_k^{-1} \\ &= \hat{\boldsymbol{\theta}}_k - a_k \left( L(\hat{\boldsymbol{\theta}}_k^+) - L(\hat{\boldsymbol{\theta}}_k^-) + \varepsilon_k^+ - \varepsilon_k^- \right) \boldsymbol{\Delta}_k^{-1}. \end{aligned} \tag{2.1}$$

Adding and subtracting $\overline{\boldsymbol{g}}(\boldsymbol{\pi}(\hat{\boldsymbol{\theta}}_k))$ to the right-hand side of eqn. (2.1), we have

$$\begin{aligned}\hat{\boldsymbol{\theta}}_{k+1} &= \hat{\boldsymbol{\theta}}_k - a_k \overline{\boldsymbol{g}}(\boldsymbol{\pi}(\hat{\boldsymbol{\theta}}_k)) \\ &\quad + a_k \left( \overline{\boldsymbol{g}}(\boldsymbol{\pi}(\hat{\boldsymbol{\theta}}_k)) - \left( L(\hat{\boldsymbol{\theta}}_k^+) - L(\hat{\boldsymbol{\theta}}_k^-) \right) \boldsymbol{\Delta}_k^{-1} \right) - a_k \left( \varepsilon_k^+ - \varepsilon_k^- \right) \boldsymbol{\Delta}_k^{-1}. \end{aligned} \tag{2.2}$$

Due to condition (vi), there exists $\Omega_1 \subseteq \Omega$, such that $P(\Omega_1) = 1$ and for any $\omega \in \Omega_1$ $\left\{ \hat{\boldsymbol{\theta}}_k(\omega) \right\}$ is a bounded sequence. Thus, there exists a subsequence $\left\{ \hat{\boldsymbol{\theta}}_{k_s}(\omega) \right\}$ and point $\boldsymbol{\theta}'(\omega)$ such that $\left\{ \hat{\boldsymbol{\theta}}_{k_s}(\omega) \right\} \to \boldsymbol{\theta}'(\omega)$ as $s \to \infty$. Then, by the recursive relationship in eqn. (2.2), we have

$$\boldsymbol{\theta}'(\omega) - \hat{\boldsymbol{\theta}}_{k_s}(\omega) = \sum_{i=k_s}^{\infty} \left( \hat{\boldsymbol{\theta}}_{i+1}(\omega) - \hat{\boldsymbol{\theta}}_i(\omega) \right)$$

$$= -\sum_{i=k_s}^{\infty} a_i \overline{\boldsymbol{g}}(\boldsymbol{\pi}(\hat{\boldsymbol{\theta}}_i(\omega))) - \sum_{i=k_s}^{\infty} a_i \left( \varepsilon_i^+(\omega) - \varepsilon_i^-(\omega) \right) \boldsymbol{\Delta}_i^{-1}(\omega)$$

$$+ \sum_{i=k_s}^{\infty} a_i \left( \overline{\boldsymbol{g}}(\boldsymbol{\pi}(\hat{\boldsymbol{\theta}}_i(\omega))) - \left( L(\hat{\boldsymbol{\theta}}_i^+(\omega)) - L(\hat{\boldsymbol{\theta}}_i^-(\omega)) \right) \boldsymbol{\Delta}_i^{-1}(\omega) \right). \tag{2.3}$$

Now let us start to discuss the terms on the right-hand side of eqn. (2.3) (suppressing $\omega$ ). Since we know $\left\{ \sum_{i=k}^{m} a_i \left( \overline{\boldsymbol{g}}(\boldsymbol{\pi}(\hat{\boldsymbol{\theta}}_i)) - \left( L(\hat{\boldsymbol{\theta}}_i^+) - L(\hat{\boldsymbol{\theta}}_i^-) \right) \boldsymbol{\Delta}_i^{-1} \right) \right\}_{m \geq k}$ is a martingale sequence, by Doob's martingale inequality (Kushner and Clark, 1978, p. 27), we have that for any $\eta > 0$

$$P\left( \sup_{m \geq k} \left\| \sum_{i=k}^{m} a_i \left( \overline{\boldsymbol{g}}(\boldsymbol{\pi}(\hat{\boldsymbol{\theta}}_i)) - \left( L(\hat{\boldsymbol{\theta}}_i^+) - L(\hat{\boldsymbol{\theta}}_i^-) \right) \boldsymbol{\Delta}_i^{-1} \right) \right\| \geq \eta \right)$$

$$\leq \eta^{-2} E \left\| \sum_{i=k}^{\infty} a_i \left( \overline{\boldsymbol{g}}(\boldsymbol{\pi}(\hat{\boldsymbol{\theta}}_i)) - \left( L(\hat{\boldsymbol{\theta}}_i^+) - L(\hat{\boldsymbol{\theta}}_i^-) \right) \boldsymbol{\Delta}_i^{-1} \right) \right\|^2 . \tag{2.4}$$

Now let us consider the right-hand side of inequality (2.4). By the definition of $\overline{\boldsymbol{g}}(\cdot)$, we have

$$\overline{\boldsymbol{g}}(\boldsymbol{\pi}(\hat{\boldsymbol{\theta}}_k)) = E\left( \left( L(\hat{\boldsymbol{\theta}}_k^+) - L(\hat{\boldsymbol{\theta}}_k^-) \right) \boldsymbol{\Delta}_k{}^{-1} \middle| \Im_k \right).$$

Then for all $i < j$, we have

$$E\left(\left(\bar{\boldsymbol{g}}(\boldsymbol{\pi}(\hat{\boldsymbol{\theta}}_i)) - (L(\hat{\boldsymbol{\theta}}_i^+) - L(\hat{\boldsymbol{\theta}}_i^-))\boldsymbol{\Delta}_i^{-1}\right)^T \left(\bar{\boldsymbol{g}}(\boldsymbol{\pi}(\hat{\boldsymbol{\theta}}_j)) - (L(\hat{\boldsymbol{\theta}}_j^+) - L(\hat{\boldsymbol{\theta}}_j^-))\boldsymbol{\Delta}_j^{-1}\right)\right)$$

$$= E\left( E\left(\left(\bar{\boldsymbol{g}}(\boldsymbol{\pi}(\hat{\boldsymbol{\theta}}_i)) - (L(\hat{\boldsymbol{\theta}}_i^+) - L(\hat{\boldsymbol{\theta}}_i^-))\boldsymbol{\Delta}_i^{-1}\right)^T \left(\bar{\boldsymbol{g}}(\boldsymbol{\pi}(\hat{\boldsymbol{\theta}}_j)) - (L(\hat{\boldsymbol{\theta}}_j^+) - L(\hat{\boldsymbol{\theta}}_j^-))\boldsymbol{\Delta}_j^{-1}\right) \middle| \Im_j, \xi_{j-1}\right)\right)$$

$$= E\left(\left(\bar{\boldsymbol{g}}(\boldsymbol{\pi}(\hat{\boldsymbol{\theta}}_i)) - (L(\hat{\boldsymbol{\theta}}_i^+) - L(\hat{\boldsymbol{\theta}}_i^-))\boldsymbol{\Delta}_i^{-1}\right)^T E\left(\left(\bar{\boldsymbol{g}}(\boldsymbol{\pi}(\hat{\boldsymbol{\theta}}_j)) - (L(\hat{\boldsymbol{\theta}}_j^+) - L(\hat{\boldsymbol{\theta}}_j^-))\boldsymbol{\Delta}_j^{-1}\right) \middle| \Im_j\right)\right)$$

$$= 0. \tag{2.5}$$

Due to conditions (i), (ii), (iv) and eqn. (2.5), we have for any $k$

$$E\left(\left\|\sum_{i=k}^{\infty} a_i \left\{\bar{\boldsymbol{g}}(\boldsymbol{\pi}(\hat{\boldsymbol{\theta}}_i)) - \left(L(\hat{\boldsymbol{\theta}}_i^+) - L(\hat{\boldsymbol{\theta}}_i^-)\right)\boldsymbol{\Delta}_i^{-1}\right\}\right\|^2\right)$$

$$= \sum_{i=k}^{\infty} a_i^2 E\left(\left\|\bar{\boldsymbol{g}}(\boldsymbol{\pi}(\hat{\boldsymbol{\theta}}_i)) - \left(L(\hat{\boldsymbol{\theta}}_i^+) - L(\hat{\boldsymbol{\theta}}_i^-)\right)\boldsymbol{\Delta}_i^{-1}\right\|^2\right)$$

$$\leq \sum_{i=k}^{\infty} a_i^2 E\left(\left(L(\hat{\boldsymbol{\theta}}_i^+) - L(\hat{\boldsymbol{\theta}}_i^-)\right)^2 \boldsymbol{\Delta}_i^{-T}\boldsymbol{\Delta}_i^{-1}\right) < \infty.$$

Thus, the inequality (2.4) can be written as

$$P\left(\sup_{m \geq k}\left\|\sum_{i=k}^{m} a_i \left(\bar{\boldsymbol{g}}(\boldsymbol{\pi}(\hat{\boldsymbol{\theta}}_i)) - \left(L(\hat{\boldsymbol{\theta}}_i^+) - L(\hat{\boldsymbol{\theta}}_i^-)\right)\boldsymbol{\Delta}_i^{-1}\right)\right\| \geq \eta\right)$$

$$\leq \eta^{-2} \sum_{i=k}^{\infty} a_i^2 E\left(\left\|\bar{\boldsymbol{g}}(\boldsymbol{\pi}(\hat{\boldsymbol{\theta}}_i)) - \left(L(\hat{\boldsymbol{\theta}}_i^+) - L(\hat{\boldsymbol{\theta}}_i^-)\right)\boldsymbol{\Delta}_i^{-1}\right\|^2\right)$$

$$< \infty,$$

from which follows that for any $n > 0$

$$P\left(\sup_{k\geq n}\left\|\sum_{i=k}^{\infty}a_i\left(\bar{\boldsymbol{g}}(\boldsymbol{\pi}(\hat{\boldsymbol{\theta}}_i))-\left(L(\hat{\boldsymbol{\theta}}_i^+)-L(\hat{\boldsymbol{\theta}}_i^-)\right)\boldsymbol{\Delta}_i^{-1}\right)\right\|\geq\eta\right)$$

$$\leq P\left(\sup_{\substack{m\geq k\\k\geq n}}\left\|\sum_{i=k}^{m}a_i\left(\bar{\boldsymbol{g}}(\boldsymbol{\pi}(\hat{\boldsymbol{\theta}}_i))-\left(L(\hat{\boldsymbol{\theta}}_i^+)-L(\hat{\boldsymbol{\theta}}_i^-)\right)\boldsymbol{\Delta}_i^{-1}\right)\right\|\geq\eta\right)$$

$$\leq\eta^{-2}\sum_{i=n}^{\infty}{a_i}^2E\left(\left\|\bar{\boldsymbol{g}}(\boldsymbol{\pi}(\hat{\boldsymbol{\theta}}_i))-\left(L(\hat{\boldsymbol{\theta}}_i^+)-L(\hat{\boldsymbol{\theta}}_i^-)\right)\boldsymbol{\Delta}_i^{-1}\right\|^2\right)$$

$$<\infty.$$

Thus, we have

$$\lim_{n\to\infty}P\left(\sup_{k\geq n}\left\|\sum_{i=k}^{\infty}a_i\left(\bar{\boldsymbol{g}}(\boldsymbol{\pi}(\hat{\boldsymbol{\theta}}_i))-\left(L(\hat{\boldsymbol{\theta}}_i^+)-L(\hat{\boldsymbol{\theta}}_i^-)\right)\boldsymbol{\Delta}_i^{-1}\right)\right\|\geq\eta\right)=0. \tag{2.6}$$

By Theorem 4.1.1 in Chung (2001) on almost sure convergence, due to eqn. (2.6), we have

$$\lim_{k\to\infty}\sum_{i=k}^{\infty}a_i\left(\bar{\boldsymbol{g}}(\boldsymbol{\pi}(\hat{\boldsymbol{\theta}}_i))-\left(L(\hat{\boldsymbol{\theta}}_i^+)-L(\hat{\boldsymbol{\theta}}_i^-)\right)\boldsymbol{\Delta}_i^{-1}\right)=\mathbf{0}\text{ a.s.}$$

By similar arguments, we know $\left\{\sum_{i=k}^{m}a_i(\varepsilon_i^+-\varepsilon_i^-)\boldsymbol{\Delta}_i^{-1}\right\}_{m\geq k}$ is a martingale. Through condition (i), (ii), (iii) and Doob's martingale inequality (Kushner and Clark, 1978, p. 27), we have

$$P\left(\sup_{m\ge k}\left\|\sum_{i=k}^{m} a_i(\varepsilon_i^+ - \varepsilon_i^-)\mathbf{\Delta}_i^{-1}\right\| \ge \eta\right) \le \eta^{-2} E\left\|\sum_{i=k}^{\infty} a_i(\varepsilon_i^+ - \varepsilon_i^-)\mathbf{\Delta}_i^{-1}\right\|^2$$

$$= \eta^{-2}\sum_{i=k}^{\infty} a_i^2 E\left\|(\varepsilon_i^+ - \varepsilon_i^-)\mathbf{\Delta}_i^{-1}\right\|^2$$

$$< \infty,$$

which indicates that for any $n > 0$

$$P\left(\sup_{k\ge n}\left\|\sum_{i=k}^{\infty} a_i(\varepsilon_i^+ - \varepsilon_i^-)\mathbf{\Delta}_i^{-1}\right\| \ge \eta\right) \le P\left(\sup_{\substack{m\ge k\\ k\ge n}}\left\|\sum_{i=k}^{m} a_i(\varepsilon_i^+ - \varepsilon_i^-)\mathbf{\Delta}_i^{-1}\right\| \ge \eta\right)$$

$$\le \eta^{-2}\sum_{i=n}^{\infty} a_i^2 E\left\|(\varepsilon_i^+ - \varepsilon_i^-)\mathbf{\Delta}_i^{-1}\right\|^2$$

$$< \infty,$$

leading to

$$\lim_{n\to\infty} P\left(\sup_{k\ge n}\left\|\sum_{i=k}^{\infty} a_i(\varepsilon_i^+ - \varepsilon_i^-)\mathbf{\Delta}_i^{-1}\right\| \ge \eta\right) = 0. \tag{2.7}$$

Thus, by Theorem 4.1.1 in Chung (2001), due to eqn. (2.7), we have

$$\lim_{k\to\infty}\sum_{i=k}^{\infty} a_i(\varepsilon_i^+ - \varepsilon_i^-)\mathbf{\Delta}_i^{-1} = \mathbf{0} \text{ a.s.}$$

Overall, there exists $\Omega_2 \subseteq \Omega$ such that $P(\Omega_2) = 1$, and for any $\omega \in \Omega_2$ we have $\sum_{i=k}^{\infty} a_i\left(\bar{\boldsymbol{g}}(\boldsymbol{\pi}(\hat{\boldsymbol{\theta}}_i(\omega))) - \left(L(\hat{\boldsymbol{\theta}}_i^+(\omega)) - L(\hat{\boldsymbol{\theta}}_i^-(\omega))\right)\mathbf{\Delta}_i^{-1}(\omega)\right) \to \mathbf{0}$ as $k \to \infty$. There also exists $\Omega_3 \subseteq \Omega$, such that $P(\Omega_3) = 1$, and for any $\omega \in \Omega_3$, $\sum_{i=k}^{\infty} a_i(\varepsilon_i^+(\omega) - \varepsilon_i^-(\omega))\mathbf{\Delta}_i^{-1}(\omega) \to \mathbf{0}$ as $k \to \infty$. In all, we have $P(\Omega_1 \cap \Omega_2 \cap \Omega_3) = 1$, and for any $\omega \in \Omega_1 \cap \Omega_2 \cap \Omega_3$, the terms in eqn. (2.3) have the results: as $s \to \infty$, we have that $\boldsymbol{\theta}'(\omega) - \hat{\boldsymbol{\theta}}_{k_s}(\omega) \to \mathbf{0}$, $\sum_{i=k_s}^{\infty} a_i \left( \bar{\boldsymbol{g}}(\boldsymbol{\pi}(\hat{\boldsymbol{\theta}}_i(\omega))) - \left( L(\hat{\boldsymbol{\theta}}_i^+(\omega)) - L(\hat{\boldsymbol{\theta}}_i^-(\omega)) \right) \boldsymbol{\Delta}_i^{-1}(\omega) \right) \to \mathbf{0}$, $\sum_{i=k_s}^{\infty} a_i \left( \varepsilon_i^+(\omega) - \varepsilon_i^-(\omega) \right) \boldsymbol{\Delta}_i^{-1}(\omega)$ $\to \mathbf{0}$, which implies that

$$\sum_{i=k_s}^{\infty} a_i \bar{\boldsymbol{g}}(\boldsymbol{\pi}(\hat{\boldsymbol{\theta}}_i(\omega))) \to \mathbf{0} \text{ as } s \to \infty. \tag{2.8}$$

Because $\hat{\boldsymbol{\theta}}_{k_s}(\omega) \to \boldsymbol{\theta}'(\omega)$, then for any $\delta > 0$, there exists $S > 0$, such that when $s > S$, $\left\| \hat{\boldsymbol{\theta}}_{k_s}(\omega) - \boldsymbol{\theta}'(\omega) \right\| < \delta$. Thus, there exists $S'$ such that when $s > S'$, all $\boldsymbol{\pi}(\hat{\boldsymbol{\theta}}_{k_s}(\omega)) \in \mathcal{M}_{\boldsymbol{\theta}'}$.

We now show that $\boldsymbol{\theta}'(\omega)$ is the optimal point $\boldsymbol{\theta}^*$. By contradiction, suppose $\boldsymbol{\theta}'(\omega)$ is not the optimal solution. Then, by condition (v), we have $\bar{\boldsymbol{g}}(\boldsymbol{m}_{\boldsymbol{\theta}'})^T (\boldsymbol{\theta}'(\omega) - \boldsymbol{\theta}^*) > 0$ for all $\boldsymbol{m}_{\boldsymbol{\theta}'} \in \mathcal{M}_{\boldsymbol{\theta}'}$ (at most $2^p$ points in $\mathcal{M}_{\boldsymbol{\theta}'}$), and it indicates that there exists a constant $\upsilon > 0$ such that for all $\boldsymbol{m}_{\boldsymbol{\theta}'} \in \mathcal{M}_{\boldsymbol{\theta}'}$, $\bar{\boldsymbol{g}}(\boldsymbol{m}_{\boldsymbol{\theta}'})^T (\boldsymbol{\theta}'(\omega) - \boldsymbol{\theta}^*) > \upsilon$. Moreover, we have $\sum_{i=k_s}^{\infty} a_i$ = $\infty$. Therefore, we have

$$\sum_{i=k_s}^{\infty} a_i \bar{\boldsymbol{g}}(\boldsymbol{\pi}(\hat{\boldsymbol{\theta}}_i(\omega)))^T (\boldsymbol{\theta}'(\omega) - \boldsymbol{\theta}^*) \to \infty \text{ as } s \to \infty.$$

But, at the same time from eqn. (2.8) we have

$$\sum_{i=k_s}^{\infty} a_i \bar{\boldsymbol{g}}(\boldsymbol{\pi}(\hat{\boldsymbol{\theta}}_i(\omega)))^T (\boldsymbol{\theta}'(\omega) - \boldsymbol{\theta}^*) \to 0 \text{ as } s \to \infty,$$

which is a contradiction. Then, for all $\omega \in \Omega_1 \cap \Omega_2 \cap \Omega_3$, the limiting points of all convergent subsequence $\left\{\hat{\boldsymbol{\theta}}_{k_s}(\omega)\right\}$ is equal to $\boldsymbol{\theta}^*$, which indicates that the bounded sequences $\left\{\hat{\boldsymbol{\theta}}_k(\omega)\right\}$ only has one cluster point $\boldsymbol{\theta}^*$ for all $\omega \in \Omega_1 \cap \Omega_2 \cap \Omega_3$. Therefore, we have $\hat{\boldsymbol{\theta}}_k$ converges to $\boldsymbol{\theta}^*$ a.s. Q.E.D.

In the following part of this section, we further discuss condition (v) in Theorem 2.1 under the case when the $\Delta_{ki}$ are independent Bernoulli random variables taking the values $\pm 1$ with probability $1/2$. Except for condition (v), all the other conditions in Theorem 2.1 are general and have similar ideas as those conditions in "standard" (continuous) SPSA. Therefore, we will only focus on the discussions of condition (v) in the following.

First let us discuss the relationship between condition (v) in Theorem 2.1 and discrete convexity. Miller (1971) is a pioneer in the early 1970s in the area of discrete convex function. Miller (1971) introduced a definition of discrete convex function and showed that the local optimal points for discrete convex function are also global optimal solutions. The definition in Miller (1971) is: the function $f: X \to \mathbb{R}$ is discretely convex if $X$ is a discrete rectangle and given $\boldsymbol{x}^1, \boldsymbol{x}^2 \in X$ and $0 \le \lambda \le 1$

$$\min_{\boldsymbol{x} \in N(\lambda \boldsymbol{x}^1 + (1-\lambda)\boldsymbol{x}^2)} f(\boldsymbol{x}) \le \lambda f(\boldsymbol{x}^1) + (1-\lambda) f(\boldsymbol{x}^2),$$

where $N(\lambda \boldsymbol{x}^1 + (1-\lambda)\boldsymbol{x}^2)$ is the discrete neighborhood of the point $\lambda \boldsymbol{x}^1 + (1-\lambda)\boldsymbol{x}^2$, and $N(\lambda \boldsymbol{x}^1 + (1-\lambda)\boldsymbol{x}^2) = \left\{\boldsymbol{x} : \boldsymbol{x} \in X, \left\|\boldsymbol{x} - (\lambda \boldsymbol{x}^1 + (1-\lambda)\boldsymbol{x}^2)\right\| < 1\right\}$ (the norm in the braces is the standard Euclidean norm). There are other definitions of discrete convex functions, such as Favati and Tardellan (1990), Murota (1998), Murota and Shioura (1999), and Fujishige and Murota (2000). Murota and Shioura (2001) show that Miller's discrete convexity contains the other classes of discrete convexities based on definitions.

Condition (v) in Theorem 2.1 (for Bernoulli $\pm 1$ case) is not a form of discrete convexity, but it has some connections to Miller's definition of convexity. When $p = 1$, discrete convex functions satisfying Miller's definition of discrete convexity also satisfy condition (v) in Theorem 2.1. However, for the higher dimensional case, the set of functions that satisfy the Miller's definition of discrete convexity does not include all functions that satisfy condition (v) in Theorem 2.1. Meanwhile, the set of functions that satisfy condition (v) also does not include all functions that satisfy Miller's discrete convexity.

In Figure 2.2, we see two simple functions where the first function satisfies condition (v) but does not satisfy Miller's definition, and the second function satisfies Miller's definition but does not satisfy condition (v). In both plots of Figure 2.2, suppose the middle point of the left unit hypercube is $\boldsymbol{m}_l$, and the middle point of the right unit hypercube is $\boldsymbol{m}_r$. In plot (a) of Figure 2.2, we know $0.5 f(\boldsymbol{x}^1) + 0.5 f(\boldsymbol{x}^2)$ = 0.75. However, $\min_{\boldsymbol{x} \in N(0.5\boldsymbol{x}^1 + 0.5\boldsymbol{x}^2)} f(\boldsymbol{x})$ = 0.8 > 0.75. Thus, the simple function in plot (a) does not satisfy Miller's definition of discrete convexity. Meanwhile, it is easy to check that condition (v) is satisfied for the function in plot (a). In plot (b) of Figure 2.2, for the right unit hypercube and the point $\boldsymbol{x}$, we have $\bar{\boldsymbol{g}}(\boldsymbol{m}_r)^T(\boldsymbol{x}-\boldsymbol{\theta}^*) = -0.2 < 0$, which means condition (v) is not satisfied. However, it is easy to check the function in plot (b) satisfies Miller's definition of discrete convexity.

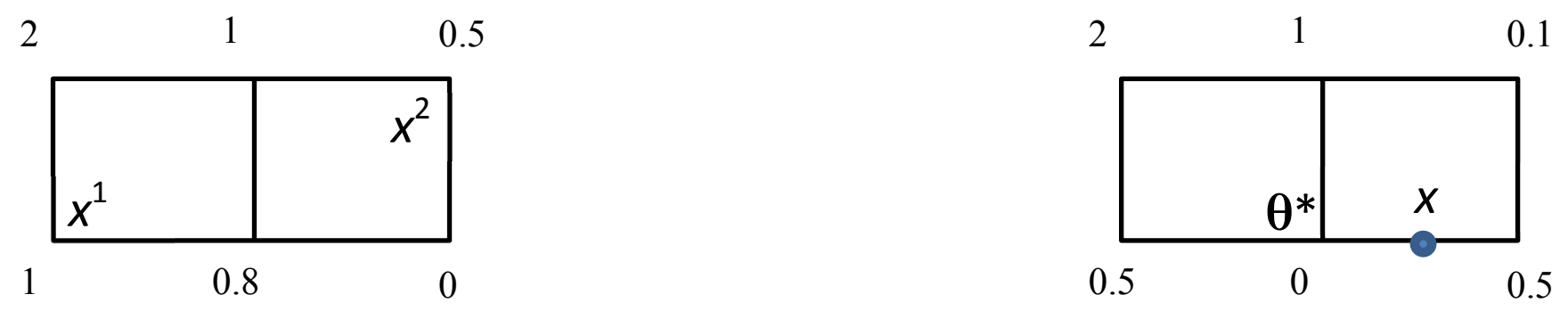

(a) (b)

**Figure 2.2** Examples of the relationship between condition (v) in Theorem 2.1 (for Bernoulli $\pm 1$ case) and Miller's definition of discrete convexity. Plot (a) is the simple function that satisfies condition (v) in Theorem 2.1, but does not satisfy Miller's definition of discrete convexity. Plot (b) is the simple function that satisfies Miller's definition of discrete convexity, but does not satisfy condition (v) in Theorem 2.1.

In the propositions below, we use three common classes of discrete functions to show that condition (v) in Theorem 2.1 is reasonable (it is the reason why these propositions are useful). Since in this part we assume that the $\Delta_{ki}$ are independent Bernoulli random variables taking the values $\pm 1$ with probability $1/2$, we know that $\boldsymbol{\theta}^+ = \boldsymbol{\pi}(\boldsymbol{\theta}) + \boldsymbol{\Delta}/2$, and $\boldsymbol{\theta}^- = \boldsymbol{\pi}(\boldsymbol{\theta}) - \boldsymbol{\Delta}/2$. For any discrete function $L$, there are an infinite number of continuous extensions $\bar{L}$ such that $\bar{L}$ is a function defined on $\mathbb{R}^p$, and for all $\boldsymbol{\theta} \in \mathbb{Z}^p$, we have $\overline{L}(\boldsymbol{\theta}) = L(\boldsymbol{\theta})$. Here we use $\overline{L}$ to restrict the class of $L$. Suppose $\boldsymbol{\theta} = \left[t_1,...,t_p\right]^T$ and the optimal solution of $L$ is $\boldsymbol{\theta}^* = \left[t_1^*,...,t_p^*\right]^T$. The three classes of discrete functions that we will discuss have very general forms of continuous extensions, which are separable function, piecewise linear function and quadratic function. The separable function is defined as $\overline{L}(\boldsymbol{\theta}) = \sum_{i=1}^{p} \overline{L}_i(t_i)$. The piecewise linear function that we consider is assumed to be linear ($\overline{L}(\boldsymbol{\theta}) = \boldsymbol{B}_{\boldsymbol{m}_\theta}^T\boldsymbol{\theta} + \boldsymbol{C}_{\boldsymbol{m}_\theta}$) within each unit hypercube (centered by $\boldsymbol{m}_\theta$), and this function may not be separable. The quadratic function is defined as $\overline{L}(\boldsymbol{\theta}) = \boldsymbol{\theta}^T\boldsymbol{A}\boldsymbol{\theta} + \boldsymbol{B}^T\boldsymbol{\theta} + \boldsymbol{C}$. Discrete functions with strictly convex separable continuous extensions are considered in Proposition 2.1 and this kind of function has also been discussed in Hill et al. (2004).

**Proposition 2.1.** Suppose $L$ is a discrete function having a continuous extension $\overline{L}$ that is a strictly convex separable function with the minimal value at the multivariate integer point $\boldsymbol{\theta}^*$. Then $L$ satisfies condition (v) in Theorem 2.1, where the $\Delta_{ki}$ are independent Bernoulli random variables taking the values $\pm 1$ with probability $1/2$.

*Proof.* Since $\boldsymbol{\theta}^*$ is the optimal point of $\overline{L}$. According to Theorem 23.4 in Rockafellar (1970), we know that the subgradient (defined in Section 2.1) always exists at all points in the domain for a continuous convex function. Since $\overline{L}$ is strictly convex, then for all $\boldsymbol{\theta} \in \mathbb{R}^p \setminus \{\boldsymbol{\theta}^*\}$ and any subgradient $g_i(t_i) \in \partial\overline{L}_i(t_i)$, we have $g_i(t_i)(t_i - t_i^*) > 0$ for $i$ = 1, …, $p$. Moreover, suppose $\boldsymbol{\Delta} = \left[\delta_1,...,\delta_p\right]^T$. For any $\boldsymbol{m}_\theta = [m_{\theta 1},...,m_{\theta p}]^T \in \mathcal{M}_\theta$, since the $\Delta_{ki}$ are independent Bernoulli random variables taking the values $\pm 1$ with probability $1/2$, we have

$$\begin{aligned}
\bar{\boldsymbol{g}}(\boldsymbol{m}_{\boldsymbol{\theta}}) &= \frac{1}{2^p}\sum\nolimits_{\boldsymbol{\Delta}}\left( L\left( \boldsymbol{m}_{\boldsymbol{\theta}} + \frac{\boldsymbol{\Delta}}{2}\right) - L\left( \boldsymbol{m}_{\boldsymbol{\theta}} - \frac{\boldsymbol{\Delta}}{2}\right)\right)\boldsymbol{\Delta}^{-1} \\
&= \frac{1}{2^p}\sum\nolimits_{\boldsymbol{\Delta}}\left( \sum_{i=1}^{p}\left( L_i\left( m_{\theta i} + \frac{\delta_i}{2}\right) - L_i\left( m_{\theta i} - \frac{\delta_i}{2}\right)\right)\right)\boldsymbol{\Delta}^{-1} \\
&= \sum_{i=1}^{p}\frac{1}{2^p}\sum\nolimits_{\boldsymbol{\Delta}}\left( L_i\left( m_{\theta i} + \frac{\delta_i}{2}\right) - L_i\left( m_{\theta i} - \frac{\delta_i}{2}\right)\right)\boldsymbol{\Delta}^{-1} \\
&= \sum_{i=1}^{p}\left( L_i\left( m_{\theta i} + \frac{1}{2}\right) - L_i\left( m_{\theta i} - \frac{1}{2}\right)\right)\boldsymbol{e}_i.
\end{aligned}$$

Then, we have $\bar{\boldsymbol{g}}(\boldsymbol{m}_{\boldsymbol{\theta}})^T(\boldsymbol{\theta}-\boldsymbol{\theta}^*)$ = $\sum_{i=1}^{p}\left(L_i(m_{\theta i}+1/2)-L_i(m_{\theta i}-1/2)\right)(t_i-t_i^*)$. The minimal point is a multivariate integer point and $\bar{L}_i$ is strictly convex. Therefore, $\bar{L}_i(m_{\theta i}+1/2)-\bar{L}_i(m_{\theta i}-1/2)$ has the same sign with the subgradient of $\bar{L}_i$ at $t_i$, indicating that $\left(\bar{L}_i(m_{\theta i}+1/2)-\bar{L}_i(m_{\theta i}-1/2)\right)(t_i-t_i^*)>0$ for all $i$ = 1, …, $p$. Thus $\bar{\boldsymbol{g}}(\boldsymbol{m}_{\boldsymbol{\theta}})^T(\boldsymbol{\theta}-\boldsymbol{\theta}^*)>0$ for all $\boldsymbol{\theta} \in \mathbb{R}\setminus\{\boldsymbol{\theta}^*\}$ and $\boldsymbol{m}_{\boldsymbol{\theta}} \in \mathcal{M}_{\boldsymbol{\theta}}$. Q.E.D.

In Proposition 2.2, we consider the loss functions that have piecewise linear continuous extension functions.

**Proposition 2.2.** Suppose $L$ is a discrete function having a continuous extension $\bar{L}$ that is strictly convex and linear in each unit hypercube. Then $L$ satisfies condition (v) in Theorem 2.1, where the $\Delta_{ki}$ are independent Bernoulli random variables taking the values $\pm 1$ with probability $1/2$.

*Remark*: Since the continuous extension $\overline{L}$ is linear in each unit hypercube, then the optimal solution of $\overline{L}$ is at a multivariate integer point $\boldsymbol{\theta}^*$.

*Proof.* For any $\boldsymbol{m}_{\boldsymbol{\theta}} \in \mathcal{M}_{\boldsymbol{\theta}}$, since the $\Delta_{ki}$ are independent Bernoulli random variables taking the values $\pm 1$ with probability $1/2$, then

$$\begin{aligned}\overline{\boldsymbol{g}}(\boldsymbol{m}_{\boldsymbol{\theta}}) &= \frac{1}{2^p}\sum\nolimits_{\boldsymbol{\Delta}}\left(L\left(\boldsymbol{m}_{\boldsymbol{\theta}} + \frac{\boldsymbol{\Delta}}{2}\right) - L\left(\boldsymbol{m}_{\boldsymbol{\theta}} - \frac{\boldsymbol{\Delta}}{2}\right)\right)\boldsymbol{\Delta}^{-1} \\ &= \frac{1}{2^p}\sum\nolimits_{\boldsymbol{\Delta}}\left(\overline{L}\left(\boldsymbol{m}_{\boldsymbol{\theta}} + \frac{\boldsymbol{\Delta}}{2}\right) - \overline{L}\left(\boldsymbol{m}_{\boldsymbol{\theta}} - \frac{\boldsymbol{\Delta}}{2}\right)\right)\boldsymbol{\Delta}^{-1} \\ &= \frac{1}{2^p}\sum\nolimits_{\boldsymbol{\Delta}}\left(\nabla\overline{L}(\boldsymbol{m}_{\boldsymbol{\theta}})^T\boldsymbol{\Delta}\right)\boldsymbol{\Delta}^{-1},\end{aligned}$$

where $\nabla\overline{L}(\boldsymbol{m}_{\boldsymbol{\theta}})$ is the gradient of $\overline{L}$ at point $\boldsymbol{m}_{\boldsymbol{\theta}}$, and where $\nabla\overline{L}(\boldsymbol{m}_{\boldsymbol{\theta}})$ is known to exist since $\overline{L}$ is differentiable at $\boldsymbol{m}_{\boldsymbol{\theta}}$. Thus,

$$\begin{aligned}\overline{\boldsymbol{g}}(\boldsymbol{m}_{\boldsymbol{\theta}})^T(\boldsymbol{\theta} - \boldsymbol{\theta}^*) &= \frac{1}{2^p}\sum\nolimits_{\boldsymbol{\Delta}}\left(\nabla\overline{L}(\boldsymbol{m}_{\boldsymbol{\theta}})^T\boldsymbol{\Delta}\right)\boldsymbol{\Delta}^{-T}(\boldsymbol{\theta} - \boldsymbol{\theta}^*) \\ &= \nabla\overline{L}(\boldsymbol{m}_{\boldsymbol{\theta}})^T\frac{1}{2^p}\sum\nolimits_{\boldsymbol{\Delta}}\boldsymbol{\Delta}\boldsymbol{\Delta}^{-T}(\boldsymbol{\theta} - \boldsymbol{\theta}^*) \\ &= \nabla\overline{L}(\boldsymbol{m}_{\boldsymbol{\theta}})^T(\boldsymbol{\theta} - \boldsymbol{\theta}^*),\end{aligned}$$

where $\boldsymbol{\Delta}^{-T} = (\boldsymbol{\Delta}^{-1})^T$. In addition, $\overline{L}$ is linear in each unit hypercube, so for any $\boldsymbol{m}_{\boldsymbol{\theta}} \in \mathcal{M}_{\boldsymbol{\theta}}$ we have $\nabla\overline{L}(\boldsymbol{m}_{\boldsymbol{\theta}}) \in \partial\overline{L}(\boldsymbol{\theta})$. Since $\overline{L}$ is a strictly convex function, then for all $\boldsymbol{\theta} \in \mathbb{R}^p \setminus \{\boldsymbol{\theta}^*\}$, we have $\nabla\overline{L}(\boldsymbol{m}_{\boldsymbol{\theta}})^T(\boldsymbol{\theta} - \boldsymbol{\theta}^*) > 0$, which indicates that $\overline{\boldsymbol{g}}(\boldsymbol{m}_{\boldsymbol{\theta}})^T(\boldsymbol{\theta} - \boldsymbol{\theta}^*) > 0$, for all $\boldsymbol{\theta} \in \mathbb{R} \setminus \{\boldsymbol{\theta}^*\}$ and $\boldsymbol{m}_{\boldsymbol{\theta}} \in \mathcal{M}_{\boldsymbol{\theta}}$. Q.E.D.

In Proposition 2.3, we consider the loss functions that have quadratic continuous extension functions. Before discussing Proposition 2.3, let us first define the strictly diagonal dominant matrix. The matrix $\boldsymbol{A} = \left[ a_{ij} \right]_{p \times p}$ is strictly diagonal dominant if for all $i$, $| a_{ii} | > \sum_{j \neq i} | a_{ij} |$.

**Proposition 2.3.** Suppose $L(\boldsymbol{\theta})$ is a discrete function having a quadratic continuous extension $\overline{L}(\boldsymbol{\theta}) = \boldsymbol{\theta}^T \boldsymbol{A}\boldsymbol{\theta} + \boldsymbol{B}^T \boldsymbol{\theta} + \boldsymbol{C}$, where $\boldsymbol{A}$ is a symmetric strictly diagonal dominant matrix with positive diagonal values. Suppose that the minimal value of $\overline{L}$ is at multivariate integer point $\boldsymbol{\theta}^*$. Then, $L$ satisfies condition (v) in Theorem 2.1, where the $\Delta_{ki}$ are independent Bernoulli random variables taking the values $\pm 1$ with probability $1/2$.

*Remark*: Because $\boldsymbol{A}$ is a symmetric strictly diagonal dominant matrix with positive diagonal values, it is known that $\boldsymbol{A}$ is positive definite. Hence $\overline{L}$ is a strictly convex function.

*Proof.* The gradient of $\overline{L}$ at point $\boldsymbol{\theta}$ is $2\boldsymbol{A}\boldsymbol{\theta} + \boldsymbol{B}$. Since $\boldsymbol{\theta}^*$ is also the optimal point of $\overline{L}$, then we have $2\boldsymbol{A}\boldsymbol{\theta}^* + \boldsymbol{B} = \boldsymbol{0}$. Furthermore, since the $\Delta_{ki}$ are independent Bernoulli random variables taking the values $\pm 1$ with probability $1/2$, then

$$\begin{aligned}\overline{\boldsymbol{g}}(\boldsymbol{m}_{\boldsymbol{\theta}}) &= \frac{1}{2^p} \sum_{\boldsymbol{\Delta}} \left( L\left( \boldsymbol{m}_{\boldsymbol{\theta}} + \frac{\boldsymbol{\Delta}}{2} \right) - L\left( \boldsymbol{m}_{\boldsymbol{\theta}} - \frac{\boldsymbol{\Delta}}{2} \right) \right) \boldsymbol{\Delta}^{-1} \\ &= \frac{1}{2^p} \sum_{\boldsymbol{\Delta}} (2\boldsymbol{m}_{\boldsymbol{\theta}}{}^T \boldsymbol{A}\boldsymbol{\Delta} + \boldsymbol{B}^T \boldsymbol{\Delta}) \boldsymbol{\Delta}^{-1}.\end{aligned}$$

Then for all $\boldsymbol{\theta} \in \mathbb{R} \setminus \{\boldsymbol{\theta}^*\}$, we have

$$
\begin{aligned}
\bar{\boldsymbol{g}}(\boldsymbol{m}_{\boldsymbol{\theta}})^T (\boldsymbol{\theta} - \boldsymbol{\theta}^*) &= \frac{1}{2^p} \sum\nolimits_{\boldsymbol{\Delta}} \left( 2\boldsymbol{m}_{\boldsymbol{\theta}}{}^T \boldsymbol{A}\boldsymbol{\Delta} + \boldsymbol{B}^T \boldsymbol{\Delta} \right) \boldsymbol{\Delta}^{-T} (\boldsymbol{\theta} - \boldsymbol{\theta}^*) \\
&= \frac{1}{2^p} \sum\nolimits_{\boldsymbol{\Delta}} \left( 2\boldsymbol{A}\boldsymbol{m}_{\boldsymbol{\theta}} + \boldsymbol{B} \right)^T \boldsymbol{\Delta}\boldsymbol{\Delta}^{-T} (\boldsymbol{\theta} - \boldsymbol{\theta}^*) \\
&= \left( 2\boldsymbol{A}\boldsymbol{m}_{\boldsymbol{\theta}} + \boldsymbol{B} \right)^T \frac{1}{2^p} \left( \sum\nolimits_{\boldsymbol{\Delta}} \boldsymbol{\Delta}\boldsymbol{\Delta}^{-T} \right) (\boldsymbol{\theta} - \boldsymbol{\theta}^*) \\
&= \left( 2\boldsymbol{A}\boldsymbol{m}_{\boldsymbol{\theta}} + \boldsymbol{B} \right)^T (\boldsymbol{\theta} - \boldsymbol{\theta}^*) \\
&= \left( 2\boldsymbol{A}\boldsymbol{m}_{\boldsymbol{\theta}} - 2\boldsymbol{A}\boldsymbol{\theta}^* \right)^T (\boldsymbol{\theta} - \boldsymbol{\theta}^*) \\
&= 2\left( \boldsymbol{m}_{\boldsymbol{\theta}} - \boldsymbol{\theta}^* \right)^T \boldsymbol{A}(\boldsymbol{\theta} - \boldsymbol{\theta}^*),
\end{aligned}
$$

where $\boldsymbol{\Delta}^{-T} = (\boldsymbol{\Delta}^{-1})^T$. Let $\boldsymbol{x} = \boldsymbol{m}_{\boldsymbol{\theta}} - \boldsymbol{\theta}^*$ and $\boldsymbol{y} = \boldsymbol{\theta} - \boldsymbol{\theta}^*$. Suppose $\boldsymbol{x} = (x_1, ..., x_p)^T$ and $\boldsymbol{y} = (y_1, ..., y_p)^T$. Because $\boldsymbol{m}_{\boldsymbol{\theta}}$ is the middle point of the unit hypercube where $\boldsymbol{\theta}$ is located in, each component of $\boldsymbol{m}_{\boldsymbol{\theta}}$ equals the sum of an integer and $\pm 0.5$, implying that for any $i, j \in \{1, ..., p\}$ we have $x_i \pm x_j$ are integers. Moreover, for all $i$, we have $-1/2 + x_i \le y_i \le 1/2 + x_i$, which indicates that $x_i y_i \ge 0$. Furthermore, since $\boldsymbol{\theta} \ne \boldsymbol{\theta}^*$, then there exists at least one $i' \in \{1, ..., p\}$ such that $x_{i'} y_{i'} > 0$. Since $\boldsymbol{A} = \left[ a_{ij} \right]_{p \times p}$ and $\boldsymbol{A}$ is a symmetric strictly diagonal dominant matrix with positive diagonal elements, then for all $i, j \in \{1, ..., p\}$, we have $a_{ii} > 0$, $a_{ij} = a_{ji}$, $| a_{ii} | > \sum_{j \ne i} | a_{ij} |$. Thus,

$$\sum\nolimits_{i=1}^{p} (a_{ii} - \sum\nolimits_{j \ne i} | a_{ij} |) x_i y_i > 0. \tag{2.9}$$

Let $s(a_{ij}) \in \{-1, 0, 1\}$ be the sign of $a_{ij}$. Then,

$2\boldsymbol{x}^T\boldsymbol{Ay}$

$$=2\sum_{i=1}^{p}\sum_{j=1}^{p}a_{ij}x_iy_j$$

$$=2\sum_{i=1}^{p}\left(a_{ii}x_iy_i+\sum_{j\neq i}a_{ij}x_iy_j\right)$$

$$=2\sum_{i=1}^{p}\left(a_{ii}x_iy_i+\sum_{j\neq i}(-|a_{ij}|)x_iy_i+\sum_{j\neq i}|a_{ij}|x_iy_i+\sum_{j\neq i}a_{ij}x_iy_j\right)$$

$$=2\sum_{i=1}^{p}\left(a_{ii}-\sum_{j\neq i}|a_{ij}|\right)x_iy_i+2\sum_{i=1}^{p}\left(\sum_{j\neq i}|a_{ij}|x_i(y_i+s(a_{ij})y_j)\right)$$

$$=2\sum_{i=1}^{p}\left(a_{ii}-\sum_{j\neq i}|a_{ij}|\right)x_iy_i+2\sum_{i=1}^{p}\left(\sum_{j>i}\big(|a_{ij}|x_i(y_i+s(a_{ij})y_j)\big)+\sum_{j<i}\big(|a_{ij}|x_i(y_i+s(a_{ij})y_j)\big)\right)$$

$$=2\sum_{i=1}^{p}\left(a_{ii}-\sum_{j\neq i}|a_{ij}|\right)x_iy_i+2\sum_{i=1}^{p}\left(\sum_{j>i}\big(|a_{ij}|x_i(y_i+s(a_{ij})y_j)\big)+\sum_{j>i}\big(|a_{ji}|x_j(y_j+s(a_{ji})y_i)\big)\right)$$

$$=2\sum_{i=1}^{p}\left(a_{ii}-\sum_{j\neq i}|a_{ij}|\right)x_iy_i+2\sum_{i=1}^{p}\left(\sum_{j>i}\big(|a_{ij}|x_i(y_i+s(a_{ij})y_j)+|a_{ji}|x_j(y_j+s(a_{ji})y_i)\big)\right)$$

$$=2\sum_{i=1}^{p}\left(a_{ii}-\sum_{j\neq i}|a_{ij}|\right)x_iy_i+2\sum_{i=1}^{p}\left(\sum_{j>i}\big(|a_{ij}|(x_i+s(a_{ij})x_j)(y_i+s(a_{ij})y_j)\big)\right).$$

In the following, we will show that for any $i, j \in \{1, \ldots, p\}$, we have the relationship $(x_i+s(a_{ij})x_j)(y_i+s(a_{ij})y_j) \geq 0$. When $x_i+s(a_{ij})x_j = 0$, it is trivial to have $(x_i+s(a_{ij})x_j)(y_i+s(a_{ij})y_j) \geq 0$. Now let us consider the case when $x_i+s(a_{ij})x_j \neq 0$. Since $-1/2+x_i \leq y_i \leq 1/2+x_i$ and $-1/2+x_j \leq y_j \leq 1/2+x_j$, then we have

$$-1+(x_i+s(a_{ij})x_j) \le y_i+s(a_{ij})y_j \le 1+(x_i+s(a_{ij})x_j) \tag{2.10}$$

for all $i, j \in \{1, \ldots, p\}$. We have shown that $x_i \pm x_j$ are integers, which indicates that $x_i+s(a_{ij})x_j$ is integer, and we have discussed the case when $x_i+s(a_{ij})x_j = 0$, so now we consider the case when $x_i+s(a_{ij})x_j$ is non-zero integer. If $x_i+s(a_{ij})x_j < 0$, then $1+(x_i+s(a_{ij})x_j) \le 0$. Therefore, by inequality (2.10), $y_i+s(a_{ij})y_j \le 0$. Similarly, if $x_i+s(a_{ij})x_j > 0$, then $-1+(x_i+s(a_{ij})x_j) \ge 0$. Therefore, by inequality (2.10), we have $y_i + s(a_{ij})y_j \ge 0$. Thus, overall we always have $(x_i+s(a_{ij})x_j)(y_i+s(a_{ij})y_j) \ge 0$.

From the arguments in the previous paragraph, we have the result that $\sum_{i=1}^{p}\sum_{j>i}\left(|a_{ij}|(x_i+s(a_{ij})x_j)(y_i+s(a_{ij})y_j)\right) \ge 0$. In addition, from the result in inequality (2.9), we have $\sum_{i=1}^{p}(a_{ii}-\sum_{j\ne i}|a_{ij}|)x_i y_i > 0$. Thus $2\boldsymbol{x}^T\boldsymbol{A}\boldsymbol{y} > 0$, which indicates that $\bar{\boldsymbol{g}}(\boldsymbol{m}_{\boldsymbol{\theta}})^T(\boldsymbol{\theta}-\boldsymbol{\theta}^*) > 0$ for all $\boldsymbol{m}_{\boldsymbol{\theta}} \in \mathcal{M}_{\boldsymbol{\theta}}$ and all $\boldsymbol{\theta} \in \mathbb{R}^p \setminus \{\boldsymbol{\theta}^*\}$. Q.E.D.

## 2.3 Case of Binary Components in θ

In this section, let us discuss the binary case, where $t_i = 0$ or 1 for all $i$. Based on the algorithm description in the Section 2.1, we know the difference between the algorithm for the unconstrained case and the algorithm for the constrained case (bounded case) is the projection $\boldsymbol{\Psi}(\boldsymbol{\theta})$ (defined in Section 2.1, mapping the value of $\boldsymbol{\theta}$ back to the feasible unit hypercubes). For the bounded case, the algorithm allows $\hat{\boldsymbol{\theta}}_k$ being outside the feasible unit hypercubes, so $\boldsymbol{\pi}(\hat{\boldsymbol{\theta}}_k)$ is defined as $\lfloor \boldsymbol{\Psi}(\hat{\boldsymbol{\theta}}_k) \rfloor + \mathbf{1}_p/2$, which makes $\boldsymbol{\pi}(\hat{\boldsymbol{\theta}}_k)$ to be the middle point of one feasible unit hypercubes. In addition, $\hat{\boldsymbol{\theta}}_k^{\pm}$ are defined as $\hat{\boldsymbol{\theta}}_k^{\pm}$ = $\left[ \boldsymbol{\Psi}(\boldsymbol{\pi}(\hat{\boldsymbol{\theta}}_k) \pm \boldsymbol{\Delta}_k/2) \right]$ ($[\cdot]$ is the round operator), which makes $\hat{\boldsymbol{\theta}}_k^{\pm}$ to be within the feasible domain. For the binary case, we only have one unit hypercube with middle point $\boldsymbol{m}$

$$\boldsymbol{m} = [0.5,...,0.5]^T = \frac{\mathbf{1}_p^T}{2}.$$

Therefore, for all $\boldsymbol{\theta} \in \mathbb{Z}^p$, $\boldsymbol{\pi}(\boldsymbol{\theta})$ = $\lfloor \boldsymbol{\Psi}(\boldsymbol{\theta}) \rfloor + \mathbf{1}_p/2$ = $\boldsymbol{m}$, and $\hat{\boldsymbol{\theta}}_k^{\pm}$ = $\left[ \boldsymbol{\Psi}(\boldsymbol{\pi}(\hat{\boldsymbol{\theta}}_k) \pm \boldsymbol{\Delta}_k/2) \right] \in$ $\{0,1\}^p$. According to the description of the modified DSPSA in Section 2.1 for the constrained case, the solution is $\left[ \boldsymbol{\Psi}(\hat{\boldsymbol{\theta}}_M) \right]$, where $M$ is the maximum number of iterations and $[\cdot]$ is the round operator. In the theorem below, we show the almost sure convergence property for the binary case.

**Theorem 2.2.** Assume $L$ is defined on $\{0,1\}^p$, and $L$ has unique minimal point $\boldsymbol{\theta}^*$. Suppose also (i) $a_k > 0$, $\lim_{k\to\infty} a_k = 0$, $\sum_{k=0}^{\infty} a_k = \infty$ and $\sum_{k=0}^{\infty} a_k^2 < \infty$; (ii) the components of $\boldsymbol{\Delta}_k$ are independently random variables and $\boldsymbol{\Delta}_k^{-T}\boldsymbol{\Delta}_k^{-1}$ is uniformly bounded in $k$; (iii) For all $k$, $E\left((\varepsilon_k^+ - \varepsilon_k^-) \mid \Im_k, \xi_k\right) = 0$ a.s., and $\mathrm{var}(\varepsilon_k^{\pm})$ is uniformly bounded in $k$; (iv) $\bar{\boldsymbol{g}}(\boldsymbol{m})^T(\boldsymbol{\theta} - \boldsymbol{\theta}^*) > 0$ for all $\boldsymbol{\theta} \in [0,1]^p \setminus \{\boldsymbol{\theta}^*\}$. Then $\left[ \boldsymbol{\Psi}\left( \lim_{k\to\infty} \hat{\boldsymbol{\theta}}_k \right) \right]$ = $\boldsymbol{\theta}^*$ a.s., where $[\cdot]$ is the round operator.

*Remarks:*

1. Condition (iv) in the Theorem 2.1 on unconstrained case is not needed here, because for the binary case, we have that the function $L$ is uniformly bounded, which implies that condition (iv) in the Theorem 2.1 is automatically true. Condition (vi) in Theorem 2.1 on the unconstrained case is also not given here, because for binary case we do not need it in the proof.

2. Even though condition (ii) here is the same as the condition (ii) in Theorem 2.1, in practical we only focus on more realistic distribution of $\boldsymbol{\Delta}_k$, where $\boldsymbol{\Delta}_k$ satisfies condition (ii) and $|\Delta_{ki}| \leq 1$ for all $k$ and $i$. For the binary case, there is only one unit hypercube, so we do not need to consider the case of $|\Delta_{ki}| > 1$. Specially, when $|\Delta_{ki}| \leq 1$, the definition of $\hat{\boldsymbol{\theta}}_k^{\pm}$ can be rewritten as $\hat{\boldsymbol{\theta}}_k^{\pm} = \left\lceil \boldsymbol{\Psi}(\boldsymbol{\pi}(\hat{\boldsymbol{\theta}}_k) \pm \boldsymbol{\Delta}_k/2) \right\rfloor = \left\lceil \boldsymbol{\pi}(\hat{\boldsymbol{\theta}}_k) \pm \boldsymbol{\Delta}_k/2 \right\rfloor$, because under the condition of $|\Delta_{ki}| \leq 1$, $\left\lceil \boldsymbol{\pi}(\hat{\boldsymbol{\theta}}_k) \pm \boldsymbol{\Delta}_k/2 \right\rfloor$ must belong to the feasible domain ($\lceil \cdot \rfloor$ is the round operator).

*Proof.* By the same arguments in the proof of Theorem 2.1, we have

$$
\begin{aligned}
\hat{\boldsymbol{\theta}}_k &= \hat{\boldsymbol{\theta}}_0 - \sum_{i=0}^{k} a_i \bar{\boldsymbol{g}}(\boldsymbol{m}) \\
&\quad + \sum_{i=0}^{k} a_i \left( \bar{\boldsymbol{g}}(\boldsymbol{m}) - \left( L(\hat{\boldsymbol{\theta}}_i^+) - L(\hat{\boldsymbol{\theta}}_i^-) \right) \boldsymbol{\Delta}_i^{-1} \right) - \sum_{i=0}^{k} a_i \left( \varepsilon_i^+ - \varepsilon_i^- \right) \boldsymbol{\Delta}_i^{-1}.
\end{aligned} \tag{2.11}
$$

Now let us consider the terms on the right-hand side of eqn. (2.11). Suppose $\boldsymbol{e}_i$ is the $p$-dimensional vector with the $i$th component being 1 and all the other components being

0. We also suppose that $\boldsymbol{\theta}^* = [t_1^*, \ldots, t_p^*]^T$. Since here for all *i*, $t_i^* \in \{0, 1\}$, then the sign of $0.5 - t_i^*$ is $1 - 2t_i^*$, which belongs to the set $\{-1, 1\}$. Then the sign vector of $\boldsymbol{m} - \boldsymbol{\theta}^*$ is $[1 - 2t_1^*, \ldots, 1 - 2t_p^*]^T$, and each component of the sign vector is none zero. Because of condition (iv), we have that

$$\bar{\boldsymbol{g}}(\boldsymbol{m})^T \left( (1 - 2t_i^*) \boldsymbol{e}_i \right) > 0,$$

for all *i*, which indicates that the *i*th component of $\bar{\boldsymbol{g}}(\boldsymbol{m})$ has the same sign as $1 - 2t_i^*$. Therefore, combining the results above and condition (i), we know that the second term on the right-hand side of eqn. (2.11) has the relationship

$$\lim_{k \to \infty} \sum_{i=0}^{k} a_i \bar{\boldsymbol{g}}(\boldsymbol{m}) = \left[ (1 - 2t_1^*) \times (\infty), \ldots, (1 - 2t_p^*) \times (\infty) \right]^T$$

(the entries on the right-hand side are $\pm\infty$).

In addition, since $\left\{ \sum_{i=0}^{k} a_i \left( \bar{\boldsymbol{g}}(\boldsymbol{m}) - \left( L(\hat{\boldsymbol{\theta}}_i^+) - L(\hat{\boldsymbol{\theta}}_i^-) \right) \boldsymbol{\Delta}_i^{-1} \right) \right\}_{k \geq 0}$ is a martingale sequence, by Doob's martingale inequality (Kushner and Clark, 1978, p. 27), we have that for any $\eta > 0$

$$P\left( \sup_{k \geq 0} \left\| \sum_{i=0}^{k} a_i \left( \bar{\boldsymbol{g}}(\boldsymbol{m}) - \left( L(\hat{\boldsymbol{\theta}}_i^+) - L(\hat{\boldsymbol{\theta}}_i^-) \right) \boldsymbol{\Delta}_i^{-1} \right) \right\| \geq \eta \right)$$

$$\leq \eta^{-2} E \left\| \sum_{i=0}^{\infty} a_i \left( \bar{\boldsymbol{g}}(\boldsymbol{m}) - \left( L(\hat{\boldsymbol{\theta}}_i^+) - L(\hat{\boldsymbol{\theta}}_i^-) \right) \boldsymbol{\Delta}_i^{-1} \right) \right\|^2. \qquad (2.12)$$

Now let us consider the right-hand side of inequality (2.12). By the definition of $\bar{\boldsymbol{g}}(\cdot)$

$$\bar{\boldsymbol{g}}(\boldsymbol{m}) = E\left(\left(L(\hat{\boldsymbol{\theta}}_i^+) - L(\hat{\boldsymbol{\theta}}_i^-)\right)\boldsymbol{\Delta}_i^{-1}\right).$$

For all $i < j$, we have

$$E\left(\left(\bar{\boldsymbol{g}}(\boldsymbol{m}) - \left(L(\hat{\boldsymbol{\theta}}_i^+) - L(\hat{\boldsymbol{\theta}}_i^-)\right)\boldsymbol{\Delta}_i^{-1}\right)^T \left(\bar{\boldsymbol{g}}(\boldsymbol{m}) - \left(L(\hat{\boldsymbol{\theta}}_j^+) - L(\hat{\boldsymbol{\theta}}_j^-)\right)\boldsymbol{\Delta}_j^{-1}\right)\right)$$

$$= E\left( E\left(\left(\bar{\boldsymbol{g}}(\boldsymbol{m}) - \left(L(\hat{\boldsymbol{\theta}}_i^+) - L(\hat{\boldsymbol{\theta}}_i^-)\right)\boldsymbol{\Delta}_i^{-1}\right)^T \left(\bar{\boldsymbol{g}}(\boldsymbol{m}) - \left(L(\hat{\boldsymbol{\theta}}_i^+) - L(\hat{\boldsymbol{\theta}}_i^-)\right)\boldsymbol{\Delta}_j^{-1}\right| \xi_{j-1}\right)\right)$$

$$= E\left(\left(\bar{\boldsymbol{g}}(\boldsymbol{m}) - \left(L(\hat{\boldsymbol{\theta}}_i^+) - L(\hat{\boldsymbol{\theta}}_i^-)\right)\boldsymbol{\Delta}_i^{-1}\right)^T E\left(\bar{\boldsymbol{g}}(\boldsymbol{m}) - \left(L(\hat{\boldsymbol{\theta}}_i^+) - L(\hat{\boldsymbol{\theta}}_i^-)\right)\boldsymbol{\Delta}_j^{-1}\right)\right)$$

$$= 0. \tag{2.13}$$

Due to conditions (i), (ii), and eqn. (2.13), we have

$$E\left(\left\|\sum_{i=0}^{\infty} a_i \left\{\bar{\boldsymbol{g}}(\boldsymbol{m}) - \left(L(\hat{\boldsymbol{\theta}}_i^+) - L(\hat{\boldsymbol{\theta}}_i^-)\right)\boldsymbol{\Delta}_i^{-1}\right\}\right\|^2\right)$$

$$= \sum_{i=0}^{\infty} a_i^2 E\left(\left\|\bar{\boldsymbol{g}}(\boldsymbol{m}) - \left(L(\hat{\boldsymbol{\theta}}_i^+) - L(\hat{\boldsymbol{\theta}}_i^-)\right)\boldsymbol{\Delta}_i^{-1}\right\|^2\right)$$

$$\leq \sum_{i=0}^{\infty} a_i^2 E\left(\left(L(\hat{\boldsymbol{\theta}}_i^+) - L(\hat{\boldsymbol{\theta}}_i^-)\right)^2 \boldsymbol{\Delta}_i^{-T}\boldsymbol{\Delta}_i^{-1}\right) < \infty.$$

Since inequality (2.12) is true for all η > 0 and from previous result $E\left(\left\|\sum_{i=0}^{\infty} a_i \left(\bar{\boldsymbol{g}}(\boldsymbol{m}) - \left(L(\hat{\boldsymbol{\theta}}_i^+) - L(\hat{\boldsymbol{\theta}}_i^-)\right)\boldsymbol{\Delta}_i^{-1}\right\|^2\right)$ is finite, we have

$$\lim_{\eta\to\infty} P\left(\sup_{k\geq 0}\left\|\sum_{i=0}^{k} a_i\left(\bar{\boldsymbol{g}}(\boldsymbol{m})-\left(L(\hat{\boldsymbol{\theta}}_i^{+})-L(\hat{\boldsymbol{\theta}}_i^{-})\right)\boldsymbol{\Delta}_i^{-1}\right)\right\|\geq \eta\right)$$

$$\leq \lim_{\eta\to\infty}\frac{E\left\|\sum_{i=0}^{\infty} a_i\left(\bar{\boldsymbol{g}}(\boldsymbol{m})-\left(L(\hat{\boldsymbol{\theta}}_i^{+})-L(\hat{\boldsymbol{\theta}}_i^{-})\right)\boldsymbol{\Delta}_i^{-1}\right)\right\|^2}{\eta^2}$$

$$=0,$$

which indicates that

$$\lim_{k\to\infty}\sum_{i=0}^{k} a_i\left(\bar{\boldsymbol{g}}(\boldsymbol{m})-\left(L\left(\boldsymbol{m}+\frac{1}{2}\boldsymbol{\Delta}_i\right)-L\left(\boldsymbol{m}-\frac{1}{2}\boldsymbol{\Delta}_i\right)\right)\boldsymbol{\Delta}_i^{-1}\right)$$

is finite almost surely. By similar arguments, we know $\left\{\sum_{i=0}^{k} a_i(\varepsilon_i^{+}-\varepsilon_i^{-})\boldsymbol{\Delta}_i^{-1}\right\}_{k\geq 0}$ is a martingale. By Doob's martingale inequality (Kushner and Clark, 1978, p. 27), we have that for any $\eta > 0$

$$P\left(\sup_{k\geq 0}\left\|\sum_{i=0}^{k} a_i(\varepsilon_i^{+}-\varepsilon_i^{-})\boldsymbol{\Delta}_i^{-1}\right\|\geq \eta\right)\leq \eta^{-2}E\left\|\sum_{i=0}^{\infty} a_i(\varepsilon_i^{+}-\varepsilon_i^{-})\boldsymbol{\Delta}_i^{-1}\right\|^2$$

$$=\eta^{-2}\sum_{i=0}^{\infty} a_i^2 E\left\|(\varepsilon_i^{+}-\varepsilon_i^{-})\boldsymbol{\Delta}_i^{-1}\right\|^2. \tag{2.14}$$

Through condition (i), (ii), (iii), we have

$$\sum_{i=0}^{\infty} a_i^2 E\left\|(\varepsilon_i^{+}-\varepsilon_i^{-})\boldsymbol{\Delta}_i^{-1}\right\|^2 < \infty.$$

Since inequality (2.14) is true for all $\eta > 0$, then we have

$$\lim_{\eta\to\infty} P\left(\sup_{k\geq 0}\left\|\sum_{i=0}^{k} a_i(\varepsilon_i^+ - \varepsilon_i^-)\boldsymbol{\Delta}_i^{-1}\right\| \geq \eta\right) \leq \lim_{\eta\to\infty} \frac{\sum_{i=0}^{\infty} a_i^2 E\left\|(\varepsilon_i^+ - \varepsilon_i^-)\boldsymbol{\Delta}_i^{-1}\right\|^2}{\eta^2} = 0,$$

which indicates that

$$\lim_{k\to\infty}\sum_{i=0}^{k} a_i\left(\varepsilon_i^+ - \varepsilon_i^-\right)\boldsymbol{\Delta}_i^{-1}$$

is finite almost surely.

Overall, as $k \to \infty$, we know that the third term and the fourth term on the right-hand side of eqn. (2.11) are finite almost surely, and the second term $-\sum_{i=0}^{k} a_i \bar{\boldsymbol{g}}(\boldsymbol{m})$ goes to $\left[(1-2t_1^*)\times(-\infty),\ldots,(1-2t_p^*)\times(-\infty)\right]^T$ (as above, components are $\pm\infty$). Thus, by eqn. (2.11), we have

$$\lim_{k\to\infty}\hat{\boldsymbol{\theta}}_k = \left[(1-2t_1^*)\times(-\infty),\ldots,(1-2t_p^*)\times(-\infty)\right]^T \text{ a.s.}$$

Because of the definition of the projection, we know that for all $i$

$$\psi_i\left((1-2t_i^*)\times(-\infty)\right) = \begin{cases} 0 & \text{if } t_i^* = 0 \\ 1-\tau & \text{if } t_i^* = 1, \end{cases}$$

which indicates that

$$\psi_i\left((1-2t_i^*)\times(-\infty)\right) = \begin{cases} t_i^* & \text{if } t_i^* = 0 \\ t_i^* - \tau & \text{if } t_i^* = 1, \end{cases}$$

where $\tau$ is a very small positive number. Hence, we have

$$\left[\boldsymbol{\Psi}\left(\left[(1-2t_1^*)\times(-\infty),...,(1-2t_p^*)\times(-\infty)\right]^T\right)\right]=\boldsymbol{\theta}^*,$$

where $[\cdot]$ is the round operator. Therefore, $\left[\boldsymbol{\Psi}\left(\lim_{k\to\infty}\hat{\boldsymbol{\theta}}_k\right)\right]=\boldsymbol{\theta}^*$ a.s. ($[\cdot]$ is the round operator) Q.E.D.

In the Theorem 2.2, except condition (iv), all the other conditions are easy to check. Thus, we only discuss the meaning of inner product condition (condition (iv)) for the binary settings. In the following we assume that the components of $\boldsymbol{\Delta}_k$ are independently Bernoulli $\pm 1$ distributed. The inner product condition is $\bar{\boldsymbol{g}}(\boldsymbol{m})^T(\boldsymbol{\theta}-\boldsymbol{\theta}^*) > 0$ for all $\boldsymbol{\theta}=[t_1,...,t_p]^T \in \{0,1\}^p \setminus \boldsymbol{\theta}^*$, where

$$\begin{aligned}\bar{\boldsymbol{g}}(\boldsymbol{m}) &= E\left(\left(L\left(\boldsymbol{m}+\frac{\boldsymbol{\Delta}}{2}\right)-L\left(\boldsymbol{m}-\frac{\boldsymbol{\Delta}}{2}\right)\right)\boldsymbol{\Delta}^{-1}\right)\\ &=\frac{1}{2^p}\sum\nolimits_{\boldsymbol{\Delta}}\left(L\left(\boldsymbol{m}+\frac{\boldsymbol{\Delta}}{2}\right)-L\left(\boldsymbol{m}-\frac{\boldsymbol{\Delta}}{2}\right)\right)\boldsymbol{\Delta}^{-1},\end{aligned}$$

because the components of $\boldsymbol{\Delta}_k$ are independently Bernoulli $\pm 1$ distributed. Suppose $\boldsymbol{\Delta}=\left[\delta_1,...,\delta_p\right]^T$. Since the optimal solution is $\boldsymbol{\theta}^*=\left[t_1^*,...,t_p^*\right]^T$, then $\boldsymbol{m}-\boldsymbol{\theta}^*$ = $\left[0.5-t_1^*,...,0.5-t_p^*\right]^T$. Hence the sign vector of $\boldsymbol{m}-\boldsymbol{\theta}^*$ is $\left[1-2t_1^*,...,1-2t_p^*\right]^T$, and each component of the sign vector is 1 or −1. Then the inner product condition $\bar{\boldsymbol{g}}(\boldsymbol{m})^T(\boldsymbol{\theta}-\boldsymbol{\theta}^*)$ > 0 is equivalent to that for all $i$

$$\bar{\boldsymbol{g}}(\boldsymbol{m})^T \left((1-2t_i^*)\boldsymbol{e}_i\right) > 0\,, \tag{2.15}$$

where $\boldsymbol{e}_i$ is the vector with the $i$th component being 1 and all other components being 0.

In addition, for each $\boldsymbol{\Delta}$, we have

$$\left(L\left(\boldsymbol{m}+\frac{\boldsymbol{\Delta}}{2}\right)-L\left(\boldsymbol{m}-\frac{\boldsymbol{\Delta}}{2}\right)\right)\boldsymbol{\Delta} = \left(L\left(\boldsymbol{m}+\frac{(-\boldsymbol{\Delta})}{2}\right)-L\left(\boldsymbol{m}-\frac{(-\boldsymbol{\Delta})}{2}\right)\right)(-\boldsymbol{\Delta})\,. \tag{2.16}$$

Suppose that $\Lambda = \left\{\boldsymbol{\Delta}\middle|\, \delta_j = \pm 1,\, j=1,\ldots,p\right\}$, $\Lambda_i = \left\{\boldsymbol{\Delta}\middle|\, \delta_i = 1-2t_i^*,\, \delta_j = \pm 1,\, j \neq i\right\}$ and $\bar{\Lambda}_i = \left\{\boldsymbol{\Delta}\middle|\, \delta_i = -1+2t_i^*,\, \delta_j = \pm 1,\, j \neq i\right\}$. Due to eqn. (2.16), we have

$$\sum\nolimits_{\boldsymbol{\Delta}\in\Lambda_i}\left(L\left(\boldsymbol{m}+\frac{\boldsymbol{\Delta}}{2}\right)-L\left(\boldsymbol{m}-\frac{\boldsymbol{\Delta}}{2}\right)\right)\boldsymbol{\Delta}^{-1} = \sum\nolimits_{\boldsymbol{\Delta}\in\bar{\Lambda}_i}\left(L\left(\boldsymbol{m}+\frac{\boldsymbol{\Delta}}{2}\right)-L\left(\boldsymbol{m}-\frac{\boldsymbol{\Delta}}{2}\right)\right)\boldsymbol{\Delta}^{-1}$$

for each $i \in \{1, \ldots, p\}$. We see that $\Lambda = \Lambda_i \cup \bar{\Lambda}_i$, then

$$\begin{aligned}
\bar{\boldsymbol{g}}(\boldsymbol{m}) &= \frac{1}{2^p}\sum\nolimits_{\boldsymbol{\Delta}\in\Lambda}\left(L\left(\boldsymbol{m}+\frac{\boldsymbol{\Delta}}{2}\right)-L\left(\boldsymbol{m}-\frac{\boldsymbol{\Delta}}{2}\right)\right)\boldsymbol{\Delta}^{-1} \\
&= \frac{1}{2^p}\sum\nolimits_{\boldsymbol{\Delta}\in\Lambda_i}\left(L\left(\boldsymbol{m}+\frac{\boldsymbol{\Delta}}{2}\right)-L\left(\boldsymbol{m}-\frac{\boldsymbol{\Delta}}{2}\right)\right)\boldsymbol{\Delta}^{-1} \\
&\quad + \frac{1}{2^p}\sum\nolimits_{\boldsymbol{\Delta}\in\bar{\Lambda}_i}\left(L\left(\boldsymbol{m}+\frac{\boldsymbol{\Delta}}{2}\right)-L\left(\boldsymbol{m}-\frac{\boldsymbol{\Delta}}{2}\right)\right)\boldsymbol{\Delta}^{-1} \\
&= \frac{2}{2^p}\sum\nolimits_{\boldsymbol{\Delta}\in\Lambda_i}\left(L\left(\boldsymbol{m}+\frac{\boldsymbol{\Delta}}{2}\right)-L\left(\boldsymbol{m}-\frac{\boldsymbol{\Delta}}{2}\right)\right)\boldsymbol{\Delta}^{-1},
\end{aligned}$$

for each $i \in \{1, \ldots, p\}$, which indicates that $\bar{\boldsymbol{g}}(\boldsymbol{m})$ can be written into different expressions. Then for each $i \in \{1, \ldots, p\}$,

$$\bar{\boldsymbol{g}}(\boldsymbol{m})^T\left((1-2t_i^*)\boldsymbol{e}_i\right)=\frac{1}{2^p}\sum\nolimits_{\boldsymbol{\Delta}\in\Lambda}\left(L\left(\boldsymbol{m}+\frac{\boldsymbol{\Delta}}{2}\right)-L\left(\boldsymbol{m}-\frac{\boldsymbol{\Delta}}{2}\right)\right)\boldsymbol{\Delta}^{-1}\left((1-2t_i^*)\boldsymbol{e}_i\right)$$

$$=\frac{2}{2^p}\sum\nolimits_{\boldsymbol{\Delta}\in\Lambda_i}\left(L\left(\boldsymbol{m}+\frac{\boldsymbol{\Delta}}{2}\right)-L\left(\boldsymbol{m}-\frac{\boldsymbol{\Delta}}{2}\right)\right)\boldsymbol{\Delta}^{-T}\left((1-2t_i^*)\boldsymbol{e}_i\right)$$

$$=\frac{2}{2^p}(1-2t_i^*)^2\sum\nolimits_{\boldsymbol{\Delta}\in\Lambda_i}\left(L\left(\boldsymbol{m}+\frac{\boldsymbol{\Delta}}{2}\right)-L\left(\boldsymbol{m}-\frac{\boldsymbol{\Delta}}{2}\right)\right).$$

Therefore, inequality (2.15) can be written as

$$(1-2t_i^*)^2\left(\sum\nolimits_{\boldsymbol{\Delta}\in\Lambda_i}L\left(\boldsymbol{m}+\frac{\boldsymbol{\Delta}}{2}\right)-\sum\nolimits_{\boldsymbol{\Delta}\in\Lambda_i}L\left(\boldsymbol{m}-\frac{\boldsymbol{\Delta}}{2}\right)\right)>0,$$

for all $i$, which is equivalent to

$$\sum\nolimits_{\boldsymbol{\Delta}\in\Lambda_i}L\left(\boldsymbol{m}+\frac{\boldsymbol{\Delta}}{2}\right)>\sum\nolimits_{\boldsymbol{\Delta}\in\Lambda_i}L\left(\boldsymbol{m}-\frac{\boldsymbol{\Delta}}{2}\right), \tag{2.17}$$

for all $i$. Let $\boldsymbol{\Delta}^* = \left[1-2t_1^*,\ldots,1-2t_p^*\right]^T$, then $\boldsymbol{\Delta}^* \in \Lambda_i$ for all $i \in \{1, \ldots, p\}$ and $\boldsymbol{\theta}^* = \boldsymbol{m}-\boldsymbol{\Delta}^*/2$. Then inequality (2.17) indicates that when the sum of loss function values at the points with the $i$th coordinate equal to $t_i^* \in \{0,1\}$ has a smaller value than the sum of the loss function values at the points with the $i$th coordinate equal to $1-t_i^* \in \{0,1\}$ for all $i \in \{1, \ldots, p\}$, the loss function for the binary case satisfies condition (iv) in Theorem 2.2.

# Chapter 3

# Rate of Convergence

In this chapter, we discuss the rate of convergence property of DSPSA. We have shown partial and preliminary results of this chapter in Wang and Spall (2013). We set up an upper bound for the finite sample performance and calculate the asymptotic performance of DSPSA in the big-$O$ sense. We also discuss the properties of the upper bounds. In the last section of this chapter, we discuss the guidelines on the choice of coefficients of the gain sequence.

## 3.1 Theoretical Analysis of Rate of Convergence

Now let us consider the rate of convergence of DSPSA. Generally, for convergent discrete stochastic algorithm, if the optimal solution is unique and all the points in the sequence $\{\hat{\boldsymbol{\theta}}_k\}$ are multivariate integer points, it is natural to use the rate at which $P(\hat{\boldsymbol{\theta}}_k \neq \boldsymbol{\theta}^*)$ going to 0 (or $P(\hat{\boldsymbol{\theta}}_k = \boldsymbol{\theta}^*)$ going to 1) as the measure of rate of convergence. However, for DSPSA, the points in the sequence $\{\hat{\boldsymbol{\theta}}_k\}$ are non-multivariate integer points, so it is natural to consider the mean square error $E\left\|\hat{\boldsymbol{\theta}}_k - \boldsymbol{\theta}^*\right\|^2$ instead. In this section, we provide an upper bound for $E\left\|\hat{\boldsymbol{\theta}}_k - \boldsymbol{\theta}^*\right\|^2$. Furthermore, we have $E\left\|\hat{\boldsymbol{\theta}}_k - \boldsymbol{\theta}^*\right\|^2$ $\geq$ $0^2 P([\hat{\boldsymbol{\theta}}_k] = \boldsymbol{\theta}^*)$ + $0.5^2 P([\hat{\boldsymbol{\theta}}_k] \neq \boldsymbol{\theta}^*)$ = $0.25 P([\hat{\boldsymbol{\theta}}_k] \neq \boldsymbol{\theta}^*)$, where $[\hat{\boldsymbol{\theta}}_k]$ is the nearest multivariate-integer point of $\hat{\boldsymbol{\theta}}_k$. Through the relationship between $E\left\|\hat{\boldsymbol{\theta}}_k - \boldsymbol{\theta}^*\right\|^2$ and $P([\hat{\boldsymbol{\theta}}_k] \neq \boldsymbol{\theta}^*)$, we can get an upper bound of $P([\hat{\boldsymbol{\theta}}_k] \neq \boldsymbol{\theta}^*)$ to compare DSPSA with other algorithms in the big-$O$ sense.

### 3.1.1 Upper Bound for Finite Sample Performance

Before discussing the asymptotic performance of DSPSA, we first consider the finite sample performance by deriving the upper bound for the mean squared error. Theorem 3.1, which provides an upper bound for $E\left\|\hat{\boldsymbol{\theta}}_k - \boldsymbol{\theta}^*\right\|^2$ that applies for all $k \geq 0$, is the main theorem on the rate of convergence analysis for DSPSA. In Corollary 3.2, we show that the upper bound can be written in another form by solving the integral, which leads to a clear asymptotic performance of DSPSA. Following the theorem statement (and preceding the proof), we offer several remarks related to the conditions.

Before proving Theorem 3.1, let us first discuss a useful lemma. In the expression of the upper bound in Theorem 3.1, we will use a function $C(\alpha)$, which is defined as below: for $0.5 < \alpha < 1$

$$C\left(\alpha\right)=\left(\exp\left(2\mu a\frac{(1+A+1)^{1-\alpha}}{1-\alpha}-2\mu a\frac{(1+A)^{1-\alpha}}{1-\alpha}\right)\right)\left(1+\frac{1}{1+A}\right)^{2\alpha},$$

and for $\alpha = 1$,

$$C(1)=\left(1+\frac{1}{1+A}\right)^{2\mu a+2},$$

where $\mu > 0$, $a > 0$ and $A \geq 0$. In Lemma 3.1, we set up a relationship for $C(\alpha)$, and the result will be used in the mathematical induction proof of Theorem 3.1.

**Lemma 3.1.** When $0.5 < \alpha < 1$, the following holds for all $k \geq 0$

$$(1+A+k)^{-2\alpha}\exp\left(\frac{2\mu a(1+A+k+1)^{1-\alpha}}{1-\alpha}\right)\Bigg/\int_k^{k+1}(1+A+k)^{-2\alpha}\exp\left(\frac{2\mu a(1+A+x)^{1-\alpha}}{1-\alpha}\right)dx$$

$$\leq C(\alpha),$$

and

$$\frac{(1+A+k+1)^{2\mu a}}{(1+A+k)^2}\Bigg/\int_k^{k+1}(1+A+x)^{2\mu a-2}dx \leq C(1).$$

*Proof.* For the case of $0.5 < \alpha < 1$, by mean value theorem for integration, there exists $\tilde{x} \in [k, k+1]$ such that

$$(1+A+k)^{-2\alpha} \exp\left(\frac{2\mu a(1+A+k+1)^{1-\alpha}}{1-\alpha}\right) \Bigg/ \int_k^{k+1} (1+A+x)^{-2\alpha} \exp\left(\frac{2\mu a(1+A+x)^{1-\alpha}}{1-\alpha}\right) dx$$

$$= (1+A+k)^{-2\alpha} \exp\left(\frac{2\mu a(1+A+k+1)^{1-\alpha}}{1-\alpha}\right) \Bigg/ (1+A+\tilde{x})^{-2\alpha} \exp\left(\frac{2\mu a(1+A+\tilde{x})^{1-\alpha}}{1-\alpha}\right)$$

$$= \left(\exp\left(\frac{2\mu a(1+A+k+1)^{1-\alpha}}{1-\alpha} - \frac{2\mu a(1+A+\tilde{x})^{1-\alpha}}{1-\alpha}\right)\right)\left(1+\frac{\tilde{x}-k}{1+A+k}\right)^{2\alpha}.$$

Next we will show that

$$(1+A+1)^{1-\alpha} - (1+A)^{1-\alpha} \geq (1+A+k+1)^{1-\alpha} - (1+A+\tilde{x})^{1-\alpha}.$$

When $k = 0$, obviously we have

$$(1+A+1)^{1-\alpha} - (1+A)^{1-\alpha} \geq (1+A+k+1)^{1-\alpha} - (1+A+\tilde{x})^{1-\alpha}.$$

When $k \geq 1$, a Taylor expansion indicates

$$(1+A+k+1)^{1-\alpha} = (1+A+\tilde{x})^{1-\alpha} + \frac{(1-\alpha)(k+1-\tilde{x})}{(1+A+\tilde{\tilde{x}})^{\alpha}},$$

where $\tilde{\tilde{x}} \in [\tilde{x}, k+1]$ , and

$$(1+A+1)^{1-\alpha} = (1+A)^{1-\alpha} + \frac{1-\alpha}{(1+A+\tilde{\tilde{\tilde{x}}})^{\alpha}},$$

where $\tilde{\tilde{\tilde{x}}} \in [0, 1]$. These two Taylor expansions indicate that

$$(1+A+k+1)^{1-\alpha}-(1+A+\tilde{x})^{1-\alpha}=\frac{(1-\alpha)(k+1-\tilde{x})}{(1+A+\tilde{\tilde{x}})^{\alpha}}, \tag{3.1}$$

$$(1+A+1)^{1-\alpha}-(1+A)^{1-\alpha}=\frac{1-\alpha}{(1+A+\tilde{\tilde{\tilde{x}}})^{\alpha}}. \tag{3.2}$$

Thus, when $k \geq 1$, due to the value of $k$, $\tilde{x}$, $\tilde{\tilde{x}}$ and $\tilde{\tilde{\tilde{x}}}$, we know that the right-hand side of eqn. (3.2) is larger than the right-hand side of eqn. (3.1), which means that the left-hand side of eqn. (3.2) is larger than the left-hand side of eqn. (3.1)

$$(1+A+1)^{1-\alpha}-(1+A)^{1-\alpha} \geq (1+A+k+1)^{1-\alpha}-(1+A+\tilde{x})^{1-\alpha}.$$

Therefore, for all $k \geq 0$, we have

$$(1+A+1)^{1-\alpha}-(1+A)^{1-\alpha} \geq (1+A+k+1)^{1-\alpha}-(1+A+\tilde{x})^{1-\alpha},$$

which demonstrates that

$$\begin{aligned}
&\left(\exp\left(\frac{2\mu a(1+A+k+1)^{1-\alpha}}{1-\alpha}-\frac{2\mu a(1+A+\tilde{x})^{1-\alpha}}{1-\alpha}\right)\right)\left(1+\frac{\tilde{x}-k}{1+A+k}\right)^{2\alpha} \\
&\leq \left(\exp\left(\frac{2\mu a(1+A+1)^{1-\alpha}}{1-\alpha}-\frac{2\mu a(1+A)^{1-\alpha}}{1-\alpha}\right)\right)\left(1+\frac{1}{1+A}\right)^{2\alpha} \\
&= C(\alpha).
\end{aligned}$$

Thus,

$$\frac{(1+A+k)^{-2\alpha}\exp\left(\frac{2\mu a(1+A+k+1)^{1-\alpha}}{1-\alpha}\right)}{\int_{k}^{k+1}(1+A+k)^{-2\alpha}\exp\left(\frac{2\mu a(1+A+x)^{1-\alpha}}{1-\alpha}\right)dx} \leq C(\alpha).$$

For the case of $\alpha$ = 1, the mean value theorem for integration indicates

$$\frac{(1+A+k+1)^{2\mu a}}{(1+A+k)^2} \Big/ \int_k^{k+1} (1+A+x)^{2\mu a-2}\,dx = \frac{(1+A+k+1)^{2\mu a}}{(1+A+k)^2} \Big/ (1+A+\tilde{y})^{2\mu a-2}$$

$$= \frac{(1+A+k+1)^{2\mu a}}{(1+A+k)^2} \Big/ \frac{(1+A+\tilde{y})^{2\mu a}}{(1+A+\tilde{y})^2}$$

$$= \left(\frac{1+A+k+1}{1+A+\tilde{y}}\right)^{2\mu a} \left(\frac{1+A+\tilde{y}}{1+A+k}\right)^2$$

$$= \left(1+\frac{k+1-\tilde{y}}{1+A+\tilde{y}}\right)^{2\mu a} \left(1+\frac{\tilde{y}-k}{1+A+k}\right)^2,$$

where $\tilde{y} \in [k, k+1]$, which implies that

$$\frac{(1+A+k+1)^{2\mu a}}{(1+A+k)^2} \Big/ \int_k^{k+1} (1+A+x)^{2\mu a-2}\,dx \le \left(1+\frac{1}{1+A}\right)^{2\mu a+2} = C(1).$$

Q.E.D.

**Theorem 3.1** Assume that $L$ is a function on $\mathbb{Z}^p$ and $L$ has an unique minimal point $\boldsymbol{\theta}^*$. Assume also (i) $a_k = a/(1+A+k)^{\alpha}$, 0.5 < $\alpha$ ≤ 1, $a$ > 0, and $A$ ≥ 0; (ii) the components of $\boldsymbol{\Delta}_k$ are independently distributed random variables and $0 < \boldsymbol{\Delta}_k^{-T}\boldsymbol{\Delta}_k^{-1} \le l < \infty$, where $l$ is independent of the sample point; (iii) for all $k$, $E\left((\varepsilon_k^+ - \varepsilon_k^-) \mid \Im_k, \xi_k\right) = 0$ a.s., and $\mathrm{var}(\varepsilon_k^{\pm})$ is uniformly bounded in $k$; (iv) $E\left(L(\hat{\boldsymbol{\theta}}_k^+) - L(\hat{\boldsymbol{\theta}}_k^-)\right)^2$ is uniformly bounded for all $k$; (v) there exists $\mu > 0$ such that $E\left((\hat{\boldsymbol{\theta}}_k - \boldsymbol{\theta}^*)^T \overline{\boldsymbol{g}}\left(\boldsymbol{\pi}(\hat{\boldsymbol{\theta}}_k)\right) - \mu(\hat{\boldsymbol{\theta}}_k - \boldsymbol{\theta}^*)^T(\hat{\boldsymbol{\theta}}_k - \boldsymbol{\theta}^*)\right) \ge 0$ for all $k$. Then

$$E\left\|\hat{\boldsymbol{\theta}}_k-\boldsymbol{\theta}^*\right\|^2$$

$$\leq\begin{cases}\exp\left(\dfrac{2\mu a(1+A)^{1-\alpha}}{1-\alpha}-\dfrac{2\mu a(1+A+k)^{1-\alpha}}{1-\alpha}\right)E\left\|\hat{\boldsymbol{\theta}}_0-\boldsymbol{\theta}^*\right\|^2 \\ +\exp\left(-\dfrac{2\mu a(1+A+k)^{1-\alpha}}{1-\alpha}\right)lba^2C(\alpha) \\ \times\displaystyle\int_0^k(1+A+x)^{-2\alpha}\exp\left(\dfrac{2\mu a(1+A+x)^{1-\alpha}}{1-\alpha}\right)dx,\quad 0.5<\alpha<1, & (3.3a) \\ \dfrac{(1+A)^{2\mu a}}{(1+A+k)^{2\mu a}}E\left\|\hat{\boldsymbol{\theta}}_0-\boldsymbol{\theta}^*\right\|^2+\dfrac{lba^2C(1)}{(1+A+k)^{2\mu a}}\displaystyle\int_0^k(1+A+x)^{2\mu a-2}dx,\quad \alpha=1, & (3.3b)\end{cases}$$

where for $0.5 < \alpha < 1$

$$C(\alpha)=\left(\exp\left(2\mu a\frac{(1+A+1)^{1-\alpha}}{1-\alpha}-2\mu a\frac{(1+A)^{1-\alpha}}{1-\alpha}\right)\right)\left(1+\frac{1}{1+A}\right)^{2\alpha},$$

and for $\alpha = 1$, we have

$$C(1)=\left(1+\frac{1}{1+A}\right)^{2\mu a+2},$$

and $b$ is a uniform upper bound for $E\left(L(\hat{\boldsymbol{\theta}}_k^+)-L(\hat{\boldsymbol{\theta}}_k^-)\right)^2 + E\left(\varepsilon_k^+-\varepsilon_k^-\right)^2$ .

*Remarks:*

1. A special case that satisfies condition (ii) is when the $\Delta_{ki}$ are independent Bernoulli random variables taking the values $\pm 1$ with probability $1/2$. Under this case, we have $\mathbf{\Delta}_k^{-T}\mathbf{\Delta}_k^{-1} = p$, which means the value of $l$ can be set as $p$ (all $l$ in inequality (3.3a) and (3.3b) can be replaced by $p$ under the special case).

2. From condition (v), we know there exists a $\mu^*$, such that for all $k$ $E\left[(\hat{\boldsymbol{\theta}}_k - \boldsymbol{\theta}^*)^T \overline{\boldsymbol{g}}\left(\boldsymbol{\pi}(\hat{\boldsymbol{\theta}}_k)\right)\right] \geq \mu^* E\left[(\hat{\boldsymbol{\theta}}_k - \boldsymbol{\theta}^*)^T (\hat{\boldsymbol{\theta}}_k - \boldsymbol{\theta}^*)\right]$. Then, for all $\mu$ satisfying $0 < \mu \leq \mu^*$, we have $E\left[(\hat{\boldsymbol{\theta}}_k - \boldsymbol{\theta}^*)^T \overline{\boldsymbol{g}}\left(\boldsymbol{\pi}(\hat{\boldsymbol{\theta}}_k)\right)\right] \geq \mu E\left[(\hat{\boldsymbol{\theta}}_k - \boldsymbol{\theta}^*)^T (\hat{\boldsymbol{\theta}}_k - \boldsymbol{\theta}^*)\right]$ for all $k$. The choice of $\mu$ is also restricted by the relationship $1 - 2\mu a_k > 0$ in the early iterations (for large $k$, the relationship $1 - 2\mu a_k > 0$ is automatically satisfied), which will be further discussed in the proof of Theorem 3.1. Furthermore, condition (v) is an analog of the definition of strongly convexity for continuous case. Thus, $\mu$ is also affected by the structure (curvature) of the loss function.

3. Condition (iv) requires that the sequence of $\{L(\hat{\boldsymbol{\theta}}_k)\}$ is stable, which means that the values of $\{L(\hat{\boldsymbol{\theta}}_k)\}$ are not extremely large in magnitude.

4. The conditions in Theorem 3.1 are not the same as these conditions in Theorem 2.1 (almost sure convergence), because in Theorem 3.1 we consider the mean square error convergence, while in Theorem 2.1 we consider almost sure convergence. The main difference is in the expression of condition (v). Also the condition (vi) in Theorem 2.1 is not needed in Theorem 3.1.

5. When $\alpha \to 1$, the form of the upper bound for 0.5 < $\alpha$ < 1 in inequality (3.3a) converges to the upper bound for the case of $\alpha$ = 1 in inequality (3.3b). We further discuss this result in Corollary 3.1 below.

6. The upper bounds in Theorem 3.1 each include an integral term. We keep the form of integral temporarily, because for the case of 0.5 < $\alpha$ < 1, it is not easy to solve the integral, so we leave it for Corollary 3.2. For the case of $\alpha$ = 1, we want to keep the form of the upper bound consistent with the form for the case of 0.5 < $\alpha$ < 1, and we also want to make the relationship between these two cases clearly. Thus, we also keep the form of integral for the case of $\alpha$ = 1 temporarily, and solve the integral in Corollary 3.2, even though the integral is not hard to solve for the case of $\alpha$ = 1.

*Proof.* The key formula for DSPSA is

$$\hat{\boldsymbol{\theta}}_{k+1} = \hat{\boldsymbol{\theta}}_k - a_k \left( y(\hat{\boldsymbol{\theta}}_k^+) - y(\hat{\boldsymbol{\theta}}_k^-) \right) \boldsymbol{\Delta}_k^{-1}. \tag{3.4}$$

Subtracting $\boldsymbol{\theta}^*$ from both sides of eqn. (3.4) and calculating the norm squared, we get

$$\begin{aligned} \left\| \hat{\boldsymbol{\theta}}_{k+1} - \boldsymbol{\theta}^* \right\|^2 = \left\| \hat{\boldsymbol{\theta}}_{k+1} - \boldsymbol{\theta}^* \right\|^2 - 2a_k (\hat{\boldsymbol{\theta}}_k - \boldsymbol{\theta}^*)^T \left( y(\hat{\boldsymbol{\theta}}_k^+) - y(\hat{\boldsymbol{\theta}}_k^-) \right) \boldsymbol{\Delta}_k^{-1} \\ + \left\| a_k \left( y(\hat{\boldsymbol{\theta}}_k^+) - y(\hat{\boldsymbol{\theta}}_k^-) \right) \boldsymbol{\Delta}_k^{-1} \right\|^2 . \end{aligned} \tag{3.5}$$

Adding and subtracting $2a_k (\hat{\boldsymbol{\theta}}_k - \boldsymbol{\theta}^*)^T \overline{\boldsymbol{g}}(\boldsymbol{\pi}(\hat{\boldsymbol{\theta}}_k))$ to the right-hand side of eqn. (3.5) and taking expectation on both sides, we have

$$E\left\|\hat{\boldsymbol{\theta}}_{k+1}-\boldsymbol{\theta}^*\right\|^2 = E\left\|\hat{\boldsymbol{\theta}}_k-\boldsymbol{\theta}^*\right\|^2 - 2a_k E\left((\hat{\boldsymbol{\theta}}_k-\boldsymbol{\theta}^*)^T\left(\left(y(\hat{\boldsymbol{\theta}}_k^+)-y(\hat{\boldsymbol{\theta}}_k^-)\right)\boldsymbol{\Delta}_k^{-1}-\bar{\boldsymbol{g}}(\boldsymbol{\pi}(\hat{\boldsymbol{\theta}}_k))\right)\right)$$

$$-2a_k E\left((\hat{\boldsymbol{\theta}}_k-\boldsymbol{\theta}^*)^T\bar{\boldsymbol{g}}(\boldsymbol{\pi}(\hat{\boldsymbol{\theta}}_k))\right)+a_k^2 E\left\|\left(y(\hat{\boldsymbol{\theta}}_k^+)-y(\hat{\boldsymbol{\theta}}_k^-)\right)\boldsymbol{\Delta}_k^{-1}\right\|^2. \tag{3.6}$$

Let us now discuss the terms on the right-hand side of eqn. (3.6). After dropping the $-2a_k$ multiplier, the second term on the right-hand side of eqn. (3.6) is:

$$\begin{aligned}
&E\left((\hat{\boldsymbol{\theta}}_k-\boldsymbol{\theta}^*)^T\left(\left(y(\hat{\boldsymbol{\theta}}_k^+)-y(\hat{\boldsymbol{\theta}}_k^-)\right)\boldsymbol{\Delta}_k^{-1}-\bar{\boldsymbol{g}}(\boldsymbol{\pi}(\hat{\boldsymbol{\theta}}_k))\right)\right)\\
&\quad= E\left((\hat{\boldsymbol{\theta}}_k-\boldsymbol{\theta}^*)^T\left(\left(L(\hat{\boldsymbol{\theta}}_k^+)-L(\hat{\boldsymbol{\theta}}_k^-)\right)\boldsymbol{\Delta}_k^{-1}-\bar{\boldsymbol{g}}(\boldsymbol{\pi}(\hat{\boldsymbol{\theta}}_k))\right)\right)+E\left((\hat{\boldsymbol{\theta}}_k-\boldsymbol{\theta}^*)^T(\varepsilon_k^+-\varepsilon_k^-)\boldsymbol{\Delta}_k^{-1}\right)\\
&\quad= E\left((\hat{\boldsymbol{\theta}}_k-\boldsymbol{\theta}^*)^T E\left(\left(L(\hat{\boldsymbol{\theta}}_k^+)-L(\hat{\boldsymbol{\theta}}_k^-)\right)\boldsymbol{\Delta}_k^{-1}-\bar{\boldsymbol{g}}(\boldsymbol{\pi}(\hat{\boldsymbol{\theta}}_k))\middle|\hat{\boldsymbol{\theta}}_k\right)\right)\\
&\qquad+E\left((\hat{\boldsymbol{\theta}}_k-\boldsymbol{\theta}^*)^T\boldsymbol{\Delta}_k^{-1}E\left((\varepsilon_k^+-\varepsilon_k^-)\mid\Im_k,\xi_k\right)\right)\\
&\quad= 0+0=0,
\end{aligned}$$

where $\Im_k$ and $\xi_k$ are defined in Section 2.2 ($\Im_k=\{\hat{\boldsymbol{\theta}}_0,\hat{\boldsymbol{\theta}}_1,...,\hat{\boldsymbol{\theta}}_k\}$, $\xi_k=\{\boldsymbol{\Delta}_0,\boldsymbol{\Delta}_1,...,\boldsymbol{\Delta}_k\}$). After dropping the $a_k^2$ multiplier, the fourth term on the right-hand side of eqn. (3.6) is

$$\begin{aligned}
E\left\|\left(y(\hat{\boldsymbol{\theta}}_k^+)-y(\hat{\boldsymbol{\theta}}_k^-)\right)\boldsymbol{\Delta}_k^{-1}\right\|^2 &= E\left(\left(y(\hat{\boldsymbol{\theta}}_k^+)-y(\hat{\boldsymbol{\theta}}_k^-)\right)^2\boldsymbol{\Delta}_k^{-T}\boldsymbol{\Delta}_k^{-1}\right)\\
&= E\left(\left(L(\hat{\boldsymbol{\theta}}_k^+)-L(\hat{\boldsymbol{\theta}}_k^-)\right)^2\boldsymbol{\Delta}_k^{-T}\boldsymbol{\Delta}_k^{-1}\right)\\
&\quad+2E\left(\left(L(\hat{\boldsymbol{\theta}}_k^+)-L(\hat{\boldsymbol{\theta}}_k^-)\right)(\varepsilon_k^+-\varepsilon_k^-)\boldsymbol{\Delta}_k^{-T}\boldsymbol{\Delta}_k^{-1}\right)\\
&\quad+E\left((\varepsilon_k^+-\varepsilon_k^-)^2\boldsymbol{\Delta}_k^{-T}\boldsymbol{\Delta}_k^{-1}\right).
\end{aligned}$$

In addition, we know that

$$
\begin{aligned}
&E\left(\left(L(\hat{\boldsymbol{\theta}}_k^+) - L(\hat{\boldsymbol{\theta}}_k^-)\right)(\varepsilon_k^+ - \varepsilon_k^-)\boldsymbol{\Delta}_k^{-T}\boldsymbol{\Delta}_k^{-1}\right) \\
&\quad = E\left(E\left(\left(L(\hat{\boldsymbol{\theta}}_k^+) - L(\hat{\boldsymbol{\theta}}_k^-)\right)(\varepsilon_k^+ - \varepsilon_k^-)\boldsymbol{\Delta}_k^{-T}\boldsymbol{\Delta}_k^{-1}\middle| \Im_k, \xi_k\right)\right) \\
&\quad = E\left(\left(L(\hat{\boldsymbol{\theta}}_k^+) - L(\hat{\boldsymbol{\theta}}_k^-)\right)\boldsymbol{\Delta}_k^{-T}\boldsymbol{\Delta}_k^{-1}E\left((\varepsilon_k^+ - \varepsilon_k^-)\middle| \Im_k, \xi_k\right)\right).
\end{aligned}
$$

By condition (iii), we have $E\left[(\varepsilon_k^+ - \varepsilon_k^-)\middle| \Im_k, \xi_k\right] = 0$, which indicates that

$$
E\left(\left(L(\hat{\boldsymbol{\theta}}_k^+) - L(\hat{\boldsymbol{\theta}}_k^-)\right)(\varepsilon_k^+ - \varepsilon_k^-)\boldsymbol{\Delta}_k^{-T}\boldsymbol{\Delta}_k^{-1}\right) = 0.
$$

Then, due to condition (ii), the fourth term on the right-hand side of eqn. (3.6) can be written as

$$
\begin{aligned}
&E\left\|\left(\left(y(\hat{\boldsymbol{\theta}}_k^+) - y(\hat{\boldsymbol{\theta}}_k^-)\right)\boldsymbol{\Delta}_k^{-1}\right)\right\|^2 \\
&\quad = E\left(\left(L(\hat{\boldsymbol{\theta}}_k^+) - L(\hat{\boldsymbol{\theta}}_k^-)\right)^2 \boldsymbol{\Delta}_k^{-T}\boldsymbol{\Delta}_k^{-1}\right) + E\left((\varepsilon_k^+ - \varepsilon_k^-)^2\boldsymbol{\Delta}_k^{-T}\boldsymbol{\Delta}_k^{-1}\right) \\
&\quad \le l\left(E\left(L(\hat{\boldsymbol{\theta}}_k^+) - L(\hat{\boldsymbol{\theta}}_k^-)\right)^2 + E(\varepsilon_k^+ - \varepsilon_k^-)^2\right). \qquad (3.7)
\end{aligned}
$$

Due to condition (iv), we have that $E\left(L(\hat{\boldsymbol{\theta}}_k^+) - L(\hat{\boldsymbol{\theta}}_k^-)\right)^2$ is uniformly bounded. By condition (iii), $E(\varepsilon_k^+ - \varepsilon_k^-)^2$ is also uniformly bounded. Therefore, there exists a positive scalar $b$ such that $E\left(L(\hat{\boldsymbol{\theta}}_k^+) - L(\hat{\boldsymbol{\theta}}_k^-)\right)^2 + E(\varepsilon_k^+ - \varepsilon_k^-)^2 \le b$, which indicates that the right-hand side of eqn. (3.7) is smaller than or equal to $lb$. After substituting the second term and the fourth term on the right-hand side of eqn. (3.6) by 0 and upper bound $lb$, we get

$$
E\left\|\hat{\boldsymbol{\theta}}_{k+1} - \boldsymbol{\theta}^*\right\|^2 \le E\left\|\hat{\boldsymbol{\theta}}_k - \boldsymbol{\theta}^*\right\|^2 - 2a_k E\left((\hat{\boldsymbol{\theta}}_k - \boldsymbol{\theta}^*)^T \overline{\boldsymbol{g}}(\boldsymbol{\pi}(\hat{\boldsymbol{\theta}}_k))\right) + a_k^2 lb.
$$

Furthermore, by condition (v), we have that for all $k$ the inductive step from $k$ to $k+1$ is

$$E\left\|\hat{\boldsymbol{\theta}}_{k+1}-\boldsymbol{\theta}^*\right\|^2 \le (1-2\mu a_k)E\left\|\hat{\boldsymbol{\theta}}_k-\boldsymbol{\theta}^*\right\|^2 + a_k^2 lb, \tag{3.8}$$

and by following all steps before inequality (3.8) in the proof, we know that the right-hand side of inequality (3.8) is positive. In addition, we will show that the term $1-2\mu a_k > 0$ for all $k$. In order to have $1-2\mu a_k > 0$, we must have $\mu < 1/2a_k$ for all $k$. Since $a_k \to 0$ as $k \to \infty$, then we will automatically have $\mu < 1/2a_k$ ($1-2\mu a_k > 0$) when $k$ is large. Thus, the choice of $\mu$ is restricted by $\mu < 1/2a_k$ in the early iterations. As we have discussed in Remark 2, by condition (v), we know that we can always pick $\mu$ small enough to make $1-2\mu a_k > 0$ for all $k$.

We see that the upper bound of the mean square error $E\left\|\hat{\boldsymbol{\theta}}_k-\boldsymbol{\theta}^*\right\|^2$ in the $k$th iteration in inequality (3.8) is composed of two terms. The first term is related to the mean square error from the previous iteration and the second term is the new error introduced in the current iteration.

In all, by recursive inequality (3.8), we have the following solution:

$$\begin{aligned}
&E\left\|\hat{\boldsymbol{\theta}}_{k+1}-\boldsymbol{\theta}^*\right\|^2 \\
&\le \prod_{i=0}^{k}(1-2\mu a_i)E\left\|\hat{\boldsymbol{\theta}}_0-\boldsymbol{\theta}^*\right\|^2 + \prod_{i=1}^{k}(1-2\mu a_i)a_0^2 lb + \ldots + (1-2\mu a_k)a_{k-1}^2 lb + a_k^2 lb \\
&= \prod_{i=0}^{k}(1-2\mu a_i)\left( E\left\|\hat{\boldsymbol{\theta}}_0-\boldsymbol{\theta}^*\right\|^2 + lb\sum_{i=0}^{k}\frac{a_i^2}{\prod_{j=0}^{i}(1-2\mu a_j)}\right).
\end{aligned}$$

From this solution, we find that the upper bound of the mean square error in the *k*th iteration is related to the mean square error of the initial guess and the cumulated errors introduced in each iteration.

Now let us prove Theorem 3.1 by the method of induction. First, let us consider the proof when $0.5 < \alpha < 1$. For the base case, $k = 0$, we have

$$\exp\left(\frac{2\mu a\left(1+A\right)^{1-\alpha}}{1-\alpha}-\frac{2\mu a\left(1+A\right)^{1-\alpha}}{1-\alpha}\right)E\left\|\hat{\boldsymbol{\theta}}_0-\boldsymbol{\theta}^*\right\|^2$$

$$+\exp\left(-\frac{2\mu a\left(1+A\right)^{1-\alpha}}{1-\alpha}\right)lba^2C(\alpha)\int_0^0\left(1+A+x\right)^{-2\alpha}\exp\left(\frac{2\mu a\left(1+A+x\right)^{1-\alpha}}{1-\alpha}\right)dx$$

$$=E\left\|\hat{\boldsymbol{\theta}}_0-\boldsymbol{\theta}^*\right\|^2. \tag{3.9}$$

Eqn. (3.9) indicates the inequality (3.3a) is true for $k = 0$.

Now, suppose the inequality (3.3a) is true for some $k$, which means

$$E\left\|\hat{\boldsymbol{\theta}}_k-\boldsymbol{\theta}^*\right\|^2 \le \exp\left(\frac{2\mu a(1+A)^{1-\alpha}}{1-\alpha}-\frac{2\mu a(1+A+k)^{1-\alpha}}{1-\alpha}\right)E\left\|\hat{\boldsymbol{\theta}}_0-\boldsymbol{\theta}^*\right\|^2$$

$$+\exp\left(-\frac{2\mu a\left(1+A+k\right)^{1-\alpha}}{1-\alpha}\right)lba^2C(\alpha)\int_0^k\left(1+A+x\right)^{-2\alpha}\exp\left(\frac{2\mu a\left(1+A+x\right)^{1-\alpha}}{1-\alpha}\right)dx.$$

Then, for the case of $k + 1$, by recursive relationship (3.8), we have

$$E\left\|\hat{\boldsymbol{\theta}}_{k+1}-\boldsymbol{\theta}^*\right\|^2 \le \left(1-2\mu a_k\right)E\left\|\hat{\boldsymbol{\theta}}_k-\boldsymbol{\theta}^*\right\|^2+a_k^2lb$$

$$=\left(1-\frac{2\mu a}{(1+A+k)^{\alpha}}\right)E\left\|\hat{\boldsymbol{\theta}}_k-\boldsymbol{\theta}^*\right\|^2+\frac{lba^2}{(1+A+k)^{2\alpha}}. \tag{3.10}$$

As we have discussed, there always exists $\mu > 0$ such that condition (v) is satisfied and $1-2\mu a_k > 0$ for all $k$. So, in inequality (3.10), we can replace $E\left\|\hat{\boldsymbol{\theta}}_k - \boldsymbol{\theta}^*\right\|^2$ with its upper bound (from the inductive hypothesis). Then, we have

$$
\begin{aligned}
&E\left\|\hat{\boldsymbol{\theta}}_{k+1} - \boldsymbol{\theta}^*\right\|^2 \\
&\quad \le \left(1-\frac{2\mu a}{(1+A+k)^{\alpha}}\right)\exp\left(\frac{2\mu a(1+A)^{1-\alpha}}{1-\alpha}-\frac{2\mu a(1+A+k)^{1-\alpha}}{1-\alpha}\right)E\left\|\hat{\boldsymbol{\theta}}_0 - \boldsymbol{\theta}^*\right\|^2 \\
&\quad\quad +\left(1-\frac{2\mu a}{(1+A+k)^{\alpha}}\right)\exp\left(-\frac{2\mu a\left(1+A+k\right)^{1-\alpha}}{1-\alpha}\right) \\
&\quad\quad \times lba^2 C(\alpha)\int_0^k \left(1+A+x\right)^{-2\alpha}\exp\left(\frac{2\mu a\left(1+A+x\right)^{1-\alpha}}{1-\alpha}\right)dx+\frac{lba^2}{(1+A+k)^{2\alpha}} \\
&\quad = \left(1-\frac{2\mu a}{(1+A+k)^{\alpha}}\right)\exp\left(-\frac{2\mu a(1+A+k)^{1-\alpha}}{1-\alpha}\right)\exp\left(\frac{2\mu a(1+A)^{1-\alpha}}{1-\alpha}\right)E\left\|\hat{\boldsymbol{\theta}}_0 - \boldsymbol{\theta}^*\right\|^2 \\
&\quad\quad +\left(1-\frac{2\mu a}{(1+A+k)^{\alpha}}\right)\exp\left(-\frac{2\mu a\left(1+A+k\right)^{1-\alpha}}{1-\alpha}\right) \\
&\quad\quad \times lba^2 C(\alpha)\int_0^k \left(1+A+x\right)^{-2\alpha}\exp\left(\frac{2\mu a\left(1+A+x\right)^{1-\alpha}}{1-\alpha}\right)dx+\frac{lba^2}{(1+A+k)^{2\alpha}}. \qquad (3.11)
\end{aligned}
$$

Now let us start to discuss the terms on the right-hand side of inequality (3.11). First we will provide an upper bound to the term

$$
\left(1-\frac{2\mu a}{(1+A+k)^{\alpha}}\right)\exp\left(-\frac{2\mu a(1+A+k)^{1-\alpha}}{1-\alpha}\right)
$$

that appears in inequality (3.11). By a second order Taylor expansion, we have

$$\exp\left(-\frac{2\mu a(1+A+k+1)^{1-\alpha}}{1-\alpha}\right)$$

$$=\exp\left(-\frac{2\mu a(1+A+k)^{1-\alpha}}{1-\alpha}\right)+\exp\left(-\frac{2\mu a(1+A+k)^{1-\alpha}}{1-\alpha}\right)\left(\frac{-2\mu a}{(1+A+k)^{\alpha}}\right)$$

$$+\frac{1}{2}\exp\left(-\frac{2\mu a(1+A+\bar{x})^{1-\alpha}}{1-\alpha}\right)\left(\frac{2\mu a\alpha}{(1+A+\bar{x})^{1+\alpha}}+\frac{4\mu^2 a^2}{(1+A+\bar{x})^{2\alpha}}\right)$$

$$\geq\left(1-\frac{2\mu a}{(1+A+k)^{\alpha}}\right)\exp\left(-\frac{2\mu a(1+A+k)^{1-\alpha}}{1-\alpha}\right),$$

where $\bar{x}\in[k,k+1]$, which means

$$\left(1-\frac{2\mu a}{(1+A+k)^{\alpha}}\right)\exp\left(-\frac{2\mu a(1+A+k)^{1-\alpha}}{1-\alpha}\right)\leq\exp\left(-\frac{2\mu a(1+A+k+1)^{1-\alpha}}{1-\alpha}\right). \quad (3.12)$$

By plugging the results of inequality (3.12) into the first and second term on the right-hand side of inequality (3.11), we have

$$E\left\|\hat{\boldsymbol{\theta}}_{k+1}-\boldsymbol{\theta}^*\right\|^2$$

$$\leq \exp\left(\frac{2\mu a\left(1+A\right)^{1-\alpha}}{1-\alpha}-\frac{2\mu a(1+A+k+1)^{1-\alpha}}{1-\alpha}\right)E\left\|\hat{\boldsymbol{\theta}}_0-\boldsymbol{\theta}^*\right\|^2$$

$$+\exp\left(-\frac{2\mu a(1+A+k+1)^{1-\alpha}}{1-\alpha}\right)lba^2C(\alpha)\int_0^k\left(1+A+x\right)^{-2\alpha}\exp\left(\frac{2\mu a\left(1+A+x\right)^{1-\alpha}}{1-\alpha}\right)dx$$

$$+\frac{lba^2}{(1+A+k)^{2\alpha}}$$

$$=\exp\left(\frac{2\mu a\left(1+A\right)^{1-\alpha}}{1-\alpha}-\frac{2\mu a(1+A+k+1)^{1-\alpha}}{1-\alpha}\right)E\left\|\hat{\boldsymbol{\theta}}_0-\boldsymbol{\theta}^*\right\|^2$$

$$+\exp\left(-\frac{2\mu a(1+A+k+1)^{1-\alpha}}{1-\alpha}\right)lba^2C(\alpha)\left\{\int_0^k(1+A+x)^{-2\alpha}\exp\left(\frac{2\mu a(1+A+x)^{1-\alpha}}{1-\alpha}\right)dx\right.$$

$$\left.+\frac{(1+A+k)^{-2\alpha}}{C(\alpha)}\exp\left(\frac{2\mu a(1+A+k+1)^{1-\alpha}}{1-\alpha}\right)\right\}, \tag{3.13}$$

For the term in the braces $\{\cdot\}$ within the last term on the right-hand side of inequality (3.13), the result in Lemma 3.1 (inequality for $C(\alpha)$) implies

$$\int_0^k(1+A+x)^{-2\alpha}\exp\left(\frac{2\mu a(1+A+x)^{1-\alpha}}{1-\alpha}\right)dx+\frac{(1+A+k)^{-2\alpha}}{C(\alpha)}\exp\left(\frac{2\mu a(1+A+k+1)^{1-\alpha}}{1-\alpha}\right)$$

$$\leq\int_0^k(1+A+x)^{-2\alpha}\exp\left(\frac{2\mu a(1+A+x)^{1-\alpha}}{1-\alpha}\right)dx+\int_k^{k+1}(1+A+x)^{-2\alpha}\exp\left(\frac{2\mu a(1+A+x)^{1-\alpha}}{1-\alpha}\right)dx$$

$$=\int_0^{k+1}(1+A+x)^{-2\alpha}\exp\left(\frac{2\mu a(1+A+x)^{1-\alpha}}{1-\alpha}\right)dx.$$

Therefore, inequality (3.13) can be written as

$$E\left\|\hat{\boldsymbol{\theta}}_{k+1}-\boldsymbol{\theta}^*\right\|^2$$

$$\leq \exp\left(\frac{2\mu a\left(1+A\right)^{1-\alpha}}{1-\alpha}-\frac{2\mu a\left(1+A+k+1\right)^{1-\alpha}}{1-\alpha}\right)E\left\|\hat{\boldsymbol{\theta}}_0-\boldsymbol{\theta}^*\right\|^2$$

$$+\exp\left(-\frac{2\mu a\left(1+A+k+1\right)^{1-\alpha}}{1-\alpha}\right)lba^2C(\alpha)\int_0^{k+1}\left(1+A+x\right)^{-2\alpha}\exp\left(\frac{2\mu a\left(1+A+x\right)^{1-\alpha}}{1-\alpha}\right)dx,$$

which indicates that the inductive $k\to k+1$ step is true for $0.5<\alpha<1$. Therefore, the inequality (3.3a) is true.

Next, let us consider the proof when $\alpha = 1$. For the base case $k = 0$, we have

$$\frac{(1+A)^{2\mu a}}{(1+A)^{2\mu a}}E\left\|\hat{\boldsymbol{\theta}}_0-\boldsymbol{\theta}^*\right\|^2+\frac{lba^2C(1)}{(1+A)^{2\mu a}}\int_0^0(1+A+x)^{2\mu a-2}dx=E\left\|\hat{\boldsymbol{\theta}}_0-\boldsymbol{\theta}^*\right\|^2. \tag{3.14}$$

Equation (3.14) indicates the inequality (3.3b) is true for $k$ = 0. Suppose the inequality (3.3b) is true for some $k$, which means

$$E\left\|\hat{\boldsymbol{\theta}}_k-\boldsymbol{\theta}^*\right\|^2\leq\frac{(1+A)^{2\mu a}}{(1+A+k)^{2\mu a}}E\left\|\hat{\boldsymbol{\theta}}_0-\boldsymbol{\theta}^*\right\|^2+\frac{lba^2C(1)}{(1+A+k)^{2\mu a}}\int_0^k(1+A+x)^{2\mu a-2}dx.$$

For the case of $k$ +1, we have that by recursive inequality (3.8),

$$E\left\|\hat{\boldsymbol{\theta}}_{k+1}-\boldsymbol{\theta}^*\right\|^2\leq\left(1-2\mu a_k\right)E\left\|\hat{\boldsymbol{\theta}}_k-\boldsymbol{\theta}^*\right\|^2+a_k^2lb$$

$$=\left(1-\frac{2\mu a}{1+A+k}\right)E\left\|\hat{\boldsymbol{\theta}}_k-\boldsymbol{\theta}^*\right\|^2+\frac{lba^2}{\left(1+A+k\right)^2}. \tag{3.15}$$

As we have discussed, there always exists $\mu > 0$ such that condition (v) is satisfied and $1-2\mu a_k > 0$ for all $k \geq 0$. Hence, we can replace $E\left\|\hat{\boldsymbol{\theta}}_k - \boldsymbol{\theta}^*\right\|^2$ with its upper bound (from the inductive hypothesis) into inequality (3.15). Then, we have

$$E\left\|\hat{\boldsymbol{\theta}}_{k+1} - \boldsymbol{\theta}^*\right\|^2 \leq \left(1 - \frac{2\mu a}{1+A+k}\right)\frac{(1+A)^{2\mu a}}{(1+A+k)^{2\mu a}} E\left\|\hat{\boldsymbol{\theta}}_0 - \boldsymbol{\theta}^*\right\|^2$$

$$+\left(1 - \frac{2\mu a}{1+A+k}\right)\frac{lba^2 C(1)}{(1+A+k)^{2\mu a}}\int_0^k (1+A+x)^{2\mu a-2}dx + \frac{lba^2}{(1+A+k)^2}. \quad (3.16)$$

By the second order Taylor expansion, we have

$$\left(\frac{1}{1+A+k+1}\right)^{2\mu a} = \left(\frac{1}{1+A+k}\right)^{2\mu a} - \frac{2\mu a}{(1+A+k)^{2\mu a+1}} + \frac{\mu a(2\mu a+1)}{(1+A+\overline{y})^{2\mu a+2}}$$

$$\geq \left(1 - \frac{2\mu a}{(1+A+k)}\right)\left(\frac{1}{1+A+k}\right)^{2\mu a}, \quad (3.17)$$

where $\overline{y} \in [k, k+1]$. Then plugging the upper bound of

$$\left(1 - \frac{2a\mu}{(1+A+k)}\right)\left(\frac{1}{1+A+k}\right)^{2\mu a}$$

in inequality (3.17) into the first and second term on the right-hand side of eqn. (3.16), we have

$$E\left\|\hat{\boldsymbol{\theta}}_{k+1}-\boldsymbol{\theta}^*\right\|^2 \le \frac{\left(1+A\right)^{2\mu a}}{\left(1+A+k+1\right)^{2\mu a}} E\left\|\hat{\boldsymbol{\theta}}_0-\boldsymbol{\theta}^*\right\|^2$$

$$+\frac{lba^2C(1)}{\left(1+A+k+1\right)^{2\mu a}}\int_0^k (1+A+x)^{2\mu a-2}dx+\frac{lba^2}{(1+A+k)^2}$$

$$=\frac{\left(1+A\right)^{2\mu a}}{\left(1+A+k+1\right)^{2\mu a}} E\left\|\hat{\boldsymbol{\theta}}_0-\boldsymbol{\theta}^*\right\|^2$$

$$+\frac{lba^2C(1)}{\left(1+A+k+1\right)^{2\mu a}}\left(\int_0^k \left(1+A+x\right)^{2\mu a-2}dx+\frac{\left(1+A+k+1\right)^{2\mu a}}{C(1)(1+A+k)^2}\right).$$

Due to Lemma 3.1 (inequality for $C(1)$), we know

$$\int_0^k \left(1+A+x\right)^{2\mu a-2}dx+\frac{\left(1+A+k+1\right)^{2\mu a}}{C(1)(1+A+k)^2} \le \int_0^k (1+A+x)^{2\mu a-2}dx+\int_k^{k+1}(1+A+x)^{2\mu a-2}dx$$

$$=\int_0^{k+1}(1+A+x)^{2\mu a-2}dx.$$

Then, we get

$$E\left\|\hat{\boldsymbol{\theta}}_{k+1}-\boldsymbol{\theta}^*\right\|^2 \le \frac{\left(1+A\right)^{2\mu a}}{\left(1+A+k+1\right)^{2\mu a}} E\left\|\hat{\boldsymbol{\theta}}_0-\boldsymbol{\theta}^*\right\|^2+\frac{lba^2C(1)}{\left(1+A+k+1\right)^{2\mu a}}\int_0^{k+1}(1+A+x)^{2\mu a-2}dx,$$

which indicates that the inductive $k \to k+1$ step is true for $\alpha = 1$. Therefore, the inequality (3.3b) is true.

In summary, we have

$$
E\left\|\hat{\boldsymbol{\theta}}_k-\boldsymbol{\theta}^*\right\|^2
$$

$$
\le\begin{cases}\exp\left(\dfrac{2\mu a\left(1+A\right)^{1-\alpha}}{1-\alpha}-\dfrac{2\mu a\left(1+A+k\right)^{1-\alpha}}{1-\alpha}\right)E\left\|\hat{\boldsymbol{\theta}}_0-\boldsymbol{\theta}^*\right\|^2 \\ +\exp\left(-\dfrac{2\mu a\left(1+A+k\right)^{1-\alpha}}{1-\alpha}\right)lba^2C(\alpha)\displaystyle\int_0^k\frac{\exp\left(\dfrac{2\mu a\left(1+A+x\right)^{1-\alpha}}{1-\alpha}\right)}{\left(1+A+x\right)^{2\alpha}}dx, \; 0.5<\alpha<1, \\ \dfrac{(1+A)^{2\mu a}}{(1+A+k)^{2\mu a}}E\left\|\hat{\boldsymbol{\theta}}_0-\boldsymbol{\theta}^*\right\|^2+\dfrac{lba^2C(1)}{(1+A+k)^{2\mu a}}\displaystyle\int_0^k(1+A+x)^{2\mu a-2}dx, \; \alpha=1,\end{cases}
$$

which completes the proof. Q.E.D.

In the following, we present Corollary 3.1, where we show the relationship between the upper bound in inequality (3.3a) and the upper bound in inequality (3.3b).

**Corollary 3.1.** Assume that the conditions of Theorem 3.1 are true. In inequality (3.3a), the upper bound for the case of $0.5 < \alpha < 1$ converges to the upper bound for the case of $\alpha = 1$ in inequality (3.3b), when $\alpha \to 1$.

*Proof.* Under the conditions of Theorem 3.1, we know that inequalities (3.3a) and (3.3b) are true. First we show that the term

$$
\exp\left(\frac{2\mu a(1+A)^{1-\alpha}}{1-\alpha}-\frac{2\mu a(1+A+x)^{1-\alpha}}{1-\alpha}\right)
$$

in inequality (3.3a) converges to $(1+A)^{2\mu a} \big/ (1+A+x)^{2\mu a}$ as $\alpha \to 1$ for all $x \geq 0$. Suppose the base of the logarithm used below is $e$ (natural logarithm). By L'Hôpital's rule, we have

$$\begin{aligned}
&\lim_{\alpha\to 1} \frac{2\mu a\left(1+A\right)^{1-\alpha} - 2\mu a\left(1+A+x\right)^{1-\alpha}}{1-\alpha} \\
&= \lim_{\alpha\to 1} \frac{-2\mu a\left(1+A\right)^{1-\alpha}\log(1+A) + 2\mu a\left(1+A+x\right)^{1-\alpha}\log(1+A+x)}{-1} \\
&= 2\mu a\log(1+A) - 2\mu a\log(1+A+x) \\
&= \log(1+A)^{2\mu a} - \log(1+A+x)^{2\mu a},
\end{aligned}$$

which indicates that

$$\lim_{\alpha\to 1} \exp\left(\frac{2\mu a(1+A)^{1-\alpha}}{1-\alpha} - \frac{2\mu a(1+A+x)^{1-\alpha}}{1-\alpha}\right) = \frac{(1+A)^{2\mu a}}{(1+A+x)^{2\mu a}}.$$

Similarly

$$\lim_{\alpha\to 1} \exp\left(\frac{2\mu a(1+A)^{1-\alpha}}{1-\alpha} - \frac{2\mu a(1+A+k)^{1-\alpha}}{1-\alpha}\right) = \frac{(1+A)^{2\mu a}}{(1+A+k)^{2\mu a}}, \tag{3.18}$$

$$\lim_{\alpha\to 1} \exp\left(\frac{2\mu a(1+A+1)^{1-\alpha}}{1-\alpha} - \frac{2\mu a(1+A)^{1-\alpha}}{1-\alpha}\right) = \frac{(1+A+1)^{2\mu a}}{(1+A)^{2\mu a}}, \tag{3.19}$$

and

$$\lim_{\alpha\to 1} \exp\left(\frac{2\mu a(1+A+x)^{1-\alpha}}{1-\alpha} - \frac{2\mu a(1+A+k)^{1-\alpha}}{1-\alpha}\right) = \frac{(1+A+x)^{2\mu a}}{(1+A+k)^{2\mu a}}. \tag{3.20}$$

By eqn. (3.19) we have

$$\lim_{\alpha\to 1} C(\alpha) = \lim_{\alpha\to 1} \exp\left( \frac{2\mu a(1+A+1)^{1-\alpha}}{1-\alpha} - \frac{2\mu a(1+A)^{1-\alpha}}{1-\alpha} \right)\left(1+\frac{1}{1+A}\right)^{2\alpha}$$

$$= \left(1+\frac{1}{1+A}\right)^{2\mu a+2}$$

$$= C(1). \tag{3.21}$$

In addition, by eqn. (3.20) we have

$$\lim_{\alpha\to 1} \exp\left( -\frac{2\mu a\left(1+A+k\right)^{1-\alpha}}{1-\alpha} \right) \int_0^k \left(1+A+x\right)^{-2\alpha} \exp\left( \frac{2\mu a\left(1+A+x\right)^{1-\alpha}}{1-\alpha} \right) dx$$

$$= \lim_{\alpha\to 1} \int_0^k \left(1+A+x\right)^{-2\alpha} \exp\left( \frac{2\mu a\left(1+A+x\right)^{1-\alpha}}{1-\alpha} - \frac{2\mu a\left(1+A+k\right)^{1-\alpha}}{1-\alpha} \right) dx$$

$$= \frac{1}{(1+A+k)^{2\mu a}} \int_0^k (1+A+x)^{2\mu a-2}\, dx. \tag{3.22}$$

Overall, by eqns. (3.18), (3.21), and (3.22), we have

$$\lim_{\alpha\to 1} \left\{ \exp\left( \frac{2\mu a\left(1+A\right)^{1-\alpha}}{1-\alpha} - \frac{2\mu a\left(1+A+k\right)^{1-\alpha}}{1-\alpha} \right) E\left\|\hat{\boldsymbol{\theta}}_0 - \boldsymbol{\theta}^*\right\|^2 \right.$$

$$\left. + \exp\left( -\frac{2\mu a\left(1+A+k\right)^{1-\alpha}}{1-\alpha} \right) lba^2 C(\alpha) \int_0^k \left(1+A+x\right)^{-2\alpha} \exp\left( \frac{2\mu a\left(1+A+x\right)^{1-\alpha}}{1-\alpha} \right) dx \right\}$$

$$= \frac{(1+A)^{2\mu a}}{(1+A+k)^{2\mu a}} E\left\|\hat{\boldsymbol{\theta}}_0 - \boldsymbol{\theta}^*\right\|^2 + \frac{lba^2 C(1)}{(1+A+k)^{2\mu a}} \int_0^k (1+A+x)^{2\mu a-2}\, dx,$$

which means that, when $\alpha \to 1$, the upper bound in inequality (3.3a) for the case of $0.5 < \alpha < 1$ converges to the upper bound in inequality (3.3b) for the case of $\alpha = 1$. Q.E.D.

## 3.1.2 Asymptotic Performance

Following the discussion of the upper bound for finite sample performance in last section, we consider the asymptotic performance of DSPSA in this section. In Corollary 3.2, we show that inequality (3.3a) and (3.3b) in Theorem 3.1 can be written in a new form by solving the integral, and from the new form we can get the rate of convergence of DSPSA in the big-$O$ sense.

**Corollary 3.2.** Assume that the conditions of Theorem 3.1 are true. Inequality (3.3a) and (3.3b) can be written as

$$E\left\|\hat{\boldsymbol{\theta}}_k - \boldsymbol{\theta}^*\right\|^2$$

$$\leq \begin{cases} \exp\left(\dfrac{2\mu a(1+A)^{1-\alpha}}{1-\alpha} - \dfrac{2\mu a(1+A+k)^{1-\alpha}}{1-\alpha}\right)\left(E\left\|\hat{\boldsymbol{\theta}}_0 - \boldsymbol{\theta}^*\right\|^2 - \dfrac{T(k,\alpha)}{(1+A)^{\alpha}}\right) \\ \quad + \dfrac{T(k,\alpha)}{(1+A+k)^{\alpha}},\ 0.5<\alpha<1, & (3.23a) \\ \dfrac{(1+A)^{2\mu a}}{(1+A+k)^{2\mu a}}\left(E\left\|\hat{\boldsymbol{\theta}}_0 - \boldsymbol{\theta}^*\right\|^2 - \dfrac{T(k,1)}{1+A}\right) + \dfrac{T(k,1)}{1+A+k},\ \alpha = 1, & (3.23b) \end{cases}$$

where we have that when $0.5 < \alpha < 1$,

$$T(k,\alpha) = \frac{lba^2 C(\alpha)}{2\mu a - \alpha\big/\left(1+A+f(k)\right)^{1-\alpha}},$$

$$f(k)=\left(\frac{\int_0^k (1+A+x)^{-2\alpha}\exp\left(\frac{2\mu a(1+A+x)^{1-\alpha}}{1-\alpha}\right)dx}{\int_0^k (1+A+x)^{-(1+\alpha)}\exp\left(\frac{2\mu a(1+A+x)^{1-\alpha}}{1-\alpha}\right)dx}\right)^{\frac{1}{1-\alpha}}-1-A,$$

$f(k)\to\infty$ as $k\to\infty$, and when $\alpha=1$

$$T(k,1)=\frac{lba^2C(1)}{2\mu a-1}.$$

*Remarks:*

1. In Corollary 3.2, we show that inequality (3.3a) and (3.3b) in Theorem 3.1 can be written into a new form by solving the integral, and we can derive the rate of convergence in the big-*O* sense through this new form easily. Thus, the values of $T(k,\alpha)$ in finite iterations are not important for us in discussing the asymptotic performance.

2. For $T(k,\alpha)$, we will show in the proof that $f(k)\to\infty$ as $k\to\infty$, which indicates that when 0.5 < α < 1, we have $T(k,\alpha)\to lba^2C(\alpha)/2\mu a$ as $k\to\infty$. In addition, when $k$ is large, the effect of the coefficient $A$ disappears in $T(k,\alpha)$, while the effect of $a$ is always present. The overall reasons that we introduce $f(k)$ here are that we can rewrite the integral into two clear big-*O* forms and when $k\to\infty$ the effect of $f(k)$ disappears in the coefficient $T(k,\alpha)$.

3. When α = 1, we can always adjust the value of $a$ a small amount to make $2\mu a\neq 1$.

*Proof.* Under the conditions of Theorem 3.1, we know that inequalities (3.3a) and (3.3b)

are true. Now we solve the integral in inequality (3.3a) and (3.3b). For the case of $0.5 < \alpha < 1$, let us consider the integral

$$2\mu a\int_0^k (1+A+x)^{-2\alpha}\exp\left(\frac{2\mu a(1+A+x)^{1-\alpha}}{1-\alpha}\right)dx .$$

Using integration by parts, we have

$$2\mu a\int_0^k (1+A+x)^{-2\alpha}\exp\left(\frac{2\mu a(1+A+x)^{1-\alpha}}{1-\alpha}\right)dx$$

$$=\int_0^k \frac{1}{(1+A+x)^{\alpha}}\,d\exp\left(\frac{2\mu a(1+A+x)^{1-\alpha}}{1-\alpha}\right)$$

$$=\frac{1}{(1+A+x)^{\alpha}}\exp\left(\frac{2\mu a(1+A+x)^{1-\alpha}}{1-\alpha}\right)\Bigg|_0^k-\int_0^k \exp\left(\frac{2\mu a(1+A+x)^{1-\alpha}}{1-\alpha}\right)d\frac{1}{(1+A+x)^{\alpha}}$$

$$=\frac{1}{(1+A+x)^{\alpha}}\exp\left(\frac{2\mu a(1+A+x)^{1-\alpha}}{1-\alpha}\right)\Bigg|_0^k+\int_0^k \alpha(1+A+x)^{-\alpha-1}\exp\left(\frac{2\mu a(1+A+x)^{1-\alpha}}{1-\alpha}\right)dx$$

$$=\frac{1}{(1+A+k)^{\alpha}}\exp\left(\frac{2\mu a(1+A+k)^{1-\alpha}}{1-\alpha}\right)-\frac{1}{(1+A)^{\alpha}}\exp\left(\frac{2\mu a(1+A)^{1-\alpha}}{1-\alpha}\right)$$

$$+\alpha\int_0^k (1+A+x)^{-2\alpha}\exp\left(\frac{2\mu a(1+A+x)^{1-\alpha}}{1-\alpha}\right)\frac{1}{(1+A+x)^{1-\alpha}}dx,$$

which implies that

$$\int_0^k (1+A+x)^{-2\alpha}\exp\left(\frac{2\mu a(1+A+x)^{1-\alpha}}{1-\alpha}\right)\left(2\mu a-\frac{\alpha}{(1+A+x)^{1-\alpha}}\right)dx$$

$$=\frac{1}{(1+A+k)^{\alpha}}\exp\left(\frac{2\mu a(1+A+k)^{1-\alpha}}{1-\alpha}\right)-\frac{1}{(1+A)^{\alpha}}\exp\left(\frac{2\mu a(1+A)^{1-\alpha}}{1-\alpha}\right). \quad (3.24)$$

By the mean value theorem for integration, for each *k*, there exists $f(k) \in [0, k]$ such that

$$\int_0^k (1+A+x)^{-2\alpha} \exp\left(\frac{2\mu a(1+A+x)^{1-\alpha}}{1-\alpha}\right)\frac{1}{(1+A+x)^{1-\alpha}}\,dx$$

$$= \frac{1}{(1+A+f(k))^{1-\alpha}}\int_0^k (1+A+x)^{-2\alpha} \exp\left(\frac{2\mu a(1+A+x)^{1-\alpha}}{1-\alpha}\right)dx, \tag{3.25}$$

which implies that

$$f(k) = \left(\frac{\displaystyle\int_0^k (1+A+x)^{-2\alpha} \exp\left(\frac{2\mu a(1+A+x)^{1-\alpha}}{1-\alpha}\right)dx}{\displaystyle\int_0^k (1+A+x)^{-(1+\alpha)} \exp\left(\frac{2\mu a(1+A+x)^{1-\alpha}}{1-\alpha}\right)dx}\right)^{\frac{1}{1-\alpha}} - 1 - A\,.$$

Combining the results of eqn. (3.24) and eqn. (3.25), we have

$$\left(2\mu a - \frac{\alpha}{(1+A+f(k))^{1-\alpha}}\right)\int_0^k (1+A+x)^{-2\alpha} \exp\left(\frac{2\mu a(1+A+x)^{1-\alpha}}{1-\alpha}\right)dx$$

$$= \frac{1}{(1+A+k)^{\alpha}}\exp\left(\frac{2\mu a(1+A+k)^{1-\alpha}}{1-\alpha}\right) - \frac{1}{(1+A)^{\alpha}}\exp\left(\frac{2\mu a(1+A)^{1-\alpha}}{1-\alpha}\right),$$

which implies that

$$\int_0^k (1+A+x)^{-2\alpha} \exp\left(\frac{2\mu a(1+A+x)^{1-\alpha}}{1-\alpha}\right)dx$$

$$= \frac{\exp\left(\frac{2\mu a(1+A+k)^{1-\alpha}}{1-\alpha}\right)\Big/(1+A+k)^{\alpha} - \exp\left(\frac{2\mu a(1+A)^{1-\alpha}}{1-\alpha}\right)\Big/(1+A)^{\alpha}}{2\mu a - \alpha\big/(1+A+f(k))^{1-\alpha}},$$

which means that we have solved the integral in inequality (3.3a). Thus, the upper bound in the inequality (3.3a) can be written as

$$\exp\left(\frac{2\mu a(1+A)^{1-\alpha}}{1-\alpha}-\frac{2\mu a(1+A+k)^{1-\alpha}}{1-\alpha}\right)E\left\|\hat{\boldsymbol{\theta}}_0-\boldsymbol{\theta}^*\right\|^2$$

$$+\exp\left(-\frac{2\mu a(1+A+k)^{1-\alpha}}{1-\alpha}\right)lba^2C(\alpha)\int_0^k(1+A+x)^{-2\alpha}\exp\left(\frac{2\mu a(1+A+x)^{1-\alpha}}{1-\alpha}\right)dx$$

$$=\exp\left(\frac{2\mu a(1+A)^{1-\alpha}}{1-\alpha}-\frac{2\mu a(1+A+k)^{1-\alpha}}{1-\alpha}\right)E\left\|\hat{\boldsymbol{\theta}}_0-\boldsymbol{\theta}^*\right\|^2+\exp\left(-\frac{2\mu a(1+A+k)^{1-\alpha}}{1-\alpha}\right)$$

$$\times lba^2C(\alpha)\left(\frac{\exp\left(\frac{2\mu a(1+A+k)^{1-\alpha}}{1-\alpha}\right)\Big/\left(1+A+k\right)^{\alpha}-\exp\left(\frac{2\mu a(1+A)^{1-\alpha}}{1-\alpha}\right)\Big/(1+A)^{\alpha}}{2\mu a-\alpha/(1+A+f(k))^{1-\alpha}}\right)$$

$$=\exp\left(\frac{2\mu a(1+A)^{1-\alpha}}{1-\alpha}-\frac{2\mu a(1+A+k)^{1-\alpha}}{1-\alpha}\right)\left(E\left\|\hat{\boldsymbol{\theta}}_0-\boldsymbol{\theta}^*\right\|^2-\frac{T(k,\alpha)}{(1+A)^{\alpha}}\right)+\frac{T(k,\alpha)}{(1+A+k)^{\alpha}},$$

which shows that the result in inequality (3.23a) is true.

In the following, we show that for the $f(k)$ in inequality (3.23a), we have $f(k)\to\infty$ as $k\to\infty$. After that, we show that inequality (3.23b) is true. By rewriting eqn. (3.25), we have

$$\int_0^k(1+A+x)^{-(1+\alpha)}\exp\left(\frac{2\mu a(1+A+x)^{1-\alpha}}{1-\alpha}\right)dx$$

$$=\left(\frac{1}{(1+A+f(k))^{1-\alpha}}\right)\int_0^k(1+A+x)^{-2\alpha}\exp\left(\frac{2\mu a(1+A+x)^{1-\alpha}}{1-\alpha}\right)dx,$$

which follows that

$$\frac{1}{(1+A+f(k))^{1-\alpha}} = \frac{\int_0^k (1+A+x)^{-(1+\alpha)} \exp\left(\frac{2\mu a(1+A+x)^{1-\alpha}}{1-\alpha}\right) dx}{\int_0^k (1+A+x)^{-2\alpha} \exp\left(\frac{2\mu a(1+A+x)^{1-\alpha}}{1-\alpha}\right) dx}.$$

By L'Hôpital's rule and the fact that $\alpha < 1$, we have

$$\lim_{k\to\infty} \frac{\int_0^k (1+A+x)^{-(1+\alpha)} \exp\left(\frac{2\mu a(1+A+x)^{1-\alpha}}{1-\alpha}\right) dx}{\int_0^k (1+A+x)^{-2\alpha} \exp\left(\frac{2\mu a(1+A+x)^{1-\alpha}}{1-\alpha}\right) dx}$$

$$= \lim_{k\to\infty} \frac{(1+A+k)^{-(1+\alpha)} \exp\left(\frac{2\mu a(1+A+k)^{1-\alpha}}{1-\alpha}\right)}{(1+A+k)^{-2\alpha} \exp\left(\frac{2\mu a(1+A+k)^{1-\alpha}}{1-\alpha}\right)}$$

$$= \lim_{k\to\infty} \frac{1}{(1+A+k)^{1-\alpha}} = 0.$$

Thus,

$$\lim_{k\to\infty} \frac{1}{(1+A+f(k))^{1-\alpha}} = 0,$$

which indicates that $f(k) \to \infty$ as $k \to \infty$. Then we have $T(k,\alpha) \to lba^2 C(\alpha)/2\mu a > 0$ as $k \to \infty$. Therefore, we have $E\left\|\hat{\boldsymbol{\theta}}_k - \boldsymbol{\theta}^*\right\|^2 = O(1/k^{\alpha})$ for the case of $0.5 < \alpha < 1$, and specifically $E\left\|\hat{\boldsymbol{\theta}}_k - \boldsymbol{\theta}^*\right\|^2 = \left(lba^2 C(\alpha)/2\mu a\right) k^{-\alpha} + o(k^{-\alpha})$.

Now let us show that when $\alpha = 1$ the inequality (3.23b) is true. For the case of $\alpha = 1$, as we have discussed in Remark 3 ($2\mu a \neq 1$) of this corollary, we have

$$\frac{(1+A)^{2\mu a}}{(1+A+k)^{2\mu a}} E\left\|\hat{\boldsymbol{\theta}}_0 - \boldsymbol{\theta}^*\right\|^2 + \frac{lba^2 C(1)}{(1+A+k)^{2\mu a}} \int_0^k (1+A+x)^{2\mu a - 2} dx$$

$$= \frac{(1+A)^{2\mu a}}{(1+A+k)^{2\mu a}} E\left\|\hat{\boldsymbol{\theta}}_0 - \boldsymbol{\theta}^*\right\|^2 + \frac{lba^2 C(1)}{(2\mu a - 1)(1+A+k)^{2\mu a}} \left((1+A+k)^{2\mu a - 1} - (1+A)^{2\mu a - 1}\right)$$

$$= \frac{(1+A)^{2\mu a}}{(1+A+k)^{2\mu a}} \left( E\left\|\hat{\boldsymbol{\theta}}_0 - \boldsymbol{\theta}^*\right\|^2 - \frac{T(k,1)}{1+A} \right) + \frac{T(k,1)}{1+A+k},$$

which shows that the result in inequality (3.23b) is true. Overall, we have shown that the inequality (3.23a) and (3.23b) are true. Q.E.D.

The second term on the right-hand side of inequality (3.23a) goes to 0 at a lower rate than the first term, so we have $E\left\|\hat{\boldsymbol{\theta}}_k - \boldsymbol{\theta}^*\right\|^2 = O(1/k^{\alpha})$ for the case of 0.5 < $\alpha$ < 1. Moreover, on the right-hand side of inequality (3.23b), the first term = $O(1/k^{2\mu a})$ and the second term = $O(1/k)$. Therefore, we need to discuss the relationship between $2\mu a$ and 1 before figuring out the asymptotic performance for $\alpha$ = 1 in the big-*O* sense. If we only consider the upper bound for the asymptotic performance, we only need to focus on the sequence for very large *k*. That is, for the asymptotic performance of $\{\hat{\boldsymbol{\theta}}_k\}$, we only need to consider the sequence $\{\hat{\boldsymbol{\theta}}_k\}_{k \geq N}$ when *N* is a very large integer value. Suppose we only focus on the exponent of the big-*O* function instead of the constant multiplier of the big-*O* function. For very large *N*, and $k \geq N$, it is known that $1/(1+A+k)$ is very small. Then, we can set *a* such that $2\mu a$ > 1, and this *a* can still make $1 - 2\mu a_k > 0$ when $k \geq N$. Therefore, asymptotically for the case of $\alpha$ = 1, we have $E\left\|\hat{\boldsymbol{\theta}}_k - \boldsymbol{\theta}^*\right\|^2 = O(1/k)$.

Based on the result of Corollary 3.2 and the discussion above, let us summarize the asymptotic performance of DSPSA in a more formal way. Through the results of inequality (3.23a) and the fact $T(k,\alpha) \to lba^2 C(\alpha)/2\mu a$, we have that for $0.5 < \alpha < 1$,

$$E\left\|\hat{\boldsymbol{\theta}}_k - \boldsymbol{\theta}^*\right\|^2 = O\left(\exp\left(\frac{2\mu a(1+A)^{1-\alpha}}{1-\alpha} - \frac{2\mu a(1+A+k)^{1-\alpha}}{1-\alpha}\right)\right) + O\left(\frac{1}{k^{\alpha}}\right).$$

When $k$ is large, we have

$$\frac{2\mu a(1+A)^{1-\alpha}}{1-\alpha} - \frac{2\mu a(1+A+k)^{1-\alpha}}{1-\alpha} \approx -\frac{2\mu a}{1-\alpha}k^{1-\alpha},$$

which indicates that for the case of $0.5 < \alpha < 1$, we have

$$E\left\|\hat{\boldsymbol{\theta}}_k - \boldsymbol{\theta}^*\right\|^2 = O\left(\exp\left(-\frac{2\mu a}{1-\alpha}k^{1-\alpha}\right)\right) + O\left(\frac{1}{k^{\alpha}}\right). \qquad (3.26)$$

The second term on the right-hand side of eqn. (3.26) goes to 0 at a lower rate than the first term (first term = $o$(second term)), so under the case of $0.5 < \alpha < 1$, $E\left\|\hat{\boldsymbol{\theta}}_k - \boldsymbol{\theta}^*\right\|^2$ = $O\left(1/k^{\alpha}\right)$. For the case of $\alpha = 1$, through the results of inequality (3.23b), we have

$$E\left\|\hat{\boldsymbol{\theta}}_k - \boldsymbol{\theta}^*\right\|^2 = O\left(\frac{1}{k^{2\mu a}}\right) + O\left(\frac{1}{k}\right).$$

By the discussion above, when $\alpha = 1$, we have $E\left\|\hat{\boldsymbol{\theta}}_k - \boldsymbol{\theta}^*\right\|^2 = O\left(1/k\right)$.

From these results, we find when $\alpha = 1$, $E\left\|\hat{\boldsymbol{\theta}}_k - \boldsymbol{\theta}^*\right\|^2$ goes to 0 at the fastest rate from the allowable $\alpha \in (0.5, 1]$. However, the usable solution we get from each iteration is

$[\hat{\boldsymbol{\theta}}_k]$ instead of $\hat{\boldsymbol{\theta}}_k$. Therefore, when $k$ is large enough, we will have $[\hat{\boldsymbol{\theta}}_k] = \boldsymbol{\theta}^*$, which indicates $\alpha = 1$ may not provide better asymptotic performance in terms of $[\hat{\boldsymbol{\theta}}_k]$. As we have discussed in Section 3.1, we have $P([\hat{\boldsymbol{\theta}}_k] \neq \boldsymbol{\theta}^*) \leq 4E\left\|\hat{\boldsymbol{\theta}}_k - \boldsymbol{\theta}^*\right\|^2$. Then by the analysis in this section on $E\left\|\hat{\boldsymbol{\theta}}_k - \boldsymbol{\theta}^*\right\|^2$, we get $P([\hat{\boldsymbol{\theta}}_k] \neq \boldsymbol{\theta}^*) = O(1/k^{\alpha})$.

# 3.2 Properties of Upper Bound on Convergence Rate

In this section, we discuss the properties of the upper bound in inequality (3.3a) and (3.3b) to see how the upper bound changes with the change of coefficients of the gain sequence $\{a_k\}$. Later we will find the properties discussed here are consistent with the numerical results shown in Chapter 4, which means the upper bounds provided in Section 3.1 are reasonable.

From inequality (3.3a) and (3.3b), we see that the upper bounds for both cases are composed of two terms: the first term is corresponding to the initial guess and the second term is corresponding to the cumulative errors introduced in the whole process of the algorithm. In the early iterations, where the value of $k$ is small, the second term is small (the integration interval of the integral is small) and the first term is more significant in the upper bound (initial guess may be far away from the optimal solution). In the later iterations, where the value of $k$ is large, the first term is smaller (due to the faster decaying speed of the coefficient of the initial guess), and the second term becomes more significant than the first term (due to the slower convergence speed of the second term). Therefore, we have Proposition 3.1 as below.

**Proposition 3.1.** The ratio of the first term on the right-hand side of inequality (3.3a,b) (the initial guess term) over the upper bound in inequality (3.3a,b) is strictly decreasing with *k*. The ratio of the second term on the right-hand side of inequality (3.3a,b) over the upper bound in inequality (3.3a,b) is strictly increasing with *k*.

*Proof.* First let us consider the upper bound in inequality (3.3a). The ratio of the first term over the upper bound (3.3a) is defined as

$$\frac{\exp\left(\dfrac{2\mu a(1+A)^{1-\alpha}}{1-\alpha}-\dfrac{2\mu a(1+A+k)^{1-\alpha}}{1-\alpha}\right)E\left\|\hat{\boldsymbol{\theta}}_0-\boldsymbol{\theta}^*\right\|^2}{\text{upper bound}}.$$

Thus, to prove the ratio is strictly decreasing with *k* is equivalent to show that the ratio of the first term over the second term in the inequality (3.3a) is strictly decreasing with *k*. The ratio is

$$\frac{\exp\left(\dfrac{2\mu a(1+A)^{1-\alpha}}{1-\alpha}-\dfrac{2\mu a(1+A+k)^{1-\alpha}}{1-\alpha}\right)E\left\|\hat{\boldsymbol{\theta}}_0-\boldsymbol{\theta}^*\right\|^2}{\exp\left(-\dfrac{2\mu a(1+A+k)^{1-\alpha}}{1-\alpha}\right)lba^2C(\alpha)\displaystyle\int_0^k(1+A+x)^{-2\alpha}\exp\left(\dfrac{2\mu a(1+A+x)^{1-\alpha}}{1-\alpha}\right)dx}$$

$$=\frac{\exp\left(\dfrac{2\mu a(1+A)^{1-\alpha}}{1-\alpha}\right)E\left\|\hat{\boldsymbol{\theta}}_0-\boldsymbol{\theta}^*\right\|^2}{lba^2C(\alpha)\displaystyle\int_0^k(1+A+x)^{-2\alpha}\exp\left(\dfrac{2\mu a(1+A+x)^{1-\alpha}}{1-\alpha}\right)dx}.$$

We see that the denominator of the ratio strictly increases with $k$ and the numerator of the ratio is independent of $k$. Therefore, the ratio strictly decreases with $k$, which indicates the ratio of the first term on the right-hand side of inequality (3.3a) over the upper bound in inequality (3.3a) strictly decreases with $k$. Therefore, the ratio of the second term on the right-hand side of inequality (3.3a) over the upper bound in inequality (3.3a) strictly increases with $k$.

Similarly for the upper bound in inequality (3.3b), we have that the ratio of the first term over the second term is

$$\frac{(1+A)^{2\mu a} E\left\|\hat{\boldsymbol{\theta}}_0 - \boldsymbol{\theta}^*\right\|^2}{lba^2 C(1)\int_0^k (1+A+x)^{2\mu a-2}\,dx},$$

which is strictly decreasing with $k$. Therefore, the ratio of the first term on the right-hand side of inequality (3.3b) over the upper bound in equality (3.3b) strictly decreases with $k$, and the ratio of the second term on the right-hand side of inequality (3.3b) over the upper bound in inequality (3.3b) strictly increases with $k$. Q.E.D.

From Proposition 3.1, we know that in the early iterations, the value of the first term (initial guess term) on the right-hand side of inequality (3.3a,b) is more significant compared to the second term. However, in the later iterations, the second term on the right-hand side of inequality (3.3a,b) becomes more significant. Thus, for different stages, the properties of the upper bound are different.

In the following two propositions, we discuss how the upper bounds change with the coefficients $\alpha$, $a$, and $A$ in different stages.

**Proposition 3.2.** The first term on the right-hand side of inequality (3.3a) (the initial guess term),

$$\exp\left(\frac{2\mu a(1+A)^{1-\alpha}}{1-\alpha}-\frac{2\mu a(1+A+k)^{1-\alpha}}{1-\alpha}\right)E\left\|\hat{\boldsymbol{\theta}}_0-\boldsymbol{\theta}^*\right\|^2, \tag{3.27}$$

is a non-decreasing function on $\alpha\in(0.5,1)$, a non-increasing function on $a\in(0,\infty)$, and a non-decreasing function on $A\in[0,\infty)$.

*Proof.* By ignoring the positive multiplier $E\left\|\hat{\boldsymbol{\theta}}_0-\boldsymbol{\theta}^*\right\|^2$, we consider the derivatives of the function

$$\exp\left(\frac{2\mu a(1+A)^{1-\alpha}}{1-\alpha}-\frac{2\mu a(1+A+k)^{1-\alpha}}{1-\alpha}\right).$$

Let us start the proof by considering the derivative on $\alpha\in(0.5,1)$. By the chain rule, we have

$$\frac{\partial}{\partial\alpha}\exp\left(\frac{2\mu a(1+A)^{1-\alpha}}{1-\alpha}-\frac{2\mu a(1+A+k)^{1-\alpha}}{1-\alpha}\right)$$
$$=\exp\left(\frac{2\mu a(1+A)^{1-\alpha}}{1-\alpha}-\frac{2\mu a(1+A+k)^{1-\alpha}}{1-\alpha}\right)\frac{\partial}{\partial\alpha}\left(\frac{2\mu a(1+A)^{1-\alpha}}{1-\alpha}-\frac{2\mu a(1+A+k)^{1-\alpha}}{1-\alpha}\right).$$

In addition,

$$\frac{\partial}{\partial\alpha}\left(\frac{2\mu a(1+A)^{1-\alpha}}{1-\alpha}-\frac{2\mu a(1+A+k)^{1-\alpha}}{1-\alpha}\right)$$

$$=\frac{2\mu a\left(-(1+A)^{1-\alpha}\log(1+A)+(1+A+k)^{1-\alpha}\log(1+A+k)\right)(1-\alpha)}{(1-\alpha)^2}$$

$$-\frac{-2\mu a\left((1+A)^{1-\alpha}-(1+A+k)^{1-\alpha}\right)}{(1-\alpha)^2}$$

$$=\frac{2\mu a\left(-(1+A)^{1-\alpha}\left(\log(1+A)^{(1-\alpha)}-1\right)+(1+A+k)^{1-\alpha}\left(\log(1+A+k)^{(1-\alpha)}-1\right)\right)}{(1-\alpha)^2},$$

where the base of the logarithm here is *e* (natural logarithm). As we know, the derivative of $x(\log x-1)$ over $x$ is $\log x$. Thus, when $x\geq 1$, $x(\log x-1)$ is an increasing function on *x*. Here $(1+A)^{1-\alpha}\leq x\leq(1+A+k)^{1-\alpha}$, where we always have $x\geq 1$, which indicates $x(\log x-1)$ is an increasing function in the interval of $\left[(1+A)^{1-\alpha},(1+A+k)^{1-\alpha}\right]$. Therefore, we have $(1+A+k)^{1-\alpha}\left(\log(1+A+k)^{(1-\alpha)}-1\right)\geq(1+A)^{1-\alpha}\left(\log(1+A)^{(1-\alpha)}-1\right)$, which indicates that

$$\frac{2\mu a\left(-(1+A)^{1-\alpha}\left(\log(1+A)^{(1-\alpha)}-1\right)+(1+A+k)^{1-\alpha}\left(\log(1+A+k)^{(1-\alpha)}-1\right)\right)}{(1-\alpha)^2}\geq 0.$$

Thus, we have

$$\frac{\partial}{\partial \alpha} \exp\left( \frac{2\mu a(1+A)^{1-\alpha}}{1-\alpha} - \frac{2\mu a(1+A+k)^{1-\alpha}}{1-\alpha} \right)$$

$$= \exp\left( \frac{2\mu a(1+A)^{1-\alpha}}{1-\alpha} - \frac{2\mu a(1+A+k)^{1-\alpha}}{1-\alpha} \right)$$

$$\times \frac{2\mu a\left( -(1+A)^{1-\alpha}\left( \log(1+A)^{1-\alpha} - 1 \right) + (1+A+k)^{1-\alpha}\left( \log(1+A+k)^{1-\alpha} - 1 \right) \right)}{(1-\alpha)^2}$$

$$\geq 0,$$

so the function (3.27) is a non-decreasing function on $\alpha \in (0.5, 1)$.

Second, let us consider the derivative over $a > 0$. We have

$$\frac{\partial}{\partial a} \exp\left( \frac{2\mu a(1+A)^{1-\alpha}}{1-\alpha} - \frac{2\mu a(1+A+k)^{1-\alpha}}{1-\alpha} \right)$$

$$= \exp\left( \frac{2\mu a(1+A)^{1-\alpha}}{1-\alpha} - \frac{2\mu a(1+A+k)^{1-\alpha}}{1-\alpha} \right)\left( \frac{2\mu}{1-\alpha}\left( (1+A)^{1-\alpha} - (1+A+k)^{1-\alpha} \right) \right)$$

$$\leq 0,$$

which indicates that the function (3.27) is a non-increasing function on $a > 0$.

Last, let us calculate the derivative over $A \geq 0$. We have

$$\frac{\partial}{\partial A} \exp\left( \frac{2\mu a(1+A)^{1-\alpha}}{1-\alpha} - \frac{2\mu a(1+A+k)^{1-\alpha}}{1-\alpha} \right)$$

$$= \exp\left( \frac{2\mu a(1+A)^{1-\alpha}}{1-\alpha} - \frac{2\mu a(1+A+k)^{1-\alpha}}{1-\alpha} \right)\left( 2\mu a\left( \frac{1}{(1+A)^{\alpha}} - \frac{1}{(1+A+k)^{\alpha}} \right) \right)$$

$$\geq 0,$$

which indicates that the function (3.27) is a non-decreasing function on $A \geq 0$. Q.E.D.

In Proposition 3.2, we see the properties of the upper bound in the early iterations on the coefficients of gain sequence. From the proof of Proposition 3.2, we find that the magnitude of the derivative over $\alpha$ is more significant compared to the magnitudes of the derivatives over $a$ and $A$ (the magnitude of the derivative over $A$ is the least significant one) in the big-$O$ sense. Therefore, based on properties of the first term of the upper bound, we prefer to pick relatively smaller $\alpha$ and tune the values of $a$ and $A$ to satisfy other requirements such as the stability of the algorithm.

From inequality (3.23a), we find that the second term ($T(k,\alpha)\big/(1+A+k)^{\alpha}$) on the right-hand side goes to 0 at a slower rate than the first term (first term = $o$(second term)). Thus, when $k$ is large, in inequality (3.23a), the second term ($T(k,\alpha)\big/(1+A+k)^{\alpha}$) on the right-hand side is the leading term in the upper bound, which indicates that discussing the properties of the second term in inequality (3.23a) in the later iterations is equivalent to discussing the properties of the upper bound in the later iterations. As we have discussed in Corollary 3.2, $T(k,\alpha)\to lbaC(\alpha)/2\mu$, so the second term on the right-hand side of inequality (3.23a) is very close to $lbaC(\alpha)/(2\mu(1+A+k)^{\alpha})$ when $k$ is large. Since it is not easy to discuss the second term on the right-hand side of inequality (3.23a) directly, we will consider its approximation $lbaC(\alpha)/(2\mu(1+A+k)^{\alpha})$ instead. In Proposition 3.3 below, we discuss the property of the function $lbaC(\alpha)/(2\mu(1+A+k)^{\alpha})$ on the coefficients of the gain sequence to see the properties of the upper bound in the later iterations.

**Proposition 3.3** The function

$$\frac{lbaC(\alpha)}{2\mu(1+A+k)^{\alpha}} \tag{3.28}$$

is strictly increasing on $a \in (0,\infty)$, and is strictly decreasing on $A \in [0,\infty)$. Function (3.28) is a strictly decreasing function on $\alpha \in (0.5, 1)$ when *k* is large enough to make

$$\frac{(2+A)^2/(1+A)^2}{(1+A+k)} < 1.$$

*Proof*. By ignoring the positive multiplier $lb/2\mu$ in function (3.28), we discuss the function $aC(\alpha)/(1+A+k)^{\alpha}$. Since

$$C(\alpha) = \left( \exp\left( \frac{2\mu a(1+A+1)^{1-\alpha}}{1-\alpha} - \frac{2\mu a(1+A)^{1-\alpha}}{1-\alpha} \right) \right) \left( 1 + \frac{1}{1+A} \right)^{2\alpha},$$

then

$$\frac{aC(\alpha)}{(1+A+k)^{\alpha}} = a\left( \exp\left( \frac{2\mu a(1+A+1)^{1-\alpha}}{1-\alpha} - \frac{2\mu a(1+A)^{1-\alpha}}{1-\alpha} \right) \right) \left( \frac{(2+A)^2/(1+A)^2}{(1+A+k)} \right)^{\alpha}.$$

After ignoring the positive multiplier (in function (3.28)) that is not related to the coefficients, we calculate the derivatives on $\alpha \in (0.5, 1)$, $a \in (0,\infty)$, and $A \in [0,\infty)$ as below:

$$\frac{\partial}{\partial\alpha}\left(\frac{aC(\alpha)}{(1+A+k)^{\alpha}}\right)$$

$$=a\left(\frac{(2+A)^2/(1+A)^2}{(1+A+k)}\right)^{\alpha}\frac{\partial}{\partial\alpha}\exp\left(\frac{2\mu a(1+A+1)^{1-\alpha}}{1-\alpha}-\frac{2\mu a(1+A)^{1-\alpha}}{1-\alpha}\right)$$

$$+a\exp\left(\frac{2\mu a(1+A+1)^{1-\alpha}}{1-\alpha}-\frac{2\mu a(1+A)^{1-\alpha}}{1-\alpha}\right)\frac{\partial}{\partial\alpha}\left(\frac{(2+A)^2/(1+A)^2}{(1+A+k)}\right)^{\alpha}, \tag{3.29}$$

$$\frac{\partial}{\partial a}\left(\frac{aC(\alpha)}{(1+A+k)^{\alpha}}\right)$$

$$=a\left(\frac{(2+A)^2/(1+A)^2}{(1+A+k)}\right)^{\alpha}\frac{\partial}{\partial a}\exp\left(\frac{2\mu a(1+A+1)^{1-\alpha}}{1-\alpha}-\frac{2\mu a(1+A)^{1-\alpha}}{1-\alpha}\right)$$

$$+\exp\left(\frac{2\mu a(1+A+1)^{1-\alpha}}{1-\alpha}-\frac{2\mu a(1+A)^{1-\alpha}}{1-\alpha}\right)\left(\frac{(2+A)^2/(1+A)^2}{(1+A+k)}\right)^{\alpha}, \tag{3.30}$$

and

$$\frac{\partial}{\partial A}\left(\frac{aC(\alpha)}{(1+A+k)^{\alpha}}\right)$$

$$=a\left(\frac{(2+A)^2/(1+A)^2}{(1+A+k)}\right)^{\alpha}\frac{\partial}{\partial A}\exp\left(\frac{2\mu a(1+A+1)^{1-\alpha}}{1-\alpha}-\frac{2\mu a(1+A)^{1-\alpha}}{1-\alpha}\right)$$

$$+a\exp\left(\frac{2\mu a(1+A+1)^{1-\alpha}}{1-\alpha}-\frac{2\mu a(1+A)^{1-\alpha}}{1-\alpha}\right)\frac{\partial}{\partial A}\left(\frac{(2+A)^2/(1+A)^2}{(1+A+k)}\right)^{\alpha}. \tag{3.31}$$

Now let us calculate all the derivatives on the right-hand side of eqns. (3.29), (3.30) and (3.31). By similar arguments as in the Proposition 3.2, we know

$$\frac{\partial}{\partial\alpha}\exp\left(\frac{2\mu a(1+A+1)^{1-\alpha}}{1-\alpha}-\frac{2\mu a(1+A)^{1-\alpha}}{1-\alpha}\right)<0,$$

$$\frac{\partial}{\partial a}\exp\left(\frac{2\mu a(1+A+1)^{1-\alpha}}{1-\alpha}-\frac{2\mu a(1+A)^{1-\alpha}}{1-\alpha}\right)>0,$$

and

$$\frac{\partial}{\partial A}\exp\left(\frac{2\mu a(1+A+1)^{1-\alpha}}{1-\alpha}-\frac{2\mu a(1+A)^{1-\alpha}}{1-\alpha}\right)<0.$$

Moreover, we also have

$$\frac{\partial}{\partial\alpha}\left(\frac{(2+A)^2/(1+A)^2}{(1+A+k)}\right)^{\alpha}=\left(\frac{(2+A)^2/(1+A)^2}{(1+A+k)}\right)^{\alpha}\log\left(\frac{(2+A)^2/(1+A)^2}{(1+A+k)}\right)<0,$$

when $k$ is large enough to make

$$\frac{(2+A)^2/(1+A)^2}{(1+A+k)}<1.$$

In addition,

$$\begin{aligned}
&\frac{\partial}{\partial A}\left(\frac{(2+A)^2/(1+A)^2}{(1+A+k)}\right)^{\alpha}\\
&=\frac{1}{(1+A+k)^{\alpha}}\frac{\partial}{\partial A}\left(1+\frac{1}{1+A}\right)^{2\alpha}+\left(1+\frac{1}{1+A}\right)^{2\alpha}\frac{\partial}{\partial A}\left(\frac{1}{(1+A+k)^{\alpha}}\right)\\
&=\frac{-2\alpha}{(1+A+k)^{\alpha}(1+A)^2}\left(1+\frac{1}{1+A}\right)^{2\alpha-1}+\left(1+\frac{1}{1+A}\right)^{2\alpha}\left(\frac{-\alpha}{(1+A+k)^{\alpha+1}}\right)\\
&<0.
\end{aligned}$$

Overall, we have calculated all terms in the derivatives on the right-hand side of eqns. (3.29), (3.30) and (3.31) and know the signs of all derivatives. Therefore, we have

$$\frac{\partial}{\partial a}\left(\frac{aC(\alpha)}{(1+A+k)^{\alpha}}\right) > 0 \text{ and } \frac{\partial}{\partial A}\left(\frac{aC(\alpha)}{(1+A+k)^{\alpha}}\right) < 0,$$

which implies that in the later iterations, we have that $(lbaC(\alpha))/(2\mu(1+A+k)^{\alpha})$ is a strictly increasing function on $a \in (0,\infty)$, and a strictly decreasing function on $A \in [0,\infty)$. We also have

$$\frac{\partial}{\partial \alpha}\left(\frac{aC(\alpha)}{(1+A+k)^{\alpha}}\right) < 0,$$

when $k$ is large enough to make

$$\frac{(2+A)^2/(1+A)^2}{(1+A+k)} < 1,$$

which implies that $(lbaC(\alpha))/(2\mu(1+A+k)^{\alpha})$ is a strictly decreasing function on $\alpha \in (0.5,1)$ when $k$ is large enough. Q.E.D.

In Proposition 3.3, we see the properties of the upper bound in inequality (3.3a) in the later iterations related to the coefficients of gain sequence. From the proof of Proposition 3.3, we find that the magnitude of the derivative over $\alpha$ is most significant compared to the magnitudes of the derivatives over $a$ and $A$ (the magnitude of the derivative over $A$ is the least significant one) in the big-$O$ sense. Therefore, based on properties of the upper bound in the later iteration, we prefer to pick relatively bigger $\alpha$ and tune the values of $a$ and $A$ to satisfy other requirements such as the stability of the algorithm.

Comparing the results of Proposition 3.2 and Proposition 3.3, we find that the effects of the coefficients of the gain sequence on the upper bound are in the opposite directions in different stages. Because different terms dominate the upper bound in different stages (finite samples stage and asymptotical stage), the good set of coefficients for better finite sample performance may not be the good set for better asymptotic performance under the criterion $E\left\|\hat{\boldsymbol{\theta}}_k - \boldsymbol{\theta}^*\right\|^2$ and vice versa.

In Appendix C, we will run some numerical experiments to see the properties of the upper bound numerically for the case when the $\Delta_{ki}$ are independent Bernoulli random variables taking the values $\pm 1$ with probability $1/2$. In Appendix C, we will see that these numerical results are consistent with the theoretical analysis of the properties of the upper bounds discussed above. Moreover, in Chapter 4, we will see that the performances of DSPSA have the same properties as the upper bound and it indicates that this upper bound is meaningful in the sense that it captures the true properties of the performances of DSPSA.

In the following, we discuss a little bit on the choice of μ and $b$, which arise in the upper bound of $E\left\|\hat{\boldsymbol{\theta}}_k - \boldsymbol{\theta}^*\right\|^2$. The values of μ and $b$ are dependent on the loss functions. From the condition (v) of Theorem 3.1, we have the relationship that $0 < \mu \leq E\left[(\hat{\boldsymbol{\theta}}_k - \boldsymbol{\theta}^*)^T \overline{\boldsymbol{g}}(\boldsymbol{\pi}(\hat{\boldsymbol{\theta}}_k))\right] \Big/ E\left[(\hat{\boldsymbol{\theta}}_k - \boldsymbol{\theta}^*)^T (\hat{\boldsymbol{\theta}}_k - \boldsymbol{\theta}^*)\right]$ for all $k \geq 0$. Furthermore, $b$ is a uniform upper bound for the value of $E\left(L(\hat{\boldsymbol{\theta}}_k^+) - L(\hat{\boldsymbol{\theta}}_k^-)\right)^2 + E(\varepsilon_k^+ - \varepsilon_k^-)^2$ for all $k \geq 0$. Therefore, the values of μ and $b$ are restricted by some iterations. The choice of μ and $b$ affect the tightness of the upper bound on $E\left\|\hat{\boldsymbol{\theta}}_k - \boldsymbol{\theta}^*\right\|^2$. In the Appendix C, we use the

separable quadratic loss function as an example to discuss how to pick the value of $\mu$ and $b$, and how these values affect the tightness of the upper bound. In these numerical tests, we consider the special case that the $\Delta_{ki}$ are independent Bernoulli random variables taking the values $\pm 1$ with probability $1/2$. From these results in Appendix C, we see the tradeoff between the tightness of the upper bound in the early iterations and the tightness of the upper bound in the later iterations through the choice of $\mu$ and $b$.

# 3.3 Choice of Gain Sequence

In this section, we discuss the implementation aspects regarding the choice of the gain sequence under the criterion of $E\left\|\hat{\boldsymbol{\theta}}_k - \boldsymbol{\theta}^*\right\|^2$ and the assumptions in Theorem 3.1. The idea in this section is similar to the idea in Appendix A, where we discuss the practical step size selection for continuous stochastic approximation algorithms. We consider the choice for both finite sample performance and asymptotic performance. We provide some guidelines here.

The choice for finite sample performance is based on the discussion of MSE $E\left\|\hat{\boldsymbol{\theta}}_k - \boldsymbol{\theta}^*\right\|^2$ in Section 3.1. Note that eqn. (3.6) can be written as

$$E\left\|\hat{\boldsymbol{\theta}}_{k+1} - \boldsymbol{\theta}^*\right\|^2 = E\left\|\hat{\boldsymbol{\theta}}_k - \boldsymbol{\theta}^*\right\|^2 - 2a_k E\left((\hat{\boldsymbol{\theta}}_k - \boldsymbol{\theta}^*)^T \bar{\boldsymbol{g}}(\boldsymbol{\pi}(\hat{\boldsymbol{\theta}}_k))\right)$$

$$+ a_k^2 E\left[\left(\left(L(\hat{\boldsymbol{\theta}}_k^+) - L(\hat{\boldsymbol{\theta}}_k^-)\right)^2 + (\varepsilon_k^+ - \varepsilon_k^-)^2\right)\boldsymbol{\Delta}_k^{-T}\boldsymbol{\Delta}_k^{-1}\right]. \quad (3.32)$$

By the condition (v) in Theorem 3.1, we have $E\left[(\hat{\boldsymbol{\theta}}_k - \boldsymbol{\theta}^*)^T \bar{\boldsymbol{g}}(\boldsymbol{\pi}(\hat{\boldsymbol{\theta}}_k))\right] \geq \mu E\left[(\hat{\boldsymbol{\theta}}_k - \boldsymbol{\theta}^*)^T (\hat{\boldsymbol{\theta}}_k - \boldsymbol{\theta}^*)\right]$. Let

$$\mu_k = \frac{E\left((\hat{\boldsymbol{\theta}}_k - \boldsymbol{\theta}^*)^T \bar{\boldsymbol{g}}(\boldsymbol{\pi}(\hat{\boldsymbol{\theta}}_k))\right)}{E\left\|\hat{\boldsymbol{\theta}}_k - \boldsymbol{\theta}^*\right\|^2} > 0\,. \tag{3.33}$$

By eqn. (3.33) and the recursive eqn. (3.32), we have

$$\begin{aligned}
&E\left\|\hat{\boldsymbol{\theta}}_{k+1} - \boldsymbol{\theta}^*\right\|^2 \\
&= \left(1 - 2a_k\mu_k\right) E\left\|\hat{\boldsymbol{\theta}}_k - \boldsymbol{\theta}^*\right\|^2 + a_k^2 E\left[\left(\left(L(\hat{\boldsymbol{\theta}}_k^+) - L(\hat{\boldsymbol{\theta}}_k^-)\right)^2 + (\varepsilon_k^+ - \varepsilon_k^-)^2\right)\boldsymbol{\Delta}_k^{-T}\boldsymbol{\Delta}_k^{-1}\right] \\
&= \prod_{j=0}^{k}(1 - 2a_j\mu_j) E\left\|\hat{\boldsymbol{\theta}}_0 - \boldsymbol{\theta}^*\right\|^2 \\
&\quad + \sum_{i=0}^{k}\prod_{j=i+1}^{k}(1 - 2a_j\mu_j) a_i^2 E\left[\left(\left(L(\hat{\boldsymbol{\theta}}_i^+) - L(\hat{\boldsymbol{\theta}}_i^-)\right)^2 + (\varepsilon_i^+ - \varepsilon_i^-)^2\right)\boldsymbol{\Delta}_i^{-T}\boldsymbol{\Delta}_i^{-1}\right].
\end{aligned} \tag{3.34}$$

There are three terms on the right-hand side of eqn. (3.32). The first term is related to the MSE in the previous iteration, the second term is negative, and the third term is positive. In order to have stable performance (the algorithm is stable means that the sequence generated by the algorithm stays within a reasonable distance of the optimum) in the early iterations, the third term of eqn. (3.32) should have small values. It can be achieved by small value of step size, which is equivalent to having the value of $a$ not too big and the value of $A$ not too small (we regard these requirements as the constraints on the gain sequence for stable performance).

Now we consider the choice of gain sequence for good finite sample performance under the constraints for stable performance. Suppose $\hat{\boldsymbol{\theta}}_0$ is far from the optimal solution $\boldsymbol{\theta}^*$, and $\mu_k$ does not change significantly on different sets of coefficients in the early iterations. Then, the first term on the right-hand side of eqn. (3.34) is the dominant term in the early iterations, so we prefer relatively bigger value of $a_k$. We can achieve relatively big value of $a_k$ by choosing small value of $\alpha$ ($\alpha = 0.501$), large value of $a$ and small value of $A$. We see that the requirements for $a$ and $A$ are on the opposite sides for stability and better finite sample performance. In order to have better finite sample performance, we need to pick $a$ relatively large and pick $A$ relatively small in the domain that provides stable performance.

Now we start to discuss a clear rule for the choice of $a$ and $A$. For some cases, the cost of noisy measurements of the loss functions might be high, so the number of iterations may be limited. In order to achieve reasonable performance for the limited number of iterations, we want $A$ to be proportional to the maximum number of allowed iterations $M$. Let $A = \eta M$, where $\eta > 0$. It indicates that for small number of allowed iterations, we prefer $A$ to be small, which can make the gain sequence large. A large step size can lead the sequence generated by DSPSA to a reasonable result even for very limited number of iterations. In addition, we want the effect of $A$ to disappear for later iterations to achieve proper decaying gain. Thus, we prefer $A$ to be at a lower value than $M$, which indicates $\eta$ = 0.1, 0.01, 0.001, ... Moreover, as we have discussed, we want $a/(1+A+i)^{\alpha}$ to be reasonably small in the early iterations to keep the stability of the algorithm, but we also want the gain step not too small in the later iterations to foster a reasonable speed of convergence. We know the effect of $A$ disappears when $k$ is large, and the effect of $a$ is always there, which follows that we prefer to pick large $A$ instead of small $a$ to keep the stability of the algorithm, because we still need $a$ not to be too small for better finite sample performance. Therefore, among all the choices of $\eta$, we prefer a larger value, which is a value such as $\eta = 0.1$, and this guideline is similar to the guideline for SPSA in continuous problems (Spall, 2003, p. 190). We find that $A = 0.1M$ satisfies all the requirements, including the stability and finite sample performance. After getting the rule of $A$, we can make the multiplication of $a_0$ and the magnitude of elements in $\hat{\boldsymbol{g}}_0(\hat{\boldsymbol{\theta}}_0)$ approximately equal to the desired change of the magnitude of $\hat{\boldsymbol{\theta}}_k$ in the early iterations, which leads to the possible choice of $a$.

Furthermore, let us discuss the choice of coefficients for better asymptotic performance. Suppose we have $E\left\|\hat{\boldsymbol{\theta}}_k - \boldsymbol{\theta}^*\right\|^2 \to 0$ as $k \to \infty$. When $k \to \infty$, the generated sequence $\{\hat{\boldsymbol{\theta}}_k\}$ bounces around the optimal solution. Asymptotically we want the sequence to bounce as little as possible, and we can achieve it by picking smaller $a_k = a/(1 + A + k)^{\alpha}$ (the discussion here on asymptotic performance not only focus on the exponent of the big-$O$ function, but also focus on the constant multiplier of the big-$O$ function). From eqn. (3.32), we have

$$E\left\|\hat{\boldsymbol{\theta}}_{k+1}-\boldsymbol{\theta}^*\right\|^2 - E\left\|\hat{\boldsymbol{\theta}}_k-\boldsymbol{\theta}^*\right\|^2 = -2a_k E\left((\hat{\boldsymbol{\theta}}_k-\boldsymbol{\theta}^*)^T \bar{\boldsymbol{g}}(\boldsymbol{\pi}(\hat{\boldsymbol{\theta}}_k))\right)$$

$$+a_k^2 E\left[\left(\left(L(\hat{\boldsymbol{\theta}}_k^+)-L(\hat{\boldsymbol{\theta}}_k^-)\right)^2+(\varepsilon_k^+-\varepsilon_k^-)^2\right)\boldsymbol{\Delta}_k^{-T}\boldsymbol{\Delta}_k^{-1}\right].$$

Since $E\left\|\hat{\boldsymbol{\theta}}_k-\boldsymbol{\theta}^*\right\|^2 \to 0$, then asymptotically we want $E\left\|\hat{\boldsymbol{\theta}}_{k+1}-\boldsymbol{\theta}^*\right\|^2 - E\left\|\hat{\boldsymbol{\theta}}_k-\boldsymbol{\theta}^*\right\|^2$ to converges to 0 as fast as possible to achieve higher convergence rate, which indicates smaller value of gain sequence is preferred asymptotically. Taking the derivatives of $a_k$ over $\alpha$, $a$ and $A$, we get

$$\frac{\partial}{\partial\alpha}\frac{a}{(1+A+k)^\alpha}=\frac{-a\log(1+A+k)}{(1+A+k)^\alpha}\le 0,$$

$$\frac{\partial}{\partial a}\frac{a}{(1+A+k)^\alpha}=\frac{1}{(1+A+k)^\alpha}\ge 0,$$

$$\frac{\partial}{\partial A}\frac{a}{(1+A+k)^\alpha}=\frac{-a\alpha}{(1+A+k)^{\alpha+1}}\le 0.$$

We find that the magnitude of the derivative over $\alpha$ is most significant compared to the magnitudes of the derivatives over $a$ and $A$ (the magnitude of the derivative over $A$ is the least significant one, so the effect of $A$ is negligible for asymptotic performance). Therefore, asymptotically we prefer relatively larger $\alpha$.

Overall, based on the discussion in this section, we know that for the selection of the coefficients we need to consider three parts: the stability of the algorithm, the finite sample performance, and the asymptotic performance. Each part prefers different set of coefficients. If we want to have stable performance and good finite sample performance, we prefer $\alpha = 0.501$, $A = 0.1M$ and make the multiplication of $a_0$ and the magnitude of elements in $\hat{\boldsymbol{g}}_0(\hat{\boldsymbol{\theta}}_0)$ approximately equal to the desired change of the magnitude of $\hat{\boldsymbol{\theta}}_k$ in the early iterations, which leads to the possible choice of *a*. If we focus more on the asymptotic performance under the criterion $E\left\|\hat{\boldsymbol{\theta}}_k - \boldsymbol{\theta}^*\right\|^2$, we prefer $\alpha = 1$, which is identical to the asymptotically optimal choice of $\alpha$ for the continuous case (e.g. Spall, 2003, Chapter 4).

# Chapter 4

# Performance of DSPSA

We have introduced the algorithm of DSPSA and discuss the convergence properties of DSPSA theoretically in the last two chapters. In this chapter, first we do some numerical experiments to check the effects of the coefficients of the gain sequence. Then, based on the guidelines of choice of coefficients discussed in Section 3.3, we test the performance of DSPSA on some general loss functions numerically.

## 4.1 Description of Loss Functions

In this chapter, we do the numerical tests on the case when the $\Delta_{ki}$ are independent Bernoulli random variables taking the values $\pm 1$ with probability $1/2$. We consider three loss functions here. The first function is a separable loss function defined on $\mathbb{Z}^p$

$$L(\boldsymbol{\theta}) = \boldsymbol{\theta}^T \boldsymbol{\theta} = \sum_{i=1}^{p} t_i^2 , \tag{4.1}$$

where $\boldsymbol{\theta} = [t_1, ..., t_p]^T$ . The second loss function is a quadratic function defined on $\mathbb{Z}^p$

$$L(\boldsymbol{\theta}) = \boldsymbol{\theta}^T \boldsymbol{D}^T \boldsymbol{D} \boldsymbol{\theta} - 2\boldsymbol{d}^T \boldsymbol{D} \boldsymbol{\theta} + \boldsymbol{d}^T \boldsymbol{d} , \tag{4.2}$$

where $\boldsymbol{D}$ is a matrix with each diagonal component being $1 + 1/p$ and all the other components being $1/p$ , and $\boldsymbol{d}$ is a vector with each component being 2. Here matrix $\boldsymbol{D}$ is a strict diagonal dominant matrix. By Theorem (6.1.10) in Horn and Johnson (1985), we know that strict diagonal dominant matrix is positive definite. Therefore, we have that $\boldsymbol{D}$ is positive definite, which indicates that the second loss functions has a strictly convex quadratic continuous extension. In addition, we see that the second loss function can be written as $L(\boldsymbol{\theta}) = (\boldsymbol{D}\boldsymbol{\theta} - \boldsymbol{d})^T (\boldsymbol{D}\boldsymbol{\theta} - \boldsymbol{d})$ , and the optimal solution is $\boldsymbol{\theta}^* = [1, ..., 1]^T \equiv \mathbf{1}_p^T$ . The third loss function is a skewed quartic function (Spall, 2003, Ex 6.6) defined on $\mathbb{Z}^p$ :

$$L(\boldsymbol{\theta}) = \boldsymbol{\theta}^T \boldsymbol{B}^T \boldsymbol{B} \boldsymbol{\theta} + 0.1 \sum_{i=1}^{p} (\boldsymbol{B}\boldsymbol{\theta})_i^3 + 0.01 \sum_{i=1}^{p} (\boldsymbol{B}\boldsymbol{\theta})_i^4 , \tag{4.3}$$

where $(\cdot)_i$ represents the $i$th component of the vector $\boldsymbol{B}\boldsymbol{\theta}$ and $p\boldsymbol{B}$ is an upper triangular matrix of 1's.

In Theorem 2.1 (almost sure convergence theorem), condition (i), (ii), (iii), (iv) and (vi) are general conditions, and these conditions are not hard to check. Then condition (v) ($\bar{\boldsymbol{g}}(\boldsymbol{m}_{\boldsymbol{\theta}})^T (\boldsymbol{\theta} - \boldsymbol{\theta}^*) > 0$ for all $\boldsymbol{m}_{\boldsymbol{\theta}} \in \mathcal{M}_{\boldsymbol{\theta}}$ and all $\boldsymbol{\theta} \in \mathbb{R}^p \setminus \{\boldsymbol{\theta}^*\}$) in Theorem 2.1 is the main focus. By Proposition 2.1 in Section 2.2, the first loss function satisfies condition (v) in Theorem 2.1. By Proposition 2.3 in Section 2.2, we know the second loss function also satisfies condition (v) in Theorem 2.1.

However, the third loss function does not satisfy condition (v) in Theorem 2.1, and we will show it in the following. In the numerical test in the Section 4.5 on the third loss function, we will test high-dimensional case, where $p = 200$. Let us consider the point $\boldsymbol{\theta}$ with all components being 0 except the $200^{\text{th}}$ component being $-1$, then the unit hypercube centered by $\boldsymbol{\pi}(\boldsymbol{\theta})$ is a vector with all components being 0.5 except the $200^{\text{th}}$ component being $-0.5$. By the definition shown in Section 2.1, we know that when the $\Delta_{ki}$ are independent Bernoulli random variables taking the values $\pm 1$ with probability $1/2$

$$\bar{\boldsymbol{g}}\left(\boldsymbol{\pi}(\boldsymbol{\theta})\right) = \frac{1}{2^p} \sum\nolimits_{\boldsymbol{\Delta}} \left( L\left(\boldsymbol{\theta} + \frac{1}{2}\boldsymbol{\Delta}\right) - L\left(\boldsymbol{\theta} - \frac{1}{2}\boldsymbol{\Delta}\right)\right)\boldsymbol{\Delta}^{-1},$$

and we can calculate the value of $\bar{\boldsymbol{g}}\left(\boldsymbol{\pi}(\boldsymbol{\theta})\right)$ directly. We find that $\bar{\boldsymbol{g}}(\boldsymbol{\pi}(\boldsymbol{\theta}))^T(\boldsymbol{\theta} - \boldsymbol{\theta}^*) \approx$ $-0.52 < 0$ (there are $2^{200}$ possible choices for $\boldsymbol{\Delta}$; we randomly pick 5,000,000 of them to do the approximation, and the final estimator is the mean value of 20 replicates), where $\boldsymbol{\theta}^* = \mathbf{0}_{200}$, and $\mathbf{0}_{200}$ is a 200-dimensional vector with all components being 0. Therefore, the condition (v) in Theorem 2.1 is not satisfied for the skewed quartic loss function defined on $\mathbb{Z}^p$. However, we will see that DSPSA still works for high-dimensional skewed quartic loss function defined on $\mathbb{Z}^p$ in Section 4.5.

## 4.2 Numerical Tests on the Effects of Coefficients of the Gain Sequence

Before doing the numerical tests on the loss functions mentioned above, we first do some simple tests to show that the discussion in Section 3.3 on the effects of coefficients of the gain sequence are consistent with the numerical results. We pick the separable loss function (4.1) as the example to do the tests. Suppose the measurement noises $\varepsilon$ are i.i.d. $N(0,1)$, the dimension $p$ = 10, the initial guess is $10\times\mathbf{1}_{10}$, and the gain step $a_k$ = $a/(k+1+A)^{\alpha}$. Let the number of iterations in each replicate be 1000 and let the number of replicates be 20. We use the sensitivity analysis for the tests, which means that in each test we only change one coefficient and fix all the rest. Suppose the base case is: $\alpha$ = 0.75, $A$ = 100, $a$ = 0.4.

In Figure 4.1, we test the effect of $a$. The values that we choose for $a$ are 0.4 and 1. From Figure 4.1, we see that after having stable performance for DSPSA, relatively larger $a$ leads to better performance in the early iterations. Meanwhile, smaller value of $a$ provides better performance in later iterations. We see a crossing point at just under 300 iterations between the line for $a$ = 0.4 and $a$ = 1. As the discussion in Section 3.3, smaller value of $a$ can help to provide stable performance for DSPSA in the early iterations, after obtaining stable performance, relatively bigger value of $a$ can work better in the early iterations, and smaller value of $a$ leads to better asymptotic performance. Therefore, the numerical result here is consistent with the discussion in Section 3.3.

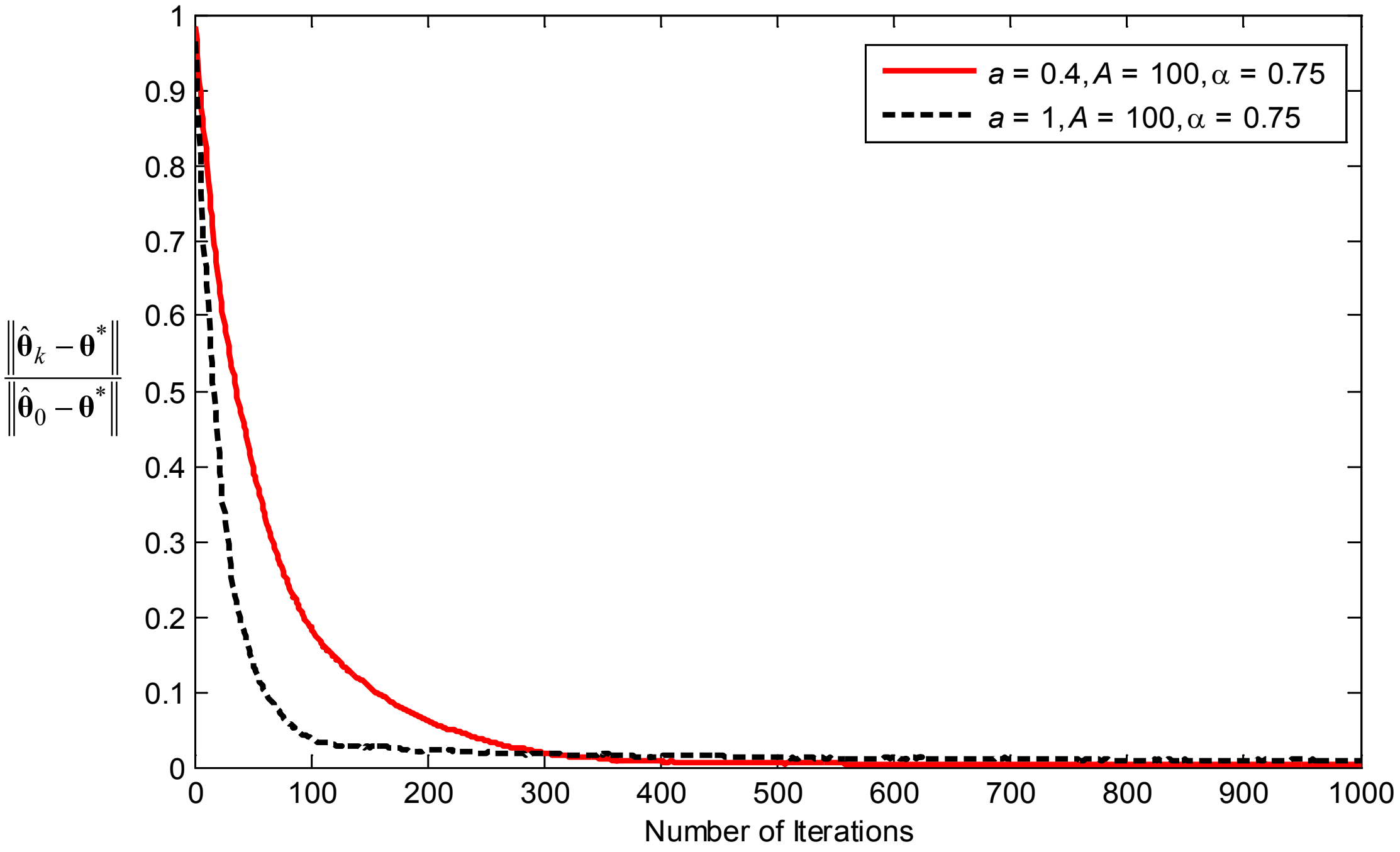

**Figure 4.1** Performance comparison with different values of $a$ for separable loss function. Under the premise of stable performance, by fixing the values of $A$ and $\alpha$, we find that the relatively larger value of $a$ leads to better performance in the early iterations; while smaller value of $a$ provides better performance asymptotically. Each curve represents the sample mean of 20 independent replicates.

In Figure 4.2, we test the effect of $A$. The values that we choose for $A$ are 10 and 100. We find that the effect of $A$ disappears in the later iterations. After having stable performance for DSPSA, the relatively smaller $A$ works better in the early iterations. As the discussion in Section 3.3, big value of $A$ can help to provide stable performance for DSPSA in the early iterations, and after obtaining stable performance, relatively smaller value of $A$ can provide better performance in the early iterations. Therefore, the numerical result here is consistent with the discussion in Section 3.3.

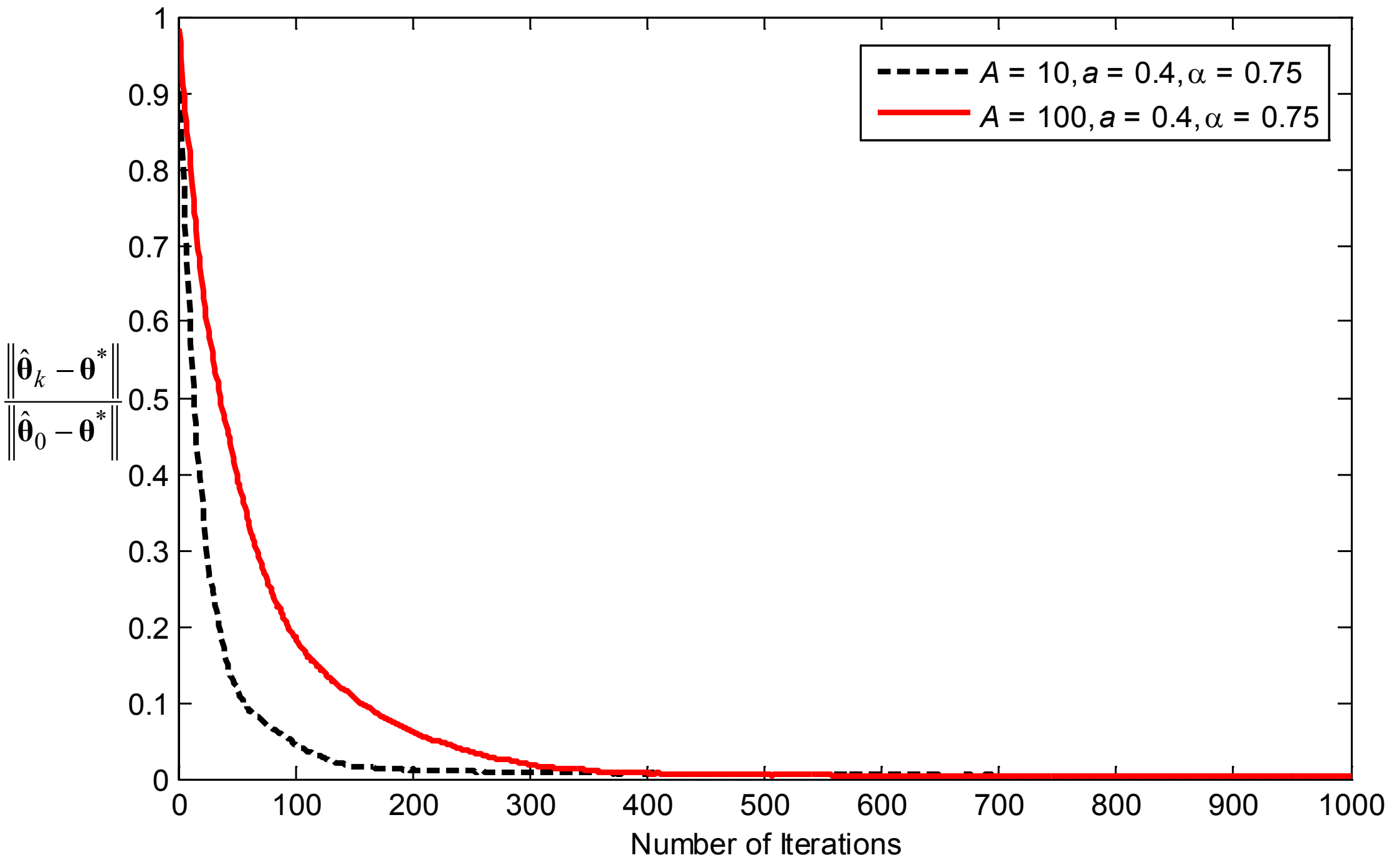

**Figure 4.2** Performance comparison with different values of $A$ for separable loss function. Under the premise of stable performance, by fixing the values of $a$ and $\alpha$, we find that the relatively smaller value of $A$ provides better performance in the early iterations; while the effect of $A$ is negligible asymptotically. Each curve represents the sample mean of 20 independent replicates.

In Figure 4.3, we test the effect of $\alpha$. The values that we choose for $\alpha$ are 0.501 and 0.75. Figure 4.3 tells that for finite sample performance in the early iterations, smaller $\alpha$ gives better performance; while for later iterations, larger α leads to better asymptotic performance. We also see a clear crossing point between the line for α = 0.501 and α = 0.75. The numerical result here is also consistent with the discussion in Section 3.3.

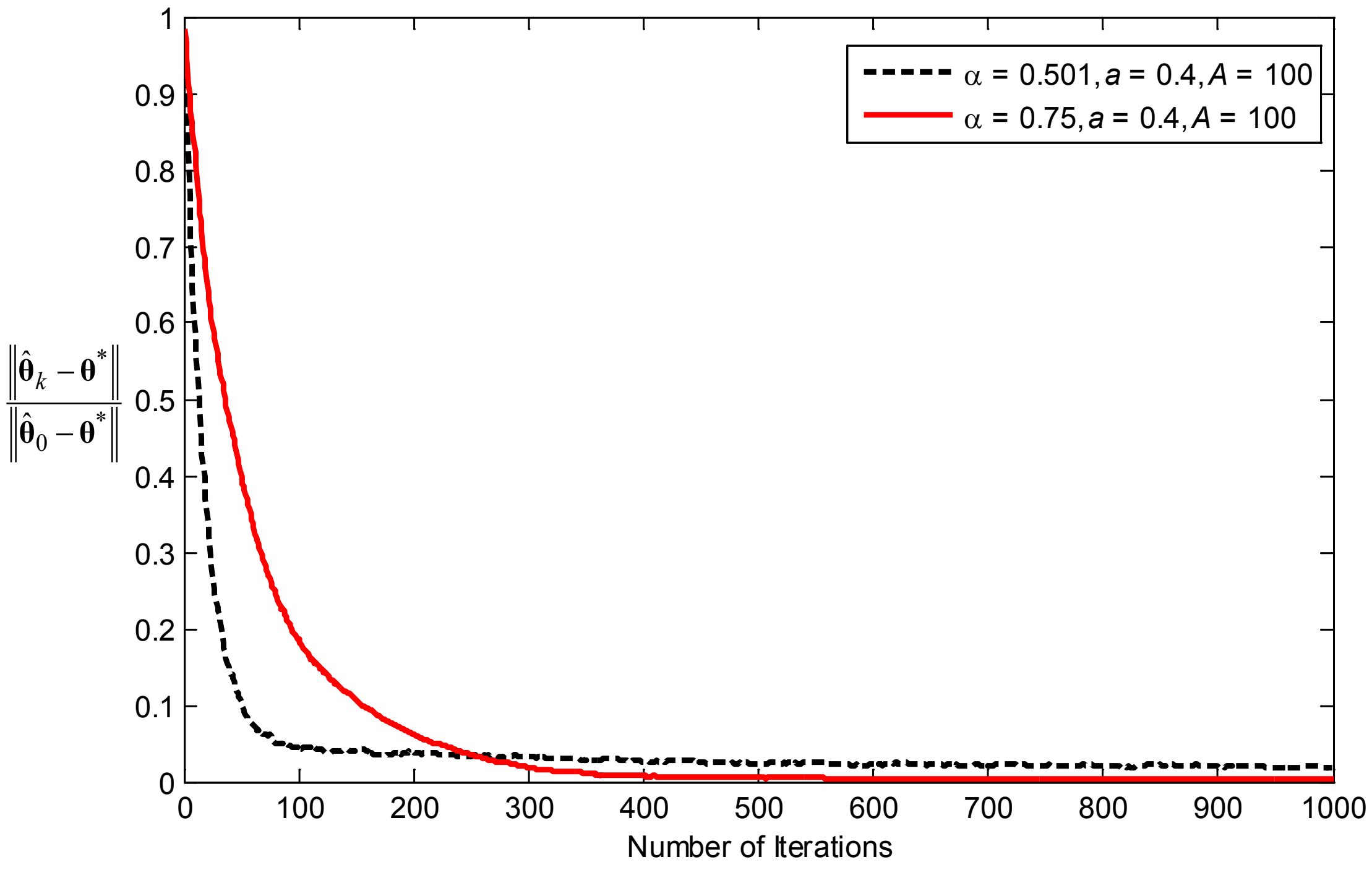

**Figure 4.3** Performance comparison with different values of α for the separable loss function. By fixing the values of $A$ and $a$, we find smaller value of α provides better performance in the early iterations; while bigger value of α provides better performance asymptotically. Each curve represents the sample mean of 20 independent replicates.

Overall, since these numerical results are consistent with the discussions in the Section 3.3, then the guidelines of the coefficients selection of the gain sequence given in Section 3.3 seem reasonable.

## 4.3 Performance of DSPSA on the Separable Loss Function

In this section, we discuss the performance of DSPSA on the first loss function (4.1) (separable loss function). The optimal solution is $\boldsymbol{\theta}^* = \mathbf{0}_p$ . We consider the high-dimensional case with $p$ = 200, and set the measurement noises $\varepsilon$ to be i.i.d. $N$(0,1). The initial guess is $10 \times \mathbf{1}_{200}$. Let the number of iterations in each replicate be 10,000, and let the number of replicates be 20. Here $a_k = a/(k+1+A)^{\alpha}$ . We pick the coefficients based on the guidelines in Section 3.3. By the guideline on $A$, we have $A = 0.1 \times 10000 = 1000$. For practical selection of $\alpha$, we choose $\alpha$ = 0.501, and for the asymptotically optimal selection of $\alpha$, we choose $\alpha$ = 1. After computing some values of $\hat{\boldsymbol{g}}_0(\hat{\boldsymbol{\theta}}_0)$, we know that the largest of the mean values of the magnitudes of the components in $\hat{\boldsymbol{g}}_0(\hat{\boldsymbol{\theta}}_0)$ is approximately 30. Suppose we want the elements of $\boldsymbol{\theta}$ move by a magnitude of 0.05 in the early iterations. Then, for the practical selection of $\alpha$, $a$ = 0.05 is according to $\left(0.05/1001^{0.501}\right) \times 30 \approx 0.05$. For the asymptotically optimal selection of $\alpha$, $a$ = 1.57 is according to $\left(1.57/1001\right) \times 30 \approx 0.05$.

Therefore, we do the numerical tests with two sets of coefficients: 1) (practical selection of $\alpha$) $\alpha = 0.501$, $A$ = 1000, $a$ = 0.05; 2) (asymptotically optimal $\alpha$) $\alpha$ = 1, $A$ = 1000, $a$ = 1.57. These two sets of coefficients provide the same gain step in the first iteration, because $a$ is picked by matching the desired change magnitude in the early iterations.

In the numerical tests below, we not only check the criterion $\left\|\hat{\boldsymbol{\theta}}_k - \boldsymbol{\theta}^*\right\| \big/ \left\|\hat{\boldsymbol{\theta}}_0 - \boldsymbol{\theta}^*\right\|$, but also check the criterion $\left\|[\hat{\boldsymbol{\theta}}_k] - \boldsymbol{\theta}^*\right\| \big/ \left\|[\hat{\boldsymbol{\theta}}_0] - \boldsymbol{\theta}^*\right\|$, where $[\cdot]$ is the round operator. (Recall that the solution provided by DSPSA in the $k$th iteration is $[\hat{\boldsymbol{\theta}}_k]$). Figures 4.4 and 4.5 show the results for the separable loss function, and under both criteria, the sequences generated by DSPSA converge to the optimal solution. We see that the results given by these two different criteria are quite similar. In Figure 4.4, we see that the first set of coefficients provides better performance in the early iterations; while the second set of coefficients provides better performance in the later iterations. In Figure 4.5, by using the round operation in the criterion, we see that for finite samples performance the first set of coefficients still leads to better performance, but the difference between the asymptotic performances on these two sets is not as large as that in Figure 4.4. The reason is that when $\hat{\boldsymbol{\theta}}_k$ is very close to the optimal solution, $[\hat{\boldsymbol{\theta}}_k]$ always equals $\boldsymbol{\theta}^*$, which makes the difference between the asymptotic performances on these two sets of coefficients smaller.

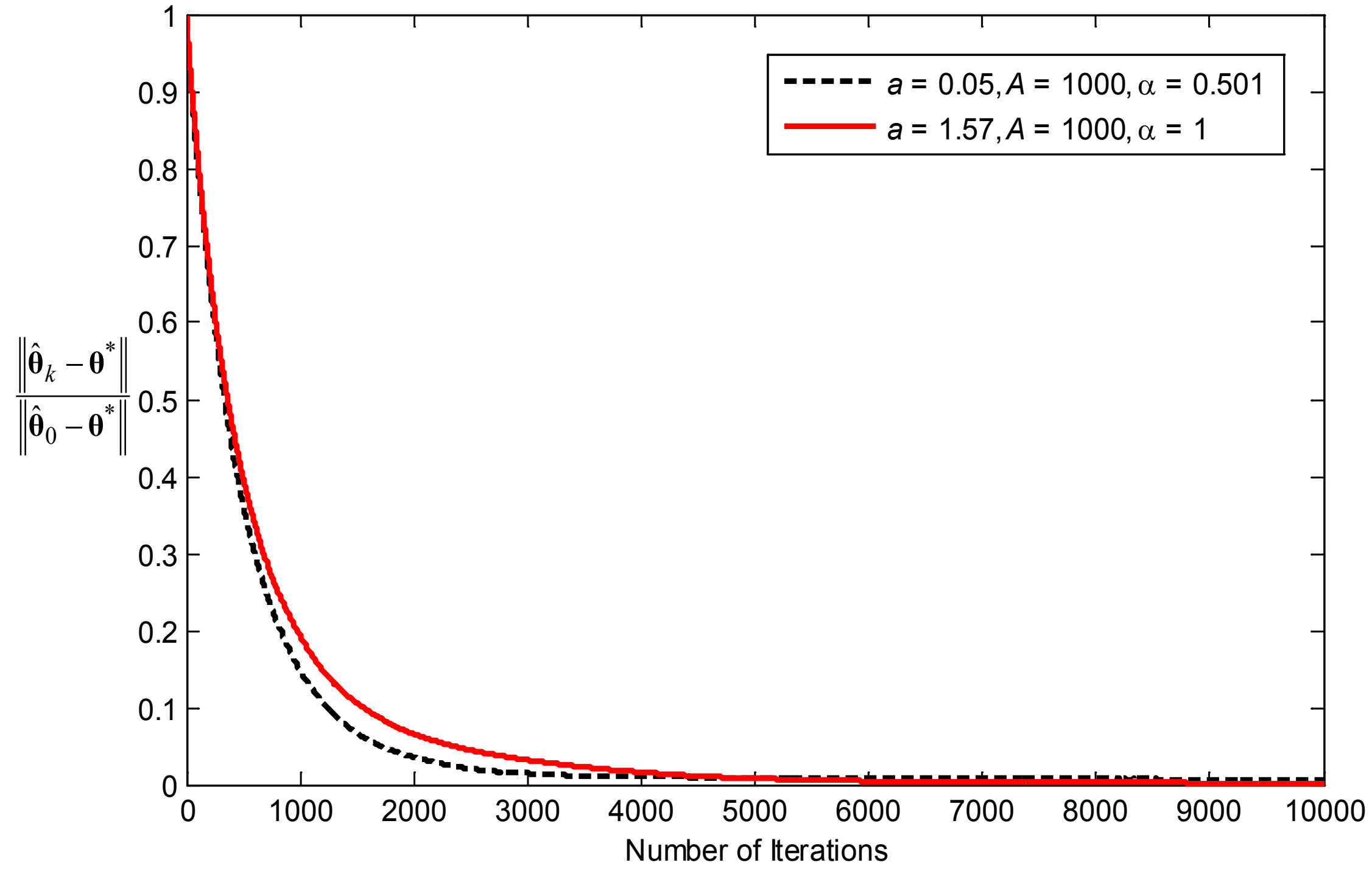

**Figure 4.4** Performance of DSPSA with two sets of coefficients for the criterion $\|\hat{\boldsymbol{\theta}}_k - \boldsymbol{\theta}^*\| / \|\hat{\boldsymbol{\theta}}_0 - \boldsymbol{\theta}^*\|$ on separable loss function. The first set of coefficients (practical selection of $\alpha$, $\alpha = 0.501$) leads to better performance in the early iterations; while in the later iterations the second set of coefficients (asymptotically optimal $\alpha$, $\alpha = 1$) provides better performance under this criterion. Each curve represents the sample mean of 20 independent replicates.

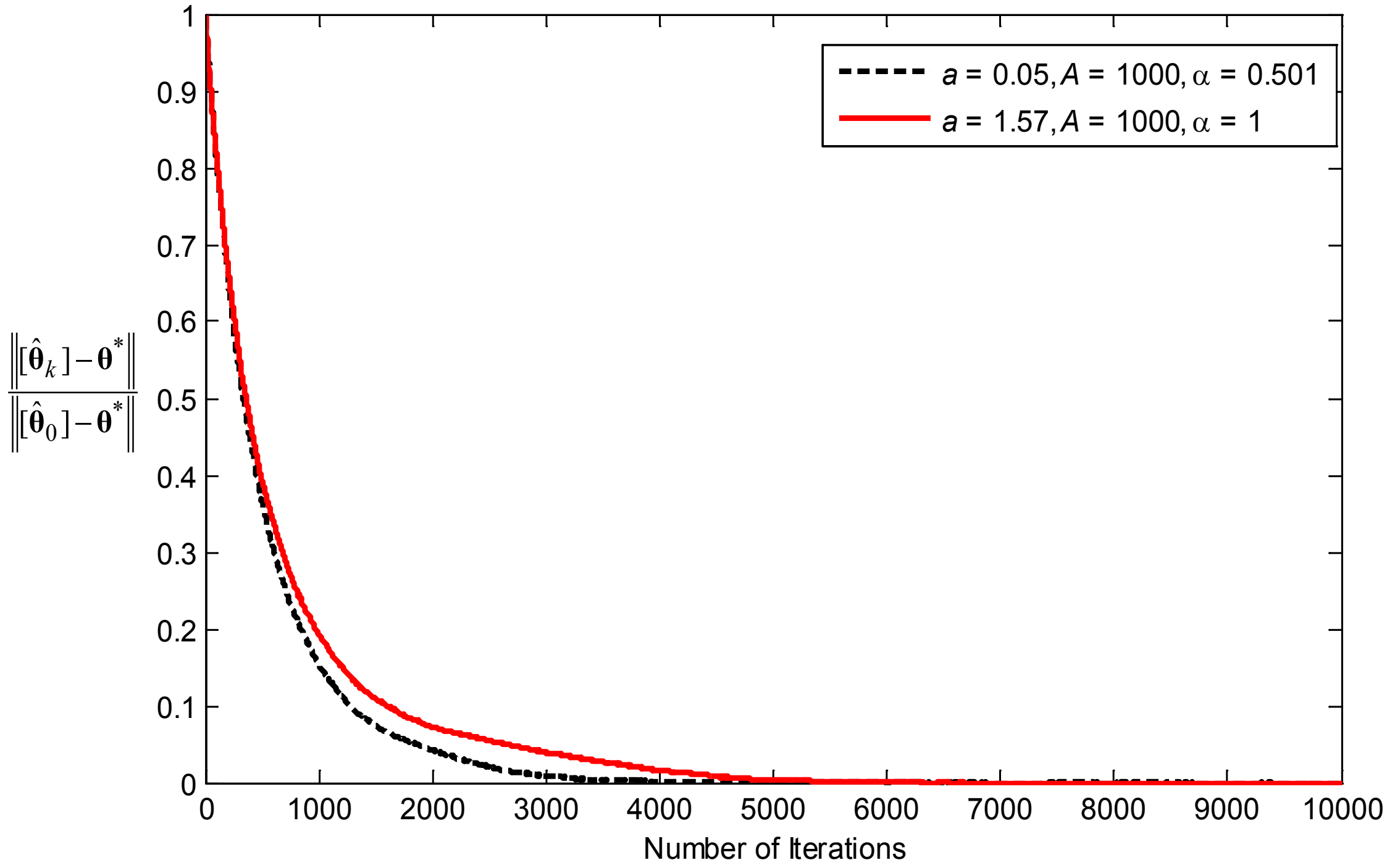

**Figure 4.5** Performance of DSPSA with two sets of coefficients for the criterion $\left\|[\hat{\boldsymbol{\theta}}_k]-\boldsymbol{\theta}^*\right\| / \left\|[\hat{\boldsymbol{\theta}}_0]-\boldsymbol{\theta}^*\right\|$ on the separable loss function. The first set of coefficients (practical selection of $\alpha$, $\alpha$ = 0.501) leads to better performance in the early iterations, and in the later iterations, the difference between the performances of DSPSA on both coefficient sets are small under this criterion. Each curve represents the sample mean of 20 independent replicates.

# 4.4 Performance of DSPSA on the Quadratic Loss Function

In this section, we discuss the performance of DSPSA on the second loss function (4.2) (quadratic loss function). The optimal solution is $\boldsymbol{\theta}^* = \mathbf{1}_p$, where $\mathbf{1}_p$ is the vector with all components being 1. We consider the high-dimensional case with $p$ = 200, and set the measurement noises $\varepsilon$ to be i.i.d. $N(0,1)$. Similar to the settings for the first loss function, we set the initial guess to be $10 \times \mathbf{1}_{200}$, the number of replicates to be 20 and the number of iterations in each replicate to be 10,000. We pick the coefficients based on the guidelines in Section 3.3. Similar to the first loss function, we set $A$ = 1000. For practical selection of $\alpha$, we choose $\alpha = 0.501$, and for the asymptotically optimal selection of $\alpha$, we choose $\alpha = 1$. After computing some values of $\hat{\boldsymbol{g}}_0(\hat{\boldsymbol{\theta}}_0)$, we know that the largest of the mean values of the magnitudes of the components in $\hat{\boldsymbol{g}}_0(\hat{\boldsymbol{\theta}}_0)$ is approximately 150. Suppose we want the elements of $\boldsymbol{\theta}$ to move by a magnitude of 0.05 in the early iterations. Then, for the practical selection of $\alpha$, $a$ = 0.01 is according to $\left(0.01/1001^{0.501}\right) \times 150 \approx 0.05$. For the asymptotically optimal selection of $\alpha$, $a$ = 0.314 is according to $\left(0.314/1001\right) \times 150 \approx 0.05$.

Therefore, we do the numerical tests with two sets of coefficients: 1) (practical selection of $\alpha$) $\alpha = 0.501$, $A$ = 1000, $a$ = 0.01; 2) (asymptotically optimal $\alpha$) $\alpha = 1$, $A$ = 1000, $a$ = 0.314. These two sets of coefficients provide the same gain step in the first iteration, because $a$ is picked by matching the desired change magnitude in the early iterations.

In the numerical tests below, we not only check the criterion $\left\|\hat{\boldsymbol{\theta}}_k - \boldsymbol{\theta}^*\right\| \big/ \left\|\hat{\boldsymbol{\theta}}_0 - \boldsymbol{\theta}^*\right\|$, but also check the criterion $\left\|[\hat{\boldsymbol{\theta}}_k] - \boldsymbol{\theta}^*\right\| \big/ \left\|[\hat{\boldsymbol{\theta}}_0] - \boldsymbol{\theta}^*\right\|$. Figures 4.6 and 4.7 show the results for the quadratic loss function. The results here are similar to the results of the separable loss function. We see that DSPSA can provide convergent sequences for the quadratic loss function under both criteria. We also see that the first set of coefficients (practical selection of $\alpha$) leads to better performance in the early iterations under both criteria. In the later iterations, the second set of coefficients (asymptotically optimal $\alpha$) provides better performance under the criterion $\left\|\hat{\boldsymbol{\theta}}_k - \boldsymbol{\theta}^*\right\| \big/ \left\|\hat{\boldsymbol{\theta}}_0 - \boldsymbol{\theta}^*\right\|$ and the difference in the performance of DSPSA on both coefficients sets are smaller under the criterion $\left\|[\hat{\boldsymbol{\theta}}_k] - \boldsymbol{\theta}^*\right\| \big/ \left\|[\hat{\boldsymbol{\theta}}_0] - \boldsymbol{\theta}^*\right\|$.

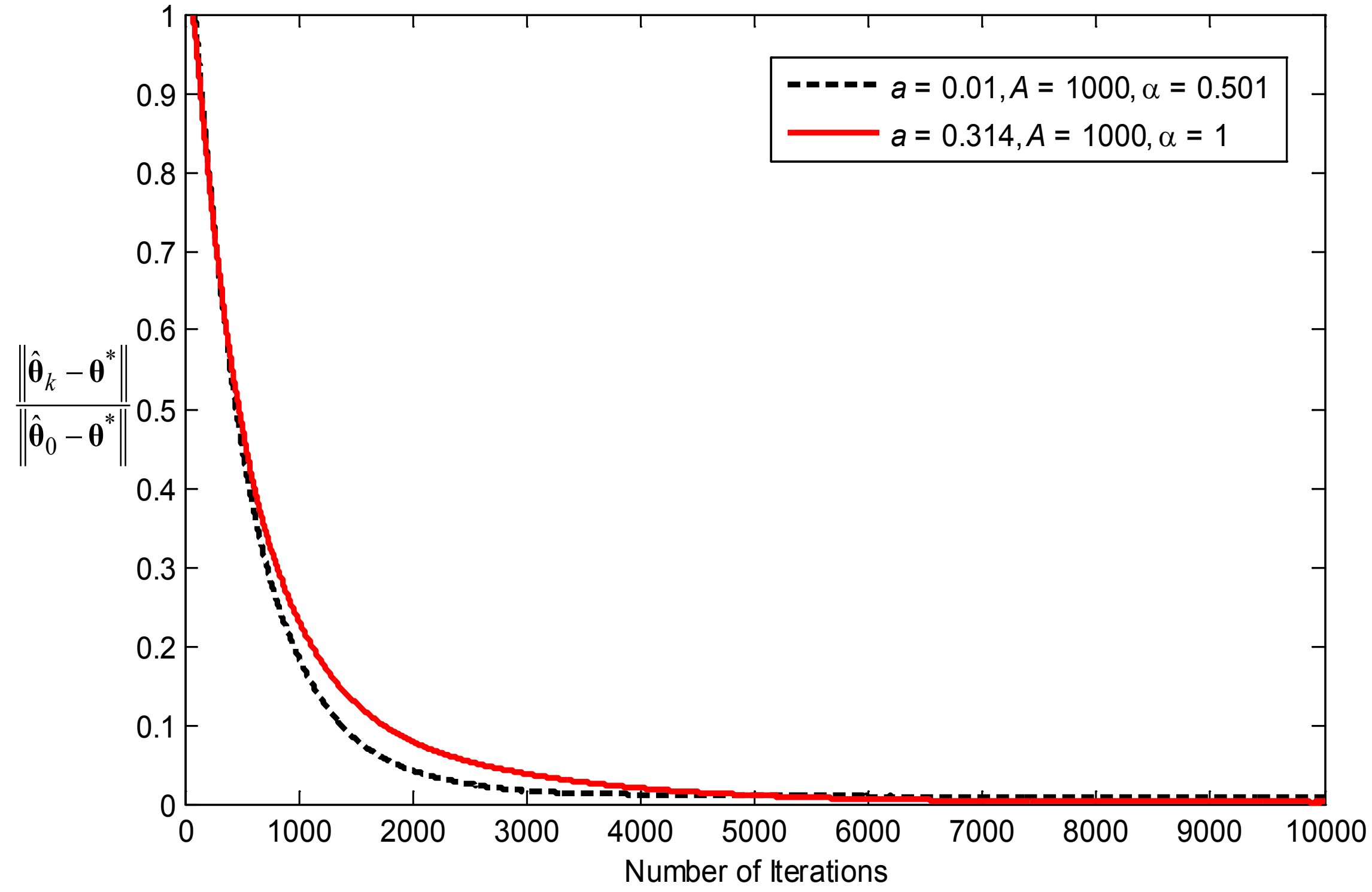

**Figure 4.6** Performance of DSPSA with two sets of coefficients for the criterion $\|\hat{\boldsymbol{\theta}}_k - \boldsymbol{\theta}^*\| / \|\hat{\boldsymbol{\theta}}_0 - \boldsymbol{\theta}^*\|$ on quadratic loss function. The first set of coefficients (practical selection of $\alpha$, $\alpha = 0.501$) leads to better performance in the early iterations; while in the later iterations the second set of coefficients (asymptotically optimal $\alpha$, $\alpha = 1$) provides better performance under this criterion. Each curve represents the sample mean of 20 independent replicates.

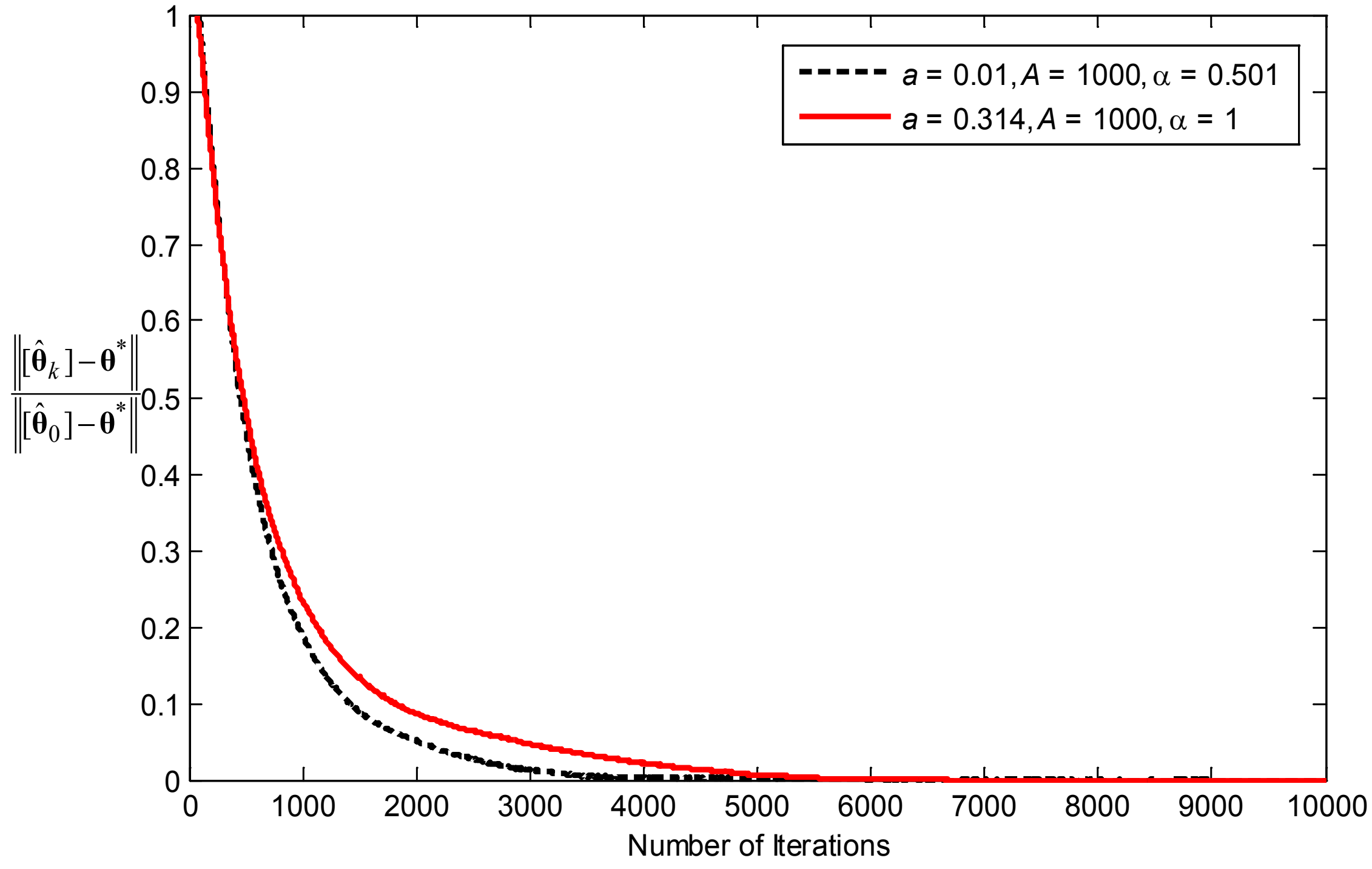

**Figure 4.7** Performance of DSPSA with two sets of coefficients for the criterion $\left\|[\hat{\boldsymbol{\theta}}_k]-\boldsymbol{\theta}^*\right\| / \left\|[\hat{\boldsymbol{\theta}}_0]-\boldsymbol{\theta}^*\right\|$ on quadratic loss function. The first set of coefficients (practical selection of $\alpha$, $\alpha$ = 0.501) leads to better performance in the early iterations, and in the later iterations, the difference between the performances of DSPSA on both coefficient sets are small under this criterion. Each curve represents the sample mean of 20 independent replicates.

# 4.5 Performance of DSPSA on the Skewed Quartic Loss Function

We discuss the performance of DSPSA on the third loss function, the skewed quartic loss function. The optimal solution is $\boldsymbol{\theta}^* = \mathbf{0}_p$. We consider the high-dimensional case with $p$ = 200, and set the measurement noises $\varepsilon$ to be i.i.d. $N(0,1)$. Similar to the settings for the last two loss functions, we set the initial guess to be $10\times\mathbf{1}_{200}$, the number of replicates to be 20, and the number of iterations in each replicate to be 10,000. We pick the coefficients based on the guidelines in Section 3.3. Similar to the last two loss functions, we set $A$ = 1000. For practical selection of $\alpha$, we choose $\alpha = 0.501$, and for the asymptotically optimal selection of $\alpha$, we choose $\alpha = 1$. After computing some values of $\hat{\boldsymbol{g}}_0(\hat{\boldsymbol{\theta}}_0)$, we see that the largest of the mean values of the magnitudes of the components in $\hat{\boldsymbol{g}}_0(\hat{\boldsymbol{\theta}}_0)$ is approximately 150. Suppose we want the elements of $\boldsymbol{\theta}$ move by a magnitude of 0.05 in the early iterations. Then, for the practical selection, $a$ = 0.01 is according to $\left(0.01/1001^{0.501}\right)\times 150 \approx 0.05$. For the asymptotically optimal selection, $a$ = 0.314 is according to $\left(0.314/1001\right)\times 150 \approx 0.05$.

Therefore, we do the numerical tests for two sets of coefficients: 1) (practical selection of $\alpha$) $\alpha = 0.501$, $A$ = 1000, $a$ = 0.01; 2) (asymptotically optimal $\alpha$) $\alpha$ = 1, $A$ = 1000, $a$ = 0.314. These two sets of coefficients provide the same gain step in the first iteration, because $a$ is picked by matching the desired change magnitude in the early iterations.

In the numerical tests below, we not only check the criterion $\left\|\hat{\boldsymbol{\theta}}_k - \boldsymbol{\theta}^*\right\| \big/ \left\|\hat{\boldsymbol{\theta}}_0 - \boldsymbol{\theta}^*\right\|$, but also check the criterion $\left\|[\hat{\boldsymbol{\theta}}_k] - \boldsymbol{\theta}^*\right\| \big/ \left\|[\hat{\boldsymbol{\theta}}_0] - \boldsymbol{\theta}^*\right\|$. Figures 4.8 and 4.9 show the performance of DSPSA for the skewed quartic loss function. Due to the special structure of the skewed quartic loss function, the sequence $\left\{\hat{\boldsymbol{\theta}}_k\right\}$ does not converge as fast to $\boldsymbol{\theta}^*$ as the first two loss functions. However, we still see the efficient convergent performance for the skewed quartic loss function by using DSPSA. We also check the performance in terms of the errors in loss function value $L([\hat{\boldsymbol{\theta}}_k])$ in Figure 4.10. We do not check the value of $L(\hat{\boldsymbol{\theta}}_k)$, because in the discrete problem, loss functions are not defined on points that are non-multivariate integer. We see that Figure 4.10 shows good convergence of DSPSA for the skewed quartic loss function in term of the errors in loss function values. Comparing the results in Figures 4.10 and 4.9, we find that the relatively slow convergence of the sequence in terms of the errors in points is due to the special structure of the loss function, where there is a large flat area near the optimal solution. From Figures 4.8 and 4.9, we find that, relative to the errors reported in Figures 4.4, 4.5, 4.6 and 4.7 based on the previous loss functions, $\hat{\boldsymbol{\theta}}_k$ is still not very close to the optimal solution $\boldsymbol{\theta}^*$ in the last iteration. Thus, we only see the finite sample performance of DSPSA for the skewed quartic loss function, and the difference between the two curves are similar in Figure 4.8 and Figure 4.9. We see that the first set of coefficients (practical selection of $\alpha$) leads to better finite sample performance under both criteria.

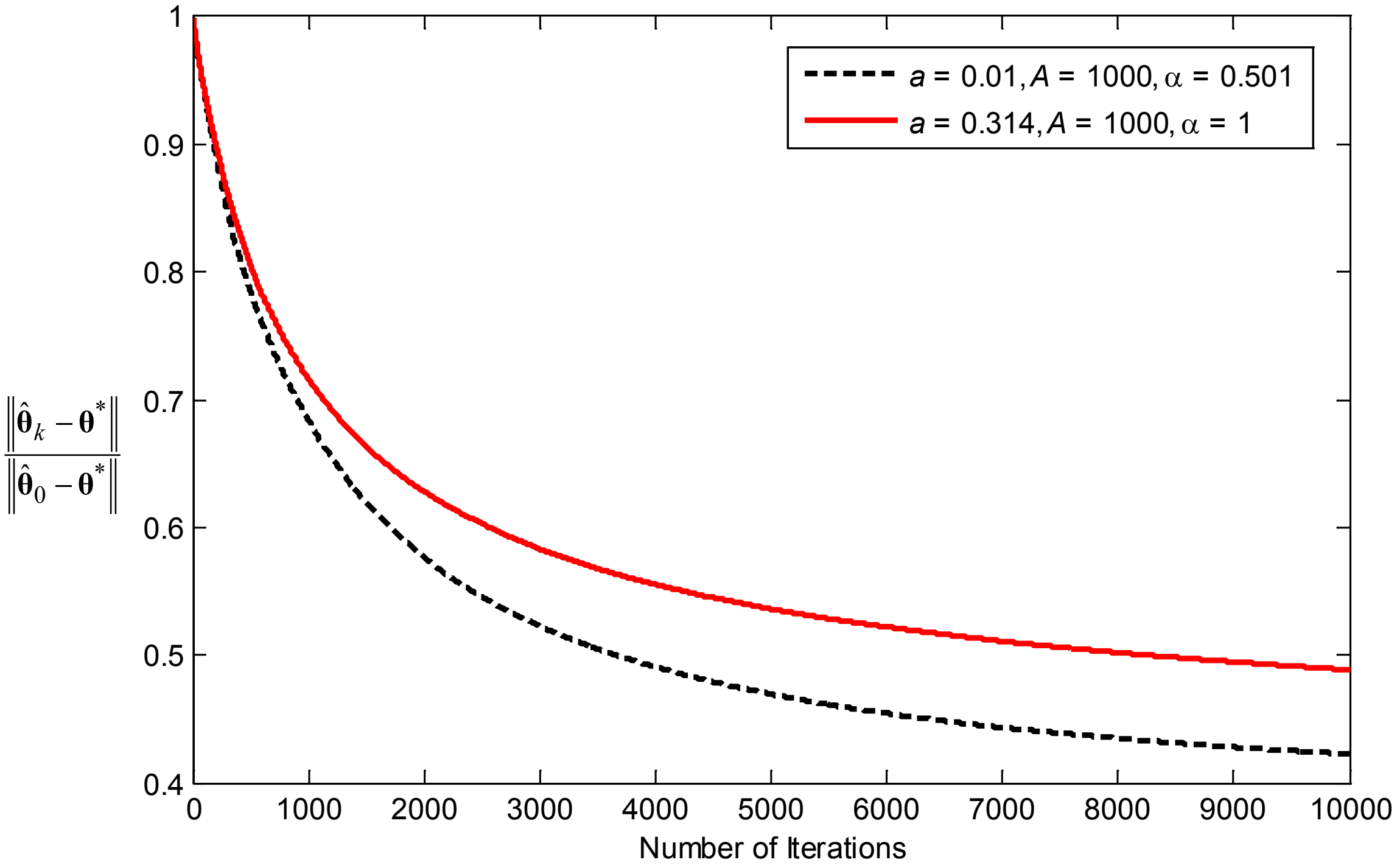

**Figure 4.8** Performance of DSPSA with two sets of coefficients for the criterion $\left\|\hat{\boldsymbol{\theta}}_k - \boldsymbol{\theta}^*\right\| / \left\|\hat{\boldsymbol{\theta}}_0 - \boldsymbol{\theta}^*\right\|$ on the skewed quartic loss function. The first set of coefficients (practical selection of $\alpha$, $\alpha$ = 0.501) leads to better finite sample performance under this criterion. Each curve represents the sample mean of 20 independent replicates.

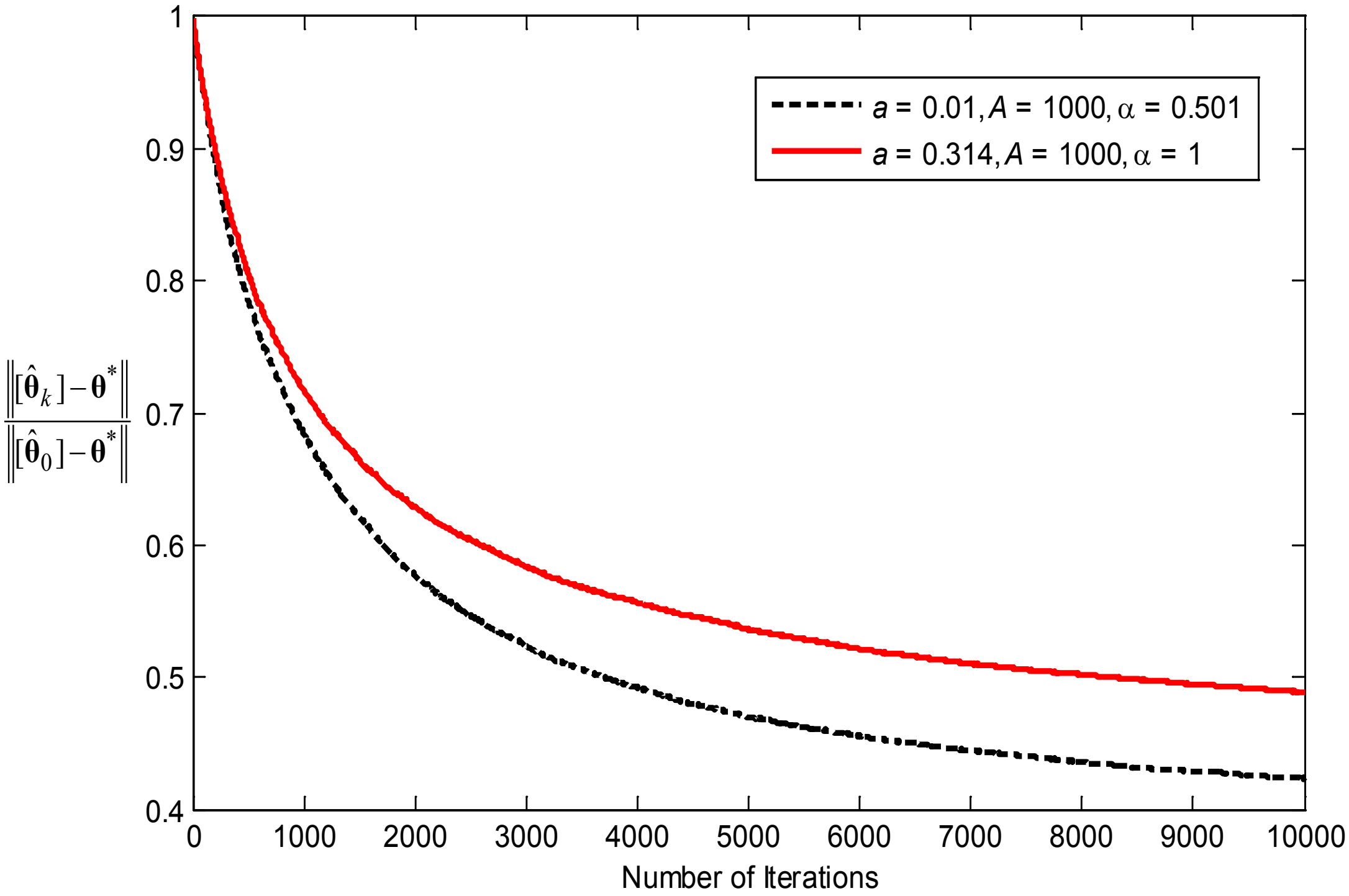

**Figure 4.9** Performance of DSPSA with two sets of coefficients for the criterion $\left\|[\hat{\boldsymbol{\theta}}_k]-\boldsymbol{\theta}^*\right\| / \left\|[\hat{\boldsymbol{\theta}}_0]-\boldsymbol{\theta}^*\right\|$ on the skewed quartic loss function. The first set of coefficients (practical selection of $\alpha$, $\alpha = 0.501$) leads to better finite sample performance under this criterion. Each curve represents the sample mean of 20 independent replicates.

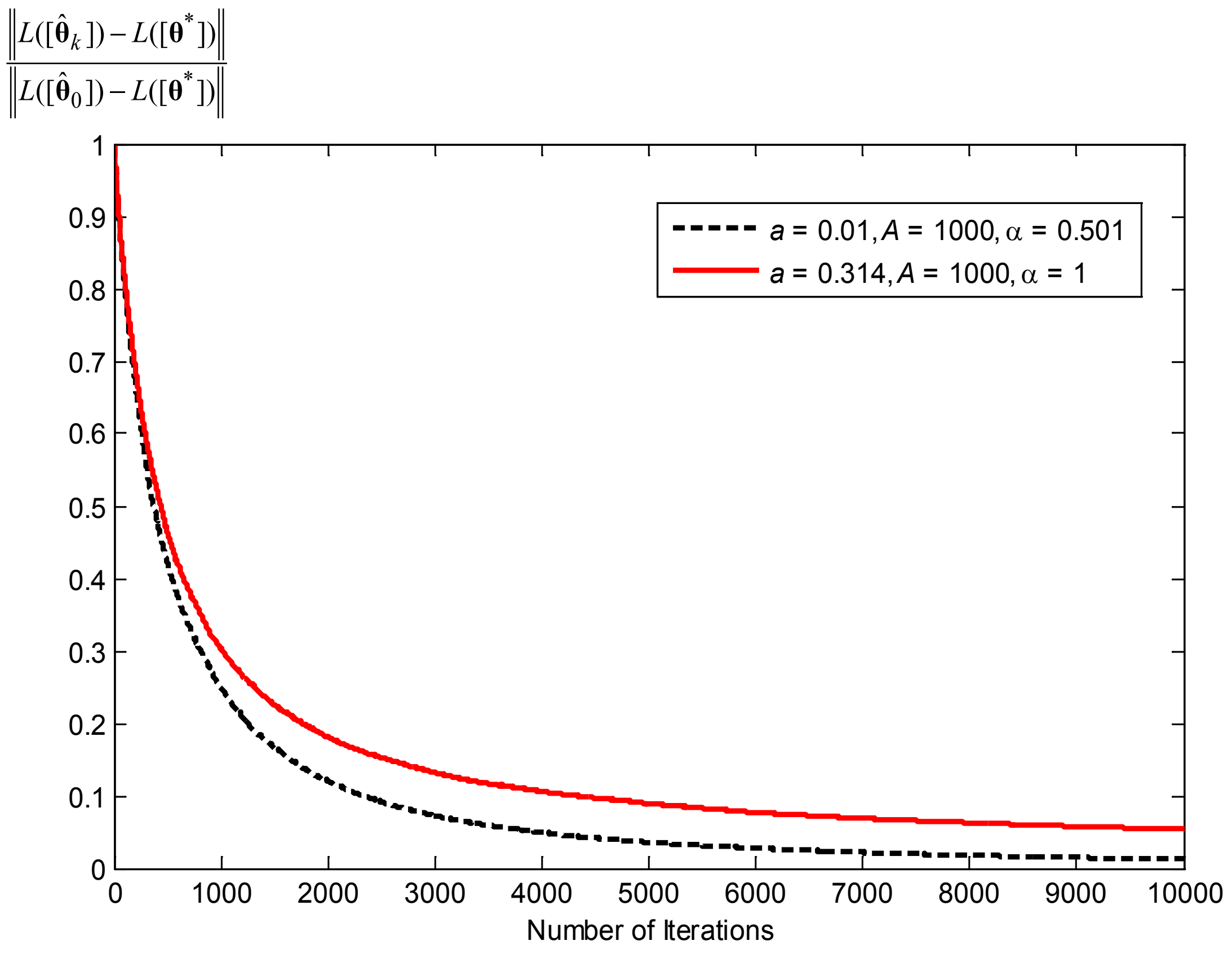

**Figure 4.10** Performance of DSPSA with two sets of coefficients for the criterion $\left(L([\hat{\boldsymbol{\theta}}_k])-L([\boldsymbol{\theta}^*])\right)\big/\left(L([\hat{\boldsymbol{\theta}}_0])-L([\boldsymbol{\theta}^*])\right)$ on the skewed quartic loss function. The first set of coefficients (practical selection of $\alpha$, $\alpha$ = 0.501) leads to better finite sample performance under this criterion. Each curve represents the sample mean of 20 independent replicates.

Overall, the first set of coefficients (practical selection of $\alpha$) can achieve better finite sample performance, and in the later iterations the second set of coefficients (asymptotically optimal $\alpha$) might do better in terms of $\left\|\hat{\boldsymbol{\theta}}_k - \boldsymbol{\theta}^*\right\| / \left\|\hat{\boldsymbol{\theta}}_0 - \boldsymbol{\theta}^*\right\|$. But after considering the round operation, we find that the better performance of the second set of coefficients in the later iterations becomes weaker. For the function with special structure, such as the skewed quartic loss function that contains a large flat area near the optimal solution, DSPSA provides good convergent performance especially in terms of the errors in loss function values.

# 4.6 Example of Non-Bernoulli Distribution for the Simultaneous Perturbation Direction

For the simultaneous perturbation direction $\boldsymbol{\Delta}$, we have mainly discussed the case when the components of $\boldsymbol{\Delta}$ are independent Bernoulli random variables taking the values $\pm 1$ with probability $1/2$. Now let us discuss a non-Bernoulli distribution, where the components of $\boldsymbol{\Delta}$ can be discrete uniformly distributed over the set $\{\pm 1, \pm 3\}$. Note that non-Bernoulli distributions have also been considered for continuous (non-discrete) problems as well (e.g. Spall, 2003, Chapter 7 and Cao, 2011).

As an example, let us explain the idea of the non-Bernoulli case on the domain of $\mathbb{Z}^2$. For example, suppose in the *k*th iteration, we have that $\hat{\boldsymbol{\theta}}_k = [0.1, 0.1]^T$ and the simulation result for $\boldsymbol{\Delta}_k$ is $\boldsymbol{\Delta}_k = [3, 1]^T$, then $\hat{\boldsymbol{\theta}}_k^+ = [2, 1]^T$ and $\hat{\boldsymbol{\theta}}_k^- = [-1, 0]^T$. In Figure 4.11, the lines radiating from $\boldsymbol{\pi}(\hat{\boldsymbol{\theta}}_k)$ represent all possible values of the perturbation vector $\boldsymbol{\Delta}_k$, the triangle point indicates the position of $\hat{\boldsymbol{\theta}}_k$, and by choosing different $\boldsymbol{\Delta}_k$, the possible values for $\hat{\boldsymbol{\theta}}_k^+$ and $\hat{\boldsymbol{\theta}}_k^-$ belong to the set of the round points. We find that this non-Bernoulli distribution of $\boldsymbol{\Delta}_k$ makes DSPSA to explore the points in a larger hypercube (not restricted within one unit hypercube). This is just one example to describe the non-Bernoulli case.

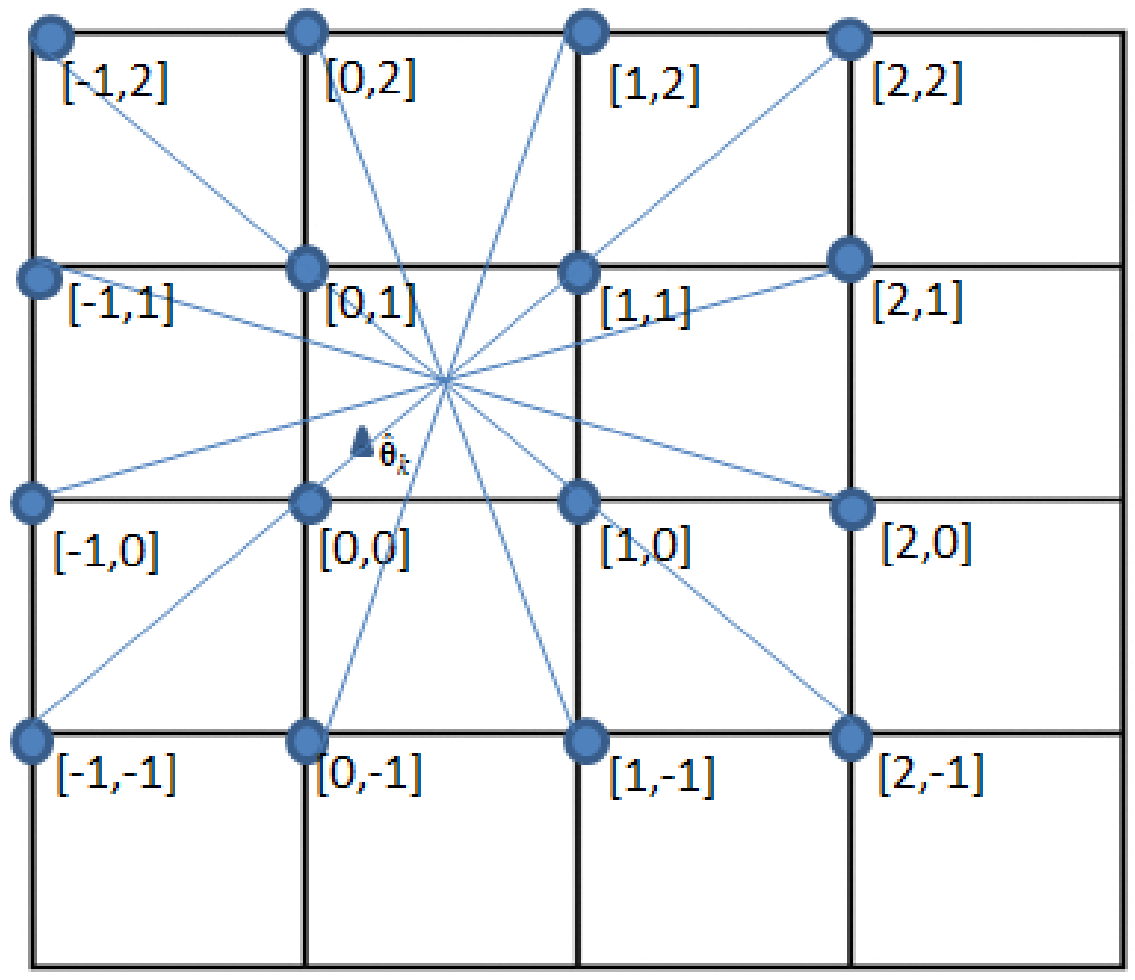

**Figure 4.11** Brief idea of the extension on $\boldsymbol{\Delta}_k$. The triangle point is $\hat{\boldsymbol{\theta}}_k$, the lines radiating from $\boldsymbol{\pi}(\hat{\boldsymbol{\theta}}_k)$ represent all possible values of the perturbation vector $\boldsymbol{\Delta}_k$, and the round points are the possible values for $\hat{\boldsymbol{\theta}}_k^+$ and $\hat{\boldsymbol{\theta}}_k^-$ based on different values of $\boldsymbol{\Delta}_k$.

We do the numerical tests by using both the non-Bernoulli distribution (the uniform $\{\pm 1, \pm 3\}$) and the Bernoulli ±1 distribution for the simultaneous perturbation vector $\mathbf{\Delta}_k$ on the 200-dimensional separable loss function. The initial guess is $10\times\mathbf{1}_p$, the number of replicates is 20, and the number of iterations in each replicate is 10,000. We see that these settings are the same as that in Section 4.3, then for the Bernoulli case we pick the same set of coefficients (practical selection), which is $\alpha = 0.501$, $A$ = 1000, $a$ = 0.05. For non-Bernoulli case, we pick the coefficients also based on the guidelines in Section 3.3. Similar to the Bernoulli case, for the non-Bernoulli case, we set $A$ = 1000. For practical selection of $\alpha$, we choose $\alpha = 0.501$. After computing some values of $\hat{\boldsymbol{g}}_0(\hat{\boldsymbol{\theta}}_0)$, we see that the largest of the mean values of the magnitudes of the components in $\hat{\boldsymbol{g}}_0(\hat{\boldsymbol{\theta}}_0)$ is approximately 45. Suppose we want the elements of $\boldsymbol{\theta}$ move by a magnitude of 0.05 in the early iterations. Then, for the practical selection, $a$ = 0.0354 is according to $\left(0.0354/1001^{0.501}\right)\times 45 \approx 0.05$.

Therefore, we do the numerical test by using the coefficients $\alpha = 0.501$, $A$ = 1000, $a$ = 0.05 for Bernoulli ±1 and $\alpha = 0.501$, $A$ = 1000, $a$ = 0.0354 for the non-Bernoulli distribution (uniform $\{\pm 1, \pm 3\}$). In Figures 4.12 and 4.13, we find by using the non-Bernoulli distribution, we also can get a convergent sequence. However, we see that the Bernoulli ±1 can provide better performance for this high-dimensional separable loss function in terms of both criteria $\left\|\hat{\boldsymbol{\theta}}_k - \boldsymbol{\theta}^*\right\| \big/ \left\|\hat{\boldsymbol{\theta}}_0 - \boldsymbol{\theta}^*\right\|$ and $\left\|[\hat{\boldsymbol{\theta}}_k] - \boldsymbol{\theta}^*\right\| \big/ \left\|[\hat{\boldsymbol{\theta}}_0] - \boldsymbol{\theta}^*\right\|$. Of course, these results do not necessarily indicate similar comparative performance using different loss functions. We further discuss the choice of distribution for the simultaneous perturbation vector $\mathbf{\Delta}$ as the future research direction in Chapter 8.

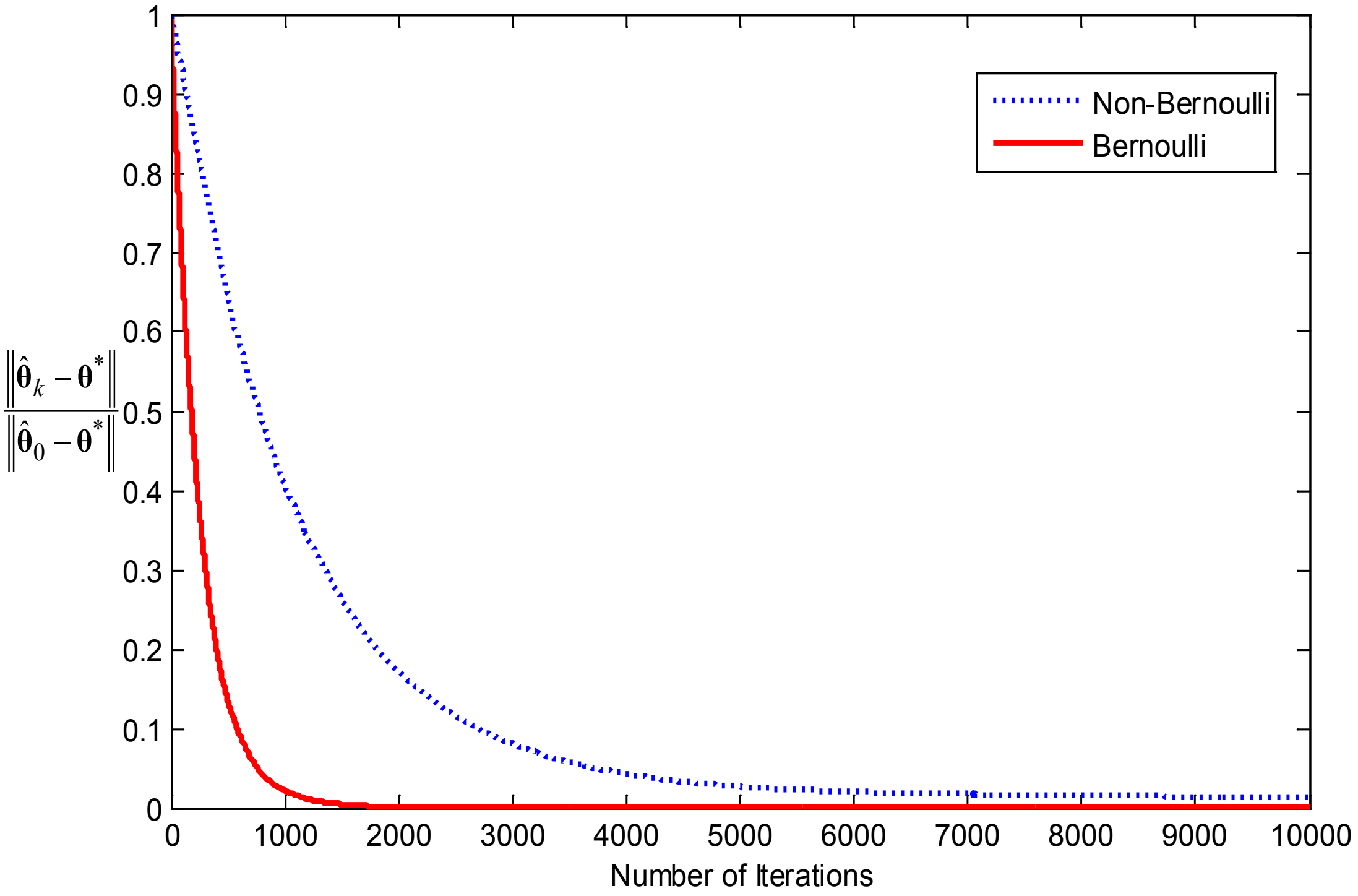

**Figure 4.12** Performance of DSPSA by using different distributions for the perturbation vector $\mathbf{\Delta}_k$ on the high-dimensional separable loss function. The non-Bernoulli distribution (uniform $\{\pm 1, \pm 3\}$) provides convergent performance. The Bernoulli $\pm 1$ provides better performance for the separable loss function in terms of the criterion $\left\|\hat{\boldsymbol{\theta}}_k - \boldsymbol{\theta}^*\right\| / \left\|\hat{\boldsymbol{\theta}}_0 - \boldsymbol{\theta}^*\right\|$. Each curve represents the sample mean of 20 independent replicates.

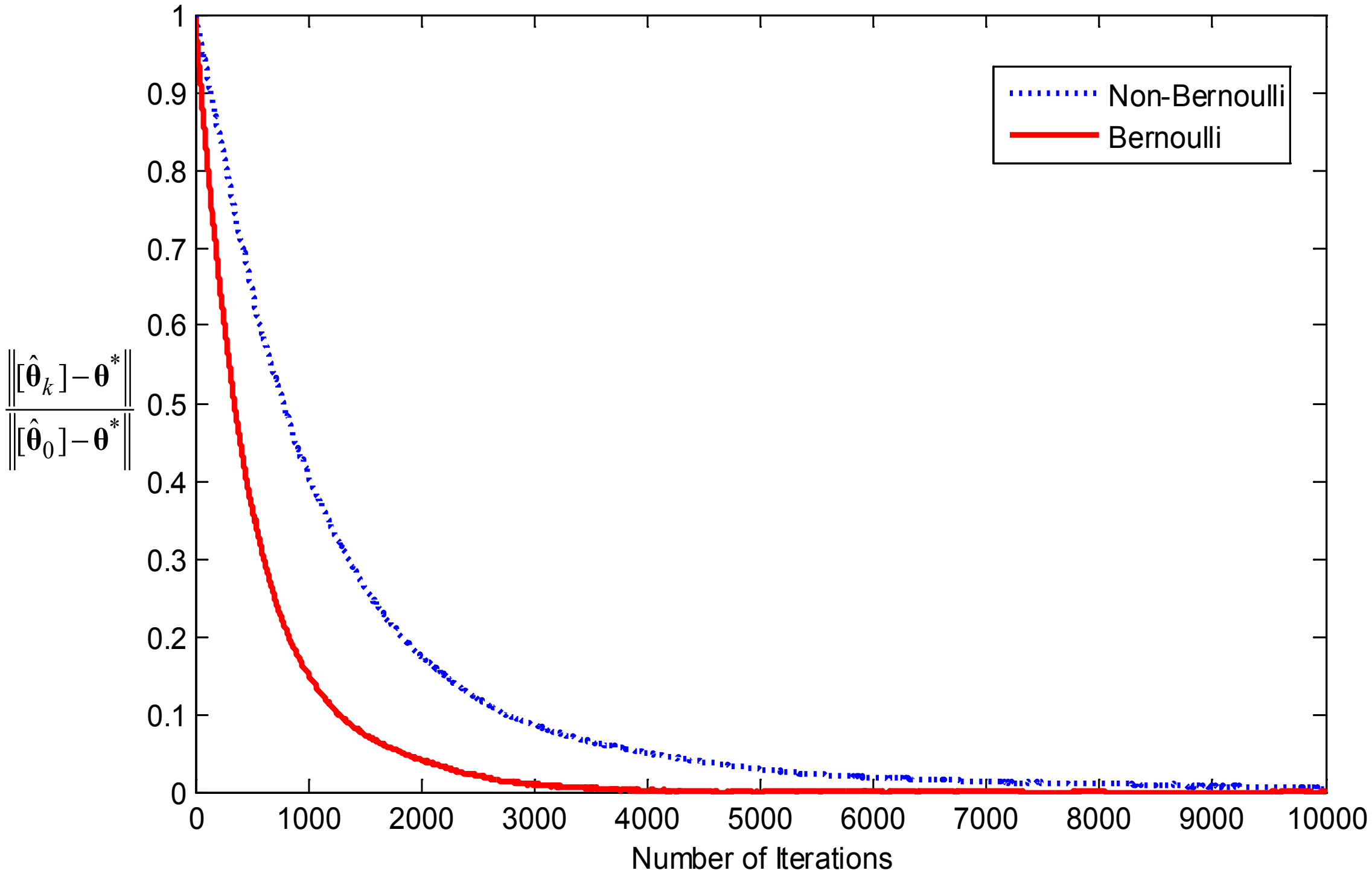

**Figure 4.13** Performance of DSPSA by using different distributions for the perturbation vector $\boldsymbol{\Delta}_k$ on the high-dimensional separable loss function. The non-Bernoulli distribution (uniform $\{\pm 1, \pm 3\}$) provides convergent performance. The Bernoulli $\pm 1$ provides better performance for the separable loss function in terms of the criterion $\|[\hat{\boldsymbol{\theta}}_k]-\boldsymbol{\theta}^*\| / \|[\hat{\boldsymbol{\theta}}_0]-\boldsymbol{\theta}^*\|$. Each curve represents the sample mean of 20 independent replicates.

# Chapter 5

# Formal Comparison of Convergence Rates of DSPSA and Two Random Search Algorithms

In this chapter, we first present some basic results of two random search type algorithms (stochastic ruler algorithm and stochastic comparison algorithm). Then, we calculate the rate of convergence in the big-$O$ sense for these two random search algorithms. At last, we do the comparison of the rate of convergence of DSPSA and two random search type algorithms theoretically.

# 5.1 Introduction

As we discussed in the Section 1.2.2, there are three main classes of algorithms designed for the discrete stochastic optimization problem: random search class, statistical class, and stochastic approximation class. Based on the literature review in Chapter 1, we see that many algorithms belong to the class of random search type algorithms. Random search type algorithms provide convergent sequences, and many theoretical results on them are available. However, the statistical class of algorithms does not provide a sequence that converges to the optimal solution. Meanwhile, the stochastic approximation class does not contain many existing algorithms, and the theoretical results of the stochastic approximation class are not well developed.

Stochastic ruler (SR) algorithm and stochastic comparison (SC) algorithm are two basic representatives of random search algorithms. Andradottir's (1999) idea, discussed in Section 1.2.2.2, is a modification of basic random search type algorithms, according to which all of the old information is stored and used to pick the current optimal point. For other random search algorithms, such as COMPASS (discussed in Section 1.2.2.2) and locally convergent random-search algorithm, they are also the extensions of basic random search type algorithms, but they are much more complicated than the basic ones. These extensions improve the performance of basic random search type algorithms for some specific kinds of problems, but at the same time the simplicities of the original algorithms are sacrificed. As we know, DSPSA is a simple algorithm, and we can also modify it to improve the performance, such as storing old information. Here we choose SR algorithm and SC algorithm as two representatives of the basic random search type algorithms, and compare them with DSPSA. In addition, as far as we know there is no generally accepted criterion for the comparisons of different discrete stochastic optimization algorithms, because different algorithms use different rate of convergence measurements. However, for each convergent algorithm of discrete stochastic optimization problem, we know $P([\hat{\boldsymbol{\theta}}_k] \neq \boldsymbol{\theta}^*) \to 0$, where $[\cdot]$ is the round operator, so we consider the rate at which $P([\hat{\boldsymbol{\theta}}_k] \neq \boldsymbol{\theta}^*)$ goes to 0 as the basis of comparisons for the three algorithms. For the SR algorithm and the SC algorithm, $\{\hat{\boldsymbol{\theta}}_k\}$ is composed of multivariate integer points, so $[\hat{\boldsymbol{\theta}}_k]$ = $\hat{\boldsymbol{\theta}}_k$. In the following analysis on the rate of convergence of the SR algorithm and SC algorithm, we consider the criterion $P(\hat{\boldsymbol{\theta}}_k \neq \boldsymbol{\theta}^*)$ directly.

# 5.2 Stochastic Ruler Algorithm

Yan and Mukai (1992) introduce the SR algorithm. The idea of this algorithm is to change the minimization problem into a particular maximization problem. A stochastic ruler $U_{u,v}$ is defined for this algorithm, and $U_{u,v}$ is a random variable uniformly distributed over the interval $[u,v]$. Yan and Mukai (1992) change the original problem of

$$\min_{\boldsymbol{\theta}\in\Theta} E(y(\boldsymbol{\theta}))$$

to the maximization problem

$$\max_{\boldsymbol{\theta}\in\Theta} P\left(y(\boldsymbol{\theta}) \leq U_{u,v}\right).$$

By Theorem 3.1 in Yan and Mukai (1992), the authors show that there exist real numbers $\bar{u}$ and $\bar{v}$, such that $\bar{u} < \bar{v}$ and for any $u < \bar{u}$ and $v > \bar{v}$, the following results hold: 1) if $E\left(y(\boldsymbol{\theta}')\right)$ < $E\left(y(\boldsymbol{\theta}'')\right)$, then $P\left(y(\boldsymbol{\theta}') \le U_{u,v}\right)$ > $P\left(y(\boldsymbol{\theta}'') \le U_{u,v}\right)$ for any $\boldsymbol{\theta}', \boldsymbol{\theta}'' \in \Theta$; 2) 0 < $P\left(y(\boldsymbol{\theta}) \le U_{u,v}\right)$ < 1 for all $\boldsymbol{\theta} \in \Theta$. The first result of Theorem 3.1 in Yan and Mukai (1992) indicates that the optimal solution $\boldsymbol{\theta}^*$ corresponding to the smallest value of $E(y(\boldsymbol{\theta}^*))$ achieves the biggest value in probability $P\left(y(\boldsymbol{\theta}) \le U_{u,v}\right)$. Here $u$ and $v$ are picked by the user, and the user can always pick very small value for $u$ and very large value for $v$ to make the condition of Theorem 3.1 in Yan and Mukai (1992) to be satisfied. In Yan and Mukai (1992), the authors indicate that if the noisy measurements of loss function are uniformly bounded, we can choose $u$ as the lower bound of $y(\boldsymbol{\theta})$ and choose $v$ as the upper bound of $y(\boldsymbol{\theta})$. If the noisy measurements of loss function are not bounded, we can still pick finite values for $u$ and $v$ to make 0 < $P\left(y(\boldsymbol{\theta}) \le U_{u,v}\right)$ < 1 for all $\boldsymbol{\theta} \in \Theta$ (the second result in Theorem 3.1 of Yan and Mukai, 1992). However, the authors do not provide a clear guideline on the choice of $u$ and $v$ for general cases.

## 5.2.1 Algorithm Description

Before describing the algorithm, let us denote $N(\boldsymbol{\theta}) \subseteq \Theta \setminus \boldsymbol{\theta}$ as a set of neighbor points of $\boldsymbol{\theta}$, excluding $\boldsymbol{\theta}$ itself. For example, for the domain of $\mathbb{Z}^p$, the local neighbor of $\boldsymbol{\theta} = \{t_1, ..., t_p\}$ may be $N(\boldsymbol{\theta}) = \left\{\boldsymbol{\theta}' = \{t_1', ..., t_p'\} \in \mathbb{Z}^p \middle| |t_i - t_i'| \le 1 \text{ for all } i\right\} \setminus \boldsymbol{\theta}$ and we call it local square-ring neighbor of $\boldsymbol{\theta}$. Meanwhile, the global square-ring neighbor of $\boldsymbol{\theta}$ may be defined as $N(\boldsymbol{\theta}) = \mathbb{Z}^p \setminus \boldsymbol{\theta}$. Also we denote $M_k$ as the maximum number of comparisons in the $k$th iteration, and $M_k$ changes with iteration. Furthermore, for any $\boldsymbol{\theta} \in \Theta$ and $\boldsymbol{\theta}' \in N(\boldsymbol{\theta})$, we define $R(\boldsymbol{\theta}, \boldsymbol{\theta}')$ as the transition probability with $R(\boldsymbol{\theta}, \boldsymbol{\theta}') > 0$ and $\sum_{\boldsymbol{\theta}' \in N(\boldsymbol{\theta})} R(\boldsymbol{\theta}, \boldsymbol{\theta}') = 1$. In particular, given a current point $\boldsymbol{\theta}$, the probability to generate $\boldsymbol{\theta}' \in N(\boldsymbol{\theta})$ as the possible candidate is $R(\boldsymbol{\theta}, \boldsymbol{\theta}')$. For example, the simplest way to define $R(\boldsymbol{\theta}, \boldsymbol{\theta}')$ is uniform distribution, i.e., $R(\boldsymbol{\theta}, \boldsymbol{\theta}') = 1/|N(\boldsymbol{\theta})|$, where $|N(\boldsymbol{\theta})|$ is the number of elements in $N(\boldsymbol{\theta})$.

The basic description of the SR algorithm is:

Step 1: Pick the initial guess $\hat{\boldsymbol{\theta}}_0$, $k = 0$.

Step 2: Given $\hat{\boldsymbol{\theta}}_k$, generate $\tilde{\boldsymbol{\theta}}_k \in N(\hat{\boldsymbol{\theta}}_k) \subseteq \Theta \setminus \hat{\boldsymbol{\theta}}_k$ according to the probability distribution $R(\hat{\boldsymbol{\theta}}_k, \tilde{\boldsymbol{\theta}}_k)$ (defined above).

Step 3: Given $\tilde{\boldsymbol{\theta}}_k$, set

$$\hat{\boldsymbol{\theta}}_{k+1} = \begin{cases} \tilde{\boldsymbol{\theta}}_k & \text{with probability } P\left(y(\tilde{\boldsymbol{\theta}}_k) \le U_{u,v} \middle| \tilde{\boldsymbol{\theta}}_k\right)^{M_k} \\ \hat{\boldsymbol{\theta}}_k & \text{with probability } 1 - P\left(y(\tilde{\boldsymbol{\theta}}_k) \le U_{u,v} \middle| \tilde{\boldsymbol{\theta}}_k\right)^{M_k} \end{cases},$$

where the values of $u$ and $v$ are picked by the user based on their experience, and the values of $u$ and $v$ are fixed during the whole process of the algorithm.

Step 4: Replace $k$ with $k$+1, and go to step 2.

The implementation of Step 3 in the SR algorithm can be accomplished by generating a noisy measurement of the loss function $y(\tilde{\boldsymbol{\theta}}_k)$, and a uniformly distributed random variable $U_{u,v}$. If $y(\tilde{\boldsymbol{\theta}}_k) > U_{u,v}$, then set $\hat{\boldsymbol{\theta}}_{k+1} = \hat{\boldsymbol{\theta}}_k$; otherwise generate another noisy measurement of loss function $y(\tilde{\boldsymbol{\theta}}_k)$ and another uniform distributed random variable $U_{u,v}$, and continue the comparison. If $y(\tilde{\boldsymbol{\theta}}_k) \le U_{u,v}$ in all $M_k$ comparisons, then we accept the candidate and set $\hat{\boldsymbol{\theta}}_{k+1} = \tilde{\boldsymbol{\theta}}_k$. Conditional on $\tilde{\boldsymbol{\theta}}_k$, these comparisons are independent. Therefore, given $\tilde{\boldsymbol{\theta}}_k$, the probability that $y(\tilde{\boldsymbol{\theta}}_k) \le U_{u,v}$ in all $M_k$ comparisons is $P\left(y(\tilde{\boldsymbol{\theta}}_k) \le U_{u,v} \middle| \tilde{\boldsymbol{\theta}}_k\right)^{M_k}$. Overall, the implementation of the algorithm is to do at most $M_k$ comparisons between the noisy observations at point $\tilde{\boldsymbol{\theta}}_k$ and the stochastic ruler in the $k$th iteration. If one of the comparisons shows that $U_{u,v}$ is smaller, then $\hat{\boldsymbol{\theta}}_{k+1} = \hat{\boldsymbol{\theta}}_k$; otherwise the current point $\hat{\boldsymbol{\theta}}_k$ will be replaced by $\tilde{\boldsymbol{\theta}}_k$. Suppose $y_i(\tilde{\boldsymbol{\theta}}_k)$ is the $i$th noisy measurement at point $\tilde{\boldsymbol{\theta}}_k$ and $U_{u,v}^{(i)}$ is the $i$th simulation result of the uniform distributed stochastic ruler, then in short we have

$$\hat{\boldsymbol{\theta}}_{k+1} = \begin{cases} \tilde{\boldsymbol{\theta}}_k & \text{if} \quad y_i(\tilde{\boldsymbol{\theta}}_k) \le U_{u,v}^{(i)} \quad \forall i = 1, \ldots, M_k \\ \hat{\boldsymbol{\theta}}_k & \text{otherwise.} \end{cases}$$

### 5.2.2 Convergence Properties

In Yan and Mukai (1992), the authors consider the case where $\Theta$ contains a finite number of points (i.e., $\Theta$ is bounded). In their paper, they show that when $M_k$ is fixed ($M_k = M$) for all iterations, which means the maximum number of comparisons does not change with the number of iteration, the sequence $\{\hat{\boldsymbol{\theta}}_k\}$ generated by the algorithm produces a stationary Markov chain. Suppose $\boldsymbol{\Pi}(M)$ is the stationary probability distribution vector conditioned on the value of $M$ when $k \to \infty$. The dimension of $\boldsymbol{\Pi}(M)$ equals the number of points in the bounded domain $\Theta$. For each $\boldsymbol{\theta} \in \Theta$, the corresponding component of $\boldsymbol{\Pi}(M)$ is

$$\Pi_{\boldsymbol{\theta}}(M) = \lim_{k\to\infty} P\left(\hat{\boldsymbol{\theta}}_k = \boldsymbol{\theta} \middle| M\right) = \frac{\left\{P\left(y(\boldsymbol{\theta}) \le U_{u,v}\right)\right\}^M}{\sum_{\boldsymbol{\theta}'\in\Theta}\left\{P\left(y(\boldsymbol{\theta}') \le U_{u,v}\right)\right\}^M}.$$

Theorem 6.1 of Yan and Mukai (1992) shows that when $M \to \infty$, the probability vector $\boldsymbol{\Pi}(M)$ converges to the probability distribution $\boldsymbol{\Pi}^*$, where each component of $\boldsymbol{\Pi}^*$ is

$$\Pi_{\boldsymbol{\theta}}^* = \lim_{k\to\infty} P(\hat{\boldsymbol{\theta}}_k = \boldsymbol{\theta}) = \begin{cases} \dfrac{1}{|\Theta^*|} & \text{if} \quad \boldsymbol{\theta} \in \Theta^* \\ 0 & \text{otherwise,} \end{cases}$$

where $\Theta^*$ is the set of optimal points, and $|\Theta^*|$ represents the number of points in $\Theta^*$.

Furthermore, in Theorem 7.2 of Yan and Mukai (1992), the authors provide a convergence theorem for the case with $M_k = \lfloor c \log_\sigma(1 + k_0 + k) \rfloor$ (not fixed at one value

*M*). Suppose $r = \min_{\boldsymbol{\theta}\in\Theta} \max_{\boldsymbol{\theta}'\in\Theta} d(\boldsymbol{\theta},\boldsymbol{\theta}')$, where $d(\boldsymbol{\theta},\boldsymbol{\theta}')$ is the distance between $\boldsymbol{\theta}$ and $\boldsymbol{\theta}'$, which is defined as the minimum length of path from $\boldsymbol{\theta}$ to $\boldsymbol{\theta}'$ (the length of path from $\boldsymbol{\theta}$ to $\boldsymbol{\theta}'$ is defined as the number of edges on the path, and $d(\boldsymbol{\theta},\boldsymbol{\theta}')$ is an integer). Here *r* describes the radius of the graph of the domain. The constraints for $\sigma$, *c*, and $k_0$ in the definition of $M_k$ are set as $\sigma \geq 1/\min_{\boldsymbol{\theta}\in\Theta} P(y(\boldsymbol{\theta}) \leq U_{u,v})$ (by conclusion (2) of Theorem 3.1 in Yan and Mukai, 1992, they show that $0 < P(y(\boldsymbol{\theta}) \leq U_{u,v}) < 1$), $0 < c \leq 1/r$, and $c \log_{\sigma}(k_0 + 1) \geq 1$. The authors show that when $M_k = \lfloor c \log_{\sigma}(1 + k_0 + k) \rfloor$, under the constraints on the coefficients ($\sigma$, *c*, and $k_0$), the sequence $\{\hat{\boldsymbol{\theta}}_k\}$ generated by the SR algorithm converges in probability

$$\lim_{k\to\infty} P(\hat{\boldsymbol{\theta}}_k \in \Theta^*) = 1.$$

In addition, suppose the vector $\tilde{\boldsymbol{\Pi}}_k(\hat{\boldsymbol{\theta}}_k)$ (with each component being $P(\hat{\boldsymbol{\theta}}_k = \boldsymbol{\theta})$) is the probability distribution for $\hat{\boldsymbol{\theta}}_k$ in the *k*th iteration for non-fixed $M_k$ case. In Theorem 8.1 of Yan and Mukai (1992), the authors show that for this inhomogenous Markov chain (i.e., time varying transition matrix), $\left\|\tilde{\boldsymbol{\Pi}}_k(\hat{\boldsymbol{\theta}}_k) - \boldsymbol{\Pi}^*\right\| \leq O(1/k^h)$ (a geometric convergence rate is for homogenous Markov chains (Spall, 2003, Theorem E.1)), where $h > 0$ is determined by the transition probabilities, domain $\Theta$, coefficients *c* and $\sigma$, and value of $\max_{\boldsymbol{\theta}\in\Theta} P\left(y(\boldsymbol{\theta}) \leq U_{u,v}\right)$.

## 5.3 Stochastic Comparison Algorithm

Gong et al. (1999) discuss the SC algorithm, and the idea is also to change the minimization problem

$$\min_{\boldsymbol{\theta} \in \Theta} E(y(\boldsymbol{\theta}))$$

to a maximization problem

$$\max_{\boldsymbol{\theta} \in \Theta} \sum_{\tilde{\boldsymbol{\theta}} \in N(\boldsymbol{\theta})} P(y(\boldsymbol{\theta}) < y(\tilde{\boldsymbol{\theta}})),$$

where $N(\boldsymbol{\theta}) = \Theta \setminus \boldsymbol{\theta}$. In Assumption 3.1 of Gong et al. (1999), the authors suppose that the noises in the measurements of the loss function are i.i.d. and have symmetric distribution with mean 0. Furthermore, in Theorem 3.1 of Gong et al. (1999), the authors show that under Assumption 3.1 of Gong et al. (1999), the optimal solution for the minimization problem can achieve highest value for the related maximization problem above.

### 5.3.1 Algorithm Description

In Gong et al. (1999), the authors only consider global square-ring neighborhood structure $N(\boldsymbol{\theta}) = \Theta \setminus \boldsymbol{\theta}$. The authors state: “While a good neighborhood structure can speed up the search process of algorithms like SR, a poor neighborhood structure can hurt performance.” Thus, Gong et al. (1999) eliminate the use of a neighborhood structure. In the SC algorithm, $M_k$ is still the maximum number of comparisons in the *k*th iteration, and $M_k$ changes with iteration. Furthermore, for any $\boldsymbol{\theta} \in \Theta$ and $\boldsymbol{\theta}' \in N(\boldsymbol{\theta})$, let $R(\boldsymbol{\theta}, \boldsymbol{\theta}')$ be the transition probability to generate $\boldsymbol{\theta}'$ as possible candidate with $R(\boldsymbol{\theta}, \boldsymbol{\theta}') > 0$ and $\sum_{\boldsymbol{\theta}' \in N(\boldsymbol{\theta})} R(\boldsymbol{\theta}, \boldsymbol{\theta}') = 1$. The basic description of the SC algorithm is quite similar to the SR algorithm:

Step 1: Pick the initial guess $\hat{\boldsymbol{\theta}}_0$, $k = 0$.

Step 2: Given $\hat{\boldsymbol{\theta}}_k$, generate $\tilde{\boldsymbol{\theta}}_k \in N(\hat{\boldsymbol{\theta}}_k) = \Theta \setminus \hat{\boldsymbol{\theta}}_k$ according to the probability distribution $R(\hat{\boldsymbol{\theta}}_k, \tilde{\boldsymbol{\theta}}_k)$ (defined above).

Step 3: Given $\tilde{\boldsymbol{\theta}}_k$, set

$$\hat{\boldsymbol{\theta}}_{k+1} = \begin{cases} \tilde{\boldsymbol{\theta}}_k & \text{with probability } P\left( y(\tilde{\boldsymbol{\theta}}_k) \le y(\hat{\boldsymbol{\theta}}_k) \middle| \tilde{\boldsymbol{\theta}}_k, \hat{\boldsymbol{\theta}}_k \right)^{M_k} \\ \hat{\boldsymbol{\theta}}_k & \text{with probability } 1 - P\left( y(\tilde{\boldsymbol{\theta}}_k) \le y(\hat{\boldsymbol{\theta}}_k) \middle| \tilde{\boldsymbol{\theta}}_k, \hat{\boldsymbol{\theta}}_k \right)^{M_k} \end{cases}.$$

Step 4: Replace $k$ with $k+1$, and go to step 2.

The implementation of Step 3 of the SC algorithm can be accomplished by first generating the pair of noisy measurements of the loss function, $y(\tilde{\boldsymbol{\theta}}_k)$ and $y(\hat{\boldsymbol{\theta}}_k)$. If $y(\tilde{\boldsymbol{\theta}}_k) > y(\hat{\boldsymbol{\theta}}_k)$, then set $\hat{\boldsymbol{\theta}}_{k+1} = \hat{\boldsymbol{\theta}}_k$; otherwise, given $\hat{\boldsymbol{\theta}}_k$ and $\tilde{\boldsymbol{\theta}}_k$, generate another independent pair of noisy measurements of loss function $y(\tilde{\boldsymbol{\theta}}_k)$ and $y(\hat{\boldsymbol{\theta}}_k)$, and continue the comparison. If in all $M_k$ comparisons $y(\tilde{\boldsymbol{\theta}}_k) \le y(\hat{\boldsymbol{\theta}}_k)$, then we accept the candidate and set $\hat{\boldsymbol{\theta}}_{k+1} = \tilde{\boldsymbol{\theta}}_k$. Overall, given $\tilde{\boldsymbol{\theta}}_k$ and $\hat{\boldsymbol{\theta}}_k$, the implementation of the algorithm is to do at most $M_k$ conditionally independent comparisons between the noisy observations at points $\tilde{\boldsymbol{\theta}}_k$ and $\hat{\boldsymbol{\theta}}_k$. If one of the comparisons indicates that $y(\hat{\boldsymbol{\theta}}_k)$ is smaller, then $\hat{\boldsymbol{\theta}}_{k+1} = \hat{\boldsymbol{\theta}}_k$; otherwise the current point $\hat{\boldsymbol{\theta}}_k$ will be replaced by $\tilde{\boldsymbol{\theta}}_k$. Suppose $y_i(\tilde{\boldsymbol{\theta}}_k)$ and $y_i(\hat{\boldsymbol{\theta}}_k)$ are the $i$th noisy measurement at point $\tilde{\boldsymbol{\theta}}_k$ and $\hat{\boldsymbol{\theta}}_k$, then in short we have

$$\hat{\boldsymbol{\theta}}_{k+1} = \begin{cases} \tilde{\boldsymbol{\theta}}_k & \text{if} \quad y_i(\tilde{\boldsymbol{\theta}}_k) \le y_i(\hat{\boldsymbol{\theta}}_k) \quad \forall i = 1, \ldots, M_k \\ \hat{\boldsymbol{\theta}}_k & \text{otherwise.} \end{cases}$$

## 5.3.2 Convergence Properties

In Gong et al. (1999), the authors consider the case where $\Theta$ contains finite points (i.e., $\Theta$ is bounded), and the authors show the convergence properties in their Theorem 5.1. They suppose $\sigma \ge 1 / \min_{\boldsymbol{\theta} \in \Theta, \boldsymbol{\theta}' \in \Theta \setminus \boldsymbol{\theta}} P(y(\boldsymbol{\theta}') < y(\boldsymbol{\theta}))$, $0 < c \le 1$, $c \log_\sigma(1 + k_0) \ge 1$. Under these assumptions, when $M_k = \lfloor c \log_\sigma(1 + k_0 + k) \rfloor$, Gong et al. (1999) show that the sequence $\{\hat{\boldsymbol{\theta}}_k\}$ generated by the SC algorithm is a strongly ergodic inhomogeneous Markov chain, and

$$\lim_{k \to \infty} P(\hat{\boldsymbol{\theta}}_k \in \Theta^*) = 1,$$

where $\Theta^*$ is the set of optimal points. However, Gong et al. (1999) do not provide a formal rate of convergence analysis, and they only do some numerical experiments to compare the performance of the SR algorithm and the SC algorithm. In the following, we present a rate of convergence analysis for both SR algorithm and SC algorithm.

## 5.4 Rate of Convergence Analysis and Comparison

In this section, we derive a rate of convergence for both the SR algorithm and the SC algorithm. From the description of Section 5.2.1 and Section 5.3.1, we see that these two algorithms have similar descriptions. Therefore, in the following, first we provide a general form of algorithm description for these two algorithms (the general form might also represent other random search type algorithms), second we discuss the rate of convergence result for the general form in Theorem 5.1, and last we consider the rate of convergence results of the SR and the SC separately in Corollary 5.1 and Corollary 5.2. We can describe the general form of these two algorithms as:

Step 1: Pick the initial guess $\hat{\boldsymbol{\theta}}_0$, $k = 0$.

Step 2: Given $\hat{\boldsymbol{\theta}}_k$, generate $\tilde{\boldsymbol{\theta}}_k \in N(\hat{\boldsymbol{\theta}}_k) \subseteq \Theta \setminus \hat{\boldsymbol{\theta}}_k$ according to the probability distribution $R(\hat{\boldsymbol{\theta}}_k, \tilde{\boldsymbol{\theta}}_k)$.

Step 3: Given $\tilde{\boldsymbol{\theta}}_k$, set

$$\hat{\boldsymbol{\theta}}_{k+1} = \begin{cases} \tilde{\boldsymbol{\theta}}_k & \text{with probability } p(\hat{\boldsymbol{\theta}}_k, \tilde{\boldsymbol{\theta}}_k)^{M_k} \\ \hat{\boldsymbol{\theta}}_k & \text{with probability } 1 - p(\hat{\boldsymbol{\theta}}_k, \tilde{\boldsymbol{\theta}}_k)^{M_k}, \end{cases}$$

where $M_k = \lfloor c \log_\sigma (1 + k_0 + k) \rfloor$, $p(\hat{\boldsymbol{\theta}}_k, \tilde{\boldsymbol{\theta}}_k) = P\left( y(\tilde{\boldsymbol{\theta}}_k) \le U_{u,v} \middle| \tilde{\boldsymbol{\theta}}_k \right)$ for the SR algorithm, and $p(\hat{\boldsymbol{\theta}}_k, \tilde{\boldsymbol{\theta}}_k) = P\left( y(\tilde{\boldsymbol{\theta}}_k) \le y(\hat{\boldsymbol{\theta}}_k) \middle| \hat{\boldsymbol{\theta}}_k, \tilde{\boldsymbol{\theta}}_k \right)$ for the SC algorithm.

Step 4: Replace $k$ with $k+1$, and go to step 2.

In the following, we state Theorem 5.1, and this theorem provides an upper bound on the rate of convergence of $P(\hat{\boldsymbol{\theta}}_k \ne \boldsymbol{\theta}^*)$. The rate of convergence describes the asymptotical performance of the algorithm, so the result is only related to loss function information in the neighbor of the optimal solution.

**Theorem 5.1.** For the general algorithm discussed above (Step 1−4) on the bounded domain $\Theta$, if (i) the sequence $\{\hat{\boldsymbol{\theta}}_k\}$ generated by the algorithm convergences to the optimal point $\boldsymbol{\theta}^*$, $\lim_{k\to\infty} P\left(\hat{\boldsymbol{\theta}}_k = \boldsymbol{\theta}^*\right) = 1$; (ii) for all $\boldsymbol{\theta} \in \Theta$ and $\boldsymbol{\theta}' \in N(\boldsymbol{\theta}) \subseteq \Theta \backslash \boldsymbol{\theta}$, we have that $\sum_{\boldsymbol{\theta}' \in N(\boldsymbol{\theta})} R(\boldsymbol{\theta}, \boldsymbol{\theta}') = 1$, $R(\boldsymbol{\theta}, \boldsymbol{\theta}') > 0$, and $0 < p(\boldsymbol{\theta}, \boldsymbol{\theta}') < 1$; and (iii) $M_k = \lfloor c \log_\sigma (1 + k_0 + k) \rfloor$, $c > 0$, $\sigma > 1$, and $c \log_\sigma (1 + k_0) \ge 1$. Then

$$k^{-\frac{c \log\left( 1 / \max_{\boldsymbol{\theta} \in N(\boldsymbol{\theta}^*)} p(\boldsymbol{\theta}^*, \boldsymbol{\theta}) \right)}{\log(\sigma)}} = O\left( P\left( \hat{\boldsymbol{\theta}}_k \ne \boldsymbol{\theta}^* \right) \right),$$

where $N(\boldsymbol{\theta}^*) \subseteq \Theta \backslash \boldsymbol{\theta}^*$ and the base of logarithm here is $e$ (natural logarithm), which means $P\left(\hat{\boldsymbol{\theta}}_k \neq \boldsymbol{\theta}^*\right)$(i.e. $1 - P(\hat{\boldsymbol{\theta}}_k = \boldsymbol{\theta}^*)$) goes to 0 at a rate not faster than the rate at which $k^{-c\log\left(1/\max_{\boldsymbol{\theta}\in N(\boldsymbol{\theta}^*)} p(\boldsymbol{\theta}^*,\boldsymbol{\theta})\right)/\log(\sigma)}$ goes to 0.

*Remarks:*

1. In Theorem 5.1, our purpose is to calculate the bound of the rate of convergence of the algorithm under the condition (i): $\lim_{k\to\infty} P\left(\hat{\boldsymbol{\theta}}_k = \boldsymbol{\theta}^*\right) = 1$. The condition for the coefficients of $M_k$ (condition (iii)) is weaker than that in the SR algorithm (discussed in Section 5.2.2) and that in the SC algorithm (conditions are discussed in Section 5.3.2), because here we have already assumed that the sequence generated by the algorithm converges to the optimal solution. However, the conditions of the coefficients in SR and SC are still needed in Corollaries 5.1 and 5.2 to show that the sequences generated by the algorithms SR and SC converge to the optimal solution.

2. The value of $k_0$ is used to make the number of comparison in each iteration not to be smaller than 1. As $k \to \infty$, the effect of $k_0$ disappears.

3. The value of $c$ affects the increase rate of $M_k$. If we pick relatively small value for $c$, then $M_k$ goes to $\infty$ at a slow speed, and if we pick relatively large value for $c$, then $M_k$ goes to $\infty$ at a fast speed. Thus, the value of $c$ affects the rate of convergence.

4. The value of $\sigma$ also affects the increase rate of $M_k$. If we pick relatively small value for $\sigma$, then $M_k$ goes to $\infty$ at a fast speed, and if we pick relatively large value for $\sigma$, then $M_k$ goes to $\infty$ at a slow speed. Thus, the value of $\sigma$ affects the rate of convergence.

*Proof.* By law of total expectation, we have

$$
\begin{aligned}
P\left(\hat{\boldsymbol{\theta}}_k \neq \boldsymbol{\theta}^*\right) = {} & P\left(\hat{\boldsymbol{\theta}}_k \neq \boldsymbol{\theta}^* \middle| \hat{\boldsymbol{\theta}}_{k-1} = \boldsymbol{\theta}^*\right) P\left(\hat{\boldsymbol{\theta}}_{k-1} = \boldsymbol{\theta}^*\right) \\
& + P\left(\hat{\boldsymbol{\theta}}_k \neq \boldsymbol{\theta}^* \middle| \hat{\boldsymbol{\theta}}_{k-1} \neq \boldsymbol{\theta}^*\right) P\left(\hat{\boldsymbol{\theta}}_{k-1} \neq \boldsymbol{\theta}^*\right) \\
= {} & P\left(\hat{\boldsymbol{\theta}}_k \neq \boldsymbol{\theta}^* \middle| \hat{\boldsymbol{\theta}}_{k-1} = \boldsymbol{\theta}^*\right)\left(1 - P\left(\hat{\boldsymbol{\theta}}_{k-1} \neq \boldsymbol{\theta}^*\right)\right) \\
& + \left(1 - P\left(\hat{\boldsymbol{\theta}}_k = \boldsymbol{\theta}^* \middle| \hat{\boldsymbol{\theta}}_{k-1} \neq \boldsymbol{\theta}^*\right)\right) P\left(\hat{\boldsymbol{\theta}}_{k-1} \neq \boldsymbol{\theta}^*\right) \\
= {} & P\left(\hat{\boldsymbol{\theta}}_k \neq \boldsymbol{\theta}^* \middle| \hat{\boldsymbol{\theta}}_{k-1} = \boldsymbol{\theta}^*\right) \\
& + \left(1 - P\left(\hat{\boldsymbol{\theta}}_k = \boldsymbol{\theta}^* \middle| \hat{\boldsymbol{\theta}}_{k-1} \neq \boldsymbol{\theta}^*\right) - P\left(\hat{\boldsymbol{\theta}}_k \neq \boldsymbol{\theta}^* \middle| \hat{\boldsymbol{\theta}}_{k-1} = \boldsymbol{\theta}^*\right)\right) P\left(\hat{\boldsymbol{\theta}}_{k-1} \neq \boldsymbol{\theta}^*\right). \qquad (5.1)
\end{aligned}
$$

In the following, we will show that the multiplier of $P\left(\hat{\boldsymbol{\theta}}_{k-1} \neq \boldsymbol{\theta}^*\right)$ on the right-hand side of eqn. (5.1) is positive when $k$ is big enough, which implies $P\left(\hat{\boldsymbol{\theta}}_k \neq \boldsymbol{\theta}^*\right)$ is a positive combination of $P\left(\hat{\boldsymbol{\theta}}_k \neq \boldsymbol{\theta}^* \middle| \hat{\boldsymbol{\theta}}_{k-1} = \boldsymbol{\theta}^*\right)$ and $P\left(\hat{\boldsymbol{\theta}}_{k-1} \neq \boldsymbol{\theta}^*\right)$. Since any positive combination of $P\left(\hat{\boldsymbol{\theta}}_k \neq \boldsymbol{\theta}^* \middle| \hat{\boldsymbol{\theta}}_{k-1} = \boldsymbol{\theta}^*\right)$ and $P\left(\hat{\boldsymbol{\theta}}_{k-1} \neq \boldsymbol{\theta}^*\right)$ goes to 0 at a rate not faster than $P\left(\hat{\boldsymbol{\theta}}_k \neq \boldsymbol{\theta}^* \middle| \hat{\boldsymbol{\theta}}_{k-1} = \boldsymbol{\theta}^*\right)$ and $P\left(\hat{\boldsymbol{\theta}}_{k-1} \neq \boldsymbol{\theta}^*\right)$ individually, it implies that $P\left(\hat{\boldsymbol{\theta}}_k \neq \boldsymbol{\theta}^*\right)$ goes to 0 at a rate not faster than the rate at which $P\left(\hat{\boldsymbol{\theta}}_k \neq \boldsymbol{\theta}^* \middle| \hat{\boldsymbol{\theta}}_{k-1} = \boldsymbol{\theta}^*\right)$ goes to 0, which means $P\left(\hat{\boldsymbol{\theta}}_k \neq \boldsymbol{\theta}^* \middle| \hat{\boldsymbol{\theta}}_{k-1} = \boldsymbol{\theta}^*\right) = O\left(P\left(\hat{\boldsymbol{\theta}}_k \neq \boldsymbol{\theta}^*\right)\right)$.

Now let us start the proof by considering the values of $P\left(\hat{\boldsymbol{\theta}}_k \neq \boldsymbol{\theta}^* \middle| \hat{\boldsymbol{\theta}}_{k-1} = \boldsymbol{\theta}^*\right)$ and $P\left(\hat{\boldsymbol{\theta}}_k = \boldsymbol{\theta}^* \middle| \hat{\boldsymbol{\theta}}_{k-1} \neq \boldsymbol{\theta}^*\right)$ separately. By the law of total expectation, we have

$$P\left(\hat{\boldsymbol{\theta}}_k \neq \boldsymbol{\theta}^* \middle| \hat{\boldsymbol{\theta}}_{k-1} = \boldsymbol{\theta}^*\right) = \sum_{\boldsymbol{\theta} \in N(\boldsymbol{\theta}^*)} P\left(\hat{\boldsymbol{\theta}}_k = \boldsymbol{\theta} \middle| \hat{\boldsymbol{\theta}}_{k-1} = \boldsymbol{\theta}^*\right)$$

$$= \sum_{\boldsymbol{\theta} \in N(\boldsymbol{\theta}^*)} R(\boldsymbol{\theta}^*, \boldsymbol{\theta}) p(\boldsymbol{\theta}^*, \boldsymbol{\theta})^{M_k}. \tag{5.2}$$

Because $M_k = \lfloor c \log_\sigma(1 + k_0 + k) \rfloor$, we have $c \log_\sigma(1 + k_0 + k) - 1 < M_k \leq c \log_\sigma(1 + k_0 + k)$, which indicates that

$$p(\boldsymbol{\theta}^*, \boldsymbol{\theta})^{c \log_\sigma(1+k_0+k)} \leq p(\boldsymbol{\theta}^*, \boldsymbol{\theta})^{M_k} < p(\boldsymbol{\theta}^*, \boldsymbol{\theta})^{c \log_\sigma(1+k_0+k)-1}. \tag{5.3}$$

We discuss the rate of convergence of $p(\boldsymbol{\theta}^*, \boldsymbol{\theta})^{M_k}$ through the lower and upper bounds in inequality (5.3). In addition, due to the assumption: $0 < p(\boldsymbol{\theta}^*, \boldsymbol{\theta}) < 1$, we have

$$\log_\sigma(1 + k_0 + k) = \frac{\log(1 + k_0 + k)}{\log(\sigma)}$$

$$= \frac{\log(1 + k_0 + k)}{\log\left(p(\boldsymbol{\theta}^*, \boldsymbol{\theta})\right)} \frac{\log\left(p(\boldsymbol{\theta}^*, \boldsymbol{\theta})\right)}{\log(\sigma)}$$

$$= \left(\log_{p(\boldsymbol{\theta}^*, \boldsymbol{\theta})}(1 + k_0 + k)\right) \log_\sigma\left(p(\boldsymbol{\theta}^*, \boldsymbol{\theta})\right),$$

and it follows that

$$p(\boldsymbol{\theta}^*, \boldsymbol{\theta})^{c \log_\sigma(1+k_0+k)} = p(\boldsymbol{\theta}^*, \boldsymbol{\theta})^{c\left(\log_{p(\boldsymbol{\theta}^*, \boldsymbol{\theta})}(1+k_0+k)\right)\log_\sigma p(\boldsymbol{\theta}^*, \boldsymbol{\theta})}$$

$$= (1 + k_0 + k)^{c \log_\sigma p(\boldsymbol{\theta}^*, \boldsymbol{\theta})}.$$

Therefore, inequality (5.3) can be written as

$$\left(1+k_0+k\right)^{c\log_\sigma p(\boldsymbol{\theta}^*,\boldsymbol{\theta})} \le p(\boldsymbol{\theta}^*,\boldsymbol{\theta})^{M_k} < \left(\frac{1}{p(\boldsymbol{\theta}^*,\boldsymbol{\theta})}\right)\left(1+k_0+k\right)^{c\log_\sigma p(\boldsymbol{\theta}^*,\boldsymbol{\theta})}. \tag{5.4}$$

Since $\sigma > 1$, $c > 0$, and $0 < p(\boldsymbol{\theta}^*,\boldsymbol{\theta}) < 1$, then we have $c\log_\sigma\left(p(\boldsymbol{\theta}^*,\boldsymbol{\theta})\right) < 0$. Therefore, $\left(1+k_0+k\right)^{c\log_\sigma p(\boldsymbol{\theta}^*,\boldsymbol{\theta})}$ goes to 0 as $k \to \infty$, implying from inequality (5.4) that $p(\boldsymbol{\theta}^*,\boldsymbol{\theta})^{M_k}$ goes to 0 as fast as $\left(1+k_0+k\right)^{c\log_\sigma p(\boldsymbol{\theta}^*,\boldsymbol{\theta})}$ goes to 0.

Furthermore, since $\sum_{\boldsymbol{\theta}\in N(\boldsymbol{\theta}^*)} R(\boldsymbol{\theta}^*,\boldsymbol{\theta}) = 1$, $R(\boldsymbol{\theta}^*,\boldsymbol{\theta}) > 0$ for all $\boldsymbol{\theta} \in N(\boldsymbol{\theta}^*)$, and $N(\boldsymbol{\theta}^*)$ contains finite points, then through eqn. (5.2) we know $P\left(\hat{\boldsymbol{\theta}}_k \ne \boldsymbol{\theta}^* \middle| \hat{\boldsymbol{\theta}}_{k-1} = \boldsymbol{\theta}^*\right)$ is a positive finite combination of all $p(\boldsymbol{\theta}^*,\boldsymbol{\theta})^{M_k}$ over $\boldsymbol{\theta} \in N(\boldsymbol{\theta}^*)$ at each $k$, which indicates that $P\left(\hat{\boldsymbol{\theta}}_k \ne \boldsymbol{\theta}^* \middle| \hat{\boldsymbol{\theta}}_{k-1} = \boldsymbol{\theta}^*\right)$ goes to 0 at the lowest rate among all rates at which $p(\boldsymbol{\theta}^*,\boldsymbol{\theta})^{M_k}$ goes to 0 for all $\boldsymbol{\theta} \in N(\boldsymbol{\theta}^*)$ (note that the multipliers in the finite combination are not dependent on $k$). The lowest rate of all $p(\boldsymbol{\theta}^*,\boldsymbol{\theta})^{M_k} \to 0$ for $\boldsymbol{\theta} \in N(\boldsymbol{\theta}^*)$ is the same as the rate at which $\left(1+k_0+k\right)^{c\log_\sigma\left(\max_{\boldsymbol{\theta}\in N(\boldsymbol{\theta}^*)} p(\boldsymbol{\theta}^*,\boldsymbol{\theta})\right)}$ goes to 0. Therefore, there exist $\tilde{c} > 0$, $\tilde{\tilde{c}} > 0$, and $\tilde{\tilde{c}} > \tilde{c}$ such that for all $k$ large enough

$$\tilde{c} < \frac{P\left(\hat{\boldsymbol{\theta}}_k \ne \boldsymbol{\theta}^* \middle| \hat{\boldsymbol{\theta}}_{k-1} = \boldsymbol{\theta}^*\right)}{k^{c\log_\sigma\left(\max_{\boldsymbol{\theta}\in N(\boldsymbol{\theta}^*)} p(\boldsymbol{\theta}^*,\boldsymbol{\theta})\right)}} < \tilde{\tilde{c}}. \tag{5.5}$$

We can derive similar arguments for $P\left(\hat{\boldsymbol{\theta}}_k = \boldsymbol{\theta}^* \middle| \hat{\boldsymbol{\theta}}_{k-1} \neq \boldsymbol{\theta}^*\right)$. Here

$$P\left(\hat{\boldsymbol{\theta}}_k = \boldsymbol{\theta}^* \middle| \hat{\boldsymbol{\theta}}_{k-1} \neq \boldsymbol{\theta}^*\right) = \sum_{\boldsymbol{\theta} \neq \boldsymbol{\theta}^*} P\left(\hat{\boldsymbol{\theta}}_k = \boldsymbol{\theta}^* \middle| \hat{\boldsymbol{\theta}}_{k-1} \neq \boldsymbol{\theta}^*, \hat{\boldsymbol{\theta}}_{k-1} = \boldsymbol{\theta}\right) P\left(\hat{\boldsymbol{\theta}}_{k-1} = \boldsymbol{\theta} \middle| \hat{\boldsymbol{\theta}}_{k-1} \neq \boldsymbol{\theta}^*\right)$$

$$= \sum_{\boldsymbol{\theta} \neq \boldsymbol{\theta}^*} P\left(\hat{\boldsymbol{\theta}}_k = \boldsymbol{\theta}^* \middle| \hat{\boldsymbol{\theta}}_{k-1} = \boldsymbol{\theta}\right) P\left(\hat{\boldsymbol{\theta}}_{k-1} = \boldsymbol{\theta} \middle| \hat{\boldsymbol{\theta}}_{k-1} \neq \boldsymbol{\theta}^*\right),$$

and for each $\boldsymbol{\theta} \neq \boldsymbol{\theta}^*$, if $\boldsymbol{\theta}^* \in N(\boldsymbol{\theta})$, then

$$P\left(\hat{\boldsymbol{\theta}}_k = \boldsymbol{\theta}^* \middle| \hat{\boldsymbol{\theta}}_{k-1} = \boldsymbol{\theta}\right) = R(\boldsymbol{\theta}, \boldsymbol{\theta}^*) p(\boldsymbol{\theta}, \boldsymbol{\theta}^*)^{M_k},$$

and if $\boldsymbol{\theta}^* \notin N(\boldsymbol{\theta})$, then $P\left(\hat{\boldsymbol{\theta}}_k = \boldsymbol{\theta}^* \middle| \hat{\boldsymbol{\theta}}_{k-1} = \boldsymbol{\theta}\right) = 0$, which indicates that

$$P\left(\hat{\boldsymbol{\theta}}_k = \boldsymbol{\theta}^* \middle| \hat{\boldsymbol{\theta}}_{k-1} \neq \boldsymbol{\theta}^*\right)$$

$$= \sum_{\boldsymbol{\theta} \neq \boldsymbol{\theta}^*,\, \boldsymbol{\theta}^* \in N(\boldsymbol{\theta})} R(\boldsymbol{\theta}, \boldsymbol{\theta}^*) p(\boldsymbol{\theta}, \boldsymbol{\theta}^*)^{M_k} P\left(\hat{\boldsymbol{\theta}}_{k-1} = \boldsymbol{\theta} \middle| \hat{\boldsymbol{\theta}}_{k-1} \neq \boldsymbol{\theta}^*\right)$$

$$= \sum_{\boldsymbol{\theta} \neq \boldsymbol{\theta}^*,\, \boldsymbol{\theta}^* \in N(\boldsymbol{\theta})} R(\boldsymbol{\theta}, \boldsymbol{\theta}^*) P\left(\hat{\boldsymbol{\theta}}_{k-1} = \boldsymbol{\theta} \middle| \hat{\boldsymbol{\theta}}_{k-1} \neq \boldsymbol{\theta}^*\right) p(\boldsymbol{\theta}, \boldsymbol{\theta}^*)^{M_k}. \tag{5.6}$$

Similar to the proof of the rate of convergence for $p(\boldsymbol{\theta}^*, \boldsymbol{\theta})^{M_k}$ in the last three paragraphs, we have $p(\boldsymbol{\theta}, \boldsymbol{\theta}^*)^{M_k} \to 0$. Furthermore, by eqn. (5.6), we know $P\left(\hat{\boldsymbol{\theta}}_k = \boldsymbol{\theta}^* \middle| \hat{\boldsymbol{\theta}}_{k-1} \neq \boldsymbol{\theta}^*\right)$ equals a nonnegative finite combination of $p(\boldsymbol{\theta}^*, \boldsymbol{\theta})^{M_k}$, implying $P\left(\hat{\boldsymbol{\theta}}_k = \boldsymbol{\theta}^* \middle| \hat{\boldsymbol{\theta}}_{k-1} \neq \boldsymbol{\theta}^*\right)$ $\to 0$.

Overall, from the proof above, we know that $P\left(\hat{\boldsymbol{\theta}}_k \neq \boldsymbol{\theta}^* \middle| \hat{\boldsymbol{\theta}}_{k-1} = \boldsymbol{\theta}^*\right) \to 0$ and $P\left(\hat{\boldsymbol{\theta}}_k = \boldsymbol{\theta}^* \middle| \hat{\boldsymbol{\theta}}_{k-1} \neq \boldsymbol{\theta}^*\right) \to 0$ as $k \to \infty$, which indicates that $1 - P\left(\hat{\boldsymbol{\theta}}_k = \boldsymbol{\theta}^* \middle| \hat{\boldsymbol{\theta}}_{k-1} \neq \boldsymbol{\theta}^*\right) - P\left(\hat{\boldsymbol{\theta}}_k \neq \boldsymbol{\theta}^* \middle| \hat{\boldsymbol{\theta}}_{k-1} = \boldsymbol{\theta}^*\right) \to 1$. Therefore, by eqn. (5.1), we know for $k$ large enough, $P\left(\hat{\boldsymbol{\theta}}_k \neq \boldsymbol{\theta}^*\right)$ is a positive combination of $P\left(\hat{\boldsymbol{\theta}}_k \neq \boldsymbol{\theta}^* \middle| \hat{\boldsymbol{\theta}}_{k-1} = \boldsymbol{\theta}^*\right)$ and $P\left(\hat{\boldsymbol{\theta}}_{k-1} \neq \boldsymbol{\theta}^*\right)$, which indicates that $P\left(\hat{\boldsymbol{\theta}}_k \neq \boldsymbol{\theta}^*\right)$ goes to 0 at a rate not faster than the rate of convergence of $P\left(\hat{\boldsymbol{\theta}}_k \neq \boldsymbol{\theta}^* \middle| \hat{\boldsymbol{\theta}}_{k-1} = \boldsymbol{\theta}^*\right)$ to 0 (shown in inequality (5.5)), which is equivalent to

$$k^{-\frac{c\log\left(1/\max_{\boldsymbol{\theta}\in N(\boldsymbol{\theta}^*)} p(\boldsymbol{\theta}^*,\boldsymbol{\theta})\right)}{\log(\sigma)}} = O\left(P\left(\hat{\boldsymbol{\theta}}_k \neq \boldsymbol{\theta}^*\right)\right),$$

Therefore, we have shown the result on the convergence rate of $P\left(\hat{\boldsymbol{\theta}}_k \neq \boldsymbol{\theta}^*\right)$. Q.E.D.

Corollaries 5.1 and 5.2 below apply the general theorem above to the special cases of the SR algorithm and the SC algorithm.

**Corollary 5.1.** For the SR algorithm, if (i) $\Theta^* = \{\boldsymbol{\theta}^*\}$; (ii) for all $\boldsymbol{\theta} \in \Theta$ and $\boldsymbol{\theta}' \in N(\boldsymbol{\theta}) \subseteq \Theta \setminus \boldsymbol{\theta}$, we have that $\sum_{\boldsymbol{\theta}'\in N(\boldsymbol{\theta})} R(\boldsymbol{\theta},\boldsymbol{\theta}')=1$, $R(\boldsymbol{\theta},\boldsymbol{\theta}') > 0$, and $0 < P\left(y(\boldsymbol{\theta}') \leq U_{u,v}\right) < 1$; and (iii) $M_k = \lfloor c\log_\sigma(1+k_0+k)\rfloor$, $0 < c \leq 1/r$, $\sigma \geq 1/\min_{\boldsymbol{\theta}\in\Theta} P(y(\boldsymbol{\theta}) < U_{u,v})$, and $c\log_\sigma(1+k_0) \geq 1$. Then

$$k^{-\frac{c\log\left(1/\max_{\boldsymbol{\theta}\in N(\boldsymbol{\theta}^*)} P\left(y(\boldsymbol{\theta})\le U_{u,v}\right)\right)}{\log(\sigma)}} = O\left(P\left(\hat{\boldsymbol{\theta}}_k \neq \boldsymbol{\theta}^*\right)\right),$$

which means $P\left(\hat{\boldsymbol{\theta}}_k \neq \boldsymbol{\theta}^*\right)$ goes to 0 at a rate not faster than the rate at which

$$k^{-\frac{c\log\left(1/\max_{\boldsymbol{\theta}\in N(\boldsymbol{\theta}^*)} P\left(y(\boldsymbol{\theta})\le U_{u,v}\right)\right)}{\log(\sigma)}}$$

goes to 0.

*Remarks:*

1. For the SR algorithm, $u$ and $v$ can be picked properly to make the assumption 0 < $P\left(y(\boldsymbol{\theta}) \le U_{u,v}\right)$ < 1 to be true (shown in Theorem 3.1 of Yan and Mukai, 1992).

2. The value of $r$ in condition (iii) is defined in Section 5.2.2.

*Proof.* For the SR algorithm and any $\boldsymbol{\theta}$, $\boldsymbol{\theta}' \in \Theta$, $p\left(\boldsymbol{\theta},\boldsymbol{\theta}'\right) = P\left(y(\boldsymbol{\theta}') \le U_{u,v}\right)$, which does not rely on the value of $\boldsymbol{\theta}$. From the Theorem 7.2 of Yan and Mukai (1992), by condition (i), (ii) and (iii) here, we have that the sequence $\left\{\hat{\boldsymbol{\theta}}_k\right\}$ generated by the SR algorithm converges in probability

$$\lim_{k\to\infty} P\left(\hat{\boldsymbol{\theta}}_k = \boldsymbol{\theta}^*\right) = 1,$$

which indicates that the condition (i) in Theorem 5.1 is true here.

Moreover, since $0 < P(y(\boldsymbol{\theta}) < U_{u,v}) < 1$, then $\sigma \geq 1/\min_{\boldsymbol{\theta}\in\Theta} P(y(\boldsymbol{\theta}) < U_{u,v}) > 1$, then the condition (iii) here implies that the condition (iii) in Theorem 5.1 is true here. Then all conditions in Theorem 5.1 are true here, and consequently we have

$$k^{-\frac{c\log\left(1/\max_{\boldsymbol{\theta}\in N(\boldsymbol{\theta}^*)} P\left(y(\boldsymbol{\theta})\leq U_{u,v}\right)\right)}{\log(\sigma)}} = O\left(P\left(\hat{\boldsymbol{\theta}}_k \neq \boldsymbol{\theta}^*\right)\right),$$

which means $P\left(\hat{\boldsymbol{\theta}}_k \neq \boldsymbol{\theta}^*\right)$ goes to 0 at a rate not faster than the rate at which

$$k^{-\frac{c\log\left(1/\max_{\boldsymbol{\theta}\in N(\boldsymbol{\theta}^*)} P\left(y(\boldsymbol{\theta})\leq U_{u,v}\right)\right)}{\log(\sigma)}}$$

goes to 0. We have shown the result on the convergence rate of $P\left(\hat{\boldsymbol{\theta}}_k \neq \boldsymbol{\theta}^*\right)$ for the SR algorithm. Q.E.D.

**Corollary 5.2.** For the SC algorithm, if (i) $\Theta^* = \{\boldsymbol{\theta}^*\}$; (ii) for all $\boldsymbol{\theta} \in \Theta$ and $\boldsymbol{\theta}' \in N(\boldsymbol{\theta}) = \Theta \setminus \boldsymbol{\theta}$, we have that $\sum_{\boldsymbol{\theta}'\in N(\boldsymbol{\theta})} R(\boldsymbol{\theta},\boldsymbol{\theta}') = 1$, $R(\boldsymbol{\theta},\boldsymbol{\theta}') > 0$, and $0 < P\left(y(\boldsymbol{\theta}') < y(\boldsymbol{\theta})\right) < 1$; (iii) $M_k = \lfloor c\log_\sigma(1+k_0+k)\rfloor$, $0 < c \leq 1$, $\sigma \geq 1/\min_{\boldsymbol{\theta}\in\Theta,\boldsymbol{\theta}'\in N(\boldsymbol{\theta})} P\left(y(\boldsymbol{\theta}') < y(\boldsymbol{\theta})\right)$, and $c\log_\sigma(1+k_0) \geq 1$. Then

$$k^{-\frac{c\log\left(1/\max_{\boldsymbol{\theta}\in N(\boldsymbol{\theta}^*)} P(y(\boldsymbol{\theta})<y(\boldsymbol{\theta}^*))\right)}{\log(\sigma)}} = O\left(P\left(\hat{\boldsymbol{\theta}}_k \neq \boldsymbol{\theta}^*\right)\right),$$

which means $P\left(\hat{\boldsymbol{\theta}}_k \neq \boldsymbol{\theta}^*\right)$ goes to 0 at a rate not faster than the rate at which

$$k^{-\frac{c\log\left(1/\max_{\boldsymbol{\theta}\in N(\boldsymbol{\theta}^*)} P(y(\boldsymbol{\theta})<y(\boldsymbol{\theta}^*))\right)}{\log(\sigma)}}$$

goes to 0.

*Remark:* For Gaussian noise case, we always have 0 < $P\big(y(\boldsymbol{\theta}')<y(\boldsymbol{\theta})\big)$ < 1 for all $\boldsymbol{\theta},\boldsymbol{\theta}'\in\Theta$, $\boldsymbol{\theta}\neq\boldsymbol{\theta}'$.

*Proof.* For the SC algorithm and any $\boldsymbol{\theta}$, $\boldsymbol{\theta}'\in\Theta$, $p\big(\boldsymbol{\theta},\boldsymbol{\theta}'\big) = P\big(y(\boldsymbol{\theta}')<y(\boldsymbol{\theta})\big)$. From the Theorem 5.1 of Gong et al. (1999), by condition (i), (ii) and (iii) here, we have that the sequence $\{\hat{\boldsymbol{\theta}}_k\}$ generated by the algorithm converges in probability

$$\lim_{k\to\infty} P\left(\hat{\boldsymbol{\theta}}_k=\boldsymbol{\theta}^*\right)=1,$$

which indicates that the condition (i) in Theorem 5.1 is true here. Moreover, since 0 < $P\big(y(\boldsymbol{\theta}')<y(\boldsymbol{\theta})\big)$< 1, then we have that σ ≥ $1/\min_{\boldsymbol{\theta}\in\Theta,\boldsymbol{\theta}'\in N(\boldsymbol{\theta})} P\big(y(\boldsymbol{\theta}')<y(\boldsymbol{\theta})\big)$ > 1, then the condition (iii) here implies that the conditions (iii) in Theorem 5.1 is true. Then all conditions in Theorem 5.1 are true here, and consequently we have

$$k^{-\frac{c\log\left(1/\max_{\boldsymbol{\theta}\in N(\boldsymbol{\theta}^*)} P(y(\boldsymbol{\theta})<y(\boldsymbol{\theta}^*))\right)}{\log(\sigma)}} = O\left(P\left(\hat{\boldsymbol{\theta}}_k\neq\boldsymbol{\theta}^*\right)\right),$$

which means $P\left(\hat{\boldsymbol{\theta}}_k\neq\boldsymbol{\theta}^*\right)$ goes to 0 at a rate not faster than the rate at which

$$k^{-\frac{c\log\left(1/\max_{\boldsymbol{\theta}\in N(\boldsymbol{\theta}^*)} P(y(\boldsymbol{\theta})<y(\boldsymbol{\theta}^*))\right)}{\log(\sigma)}}$$

goes to 0. We have shown the result on the convergence rate of $P\left(\hat{\boldsymbol{\theta}}_k \neq \boldsymbol{\theta}^*\right)$ for the SC algorithm. Q.E.D.

Let

$$\gamma_{SR} = \frac{c\log\left(1\big/\max_{\boldsymbol{\theta}\in N(\boldsymbol{\theta}^*)} P\left(y(\boldsymbol{\theta}) \leq U_{u,v}\right)\right)}{\log(\sigma)},$$

$$\gamma_{SC} = \frac{c\log\left(1\big/\max_{\boldsymbol{\theta}\in N(\boldsymbol{\theta}^*)} P(y(\boldsymbol{\theta}) < y(\boldsymbol{\theta}^*))\right)}{\log(\sigma)},$$

and now we discuss some properties of the exponent $\gamma_{SR}$ and $\gamma_{SC}$. For the SR algorithm, we have $\sigma \geq 1\big/\min_{\boldsymbol{\theta}\in\Theta} P(y(\boldsymbol{\theta}) < U_{u,v})$ and $0 < c \leq 1/r$, so the exponent can be bounded as

$$0 < \gamma_{SR} \leq \frac{c\log\left(1\big/\max_{\boldsymbol{\theta}\in N(\boldsymbol{\theta}^*)} P\left(y(\boldsymbol{\theta}) \leq U_{u,v}\right)\right)}{\log\left(1\big/\min_{\boldsymbol{\theta}\in\Theta} P\left(y(\boldsymbol{\theta}) \leq U_{u,v}\right)\right)}$$

Since $\min_{\boldsymbol{\theta}\in\Theta} P\left(y(\boldsymbol{\theta}) \leq U_{u,v}\right) \leq \max_{\boldsymbol{\theta}\in N(\boldsymbol{\theta}^*)} P\left(y(\boldsymbol{\theta}) \leq U_{u,v}\right)$, and $0 < c \leq 1/r \leq 1$, then we have that the range of the exponent $\gamma_{SR}$ is (0,1]. Similarly for the SC algorithm, we have $\sigma \geq 1\big/\min_{\boldsymbol{\theta}\in\Theta,\boldsymbol{\theta}'\in\Theta\backslash\boldsymbol{\theta}} P\left(y(\boldsymbol{\theta}') < y(\boldsymbol{\theta})\right)$ and $0 < c \leq 1$, then the exponent has lower and upper bounds as

$$0 < \gamma_{SC} \leq \frac{c\log\left(1\big/\max_{\boldsymbol{\theta}\in N(\boldsymbol{\theta}^*)} P(y(\boldsymbol{\theta}) < y(\boldsymbol{\theta}^*))\right)}{\log\left(1\big/\min_{\boldsymbol{\theta}\in\Theta,\boldsymbol{\theta}'\in N(\boldsymbol{\theta})} P\left(y(\boldsymbol{\theta}') < y(\boldsymbol{\theta})\right)\right)} \leq 1.$$

Therefore, we have that the range of the exponent $\gamma_{SC}$ is (0,1].

From the discussion in Section 3.1.2, we have $P([\hat{\boldsymbol{\theta}}_k] \neq \boldsymbol{\theta}^*) = O(k^{-\alpha})$. Overall, we have the comparative rates of convergence for DSPSA, SR, and SC summarized in Table 5.1. For SR and SC, we have $[\hat{\boldsymbol{\theta}}_k] = \hat{\boldsymbol{\theta}}_k$. In Table 5.1, we have $P([\hat{\boldsymbol{\theta}}_k] \neq \boldsymbol{\theta}^*) = O(k^{-\alpha})$, $0.5 < \alpha \leq 1$ for DSPSA, $k^{-\gamma_{SR}} = O(P([\hat{\boldsymbol{\theta}}_k] \neq \boldsymbol{\theta}^*))$, $0 < \gamma_{SR} \leq 1$ for SR, and $k^{-\gamma_{SC}} = O(P([\hat{\boldsymbol{\theta}}_k] \neq \boldsymbol{\theta}^*))$, $0 < \gamma_{SC} \leq 1$ for SC, which indicates that SR and SC cannot achieve higher rate of convergence than DSPSA.

**Table 5.1** Analysis of rates of convergence to 0 for $P([\hat{\boldsymbol{\theta}}_k] \neq \boldsymbol{\theta}^*)$ in DSPSA, stochastic ruler, and stochastic comparison. For SR and SC, we have $[\hat{\boldsymbol{\theta}}_k] = \hat{\boldsymbol{\theta}}_k$.

| Method Name | Analysis of Rate of Convergence of $P([\hat{\boldsymbol{\theta}}_k] \neq \boldsymbol{\theta}^*)$ |
|---|---|
| DSPSA | $P([\hat{\boldsymbol{\theta}}_k] \neq \boldsymbol{\theta}^*) = O(k^{-\alpha})$<br>$0.5 < \alpha \leq 1$ |
| SR | $k^{-\gamma_{SR}} = O(P([\hat{\boldsymbol{\theta}}_k] \neq \boldsymbol{\theta}^*))$, $\gamma_{SR} = \dfrac{c \log\left(1 \big/ \max_{\boldsymbol{\theta} \in N(\boldsymbol{\theta}^*)} P\left(y(\boldsymbol{\theta}) \leq U_{u,v}\right)\right)}{\log(\sigma)}$<br>$\sigma \geq 1 \big/ \min_{\boldsymbol{\theta} \in \Theta} P(y(\boldsymbol{\theta}) < U_{u,v})$, $0 < c \leq 1/r$,<br>$0 < \gamma_{SR} \leq 1$ |
| SC | $k^{-\gamma_{SC}} = O(P([\hat{\boldsymbol{\theta}}_k] \neq \boldsymbol{\theta}^*))$, $\gamma_{SC} = \dfrac{c \log\left(1 \big/ \max_{\boldsymbol{\theta} \in N(\boldsymbol{\theta}^*)} P(y(\boldsymbol{\theta}) < y(\boldsymbol{\theta}^*))\right)}{\log(\sigma)}$<br>$\sigma \geq 1 \big/ \min_{\boldsymbol{\theta} \in \Theta, \boldsymbol{\theta}' \in \Theta \backslash \boldsymbol{\theta}} P\left(y(\boldsymbol{\theta}') < y(\boldsymbol{\theta})\right)$, $0 < c \leq 1$,<br>$0 < \gamma_{SC} \leq 1$ |

The forms of the rates of convergence of the SR algorithm and the SC algorithm in Table 5.1 seem complicated, so now let us use two very simple examples to see the rates clearly. Let us consider two-dimensional separable quadratic function $t_1^2 + t_2^2$ in the domain of $\{-1, 0, 1\}^2$ and suppose the noises $\varepsilon$ are i.i.d. $N(0,1)$. For the SR algorithm, based on the function structure and our experience, we let $u = 0$ and $v = 2$. Then, we have that $\max_{\boldsymbol{\theta} \in N(\boldsymbol{\theta}^*)} P\left(y(\boldsymbol{\theta}) \le U_{u,v}\right) = P(y(1,0) \le U_{0,2}) \approx 0.5$, and $\sigma \ge 1 \big/ \min_{\boldsymbol{\theta} \in \Theta} P\left(y(\boldsymbol{\theta}) \le U_{u,v}\right) = 1 \big/ P(y(1,1) \le U_{0,2}) \approx 5.13$. Moreover, by the definition of $r$ in Section 5.2.2, we have $r = 1$, then $0 < c \le 1/r = 1$. Thus, the exponent $\gamma_{SR}$ is smaller than or equal to 0.4244. Similarly, for the SC algorithm, $\max_{\boldsymbol{\theta} \in N(\boldsymbol{\theta}^*)} P\left(y(\boldsymbol{\theta}) < y(\boldsymbol{\theta}^*)\right) = P(y(1,0) \le y(0,0)) \approx 0.24$, $\sigma \ge 1 \big/ \min_{\boldsymbol{\theta} \in \Theta, \boldsymbol{\theta}' \in N(\boldsymbol{\theta})} P(y(\boldsymbol{\theta}') < y(\boldsymbol{\theta})) = 1 \big/ P(y(1,1) \le y(0,0)) \approx 12.82$. Thus, the exponent $\gamma_{SC}$ is smaller than or equal to 0.561. We see that for this simple loss function, the values of the exponent for both SR and SC can not achieve the upper bound 1.

The other loss function is a two-dimensional function: $\max(|t_1|, |t_2|)$ defined in the domain of $\{-1, 0, 1\}^2$ and suppose the noises $\varepsilon$ are i.i.d. $N(0,1)$. The loss function value equals 0 at the point (0, 0) and equals 1 at all the other points. For the SR algorithm, $\max_{\boldsymbol{\theta} \in N(\boldsymbol{\theta}^*)} P\left(y(\boldsymbol{\theta}) \le U_{u,v}\right) = \min_{\boldsymbol{\theta} \in \Theta} P\left(y(\boldsymbol{\theta}) \le U_{u,v}\right)$ for this simple loss function. Moreover, by the definition of $r$ in Section 5.2.2, we have $r = 1$, then $0 < c \le 1/r = 1$. Thus, the exponent $\gamma_{SR}$ can achieve the upper bound 1. Similarly for the SC algorithm, we have that $\max_{\boldsymbol{\theta} \in N(\boldsymbol{\theta}^*)} P\left(y(\boldsymbol{\theta}) < y(\boldsymbol{\theta}^*)\right) = \min_{\boldsymbol{\theta} \in \Theta, \boldsymbol{\theta}' \in N(\boldsymbol{\theta})} P\left(y(\boldsymbol{\theta}') < y(\boldsymbol{\theta})\right)$ for this simple loss function. Thus, the exponent $\gamma_{SC}$ can achieve the upper bound 1. We see that for this simple loss function, the values of the exponent for both SR and SC can achieve the upper bound 1.

Moreover, for the SR algorithm, it has other properties as following: 1) If the feasible region is very large, then we have a large value for $r$, which leads to a small value for $c$. 2) The choice of $\sigma$ is heavily based on the structure of the loss function. 3) The choices of $u$ and $v$ are not trivial. These coefficients directly affect the rate of convergence of the algorithm, but it is hard to control them. For the SC algorithm, we have similar arguments as the SR algorithm. The choice of $\sigma$ is still heavily based on the structure of the loss function. It is also hard to figure out the value of $\sigma$, since the structure of the loss function may not be clear. Compared with the SR algorithm and the SC algorithm, the rate of convergence of DSPSA in the big-$O$ sense is just based on the choice of $\alpha$, which is independent on the structure of the loss function. In addition, as we have discussed, for $\gamma_{SR}$ and $\gamma_{SC}$, we have $0 < \gamma_{SR} \leq 1$ and $0 < \gamma_{SC} \leq 1$, and for $\alpha$, we have $0.5 < \alpha \leq 1$. Therefore, from Table 5.1, we know that SR and SC can not achieve higher rate of convergence than DSPSA in the big-$O$ sense. But the real performance of DSPSA is still based on the structure of the loss function, since the constant multiplier of $k^{-\alpha}$ is different for various kinds of loss functions.

# Chapter 6

# Numerical Comparison of DSPSA and Random Search Algorithms

In Chapter 5, we have introduced the SR algorithm and the SC algorithm and discussed some convergence properties of both algorithms. We have also compared the rate of convergence of DSPSA, SR and SC theoretically in the big-$O$ sense. In this chapter, we present numerical experiments on the comparisons of DSPSA, SR and SC. First we discuss the choice of coefficients in the maximum number of comparisons in each iteration for both the SR algorithm and the SC algorithm. Then, we consider the comparisons of the three algorithms for the separable loss function and the skewed quartic loss function.

# 6.1 Choice of Coefficients in the Maximum Number of Comparisons for Random Search Algorithms

In both Yan and Mukai (1992) and Gong et al. (1999), the authors do not discuss the choice of coefficients in $M_k$ (maximum number of comparisons in the *k*th iteration) for the SR algorithm and the SC algorithm. Therefore, before running the SR algorithm and the SC algorithm, we need to determine the values of the coefficients of $M_k$ for both algorithms. For the SR algorithm, besides the coefficients, we also need to pick the lower and upper bounds of the uniform distribution for the stochastic ruler.

In both algorithms $M_k = \lfloor c \log_{\sigma}(1 + k_0 + k) \rfloor$. In order to set the value for $M_k$, we need to determine the values of the coefficients *c*, σ and $k_0$. In Yan and Mukai (1992), the authors require that in the SR algorithm $\sigma \geq 1 / \min_{\boldsymbol{\theta} \in \Theta} P(y(\boldsymbol{\theta}) < U_{u,v})$, $0 < c \leq 1/r$, and $c \log_{\sigma}(1 + k_0) \geq 1$. In Gong et al. (1999), the authors require that in the SC algorithm $\sigma \geq 1 / \min_{\boldsymbol{\theta} \in \Theta, \boldsymbol{\theta}' \in \Theta \backslash \boldsymbol{\theta}} P\left(y(\boldsymbol{\theta}') < y(\boldsymbol{\theta})\right)$, $0 < c \leq 1$, $c \log_{\sigma}(1 + k_0) \geq 1$. These constraints are sufficient (not necessary) conditions that are mainly used for the proof of convergence of the algorithms, and we call these sufficient conditions as the "coefficients constraints." In the numerical tests of Gong et al. (1999), the authors do not obey the coefficients constraints. In the following, we discuss both the coefficients set that satisfies the constraints and the coefficients set that does not satisfy the constraints. Later we will see that after relaxing the coefficients constraints, we can yield better performance. Since we want to compare the SR algorithm and SC algorithm with DSPSA numerically, we try to find the coefficients sets that are as good as possible for the SR algorithm and the SC algorithm. In Yan and Mukai (1992) and Gong et al. (1999) the authors do not provide any guidelines on the coefficients selection, so it is very hard to find the optimal coefficients. Therefore, we just try our best to pick some reasonable sets of coefficients for SR and SC to do the comparisons.

First we use a simple numerical example discussed in Section 5.4 to see how the performance of both SR and SC can be improved by relaxing the coefficients constraints. The loss function is a two-dimensional separable function $L(\boldsymbol{\theta}) = t_1^2 + t_2^2$ in the domain of $\{-1, 0, 1\}^2$. Suppose the noises $\varepsilon$ are i.i.d. $N(0,1)$. Here for both SR and SC algorithms, the neighborhood of $\boldsymbol{\theta}$ is $N(\boldsymbol{\theta}) = \Theta \setminus \boldsymbol{\theta}$.

For the SR algorithm, based on the loss function, we set $u$ = 0 and $v$ = 2. Under the coefficients constraints, we require that $\sigma \geq 1/\min_{\boldsymbol{\theta}\in\Theta} P(\boldsymbol{\theta},u,v) = 1/P(y(1,1) \leq U_{0,2}) \approx$ 5.13. We also need to have $0 < c \leq 1/r$, $c\log_{\sigma}(1+k_0) \geq 1$. Here by the definition of $r$ in Section 5.2.2, we have the radius of the domain $r$ = 1, so $0 < c \leq 1$. The initial guess is set as $\hat{\boldsymbol{\theta}}_0$= $[1, 1]^T$ and the optimal solution is $[0, 0]^T$. The number of replicates is 20, and the number of noisy measurements in each replicate is 20,000. The results of the tuning process are shown in Table 6.1. We list different sets of coefficients that satisfy the coefficients constraints and the final errors are shown in terms of both sample mean of $\left\|\hat{\boldsymbol{\theta}}_q - \boldsymbol{\theta}^*\right\| / \left\|\hat{\boldsymbol{\theta}}_0 - \boldsymbol{\theta}^*\right\|$ and sample mean of $\left|L(\hat{\boldsymbol{\theta}}_q) - L(\boldsymbol{\theta}^*)\right| / \left|L(\hat{\boldsymbol{\theta}}_0) - L(\boldsymbol{\theta}^*)\right|$, where the sample mean is the arithmetic mean of the observed values across 20 replicates and $q$ is the number of iterations when the SR algorithm hits 20,000 noisy measurements in each replicate ($q$ is a random variable). Overall, the final errors are not 0 for all sets of coefficients, and the first coefficients set provides the smallest final error in terms of both criteria among all sets in Table 6.1.

**Table 6.1** Performance of the SR algorithm with different sets of coefficients that satisfy the coefficients constraints for the two-dimensional separable loss function. The first set of coefficients provides the smallest final error in terms of both sample mean of $\left\|\hat{\boldsymbol{\theta}}_q - \boldsymbol{\theta}^*\right\| / \left\|\hat{\boldsymbol{\theta}}_0 - \boldsymbol{\theta}^*\right\|$ and sample mean of $\left|L(\hat{\boldsymbol{\theta}}_q) - L(\boldsymbol{\theta}^*)\right| / \left|L(\hat{\boldsymbol{\theta}}_0) - L(\boldsymbol{\theta}^*)\right|$ among all sets in the table. The number of replicates is 20 and the number of noisy measurements of loss function in each replicate is 20,000.

| $c$ | $\sigma$ | $k_0$ | Sample Mean of $\dfrac{\left\|\hat{\boldsymbol{\theta}}_q - \boldsymbol{\theta}^*\right\|}{\left\|\hat{\boldsymbol{\theta}}_0 - \boldsymbol{\theta}^*\right\|}$ | Sample Mean of $\dfrac{\left\|L(\hat{\boldsymbol{\theta}}_q) - L(\boldsymbol{\theta}^*)\right\|}{\left\|L(\hat{\boldsymbol{\theta}}_0) - L(\boldsymbol{\theta}^*)\right\|}$ |
|---|---|---|---|---|
| 1 | 5.13 | 100 | 0.1061 | 0.075 |
| 0.8 | 5.13 | 100 | 0.2475 | 0.175 |
| 1 | 10 | 100 | 0.3889 | 0.275 |
| 1 | 5.13 | 10 | 0.2475 | 0.175 |
| 1 | 5.13 | 1000 | 0.3182 | 0.225 |

Now we start to consider the case where we may violate the coefficients constraints, $\sigma \geq 1/\min_{\boldsymbol{\theta}\in\Theta} P(y(\boldsymbol{\theta}) < U_{u,v})$ and $0 < c \leq 1/r$ (the other constraint $c\log_\sigma(1+k_0) \geq 1$, cannot be relaxed, because this constraint makes $M_k$ to be at least 1 in each iteration). First we relax the constraints, $0 < c \leq 1/r$ and $\sigma \geq 1/\min_{\boldsymbol{\theta}\in\Theta} P(y(\boldsymbol{\theta}) < U_{u,v})$ into the constraints, $c > 0$ and $\sigma > 1$. We call the constraints, $\sigma > 1$, $c > 0$, and $c\log_\sigma(1+k_0) \geq 1$ as the "relaxed coefficients constraints." It is clear that the relaxed coefficients constraints are weaker than the coefficients constraints, which means the coefficients that satisfy the coefficients constraints also satisfy the relaxed coefficient constraints. Since $M_k = \lfloor c\log_\sigma(1+k_0+k) \rfloor$, we have

$$M_k = \left\lfloor \frac{c}{\log(\sigma)} \log(1+k_0+k) \right\rfloor.$$

We see that both $c$ and $\sigma$ affect the rate at which $M_k$ increases and $k_0$ affects the value of $M_k$ in the early iterations. Many sets of $c$ and $\sigma$ can provide the same value of $c/\log(\sigma)$. For example, we see both $c=1$, $\sigma=4$ and $c=0.5$, $\sigma=2$ generate the same value for $c/\log(\sigma) = 1/\log 4$. Thus, we can set $c=1$ and only change the value of $\sigma$. In Table 6.2, we try different sets of coefficients that satisfy the relaxed coefficients constraints. From Table 6.2, we see that when we relax the coefficients constraints, the performance of the SR algorithm has been improved significantly. The final errors in terms of both sample mean of $\left\|\hat{\boldsymbol{\theta}}_q - \boldsymbol{\theta}^*\right\| / \left\|\hat{\boldsymbol{\theta}}_0 - \boldsymbol{\theta}^*\right\|$ and sample mean of $\left|L(\hat{\boldsymbol{\theta}}_q) - L(\boldsymbol{\theta}^*)\right| / \left|L(\hat{\boldsymbol{\theta}}_0) - L(\boldsymbol{\theta}^*)\right|$ are 0, where the sample mean is the arithmetic mean of the observed values across 20 replicates

and $q$ is the number of iterations when the SR algorithm hits 20,000 noisy measurements in each replicate.

**Table 6.2** Performance of the SR algorithm with different sets of coefficients that satisfy the relaxed coefficients constraints for the two-dimensional separable loss function. When we relax the coefficients constraints, the performance of the SR algorithm has been improved significantly. The number of replicates is 20 and the number of noisy measurements of the loss function in each replicate is 20,000.

| $c$ | $\sigma$ | $k_0$ | Sample Mean of $\frac{\left\|\hat{\boldsymbol{\theta}}_q - \boldsymbol{\theta}^*\right\|}{\left\|\hat{\boldsymbol{\theta}}_0 - \boldsymbol{\theta}^*\right\|}$ | Sample mean of $\frac{\left|L(\hat{\boldsymbol{\theta}}_q) - L(\boldsymbol{\theta}^*)\right|}{\left|L(\hat{\boldsymbol{\theta}}_0) - L(\boldsymbol{\theta}^*)\right|}$ |
|---|---|---|---|---|
| 1 | 1.5 | 100 | 0 | 0 |
| 1 | 2 | 100 | 0 | 0 |
| 1 | 1.5 | 10 | 0 | 0 |
| 1 | 1.5 | 1000 | 0 | 0 |

We pick the first set of coefficients in both Tables 6.1 and 6.2 and compare the performances of SR based on these two sets of coefficients. In Figure 6.1, we see that the performance based on the coefficients from Table 6.2 provides much better performance than that from Table 6.1, which indicates that relaxed coefficient constraints leads to significant better convergent performance. Compared with the coefficients constraints, the relaxed coefficients constraints allow us to pick smaller value for $\sigma$, which indicates that $M_k$ can increase faster. Thus, the relaxed constraints can handle the noise better for this two-dimensional separable loss function.

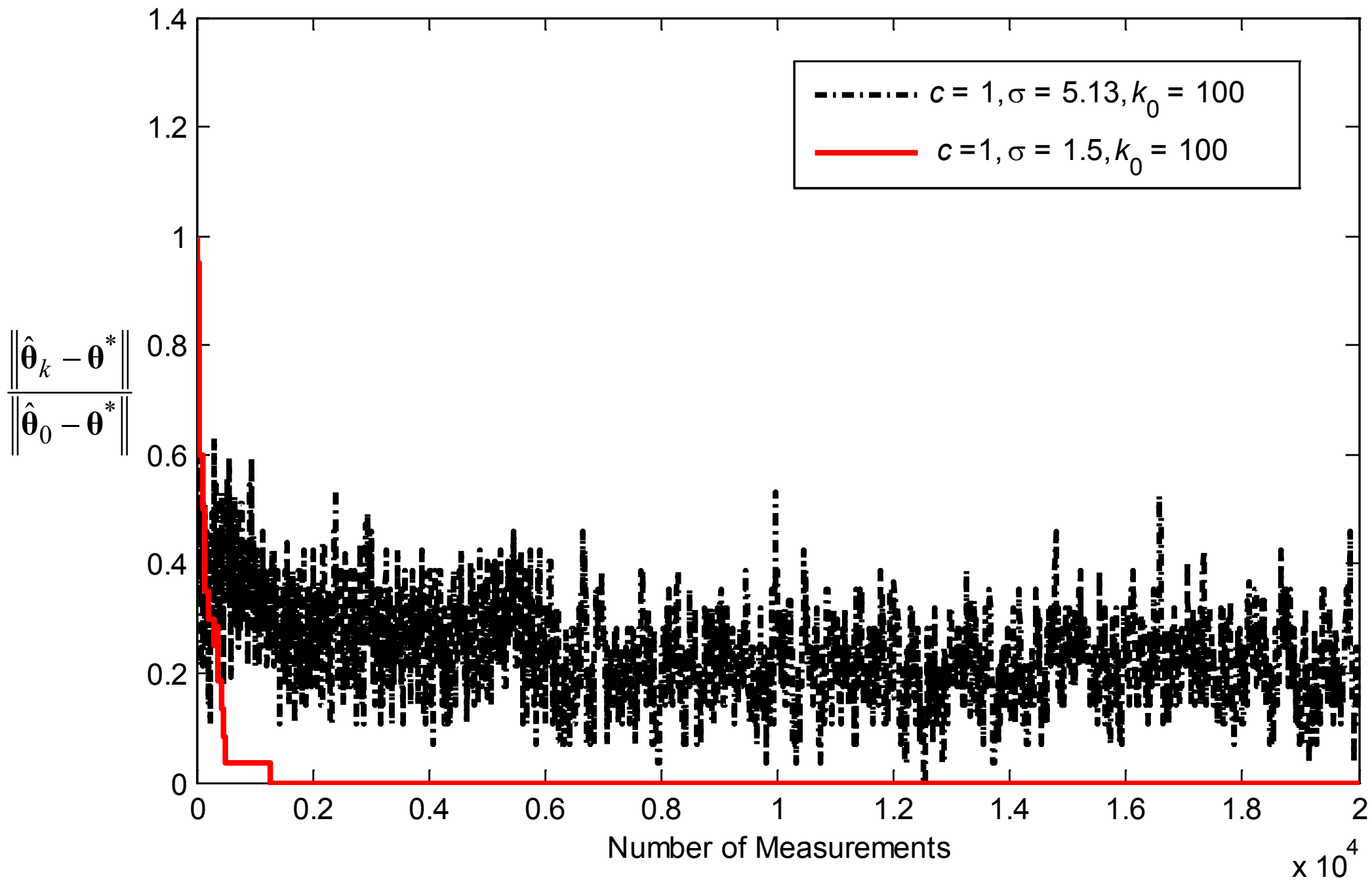

**Figure 6.1** Comparison for the SR algorithm between the best set of coefficients in Table 6.1 that satisfies the coefficients constraints (dashed line) and the set in Table 6.2 that satisfies the relaxed coefficients constraints (solid line). After relaxing the coefficients constraints, the coefficients set provides better performance for the two-dimensional separable loss function. Each curve represents the sample mean of 20 independent replicates.

For the SC algorithm, the sufficient conditions are $\sigma \geq 1/\min_{\theta\in\Theta,\theta'\in\Theta\backslash\theta} P(y(\theta') < y(\theta))$ $= 1/P(y(1,1) \leq y(0,0)) \approx 12.82$, $0 < c \leq 1$, $c\log_\sigma(1+k_0) \geq 1$. The initial guess is set as $\hat{\boldsymbol{\theta}}_0 = [1, 1]^T$. The number of replicates is 20, and the number of noisy measurements in each replicate is 20,000. We do the tuning process in Table 6.3. We list different sets of coefficients that satisfy the coefficients constraints and show the final errors in terms of both the sample mean of $\left\|\hat{\boldsymbol{\theta}}_q - \boldsymbol{\theta}^*\right\| / \left\|\hat{\boldsymbol{\theta}}_0 - \boldsymbol{\theta}^*\right\|$ and the sample mean of $\left|L(\hat{\boldsymbol{\theta}}_q) - L(\boldsymbol{\theta}^*)\right| / \left|L(\hat{\boldsymbol{\theta}}_0) - L(\boldsymbol{\theta}^*)\right|$, where the sample mean is the arithmetic mean of the observed values from 20 replicates and $q$ is the number of iterations when the SC algorithm hits 20,000 noisy measurements in each replicate. Overall, the final errors in terms of both criteria above are not 0 for all sets, and the first set provides the smallest error in terms of both criteria.

**Table 6.3** Performance of the SC algorithm with different sets of coefficients that satisfy the coefficients constraints for the two-dimensional separable loss function. The first set of coefficients provides the smallest final error among all sets in the table in terms of both sample mean of $\left\|\hat{\boldsymbol{\theta}}_q - \boldsymbol{\theta}^*\right\| \Big/ \left\|\hat{\boldsymbol{\theta}}_0 - \boldsymbol{\theta}^*\right\|$ and sample mean of $\left|L(\hat{\boldsymbol{\theta}}_q) - L(\boldsymbol{\theta}^*)\right| \Big/ \left|L(\hat{\boldsymbol{\theta}}_0) - L(\boldsymbol{\theta}^*)\right|$. The number of replicates is 20 and the number of noisy measurements of the loss function in each replicate is 20,000.

| $c$ | $\sigma$ | $k_0$ | Sample Mean of $\dfrac{\left\|\hat{\boldsymbol{\theta}}_q - \boldsymbol{\theta}^*\right\|}{\left\|\hat{\boldsymbol{\theta}}_0 - \boldsymbol{\theta}^*\right\|}$ | Sample mean of $\dfrac{\left|L(\hat{\boldsymbol{\theta}}_q) - L(\boldsymbol{\theta}^*)\right|}{\left|L(\hat{\boldsymbol{\theta}}_0) - L(\boldsymbol{\theta}^*)\right|}$ |
|---|---|---|---|---|
| 1 | 13 | 2000 | 0.0354 | 0.025 |
| 0.8 | 13 | 2000 | 0.2061 | 0.175 |
| 1 | 20 | 2000 | 0.1061 | 0.075 |
| 1 | 13 | 200 | 0.1768 | 0.125 |
| 1 | 13 | 20000 | 0.1061 | 0.075 |

Similarly, we discuss the case when we relax the coefficients constraints. We call the constraints: $\sigma > 1$, $c > 0$, and $c\log_{\sigma}(1+k_0) \geq 1$ as the "relaxed coefficients constraints." By the same arguments as the SR algorithm, we keep the condition $c\log_{\sigma}(1+k_0) \geq 1$, set $c = 1$, and only change the value of $\sigma$. In Table 6.4, we do the tuning process. We list different sets of coefficients that satisfy the relaxed coefficients constraints, and show the final errors in terms of both the sample mean of $\left\|\hat{\boldsymbol{\theta}}_q - \boldsymbol{\theta}^*\right\| \Big/ \left\|\hat{\boldsymbol{\theta}}_0 - \boldsymbol{\theta}^*\right\|$ and the sample mean

of $\left| L(\hat{\boldsymbol{\theta}}_q) - L(\boldsymbol{\theta}^*) \right| \big/ \left| L(\hat{\boldsymbol{\theta}}_0) - L(\boldsymbol{\theta}^*) \right|$, where the sample mean is the arithmetic mean of the observed values across 20 replicates and $q$ is the number of iterations when the SC algorithm hits 20,000 noisy measurements in each replicate. From Table 6.4, we see that when we relax the coefficients constraints, the performance of the SC algorithm has been improved significantly. The final errors are all 0 in terms of both criteria.

**Table 6.4** Performance of the SC algorithm with different sets of coefficients after relaxing the coefficients constraints for the two-dimensional separable loss function. When we relax the sufficient conditions, the performance of the SC algorithm has been improved significantly. The final errors are all 0 in terms of both criteria. The number of replicates is 20 and the number of noisy measurements of the loss function in each replicate is 20,000.

| $c$ | $\sigma$ | $k_0$ | Sample mean of $\dfrac{\left\| \hat{\boldsymbol{\theta}}_q - \boldsymbol{\theta}^* \right\|}{\left\| \hat{\boldsymbol{\theta}}_0 - \boldsymbol{\theta}^* \right\|}$ | Sample mean of $\dfrac{\left| L(\hat{\boldsymbol{\theta}}_q) - L(\boldsymbol{\theta}^*) \right|}{\left| L(\hat{\boldsymbol{\theta}}_0) - L(\boldsymbol{\theta}^*) \right|}$ |
|---|---|---|---|---|
| 1 | 2 | 10 | 0 | 0 |
| 1 | 3 | 10 | 0 | 0 |
| 1 | 2 | 2 | 0 | 0 |
| 1 | 2 | 100 | 0 | 0 |

We pick the first set of coefficients in both Tables 6.3 and 6.4 and compare the performance of SC based on these two sets of coefficients. In Figure 6.2, we see that the performance based on the coefficients from Table 6.4 is much better than that from Table 6.3. Compared with the coefficients constraints, the relaxed coefficients constraints allow us to pick smaller value for $\sigma$, which indicates that $M_k$ can increase faster. Thus, we can handle the noise better.

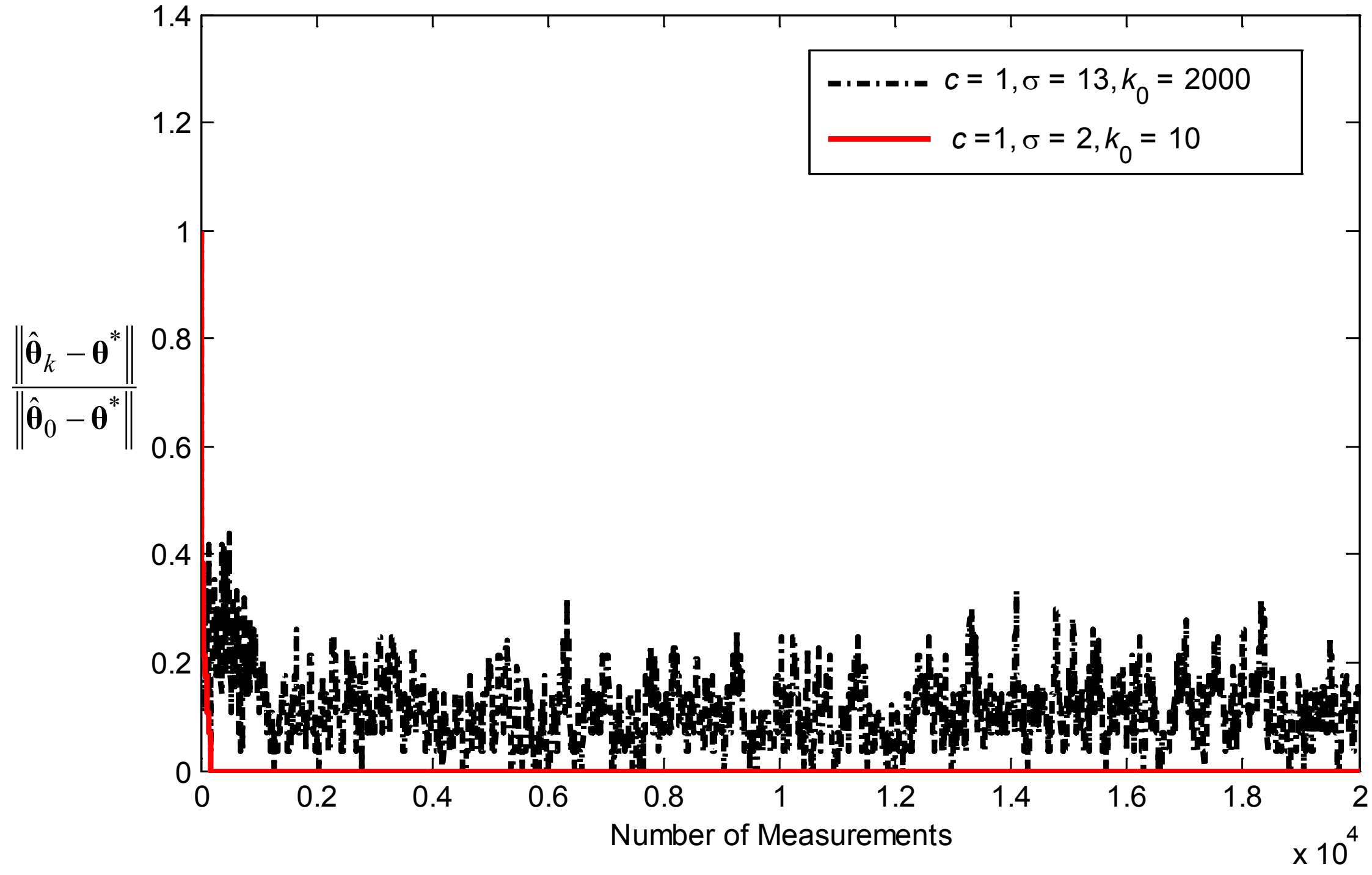

**Figure 6.2** Comparison for the SC algorithm between the best set of coefficients in Table 6.3 that satisfies the coefficients constraints (dash line) and the set in Table 6.4 that satisfies the relaxed coefficients constraints (solid line). After relaxing the coefficients constraints, the coefficients set provides better performance for the two-dimensional separable loss function. Each curve represents the sample mean of 20 independent replicates.

Overall, when we relax the coefficients constraints, the performances of both algorithms (SR and SC) have been improved significantly for this two-dimensional separable loss function. In Table 6.5, we summarize the coefficients constraints and the relaxed coefficients constraints. As we have discussed, the relaxed coefficients constraints are weaker than the coefficients constraints. Based on the analysis above on the effects of $\sigma$ and $c$ on $M_k$, we can always make $c$ = 1, and only change the value of $\sigma$.

**Table 6.5** The coefficients constraints and relaxed coefficients constraints for both SR and SC. We see that the relaxed coefficients constraints are weaker than the coefficients constraints (the maximum number of comparisons in the $k$th iteration $M_k$ = $\lfloor c\log_\sigma(1+k_0+k)\rfloor$).

| | Coefficients Constraints | Relaxed Coefficients Constraints |
|---|---|---|
| SR | $\sigma \geq 1/\min_{\boldsymbol{\theta}\in\Theta} P(y(\boldsymbol{\theta}) < U_{u,v})$, $0 < c \leq 1/r$, $c\log_\sigma(1+k_0) \geq 1$. | $\sigma > 1$, $c > 0$ $c\log_\sigma(1+k_0) \geq 1$. |
| SC | $\sigma \geq 1/\min_{\boldsymbol{\theta}\in\Theta,\boldsymbol{\theta}'\in\Theta\backslash\boldsymbol{\theta}} P\left(y(\boldsymbol{\theta}') < y(\boldsymbol{\theta})\right)$, $0 < c \leq 1$, $c\log_\sigma(1+k_0) \geq 1$. | $\sigma > 1$, $c > 0$ $c\log_\sigma(1+k_0) \geq 1$. |

In the numerical experiments in the next section, we tune the coefficients sets based on the relaxed coefficients constraints. Our goal for the relaxation of the constraints is to improve the performance of the SR algorithm and the SC algorithm as much as we can.

## 6.2 Numerical Experiments for Comparisons Based on Relaxed Coefficients Constraints of Stochastic Ruler and Stochastic Comparison

In this section, we do numerical tests to compare the performance of the SR algorithm, the SC algorithm and the DSPSA algorithm. For the SR algorithm and the SC algorithm, we relax the coefficients constraints on these algorithm coefficients. We tune the coefficients to pick an appropriate set for each algorithm, then we use these sets to do the comparisons of SR, SC and DSPSA algorithms.

The sequence generated by SR and SC is composed of multivariate integer points, so for SR and SC, $[\hat{\boldsymbol{\theta}}_k] = \hat{\boldsymbol{\theta}}_k$. However, the sequence generated by DSPSA is composed of non-multivariate integer points. Therefore, we compare the sequence $\{\hat{\boldsymbol{\theta}}_k\}$ (equal to $\{[\hat{\boldsymbol{\theta}}_k]\}$) from SR and SC with the sequence $\{[\hat{\boldsymbol{\theta}}_k]\}$ from DSPSA.

## 6.2.1 Numerical Experiments on Two-Dimensional Separable Loss Function over a Small Domain

In Section 6.1, we see that for a two-dimensional separable loss function with small domain, the performances of both SR and SC are good when we pick some appropriate sets of coefficients that satisfy the relaxed coefficients constraints. In the following, first we compare the three algorithms, SR, SC, and DSPSA, for the same two-dimensional separable loss function with the same small domain.

In Section 6.1, we figured out the appropriate values of coefficients for the SR algorithm and the SC algorithm, so we use them directly in the comparison here. Since all sets of coefficients in Tables 6.2 and 6.4 provide final errors being 0 (in terms of sample mean of both $\left\|\hat{\boldsymbol{\theta}}_q - \boldsymbol{\theta}^*\right\| \big/ \left\|\hat{\boldsymbol{\theta}}_0 - \boldsymbol{\theta}^*\right\|$ and $\left|L(\hat{\boldsymbol{\theta}}_q) - L(\boldsymbol{\theta}^*)\right| \big/ \left|L(\hat{\boldsymbol{\theta}}_0) - L(\boldsymbol{\theta}^*)\right|$), we try all of these sets and pick the one that provides the most efficient convergent performance. Therefore, for both the SR algorithm and the SC algorithm, we pick the first set of coefficients in both Tables 6.2 and 6.4. For DSPSA, based on the selection guidelines in Section 3.3, we set $\alpha = 0.501$, $A = 1000$, $a = 1$. All three algorithms converge effectively to 0 error with 20,000 noisy measurements, so in order to see their performance clearly, we use fewer noisy measurements in each replicate. Also we increase the number of replicates to reduce the effect of the noise to see a clearer comparison. We set the number of replicates to be 200, and the number of noisy measurements in each replicate to be 2500. In Figure 6.3 we provide the comparison result for this two-dimensional separable loss function.

We see that the error is ultimately 0 for all the three algorithms for the two-dimensional separable loss function with small domain.

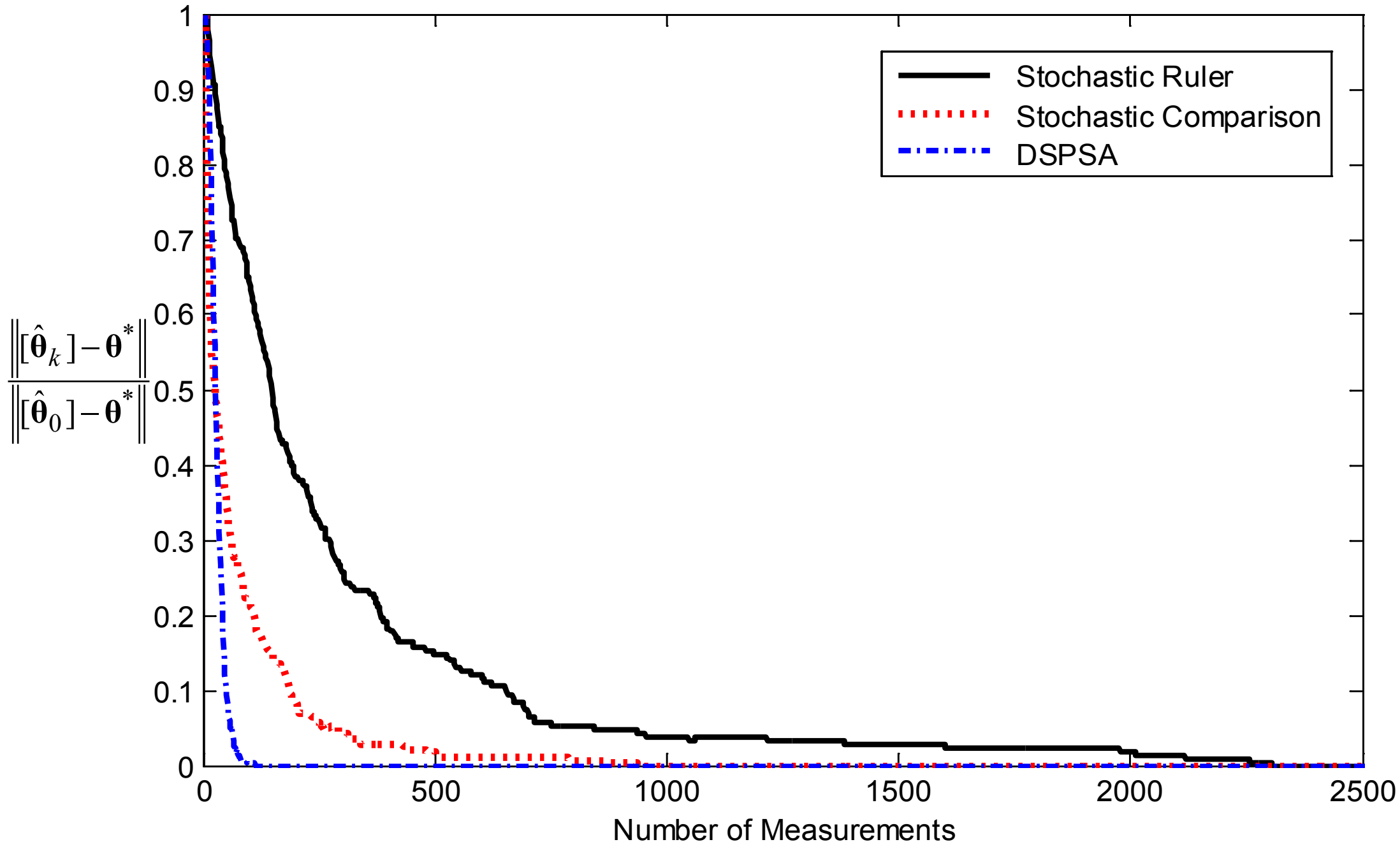

**Figure 6.3** Comparisons of SR, SC and DSPSA for the simple two-dimensional separable loss function with small domain. All the three algorithms converge to 0. Each curve represents the sample mean of 200 independent replicates.

## 6.2.2 Numerical Experiments on High-Dimensional Loss Functions over Large Domain

Let us consider two loss functions. The first loss function is a more general separable loss function defined on $\mathbb{Z}^p$:

$$L(\boldsymbol{\theta}) = \sum_{i=1}^{p} t_i^2 \ .$$

The optimal solution of the separable loss function is $\boldsymbol{\theta}^* = \mathbf{0}_p$ , where $\mathbf{0}_p$ is a $p$-dimensional vector with all components being 0. The second loss function is a skewed quartic function (Spall, 2003, Ex 6.6) defined on $\mathbb{Z}^p$:

$$L(\boldsymbol{\theta}) = \boldsymbol{\theta}^T \boldsymbol{B}^T \boldsymbol{B}\boldsymbol{\theta} + 0.1\sum_{i=1}^{p} (\boldsymbol{B}\boldsymbol{\theta})_i^3 + 0.01\sum_{i=1}^{p} (\boldsymbol{B}\boldsymbol{\theta})_i^4 \ ,$$

where $(\boldsymbol{B}\boldsymbol{\theta})_i$ represents the $i$th component of the vector $\boldsymbol{B}\boldsymbol{\theta}$ and $p\boldsymbol{B}$ is an upper triangular matrix of 1's. The optimal solution of the skewed quartic loss function is $\boldsymbol{\theta}^* = \mathbf{0}_p$.

We consider the high-dimensional case $p$ = 200 for both loss functions, since we have already discussed low-dimensional small domain problem. We set the domain to be $\{-10, \ldots, -1, 0, 1, \ldots, 10\}^{200}$, hence the search space has on the order of $10^{264}$ elements. The measurement noises $\varepsilon$ are i.i.d. $N(0,1)$. The initial guess is set as $\hat{\boldsymbol{\theta}}_0 = 10 \times \mathbf{1}_{200}$ for each loss function, where $\mathbf{1}_{200}$ is a 200-dimensional vector with all components being 1. Moreover, for the numerical comparisons, which will show in Section 6.2.2.2, we do 20 replicates and in each replicate the number of noisy measurements is 20,000.

For the SR algorithm and the SC algorithm, in the next subsection, we use a tuning process to pick the relatively good coefficients for them. For DSPSA, we have discussed the practical guidelines for coefficients selection in Section 3.3, so we follow these guidelines to pick the coefficients, and we get $\alpha$ = 0.501, $A$ = 1000, $a$ = 0.05 for the separable loss function, and $\alpha$ = 0.501, $A$ = 1000, $a$ = 0.01 for the skewed quartic loss function. From numerical results on DSPSA, we have that the final error in terms of both sample mean of $\left\|[\hat{\boldsymbol{\theta}}_{10,000}]-\boldsymbol{\theta}^*\right\| \big/ \left\|[\hat{\boldsymbol{\theta}}_0]-\boldsymbol{\theta}^*\right\|$ and $\left|L([\hat{\boldsymbol{\theta}}_{10,000}])-L(\boldsymbol{\theta}^*)\right| \big/ \left|L([\hat{\boldsymbol{\theta}}_0])-L(\boldsymbol{\theta}^*)\right|$ are 0 and 0 for the separable loss function, and 0.4242 and 0.013 for the skewed quartic loss function (10,000 iterations indicates 20,000 noisy measurements in DSPSA). The sequences $\{\hat{\boldsymbol{\theta}}_k\}$ generated by DSPSA may not be composed of multivariate integer points, and the sequences $\{\hat{\boldsymbol{\theta}}_k\}$ generated by the SR algorithm and the SC algorithms are all composed of multivariate integer points, so for DSPSA we use the criteria $\left\|[\hat{\boldsymbol{\theta}}_{10,000}]-\boldsymbol{\theta}^*\right\| \big/ \left\|[\hat{\boldsymbol{\theta}}_0]-\boldsymbol{\theta}^*\right\|$ and $\left|L([\hat{\boldsymbol{\theta}}_{10,000}])-L(\boldsymbol{\theta}^*)\right| \big/ \left|L([\hat{\boldsymbol{\theta}}_0])-L(\boldsymbol{\theta}^*)\right|$ to measure the final errors.

#### 6.2.2.1 Coefficient Selection on High-Dimensional Loss Functions over Large Domain for Stochastic Ruler and Stochastic Comparison

In this subsection, we consider coefficient selection on high-dimensional cases for SR and SC. Under these cases, the domain will be significantly larger than the domain $\{-1,0,1\}^2$ in the simple example in Section 6.2.1. For the SR algorithm, Yan and Mukai (1992) consider the general neighborhood structure $N(\boldsymbol{\theta}) \subseteq \Theta \setminus \boldsymbol{\theta}$. Therefore, we will try both the global square-ring neighbor $N(\boldsymbol{\theta}) = \Theta \setminus \boldsymbol{\theta}$ and the local square-ring neighbor $N(\boldsymbol{\theta})$ = $\left\{\boldsymbol{\theta}' \in \mathbb{Z}^p \middle| |t_i - t_i'| \le 1 \text{ for all } i\right\} \setminus \boldsymbol{\theta}$ , and pick the one that achieves the better performance. For the SC algorithm, Gong et al (1999) only consider global square-ring neighbor $N(\boldsymbol{\theta}) = \Theta \setminus \boldsymbol{\theta}$, and they think that even though a good neighborhood structure can speed up the convergence, a poor neighborhood can hurt the speed. But in order to make the SC algorithm to provide as good performance as we can for the specific loss functions, we will try both the global square-ring neighbor and the local square-ring neighbor, and pick the one that achieves the better performance.

For the sequence $\{M_k\}$, we have four strategies, which apply to both SR and SC: 1) small $M_0$, low increasing rate; 2) small $M_0$, high increasing rate; 3) big $M_0$, low increasing rate; 4) big $M_0$, high increasing rate. There are possible problems with each strategy. A small $M_0$ may make the algorithm to accept new points too easy in the early iterations, which can easily lead the algorithm to bad points in the early iterations. A big $M_0$ may make the algorithm too hard to accept new points in the early iterations, which may make $\hat{\boldsymbol{\theta}}_k$ to be stuck at the initial guess in the early iterations. Moreover, the value of $M_k$ is increasing in the whole process, big $M_0$ also indicates big $M_k$ in the later iterations, which may make the whole sequence $\left\{\hat{\boldsymbol{\theta}}_k\right\}$ to be harder to leave the initial guess. A low increasing rate of $M_k$ may make the algorithm to continue its good/bad performance of the early iterations. A high increasing rate of $M_k$ may make the probability to accept new points decrease too fast. We use a numerical tuning process to find appropriate coefficients sets in the numerical tests below.

### 6.2.2.1.1 Tuning Process for the Stochastic Ruler Algorithm

First, we do the tuning process to decide the coefficients and the neighborhood structure that will be used for the SR algorithm. As we have discussed in Section 6.1, we can pick the coefficients based on the relaxed coefficients constraints, and set $c$ = 1. Tables 6.6 and 6.7 show the results of tuning process for the separable loss function and the skewed quartic loss function.

From Tables 6.6 and 6.7, we find that for the high-dimensional separable loss function and skewed quartic loss function, the best sets of coefficients in local square-ring neighborhood structure and global square-ring neighborhood structure lead to similar level of final errors in terms of both sample mean of $\left\|\hat{\boldsymbol{\theta}}_q - \boldsymbol{\theta}^*\right\| \Big/ \left\|\hat{\boldsymbol{\theta}}_0 - \boldsymbol{\theta}^*\right\|$ and sample mean of $\left|L(\hat{\boldsymbol{\theta}}_q) - L(\boldsymbol{\theta}^*)\right| \Big/ \left|L(\hat{\boldsymbol{\theta}}_0) - L(\boldsymbol{\theta}^*)\right|$ , where the sample mean is the arithmetic mean of the observed values from 20 replicates and $q$ is the number of iterations when the SR algorithm hits 20,000 noisy measurements in each replicate. In both Tables 6.6 and 6.7, the first set of coefficients in both neighborhood structures is the base set. Based on the base set no matter how we enlarge/shrink the range of the uniform distribution of the stochastic ruler, how we increase/decrease the growing rate of $M_k$ (maximum number of comparisons), and how we increase/decrease the value of $M_0$, we cannot improve the performance of the SR algorithm significantly relative to the performance with the base set. We see that no matter how the coefficients are adjusted, the final errors from the SR algorithm are always much larger than the final errors from DSPSA (discussed at the third paragraph of Section 6.2.2) in terms of the sample mean of $\left\|\hat{\boldsymbol{\theta}}_q - \boldsymbol{\theta}^*\right\| \Big/ \left\|\hat{\boldsymbol{\theta}}_0 - \boldsymbol{\theta}^*\right\|$ for both loss functions. For the high-dimensional separable loss function, from Table 6.6, the final errors from SR in terms of sample mean of $\left|L(\hat{\boldsymbol{\theta}}_q) - L(\boldsymbol{\theta}^*)\right| \Big/ \left|L(\hat{\boldsymbol{\theta}}_0) - L(\boldsymbol{\theta}^*)\right|$ are still much larger than the final errors from DSPSA. However, for the high-dimensional skewed quartic loss function, from Table 6.7, we see that the final errors from SR in terms of sample mean of $\left|L(\hat{\boldsymbol{\theta}}_q) - L(\boldsymbol{\theta}^*)\right| \Big/ \left|L(\hat{\boldsymbol{\theta}}_0) - L(\boldsymbol{\theta}^*)\right|$ are very close to 0 relative to the error from DSPSA, because the skewed quartic loss function has a special skewed and twisted shape, and contains a large flat area.

**Table 6.6** Performance of the SR algorithm on high-dimensional ($p$ = 200) separable loss function with different sets of algorithm coefficients. Global neighborhood structure and local neighborhood structure provide similar level of final errors in terms of both sample mean of $\left\|\hat{\boldsymbol{\theta}}_q - \boldsymbol{\theta}^*\right\| \Big/ \left\|\hat{\boldsymbol{\theta}}_0 - \boldsymbol{\theta}^*\right\|$ and sample mean of $\left|L(\hat{\boldsymbol{\theta}}_q) - L(\boldsymbol{\theta}^*)\right| \Big/ \left|L(\hat{\boldsymbol{\theta}}_0) - L(\boldsymbol{\theta}^*)\right|$ for the appropriate sets. The number of replicates is 20 and the number of noisy measurements of the loss function is 20,000.

| Neighbor Structure | $u$ | $v$ | $\sigma$ | $k_0$ | Sample mean of $\dfrac{\left\|\hat{\boldsymbol{\theta}}_q - \boldsymbol{\theta}^*\right\|}{\left\|\hat{\boldsymbol{\theta}}_0 - \boldsymbol{\theta}^*\right\|}$ | Sample mean of $\dfrac{\left|L(\hat{\boldsymbol{\theta}}_q) - L(\boldsymbol{\theta}^*)\right|}{\left|L(\hat{\boldsymbol{\theta}}_0) - L(\boldsymbol{\theta}^*)\right|}$ |
|---|---|---|---|---|---|---|
| Global | 0 | 6500 | 100 | 200 | 0.5503 | 0.3030 |
| Global | −500 | 6500 | 100 | 200 | 0.5547 | 0.3077 |
| Global | 500 | 6500 | 100 | 200 | 0.5536 | 0.3065 |
| Global | 0 | 6000 | 100 | 200 | 0.9076 | 0.8579 |
| Global | 0 | 7000 | 100 | 200 | 0.5648 | 0.3191 |
| Global | 0 | 6500 | 200 | 200 | 0.5545 | 0.3076 |
| Global | 0 | 6500 | 20 | 200 | 0.8226 | 0.7239 |
| Global | 0 | 6500 | 100 | 100 | 0.5551 | 0.3082 |
| Global | 0 | 6500 | 100 | 1000 | 0.5517 | 0.3045 |
| Local | 0 | 19000 | 100 | 200 | 0.5946 | 0.3538 |
| Local | −1000 | 19000 | 100 | 200 | 0.6073 | 0.3691 |
| Local | 1000 | 19000 | 100 | 200 | 0.6001 | 0.3603 |
| Local | 0 | 18000 | 100 | 200 | 1 | 1 |

Table 6.6, continued

| Neighbor Structure | $u$ | $v$ | $\sigma$ | $k_0$ | Sample mean of $\frac{\lVert \hat{\boldsymbol{\theta}}_q - \boldsymbol{\theta}^* \rVert}{\lVert \hat{\boldsymbol{\theta}}_0 - \boldsymbol{\theta}^* \rVert}$ | Sample mean of $\frac{\lvert L(\hat{\boldsymbol{\theta}}_q) - L(\boldsymbol{\theta}^*) \rvert}{\lvert L(\hat{\boldsymbol{\theta}}_0) - L(\boldsymbol{\theta}^*) \rvert}$ |
|---|---|---|---|---|---|---|
| Local | 0 | 20000 | 100 | 200 | 0.6064 | 0.3680 |
| Local | 0 | 19000 | 200 | 200 | 0.6053 | 0.3666 |
| Local | 0 | 19000 | 20 | 200 | 0.6008 | 0.3612 |
| Local | 0 | 19000 | 100 | 100 | 0.6080 | 0.3699 |
| Local | 0 | 19000 | 100 | 1000 | 0.6101 | 0.3726 |

**Table 6.7** Performances of the SR algorithm on high-dimensional ($p$ = 200) skewed quartic loss function with different sets of coefficients. Global neighborhood structure and local neighborhood structure can provide similar level of final errors in terms of the sample mean of $\lVert \hat{\boldsymbol{\theta}}_q - \boldsymbol{\theta}^* \rVert / \lVert \hat{\boldsymbol{\theta}}_0 - \boldsymbol{\theta}^* \rVert$ for the appropriate sets. The final errors in terms of sample mean of $\lVert \hat{\boldsymbol{\theta}}_q - \boldsymbol{\theta}^* \rVert / \lVert \hat{\boldsymbol{\theta}}_0 - \boldsymbol{\theta}^* \rVert$ are not close to 0 relative to the final errors in terms of sample mean of $\lvert L(\hat{\boldsymbol{\theta}}_q) - L(\boldsymbol{\theta}^*) \rvert / \lvert L(\hat{\boldsymbol{\theta}}_0) - L(\boldsymbol{\theta}^*) \rvert$. The number of replicates is 20 and the number of noisy measurements of loss function is 20,000.

| Neighbor Structure | $u$ | $v$ | $\sigma$ | $k_0$ | Sample mean of $\frac{\lVert \hat{\boldsymbol{\theta}}_q - \boldsymbol{\theta}^* \rVert}{\lVert \hat{\boldsymbol{\theta}}_0 - \boldsymbol{\theta}^* \rVert}$ | Sample mean of $\frac{\lvert L(\hat{\boldsymbol{\theta}}_q) - L(\boldsymbol{\theta}^*) \rvert}{\lvert L(\hat{\boldsymbol{\theta}}_0) - L(\boldsymbol{\theta}^*) \rvert}$ |
|---|---|---|---|---|---|---|
| Global | 0 | 1 | 10 | 15 | 0.5953 | 0.000068 |

Table 6.7, continued

| Neighbor Structure | $u$ | $v$ | $\sigma$ | $k_0$ | Sample mean of $\dfrac{\lVert\hat{\boldsymbol{\theta}}_q-\boldsymbol{\theta}^*\rVert}{\lVert\hat{\boldsymbol{\theta}}_0-\boldsymbol{\theta}^*\rVert}$ | Sample mean of $\dfrac{\lvert L(\hat{\boldsymbol{\theta}}_q)-L(\boldsymbol{\theta}^*)\rvert}{\lvert L(\hat{\boldsymbol{\theta}}_0)-L(\boldsymbol{\theta}^*)\rvert}$ |
|---|---|---|---|---|---|---|
| Global | −1 | 1 | 10 | 15 | 0.6842 | 0.2 |
| Global | 0.5 | 1 | 10 | 15 | 0.5952 | 0.000069 |
| Global | 0 | 2 | 10 | 15 | 0.599 | 0.000077 |
| Global | 0 | 0.5 | 10 | 15 | 0.6398 | 0.1 |
| Global | 0 | 1 | 15 | 15 | 0.6040 | 0.000070 |
| Global | 0 | 1 | 5 | 15 | 0.7423 | 0.35 |
| Global | 0 | 1 | 10 | 10 | 0.6036 | 0.000063 |
| Global | 0 | 1 | 10 | 20 | 0.5964 | 0.000069 |
| Local | 0 | 15000 | 100 | 200 | 0.6038 | 0.000840 |
| Local | −1000 | 15000 | 100 | 200 | 0.6090 | 0.000937 |
| Local | 1000 | 15000 | 100 | 200 | 0.6032 | 0.000945 |
| Local | 0 | 14000 | 100 | 200 | 0.9208 | 0.8003 |
| Local | 0 | 16000 | 100 | 200 | 0.6014 | 0.000995 |
| Local | 0 | 15000 | 200 | 200 | 0.6034 | 0.0014 |
| Local | 0 | 15000 | 20 | 200 | 0.6017 | 0.0012 |
| Local | 0 | 15000 | 100 | 100 | 0.5974 | 0.0018 |
| Local | 0 | 15000 | 100 | 1000 | 0.6054 | 0.0014 |

### 6.2.2.1.2 Tuning Process for the Stochastic Comparison Algorithm

For the SC algorithm, as we have discussed in Section 6.1, we can pick the coefficients based on the relaxed coefficients constraints, and set $c$ = 1. In Gong et al (1999), the authors only consider global square-ring neighborhood $N(\boldsymbol{\theta}) = \Theta \setminus \boldsymbol{\theta}$. But here the search space has on the order of $10^{264}$ points. In order to have possible better performance for the SC algorithm, we also try the local neighborhood structure. Thus, we do the tuning process not only for the global square-ring neighborhood structure, but also for the local square-ring neighborhood structure.

Tables 6.8 and 6.9 show the results of tuning process for the high-dimensional separable loss function and skewed quartic loss function. In both Tables 6.8 and 6.9, the first set in both neighborhood structures is the base set. Based on the base set, no matter how we increase/decrease the growing rate of $M_k$, and how we increase/decrease the value of $M_0$, we cannot improve the performance of the SC algorithm significantly relative to the performance with the base set. We see that no matter how the coefficients are adjusted, the final errors from the SC algorithm are always larger than the final errors from DSPSA (discussed at the third paragraph of Section 6.2.2) in terms of the sample mean of $\left\|\hat{\boldsymbol{\theta}}_q - \boldsymbol{\theta}^*\right\| \big/ \left\|\hat{\boldsymbol{\theta}}_0 - \boldsymbol{\theta}^*\right\|$ for both loss functions. For the high-dimensional separable loss function, the final errors from SC are still larger than the final error from DSPSA in terms of the sample mean of $\left|L(\hat{\boldsymbol{\theta}}_q) - L(\boldsymbol{\theta}^*)\right| \big/ \left|L(\hat{\boldsymbol{\theta}}_0) - L(\boldsymbol{\theta}^*)\right|$. However, for the high-dimensional skewed quartic loss function, from Table 6.9, we see that the final errors

from SC in terms of sample mean of $\left|L(\hat{\boldsymbol{\theta}}_q) - L(\boldsymbol{\theta}^*)\right| \Big/ \left|L(\hat{\boldsymbol{\theta}}_0) - L(\boldsymbol{\theta}^*)\right|$ are already very close to 0 relative to the final errors in terms of sample mean of $\left\|\hat{\boldsymbol{\theta}}_q - \boldsymbol{\theta}^*\right\| \Big/ \left\|\hat{\boldsymbol{\theta}}_0 - \boldsymbol{\theta}^*\right\|$, because the skewed quartic loss function has a special skewed and twisted shape, and contains a large flat area.

From Table 6.8, we find that for the high-dimensional separable loss function, the best set of coefficients for the local square-ring neighborhood structure can lead to significantly better performance than the best sets of coefficients for the global square-ring neighborhood structure. Thus, the performance of the SC algorithm can be improved for the separable loss function by using the local neighborhood structure. For the high-dimensional skewed quartic loss function, from Table 6.9, we see that the appropriate sets of coefficients for both local square-ring neighborhood structure and global square-ring neighborhood structure can lead the sequence to the points with the similar level of final errors in terms of both the sample mean of $\left\|\hat{\boldsymbol{\theta}}_q - \boldsymbol{\theta}^*\right\| \Big/ \left\|\hat{\boldsymbol{\theta}}_0 - \boldsymbol{\theta}^*\right\|$ and sample mean of $\left|L(\hat{\boldsymbol{\theta}}_q) - L(\boldsymbol{\theta}^*)\right| \Big/ \left|L(\hat{\boldsymbol{\theta}}_0) - L(\boldsymbol{\theta}^*)\right|$. Therefore, compared with the global neighborhood structure, the local neighborhood structure does not improve the performance of SC for the skewed quartic loss function.

**Table 6.8** Performance of the SC algorithm on high-dimensional ($p$ = 200) separable loss function with different sets of coefficients. Local neighborhood structure can lead to better performance in terms of both the sample mean of $\left\|\hat{\boldsymbol{\theta}}_q - \boldsymbol{\theta}^*\right\| \big/ \left\|\hat{\boldsymbol{\theta}}_0 - \boldsymbol{\theta}^*\right\|$ and the sample mean of $\left|L(\hat{\boldsymbol{\theta}}_q) - L(\boldsymbol{\theta}^*)\right| \big/ \left|L(\hat{\boldsymbol{\theta}}_0) - L(\boldsymbol{\theta}^*)\right|$ than global neighborhood structure. The number of replicates is 20 and the number of noisy measurements of loss function is 20,000 in each replicate.

| Neighbor Structure | $k_0$ | $\sigma$ | Sample mean of $\dfrac{\left\|\hat{\boldsymbol{\theta}}_q - \boldsymbol{\theta}^*\right\|}{\left\|\hat{\boldsymbol{\theta}}_0 - \boldsymbol{\theta}^*\right\|}$ | Sample mean of $\dfrac{\left|L(\hat{\boldsymbol{\theta}}_q) - L(\boldsymbol{\theta}^*)\right|}{\left|L(\hat{\boldsymbol{\theta}}_0) - L(\boldsymbol{\theta}^*)\right|}$ |
|---|---|---|---|---|
| Global | 100 | 10 | 0.5279 | 0.2787 |
| Global | 1000 | 10 | 0.5312 | 0.2822 |
| Global | 10 | 10 | 0.5293 | 0.2803 |
| Global | 100 | 100 | 0.5305 | 0.2815 |
| Global | 100 | 5 | 0.5299 | 0.2809 |
| Local | 100 | 10 | 0.1982 | 0.0393 |
| Local | 1000 | 10 | 0.1995 | 0.0398 |
| Local | 10 | 10 | 0.2008 | 0.0403 |
| Local | 100 | 100 | 0.1985 | 0.0394 |
| Local | 100 | 5 | 0.1997 | 0.0399 |

**Table 6.9** Performance of the SC algorithm on high-dimensional ($p$ = 200) skewed quartic loss function with different sets of coefficients. Global neighborhood structure and local neighborhood structure provide similar level of final errors in terms of both sample mean of $\|\hat{\boldsymbol{\theta}}_q - \boldsymbol{\theta}^*\| / \|\hat{\boldsymbol{\theta}}_0 - \boldsymbol{\theta}^*\|$ and sample mean of $|L(\hat{\boldsymbol{\theta}}_q) - L(\boldsymbol{\theta}^*)| / |L(\hat{\boldsymbol{\theta}}_0) - L(\boldsymbol{\theta}^*)|$ for the appropriate sets. The final errors in terms of sample mean of $\|\hat{\boldsymbol{\theta}}_q - \boldsymbol{\theta}^*\| / \|\hat{\boldsymbol{\theta}}_0 - \boldsymbol{\theta}^*\|$ are not close to 0 relative to the final errors in terms of sample mean of $|L(\hat{\boldsymbol{\theta}}_q) - L(\boldsymbol{\theta}^*)| / |L(\hat{\boldsymbol{\theta}}_0) - L(\boldsymbol{\theta}^*)|$ . The number of replicates is 20 and the number of noisy measurements of loss function in each replicate is 20,000.

| Neighbor Structure | $\sigma$ | $k_0$ | Sample mean of $\dfrac{\|\hat{\boldsymbol{\theta}}_q - \boldsymbol{\theta}^*\|}{\|\hat{\boldsymbol{\theta}}_0 - \boldsymbol{\theta}^*\|}$ | Sample mean of $\dfrac{\lvert L(\hat{\boldsymbol{\theta}}_q) - L(\boldsymbol{\theta}^*)\rvert}{\lvert L(\hat{\boldsymbol{\theta}}_0) - L(\boldsymbol{\theta}^*)\rvert}$ |
|---|---|---|---|---|
| Global | 100 | 10 | 0.5951 | 0.000057415 |
| Global | 1000 | 10 | 0.6010 | 0.000069153 |
| Global | 10 | 10 | 0.5978 | 0.000071960 |
| Global | 100 | 5 | 0.5996 | 0.000064031 |
| Global | 100 | 100 | 0.6109 | 0.0001065 |
| Local | 100 | 10 | 0.6060 | 0.000073827 |
| Local | 1000 | 10 | 0.6060 | 0.000074125 |
| Local | 10 | 10 | 0.6057 | 0.000074743 |
| Local | 100 | 5 | 0.5971 | 0.000099406 |
| Local | 100 | 100 | 0.5969 | 0.000095786 |

### 6.2.2.2 Numerical Comparisons of Stochastic Ruler Algorithm, Stochastic Comparison Algorithm, and DSPSA algorithm

Based on the discussions above, in the numerical comparison below, we pick local square-ring neighbor for both SR and SC algorithms for the high-dimensional separable loss function, and pick global square-ring neighbor for both SR and SC algorithms for the high-dimensional skewed quartic loss function. We do 20 replicates, and in each replicate the number of noisy measurements is 20,000.

Now let us first consider the numerical comparison for the high-dimensional separable loss function. For the SR algorithm, based on the results of tuning process in the Table 6.6, we pick the first set of the coefficients for the local neighborhood structure $\sigma = 100$, $c = 1$, $k_0 = 200$, $u = 0$, $v = 19000$. For the SC algorithm, based on the results of tuning process in the Table 6.8, we pick the first set of the coefficients for the local neighborhood structure $\sigma = 10$, $c = 1$, $k_0 = 100$. For DSPSA, as we have discussed, we pick the coefficients $\alpha = 0.501$, $A = 1000$, $a = 0.05$. In Figure 6.4 we have the comparison results of the three algorithms (SR, SC, and DSPSA) for the high-dimensional separable loss function. The sequences $\{\hat{\boldsymbol{\theta}}_k\}$ generated by DSPSA may not be composed of multivariate integer points, and the sequences $\{\hat{\boldsymbol{\theta}}_k\}$ generated by the SR algorithm and the SC algorithms are composed of multivariate integer points, so in the numerical comparison below, we use the criteria of the sample mean of $\left\|[\hat{\boldsymbol{\theta}}_k] - \boldsymbol{\theta}^*\right\| \big/ \left\|[\hat{\boldsymbol{\theta}}_0] - \boldsymbol{\theta}^*\right\|$ and the sample mean of $\left|L([\hat{\boldsymbol{\theta}}_k]) - L(\boldsymbol{\theta}^*)\right| \big/ \left|L([\hat{\boldsymbol{\theta}}_0]) - L(\boldsymbol{\theta}^*)\right|$.

From Figure 6.4, we see that DSPSA has significant better convergence efficiency than the other two algorithms for the high-dimensional separable loss function. However, the SC algorithm has better performance than the other two algorithms in the early iterations. In Figure 6.4, the SR algorithm gets stuck after the error decrease by around 40%. For the SR algorithm, the lower and upper bounds of the uniform distribution of the stochastic ruler are fixed, so when the error decreases to some level, the upper bound may be too large. Thus, the points accepted after that may not be better points. But, when we reduce the value of the upper bound, we may come across another problem. Since the number of iterations is limited, a smaller value of the upper bound may reduce the probability of accepting new points, which may hurt the performance in the early iterations. Moreover, the local square-ring neighborhood structure does not help SR to improve the convergence speed (shown in Table 6.6), because in SR the loss function value at the candidate point is not compared with current point but compared with the stochastic ruler, and as we have discussed the upper bound of the stochastic ruler may hurt the efficiency of SR significantly. The SC algorithm provides much better performance than the SR algorithm in this case. For the SC algorithm, the local square-ring neighborhood structure does help to improve the convergence speed significantly, because the local square-ring neighborhood contains fewer points than the global square-ring neighborhood, and fewer neighbor points make SC to accept new point more efficiently in the comparison between current point and candidate point. Thus, the performance by using local square-ring neighborhood is significantly better than by using global square-ring neighborhood for SC, which can be seen in Table 6.8.

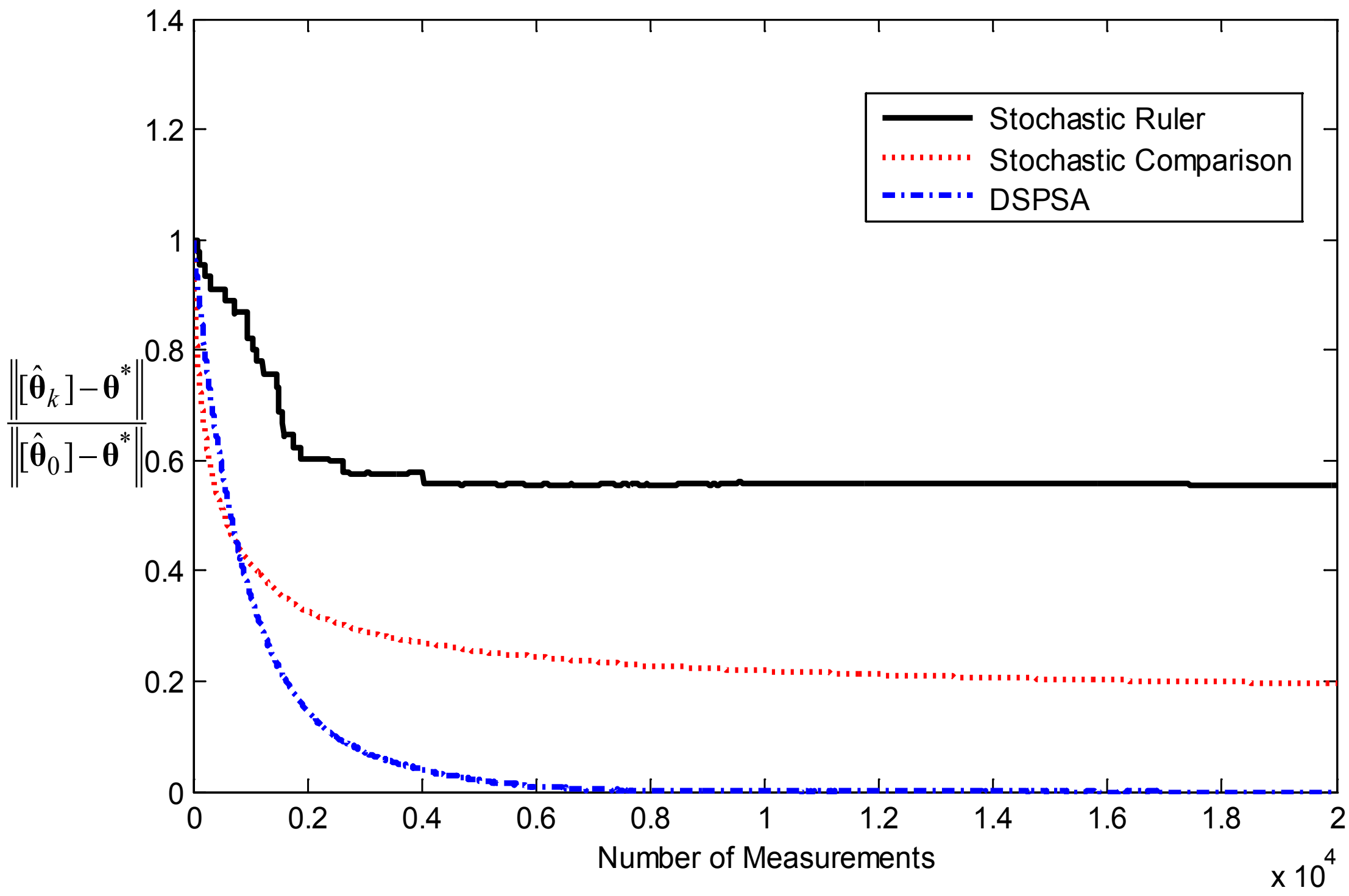

**Figure 6.4** Comparison of SR, SC and DSPSA in terms of sample mean of $\|[\hat{\boldsymbol{\theta}}_k]-\boldsymbol{\theta}^*\| / \|[\hat{\boldsymbol{\theta}}_0]-\boldsymbol{\theta}^*\|$ for high-dimensional ($p$ = 200) separable function. DSPSA provides better convergence rate than the other two algorithms. However, the SC algorithm has better performance than the other two algorithms in the early iterations. Each curve represents the sample mean of 20 independent replicates.

Next we consider the numerical comparison for the high-dimensional skewed quartic loss function. For the SR algorithm, based on the results of tuning process in the Table 6.7, we pick the first set of the coefficients for the global neighborhood structure $\sigma = 10$, $c = 1$, $k_0 = 15$, $u = 0$, $v = 1$. For the SC algorithm, based on the results of tuning process in the Table 6.9, we pick the first set of the coefficients for the global neighborhood structure $\sigma = 10$, $c = 1$, $k_0 = 100$. For DSPSA, as we discussed, we pick the coefficients $\alpha = 0.501$, $A = 1000$, $a = 0.01$. In Figures 6.5 and 6.6, we have the comparison results of the three algorithms (SR, SC and DSPSA) for high-dimensional skewed quartic loss function.

From Figure 6.5 and 6.6, we see that even though DSPSA performs much better than the other algorithms in terms of sample mean of $\left\| [\hat{\boldsymbol{\theta}}_k] - \boldsymbol{\theta}^* \right\| \big/ \left\| [\hat{\boldsymbol{\theta}}_0] - \boldsymbol{\theta}^* \right\|$, the SR algorithm and the SC algorithm provide better performance in terms of sample mean of $\left| L([\hat{\boldsymbol{\theta}}_k]) - L(\boldsymbol{\theta}^*) \right| \big/ \left| L([\hat{\boldsymbol{\theta}}_0]) - L(\boldsymbol{\theta}^*) \right|$, because the skewed quartic loss function has a special skewed and twisted shape, and contains a large flat area.

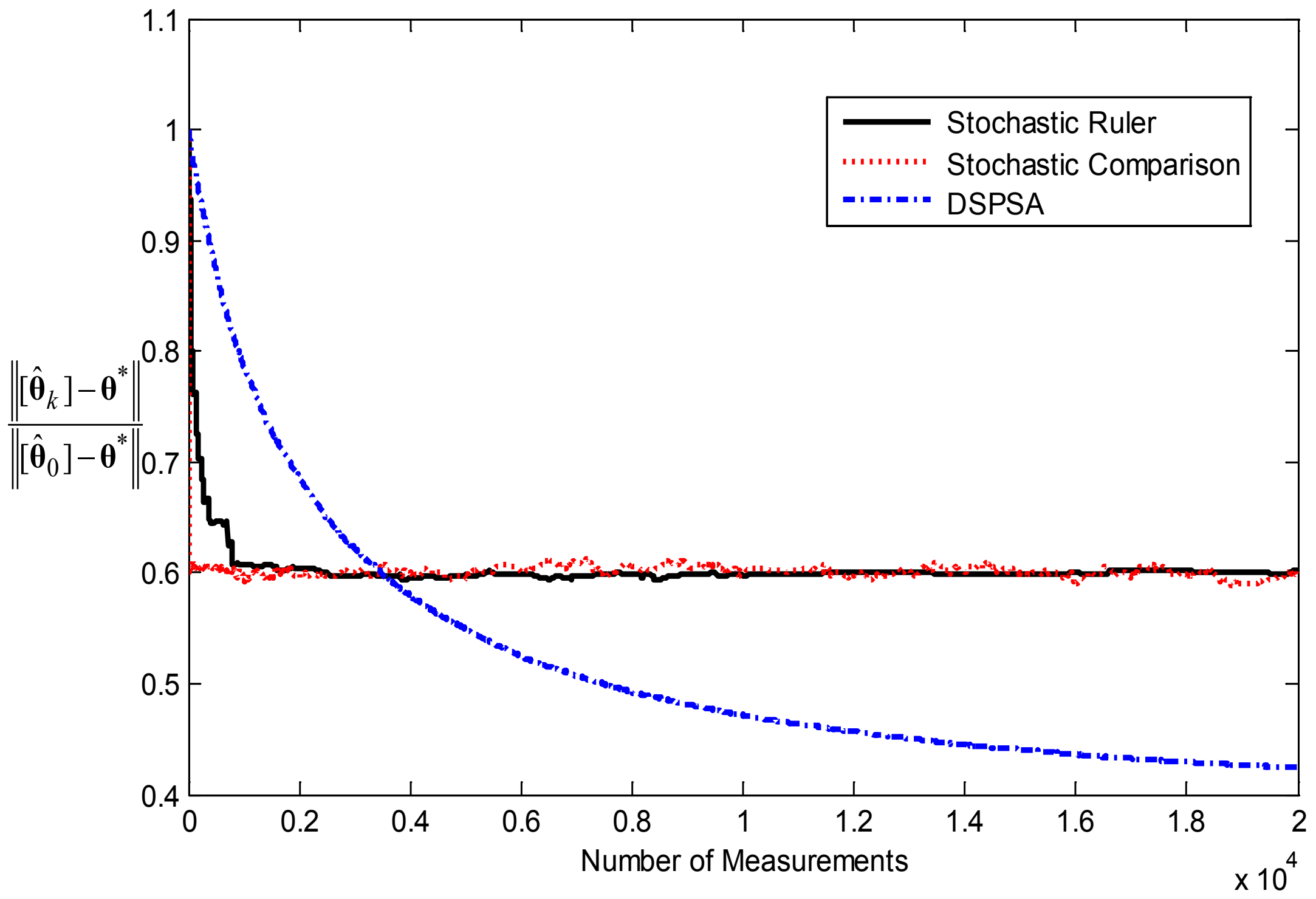

**Figure 6.5** Comparison of SR, SC and DSPSA in terms of sample mean of $\|[\hat{\boldsymbol{\theta}}_k]-\boldsymbol{\theta}^*\| / \|[\hat{\boldsymbol{\theta}}_0]-\boldsymbol{\theta}^*\|$ for high-dimensional ($p$ = 200) skewed quartic loss function. DSPSA provides better convergence rate than the other two algorithms. Each curve represents the sample mean of 20 independent replicates.

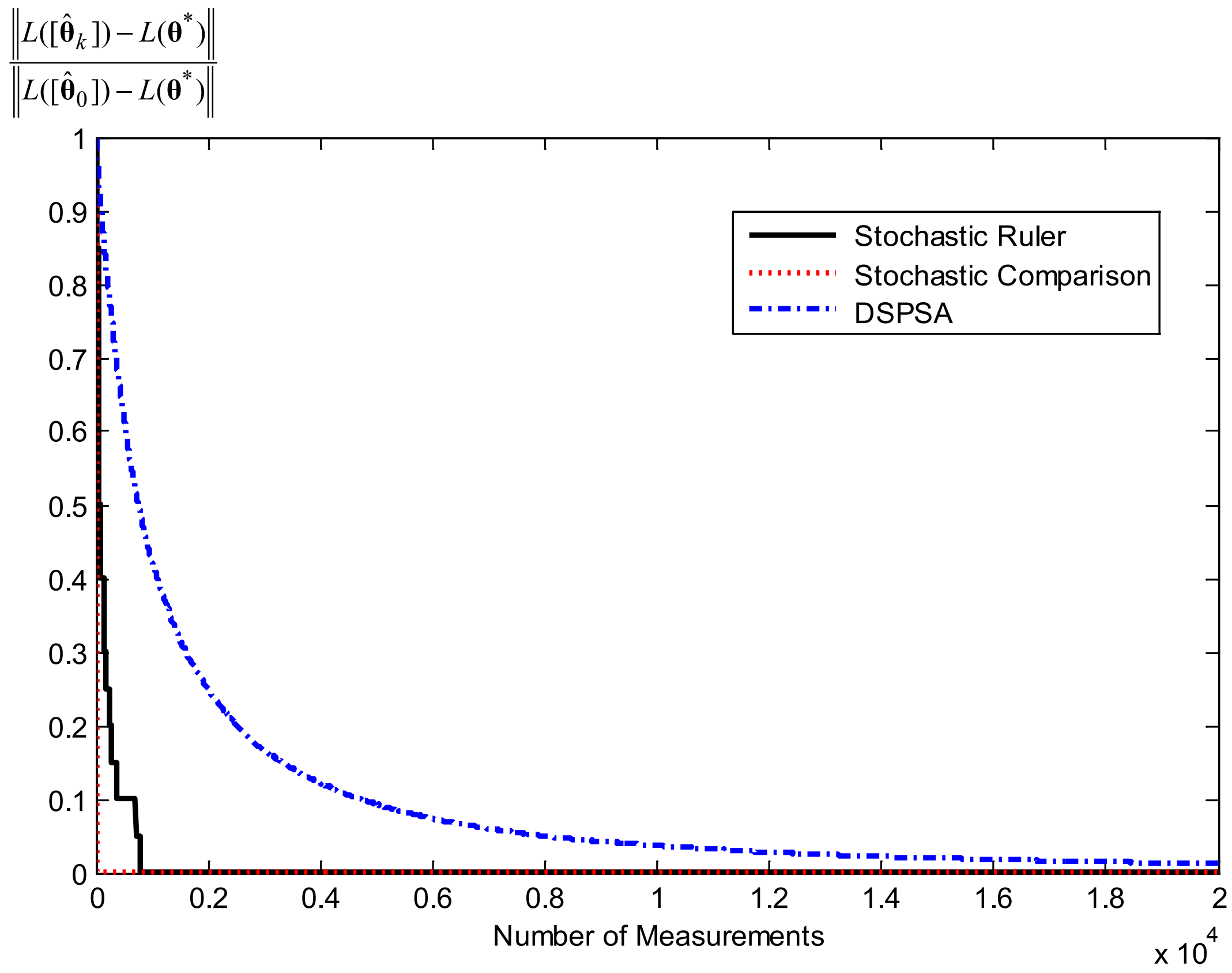

**Figure 6.6** Comparison of SR, SC and DSPSA in terms of sample mean of $\left|L([\hat{\boldsymbol{\theta}}_k]) - L(\boldsymbol{\theta}^*)\right| / \left|L([\hat{\boldsymbol{\theta}}_0]) - L(\boldsymbol{\theta}^*)\right|$ for high-dimensional ($p$ = 200) skewed quartic loss function. The sequences generated by SR and SC converge to 0 faster than DSPSA. Each curve represents the sample mean of 20 independent replicates.

In all, even though the SR is a convergent algorithm, determining the values of lower and upper bounds for the uniform distribution of the stochastic ruler is not easy. Moreover, for both SR and SC, the choices of $\sigma$, $c$ and $k_0$ are hard to determine. Therefore, the efficiencies of these two algorithms are hard to control. However, for DSPSA, we have already figured out the guidelines of coefficients selections. In addition, from these numerical results, we see that even though the SR algorithm and the SC algorithm provide bad performance in terms of sample mean of $\left\|[\hat{\boldsymbol{\theta}}_k]-\boldsymbol{\theta}^*\right\| \big/ \left\|[\hat{\boldsymbol{\theta}}_0]-\boldsymbol{\theta}^*\right\|$ for the high-dimensional loss functions, the performance in terms of sample mean of $\left|L([\hat{\boldsymbol{\theta}}_k])-L(\boldsymbol{\theta}^*)\right| \big/ \left|L([\hat{\boldsymbol{\theta}}_0])-L(\boldsymbol{\theta}^*)\right|$ seems good for some special shape of loss function (the skewed quartic loss function). For DSPSA, the performance in terms of both sample mean of $\left\|[\hat{\boldsymbol{\theta}}_k]-\boldsymbol{\theta}^*\right\| \big/ \left\|[\hat{\boldsymbol{\theta}}_0]-\boldsymbol{\theta}^*\right\|$ and sample mean of $\left|L([\hat{\boldsymbol{\theta}}_k])-L(\boldsymbol{\theta}^*)\right| \big/ \left|L([\hat{\boldsymbol{\theta}}_0])-L(\boldsymbol{\theta}^*)\right|$ are quite reasonable and stable. Furthermore, no comparisons of noisy loss function values is involved in the whole process of DSPSA, which might make DSPSA work better to handle the noise than the SR algorithm and the SC algorithm in other problem settings.

# Chapter 7

# Application of DSPSA in Resource Allocation in Public Health

In this chapter, we consider the application of DSPSA towards developing optimal public health strategies for containing the spread of influenza given limited societal resources. We use DSPSA to do the simulation based optimization to solve the optimal intervention method for H1N1 to achieve the minimal loss to the economy. The loss to the economy includes the cost related to the interventions and the cost induced by people who become infected. In the following, we introduce open source software for intervention strategies: FluTE, and based on the simulation results of FluTE we use DSPSA to determine the optimal intervention strategy instead of just doing sensitivity analysis on the effect of each intervention method.

# 7.1 Background

Seasonal influenza epidemics and worldwide epidemics (pandemics) have caused many deaths and much economic loss throughout history. As recently as 2009, the H1N1 virus was responsible for millions of confirmed infections and thousands of deaths throughout the world. H1N1 has also caused billions of dollars of loss to the economy. There are many intervention methods, including vaccination, antiviral courses, school closure, isolation of ascertained people, quarantine of family members of symptomatic people, and the use of masks and alcohol-based hand gels. In this chapter, we use the DSPSA algorithm to determine the best combination of interventions to achieve the least economic loss to society.

Due to the novelty of the H1N1 virus strain, most people are not immune to it. New vaccines for H1N1 need time to be tested and produced. Generally, the production of vaccines for each year adheres to the following schedule: In February, the World Health Organization selects the three virus strains to be used in the fall vaccination program. Then, the production of vaccines for each strain takes place in 11-day old embryonated eggs. The clinical trials to evaluate vaccines on people start in May and June, and the filling and packaging are done in July and August. The vaccines ship out in September and vaccinations begin to take place in October. The immunity takes two weeks to build up after the vaccination. The schedule of vaccine manufacturing is tight, and many disruptions may happen, such as an insufficient number of eggs and/or a delay in picking the strains.

Generally, the number of vaccines is not sufficient to cover 100% of the population. Also, different vaccines have different restrictions for different people. Due to a shortage of vaccines, we need to determine which groups of people should be given priority to receive the vaccination, so as to achieve better intervention results. Antiviral agents can also be used for treatment of infected people and protection of susceptible people. However, unlike the vaccines, the antiviral agents are effective only during the period in which they are being taken. A single course of antiviral agents is enough for ten days of protection for susceptible people or five days of treatment for infected people. In addition, if school closure becomes necessary, the cost is significant, because it takes additional teaching hours to make up missed classes, and often parents are forced to take time off from work.

Vaccinations, antiviral agents, and school closure are the three main methods that can help to control an epidemic and reduce costs incurred from hospitalizations and treatment in the ICU. Policymakers want to find the optimal intervention strategy to achieve the best result in terms of loss to society. Generally, researchers set up a base case, such as no intervention, and then they do the sensitivity analysis by adding some interventions (e.g. Halder et al., 2011, Khazeni et al., 2009). However, none of these papers solve this problem by using an optimization algorithm. Moreover, due to the noise in the loss function, a stochastic optimization algorithm is needed. Thus, we will use DSPSA to find the best combination of interventions to achieve the least economic loss to society.

# 7.2 Introduction of the Simulator (FluTE)

We use software available online to simulate the spread of the influenza virus. The simulation outputs are used to calculate the noisy measurements of the objective loss function value. The software is called FluTE and it is freely available at www.cs.unm.edu/~dlchao/flute/. This software is created based on a new stochastic model of the epidemic within a large population, and it is used to help policymakers prepare for future influenza seasonal epidemics or pandemics (Chao et al., 2010, Chao et al., 2011). In this new stochastic model, the infection processes are based on real historical influenza data. The model is calibrated to the data gathered from the 1957/1958 Asian (H2N2) and the 2009 pandemic (H1N1). The calibration involves tuning the contact probabilities of two individuals to make the final age-specific illness attack rates similar to the data in the 1957/1958 Asian (H2N2) and the 2009 pandemic (H1N1).

We will now briefly discuss the basic structure and assumptions of this simulator. Basically, FluTE is based on a person-by-person model that simulates the spread of influenza among a large population. The structure of the synthetic population is based on real communities in the United States, and the synthetic social network generated in the model is constructed according to the realistic contact networks. The model of the transmission of the diseases and the infectious process is based on the natural history of influenza. After setting up the population and the way influenza spreads, Chao et al. (2010) discuss the simulation of interventions. Among many interventions, Chao et al. (2010) mainly consider three ways, and we discuss them as follows:

The first way is the primary pharmaceutical intervention: vaccination. Vaccination can help people to reduce the probability of becoming infected, the probability of becoming ill if infected, and the probability of transmission of the infection. These probabilities are the key properties of the efficiency for different vaccines. Generally, vaccines need two weeks to reach maximum efficacy.

The second type of intervention is antiviral agents. People who take antiviral agents can reduce the probability of susceptibility, the probability of becoming ill if infected, and the probability of transmission during the period when the antiviral agents are taken.

The third type of intervention is to close school. School closure reduces the contacts within school, but it also increases daytime contact within family members and the neighborhood.

## 7.3 Input Parameters of FluTE and Loss Function

The inputs and outputs of the software FluTE are written in text files. There are many possible input parameters (shown in the README file of FluTE). We list some of the important ones in Table 7.1:

**Table 7.1** Key Input Parameters for FluTE.

| Parameter name: Parameter type [dimension of the parameter] | Meanings |
|---|---|
| `datafile`: string | Input data file names (e.g. seattle, la, usa) |
| `R0`: real | Basic reproductive number: the average number of people an infected individual infects in susceptible population |
| `preexistingimmunitybyage`: real[5] | Vector of 5 real numbers representing the fractions of individuals in each age group with pre-existing immunity: all preschoolers (0 – 4 years), all school-age children (5 – 18 years), all young adults (19 – 29 years), all older adults (30 – 64 years), and all elderly (65 + years) |
| `vaccinationfraction`: real | Fraction of assigned people to get vaccination |
| `runlength`: integer | Number of simulated days |
| `vaccinepriorities`: integer [13] | A vector represents the vaccine priority for 13 categories of people. 0 indicates no vaccination, 1 indicates highest priority, 2 indicates next-highest priority, etc. The categories are: essential workforce, pregnant women, members of families containing infants, high risk preschoolers, high risk school-age children, high risk young adults, high risk older adults, high risk elderly, all preschoolers, all school-age children, all young adults, all older adults, and all elderly. |

| Parameter name: Parameter type [dimension of the parameter] | Meanings |
|---|---|
| vaccinedoses: integer[2] | The vaccine ID followed by the number of vaccine doses available at the beginning of the simulation |
| vaccineproduction: integer [runlength+1] | Vaccine ID followed by the number of doses that become available each day |
| vaccinedata: integer, real[3], real[6], Boolean | Vaccine ID followed by the vaccine efficiency parameters and restrictions of the vaccines on 5 groups of people and pregnant women |
| vaccinebuildup: integer, integer, real [29] | Vaccine ID followed by the day that the boost should be given and how the efficacy of the vaccine is built up over the 29 days after the vaccine is given |
| vaccineefficacybyage: real [5] | A vector indicates the vaccine efficacy for each group. |
| AVEs:real | Antiviral agents efficacy for susceptible people |
| AVEi:real | Antiviral agents efficacy for infected people |
| AVEp:real | Antiviral agents efficacy for illness given infected |
| responsedelay: integer | Number of days before initiating reactive strategies |
| ascertainmentdelay: integer | Number of days to be taken to ascertain a symptomatic individual |
| ascertainmentfraction: real | Fraction of symptomatic individuals who can be |

| Parameter name: Parameter type [dimension of the parameter] | Meanings |
|---|---|
| | ascertained |
| `essentialfraction`: real | Fraction of working-age adults that are essential workforce |
| `pregnantfraction`: real [5] | Fraction of people who are pregnant in each of five groups |
| `highriskfraction`: real [5] | Fraction of people who are at high risk from influenza in each of five groups |
| `seedinfected`: integer | Number of people at time 0 to infect the whole population |
| `antiviralpolicy`: string | Represents the policy by which people get antiviral agents. Possible values include: “none”, “treatmentonly” (treat the ascertained people only), “HHTAP” (all family members get antiviral agents if one member is ascertained), “HHTAP100” (special option for Los Angeles county) |
| `schoolclosuredays`: integer | Number of days to close schools |

Some of these inputs are related to the properties of vaccines and H1N1, and some are assigned as the variables of the optimization problem. In the optimization problem, we can change the values of the vaccination fraction, the vaccination priorities, the antiviral policy, and the time of school closure. The outputs of the simulator contain information including the number of communities, the vaccines used, the antivirals used, the total number of symptomatic individuals by age, etc. Based on these outputs, we can set up the basic formula of the loss function of the optimization problem. We know that the cost of intervention and the cost of hospitalization change in opposite directions, which means that more money used for interventions can reduce the cost of hospitalization, while less money used for interventions may increase the cost of hospitalization. The objective function is to find out the best tradeoff point to minimize the total cost of the epidemic to the economy. The loss function $L(\boldsymbol{\theta})$ is defined as:

vaccination cost + antivirals cost +school closure cost
+hospitalization cost + ICU cost+death cost,

and $\boldsymbol{\theta}$ is the vector of variables that contains the input parameters related to the intervention strategies, which will be discussed in detail in the next section. The noisy measurements of the loss function can be calculated, all based on the outputs of the simulator.

We now discuss in detail the cost components of $L$ in the above. The total cost of vaccination is composed of two parts: one part is related to the cost of the vaccine itself, and the other part is related to the side effects of the vaccine. The cost per vaccine includes production costs, administration fees, and the cost of patient time. According to the Khazeni et al. (2009) paper, the cost per vaccine is approximately \$30. Khazeni et al. (2009) also show that about 0.001% of the population could have severe adverse reactions to the vaccines, which results in hospitalization (the cost of hospitalization is discussed below).

The total cost of antiviral agents is equal to the antiviral cost per course times the number of courses used. The antiviral cost per course includes the production cost and dispensing cost, which is approximately $56 per course (Halder et al., 2011).

The school closure cost is equal to the average cost per day per student times the number of students and the number of school closure days. In Halder et al. (2011), they show that the average cost for school closure is around $20 per day per student. But this value only contains the cost to make up missed classes and does not include the cost to the parents, since they need to take time off from work to stay at home to care for their infected children. Based on a report by the Brookings Institute (http://www.pbs.org/newshour/updates/health/july-dec09/flu-costs_10-08.html), the average cost of closing school is between $35 and $157 per student per day. Also, Araz et al. (2012) show that on average, parents need to take 2.5 working days off per week to take care of students if school is closed. Halder et al. (2011) indicate that the average wage per week for a person is $836. Thus, the cost to take care of students per day is around $80 dollars. Adding $80 dollars to the cost of making up missed classes ($20), this comes to around $100 in total per student per day of school closure. In FluTE, each community (500-3000 individuals) has two elementary schools, one middle school and one high school, which have 79, 128 and 155 students, respectively. Thus, the total cost for one community of school closure is around $181,000 per week.

The cost of hospitalization or ICU is calculated as the expected value of the expense for people who require hospital or ICU care. Khazeni et al. (2009) show 3.3% of symptomatic individuals need five days of hospitalization, and among these individuals, 10% need ten days in the ICU. Chao et al. (2011) provide a different percentage on Web Table 7, where high risk people might have a higher rate of hospitalization, but they do not consider the ICU care. We follow Chao et al. (2011)'s table for hospitalization rate, and also assume 10% ICU care, as shown by Khazeni et al. (2009). Khazeni et al. (2009) also show that the average cost of medical hospitalization and ICU per day per person are \$1,830 and \$3,739, respectively, and we calculate the total cost of hospitalization or ICU based on these two values.

The cost of mortality is calculated as the expected value of future earnings. Molinari et al. (2007) list the cost of mortality for different groups of people, and we use their results in our problem. In addition, Chao et al. (2011) indicate that high risk people could have a higher mortality rate, and we also follow their result.

# 7.4 Practical Implementation of DSPSA by Using FluTE

We want to solve the optimization problem by figuring out the optimal value for the vaccine fraction ($F$), vaccine priorities ($\boldsymbol{P}$), antiviral policy ($A$) and length of school closure ($S$). We see that there are 1+13+1+1 = 16 variables for which we need to solve. Here we assume $F \in \{0, 0.1, 0.2, \ldots, 1\}$ ($F$ may not be integer, then in DSPSA, we map the value to integer by multiplying 10, and in the input of the simulator we map the value back to the domain of $F$ through dividing 10), $\boldsymbol{P} \in \{0, 1, 2, 3\}^{13}$, $A \in \{0, 1, 2, 3\}$ (0:

“none”, 1: “treatmentonly”, 2: “HHTAP100”, 3: “HHTAP”), and $S \in \{0, 1, \ldots, \lfloor \text{runlength}/7 \rfloor\}$ (We measure the school closing time in terms of weeks, and in the input of the simulator we map the value of weeks back to days by multiplying 7). Then we set $\boldsymbol{\theta} = [10F, \boldsymbol{P}^T, A, S]^T$ for the optimization problem $\min_{\boldsymbol{\theta}} L(\boldsymbol{\theta})$.

However, from the output of the simulator, we cannot get the value of $L(\boldsymbol{\theta})$ because of the noise in the simulation. There are many randomness parts in the simulation. For example, whether an individual is asymptomatic after infection is random. Whether an individual has contact with an infected person on a specific day is random. For a vaccine dose, among the people who are in the same level of priority, the people to be picked are random. Due to many randomness issues in the simulation, the outputs of the simulator involve a lot of noise. Thus, we can only get the noisy measurements $y(\boldsymbol{\theta})$ of the loss function. The levels of the noises compared with the loss function values are different with different inputs and settings of the simulator. In Table 7.3, we will see the noise level for the real problem based on different inputs.

From the description of DSPSA in Section 2.1, we know that there are two main parts of the algorithm: the perturbation part (generating perturbation direction $\boldsymbol{\Delta}_k$), and the calculation part (calculating the gradient-like vector and $\hat{\boldsymbol{\theta}}_{k+1}$). In this chapter, for $\boldsymbol{\Delta}_k$, we always use the distribution where the $\Delta_{ki}$ are independent Bernoulli random variables taking the values $\pm 1$ with probability $1/2$. Figure 7.1 shows how DSPSA can be used to solve the discrete stochastic optimization problem through the simulation results from FluTE. We see that the role of FluTE is to generate a simulation result of the loss function. The simulation takes a long time when the population is very large or when the average number of people that an infected individual infects in a susceptible population is big.

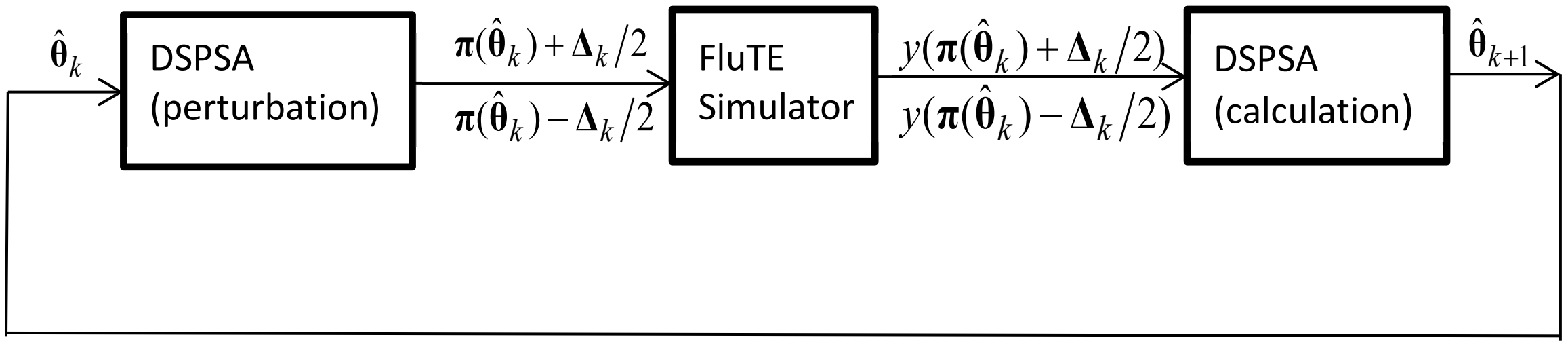

**Figure 7.1** The process used to solve the epidemic problem by using DSPSA. The simulator generates outputs that are used to calculate the noisy measurements of the loss function.

In DSPSA, we need to determine the coefficients $\alpha$, $a$ and $A$. By the guidelines of the coefficients selection in Section 3.3, we pick $\alpha$ = 0.501, $A$ = 0.1 × (number of allowed iterations). For the choice of $a$, we pick $a$ such that the multiplication of $a_0$ and approximated value of $\hat{\boldsymbol{g}}_0(\hat{\boldsymbol{\theta}}_0)$ is equal to the desired change magnitude in the early iterations.

# 7.5 Solution in Extreme Cases of Free and Very Expensive Interventions

Before presenting the main numerical experiments in Section 7.6, we do a test for two extreme cases: 1) free interventions 2) expensive interventions (no hospitalization, ICU and mortality cost) as sanity checks.

For these two extreme cases, we assume the number of available vaccines and antiviral agents is infinite, so that for the vaccine priority of any group we only need to choose between 0 and 1 (without/with vaccination), which implies that for these two extreme cases, $\boldsymbol{P} \in \{0, 1\}^{13}$. In addition, there are a finite number of choices for the vaccine fraction, the antivirals policy, and the weeks of school closure. Thus, the feasible domain is related to a finite number of unit hypercubes. Even though in Chapter 2, we only discussed the binary case in Theorem 2.2, we think that for the case of finite number of unit hypercubes, similar results can be derived. Therefore, before running the numerical tests, we first check the inner product condition ($\bar{\boldsymbol{g}}(\boldsymbol{m}_{\boldsymbol{\theta}})^T(\boldsymbol{\theta}-\boldsymbol{\theta}^*) > 0$ for all $\boldsymbol{m}_{\boldsymbol{\theta}} \in \mathcal{M}_{\boldsymbol{\theta}}$ and all $\boldsymbol{\theta} \in \Theta \setminus \{\boldsymbol{\theta}^*\}$) for these extreme cases.

We have discussed the inner product condition for the binary case in Section 2.3. Thus, here we can check the inner product condition for the extreme cases by using a similar idea.

For any point $\boldsymbol{\theta}' = \{t_1', ..., t_p'\} \neq \boldsymbol{\theta}^* = \{t_1^*, ..., t_p^*\}$ and $\boldsymbol{m}_{\boldsymbol{\theta}'} \in \mathcal{M}_{\boldsymbol{\theta}'}$, suppose $l_i^{(\boldsymbol{m}_{\boldsymbol{\theta}'})}$ and $u_i^{(\boldsymbol{m}_{\boldsymbol{\theta}'})}$ are the lower and upper bound in the *i*th coordinate of the unit hypercube centered with $\boldsymbol{m}_{\boldsymbol{\theta}'}$, respectively. Let $\mathrm{Z}_i = \left\{\boldsymbol{\theta} \middle| t_i = l_i^{(\boldsymbol{m}_{\boldsymbol{\theta}'})} \text{ and } t_j \in \{l_j^{(\boldsymbol{m}_{\boldsymbol{\theta}'})}, u_j^{(\boldsymbol{m}_{\boldsymbol{\theta}'})}\} \text{ for } j \neq i\right\}$ and $\overline{\mathrm{Z}}_i = \left\{\boldsymbol{\theta} \middle| t_i = u_i^{(\boldsymbol{m}_{\boldsymbol{\theta}'})} \text{ and } t_j \in \{l_j^{(\boldsymbol{m}_{\boldsymbol{\theta}'})}, u_j^{(\boldsymbol{m}_{\boldsymbol{\theta}'})}\} \text{ for } j \neq i\right\}$.

When the interventions are free, adding any intervention decreases the number of symptomatic persons, which reduces the total cost. Therefore,

$$\sum_{\boldsymbol{\theta} \in \mathrm{Z}_i} L(\boldsymbol{\theta}) > \sum_{\boldsymbol{\theta} \in \overline{\mathrm{Z}}_i} L(\boldsymbol{\theta}), \tag{7.1}$$

for all $i = 1, ..., 16$. By similar arguments of Section 2.3, we know that inequality (7.1) indicates that

$$\overline{\boldsymbol{g}}(\boldsymbol{m}_{\boldsymbol{\theta}'})^T \boldsymbol{e}_i < 0, \tag{7.2}$$

for all $i = 1, ..., 16$. Since the optimal solution is to vaccinate 100% of the population, take HHTAP as the antivirals policy, and close school for all weeks, we have

$$\boldsymbol{\theta}' - \boldsymbol{\theta}^* \leq \boldsymbol{0}, \tag{7.3}$$

which means $t_i' - t_i^* \leq 0$ for all *i*. Since $\boldsymbol{\theta}' \neq \boldsymbol{\theta}^*$, there exist $\tilde{i}$ such that $t_{\tilde{i}}' - t_{\tilde{i}}^* < 0$. Combining the results of inequality (7.2) and (7.3), we have

$$\overline{\boldsymbol{g}}(\boldsymbol{m}_{\boldsymbol{\theta}'})^T (\boldsymbol{\theta}' - \boldsymbol{\theta}^*) > 0. \tag{7.4}$$

Thus, inequality (7.4) shows that the inner product condition is satisfied for the free intervention case.

When the interventions are expensive (the cost of hospitalization, ICU, and mortality is 0), adding any intervention increases the total cost. Therefore, we have that for any point $\boldsymbol{\theta}' \neq \boldsymbol{\theta}^*$

$$\sum_{\boldsymbol{\theta}\in Z_i} L(\boldsymbol{\theta}) < \sum_{\boldsymbol{\theta}\in \bar{Z}_i} L(\boldsymbol{\theta})$$

for all $i = 1, \ldots, 16$ for the expensive interventions case. By similar arguments of Section 2.3, we have

$$\bar{\boldsymbol{g}}(\boldsymbol{m}_{\boldsymbol{\theta}'})^T \boldsymbol{e}_i > 0, \tag{7.5}$$

for all $i = 1, \ldots, 16$. Since the optimal solution is to have no intervention, then we have

$$\boldsymbol{\theta}' - \boldsymbol{\theta}^* \geq \boldsymbol{0}, \tag{7.6}$$

which means $t_i' - t_i^* \geq 0$ for all *i*. Since $\boldsymbol{\theta}' \neq \boldsymbol{\theta}^*$, there exist $\tilde{\tilde{i}}$ such that $t_{\tilde{\tilde{i}}}' - t_{\tilde{\tilde{i}}}^* > 0$. Combining the results of inequalities (7.5) and (7.6), we have $\bar{\boldsymbol{g}}(\boldsymbol{m}_{\boldsymbol{\theta}'})^T (\boldsymbol{\theta}' - \boldsymbol{\theta}^*) > 0$, which indicates that the inner product condition holds true for the expensive intervention extreme case.

Overall, the inner product condition holds true for both extreme cases. In addition, due to the characteristics of the simulator and the loss function, we see that the other conditions in Theorem 2.2 are obviously true here. Thus, DSPSA will converge to $\boldsymbol{\theta}^*$ in both extreme cases.

For these two extreme cases, we consider a population with 20,000 people (with 10,000 working-age individuals, among whom 7,390 are employed). We set the run length (`runlength`) to be 28, which indicates that the software runs the simulation for 28 days (4 weeks), and let the number of initial infected people be 20. We set `R0` = 2 (see Table 7.1), which indicates that on average one infected person transmits the virus to 2 people. As we mentioned in the second paragraph of this section, we assume the availability of vaccines and antiviral agents is infinite in the extreme cases, so that for the vaccine priority of any group we only need to choose between 0 and 1 (without/with vaccination). Moreover, we assume that the vaccine builds up immunity immediately. In addition, we also assume that the data on the percentage of high risk/pregnant people and the data on pre-existing immunity are the same as that of Chao et al. (2011). In the following numerical tests, we plot the graph in terms of loss function value and check the final intervention policy generated by DSPSA. Because we do not know the exact form of the loss function, we can only plot the noisy measurements of the loss function. In the following paragraphs, we present the results of the numerical tests for the extreme cases.

For the free interventions case, as we have discussed, the optimal solution is to use all possible interventions: vaccinating 100% of all people, giving all members of ascertained people antiviral agents (HHTAP), and closing schools for the full four weeks period of interest. We set the number of iterations at 1,000. By the guidelines of coefficients selection in Section 3.3, we pick $\alpha = 0.501$, $a = 0.00005$, $A = 100$. The initial policy for use in starting the algorithm is the case of no intervention expressed as {0, 0, 0, 0, 0, 0, 0, 0, 0, 0, 0, 0, 0, 0, 0, 0, 0, 0, 0, 0, 0, 0}. In order to see the performance of DSPSA, we pick this really bad initial guess. By using DSPSA, we can get the sequence to converge to the optimal solution of all possible interventions. The expression of the policy that we get from DSPSA is {10, 1, 1, 1, 1, 1, 1, 1, 1, 1, 1, 1, 1, 1, 3, 4}, which indicates that we can use DSPSA to achieve the optimal policy. Figure 7.2 shows the values of the noisy measurements of the loss function for each iteration. The optimal loss function value is not 0, because even if all of the intervention methods are used, there are still some people who may get infected. We find that in the first several iterations the noisy loss function values decrease a lot, because the value of school closure time increases from 0 to 4 weeks, and the value of antiviral agents policy changes from none to HHTAP. But after these iterations, the magnitude of the vibration of the noisy loss function seems significant compared to the noisy loss function values.

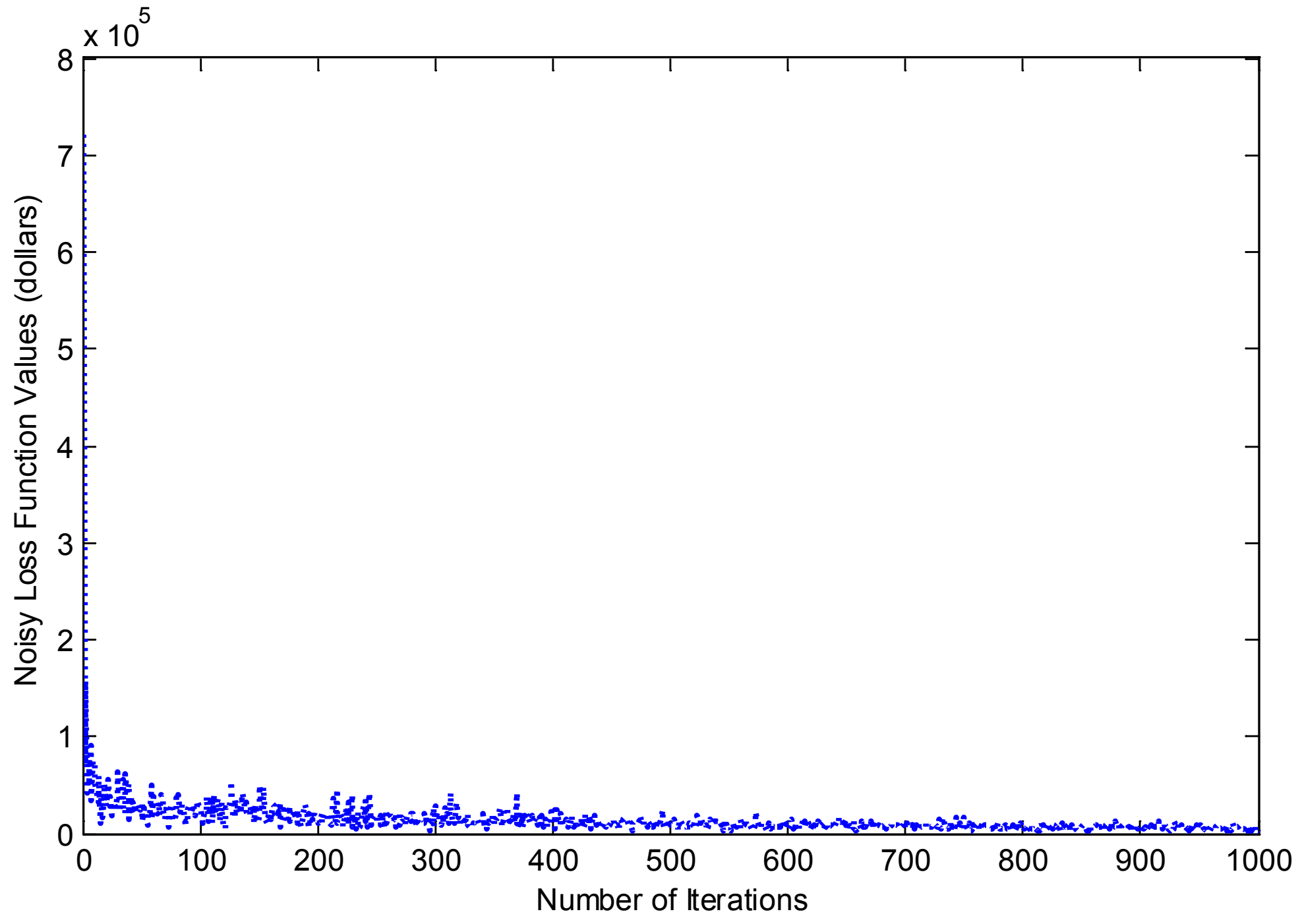

**Figure 7.2** Performance of DSPSA on free intervention case. The sequence generated by DSPSA converges to the optimal policy .

For the expensive intervention case, we set the hospitalization fee, the ICU fee and mortality cost related to the infection to be 0, and set the intervention fee to be nonzero. We know that the solution is no intervention. We set the number of iterations at 1,000. By the guidelines of coefficients selection in Section 3.3, we pick $\alpha$ = 0.501, $a$ = 0.000005, $A$ = 100. In order to see the performance of DSPSA clearly, we pick the initial policy to be using all possible interventions expressed as {10, 1, 1, 1, 1, 1, 1, 1, 1, 1, 1, 1, 1, 1, 3, 4}, which is a really bad initial guess under this extreme case. By using DSPSA, we can reach the optimal solution of no intervention. The expression of the policy that we get from DSPSA is {0, 0, 0, 0, 0, 0, 0, 0, 0, 0, 0, 0, 0, 0, 0, 0}, which indicates that we can use DSPSA to achieve the optimal policy. Figure 7.3 shows the value of noisy measurements of the loss function in each iteration. The optimal loss function value is 0, since the costs of hospitalization, ICU and mortality are all 0. We find that in the first several iterations the noisy loss function values decrease a lot, because the value of school closure time decreases from 4 weeks to none. The optimal solution is no intervention, with all costs related to the infected people (hospitalization cost, ICU cost and morality cost) being 0, so the noisy measurements of the loss function on the optimal solution are always 0.

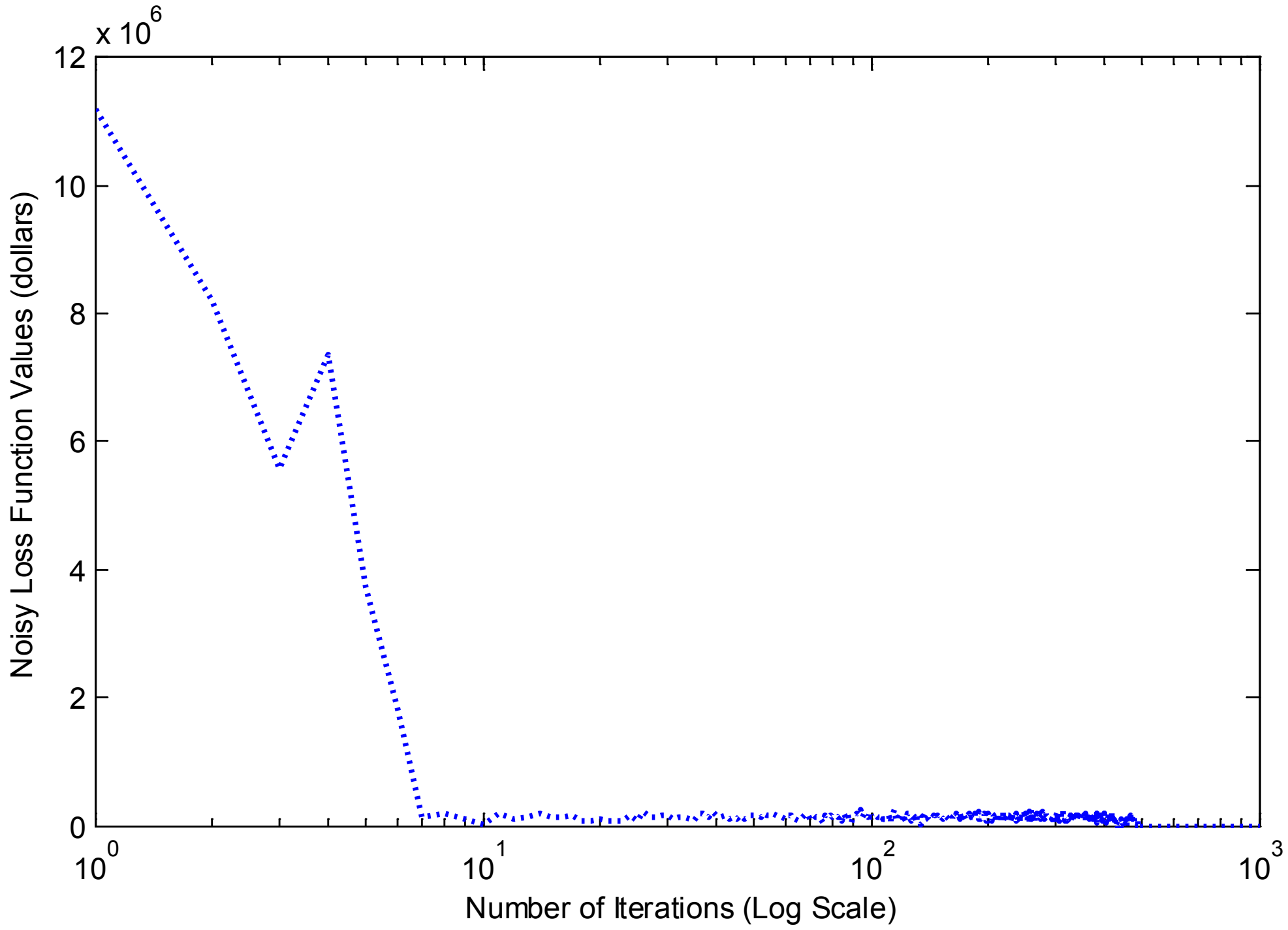

**Figure 7.3** Performance of DSPSA on expensive intervention case. The sequence generated by DSPSA converges to the optimal policy.

# 7.6 Solution for Real Problem

Now we start to solve the real problem. However, for the real problem, the loss function does not have the monotonicity property (discussed in Section 7.5) as in the extreme cases. The procedure to solve the real problem is as following. First, let us set up the data file. We pick 20 tracts of Los Angles (LA) containing a total population of 100,096, and we set this as a synthetic city having similar census properties to the whole nation. There are 9 different vaccines produced by 5 manufacturers and the detailed availability information is shown in Web Table 3 of Chao et al. (2011). This Web Table 3 provides the US H1N1 vaccine supply information, and the supply for the synthetic city here is assumed to be proportional to the national supply by population. Chao et al. (2011) show that the vaccine would be deployed in the LA counties 9 days later. We start the simulation on September 1, 2009, and the total length of simulation days is set to be 175 days (25 weeks). Although, in reality, the vaccines start to be available near the peak of the infections, we are also interested in the situation when the vaccines are available earlier than in the "real" situation because we want to check the effect if the vaccines are available earlier. Therefore, we also do simulations for "early vaccination," in which case we assume that the vaccines are available 30 days earlier than in the "real" situation. We compare the performance for these two situations (early vaccination case and late vaccination case). The key parameters are set as below (except for the parameter "`seedinfected`", all other parameters are consistent with the paper of Chao et al., 2011):

**Table 7.2** Key parameters for the real problem. These parameters are for both early vaccination case and late vaccination case.

| Parameter name | Values |
|---|---|
| `R0` | 1.3 (This value is most consistent with the cumulative illness levels of H1N1.) |
| `preexistingimmunitybyage` | {0.056, 0.0658, 0.0173, 0.0095, 0.0004} |
| `runlength` | 140 |

| Parameter name | Values |
|---|---|
| `vaccinedoses` | {0, 0} {1, 0} {2, 0} {3, 0} {4, 0} {5, 0} {6, 0} {7, 0}{8, 0} |
| `vaccinedata` | {0, 0.4, 0.4, 0.67, 1, 0, 1, 1, 1, 1, 1},<br>{1, 0.4, 0.4, 0.67, 1, 1, 0, 0, 0, 0, 0},<br>{2, 0.4, 0.4, 0.67, 1, 1, 0, 0, 0, 0, 1},<br>{3, 0.4, 0.4, 0.67, 1, 1, 0, 0, 0, 0, 0},<br>{4, 0.4, 0.4, 0.67, 1, 1, 0, 0, 0, 0, 1},<br>{5, 0.4, 0.4, 0.67, 1, 1, 0, 0, 0, 0, 0},<br>{6, 0.4, 0.4, 0.67, 1, 1, 0, 0, 0, 0, 1},<br>{7, 0.4, 0.4, 0.67, 1, 1, 1, 0, 0, 0, 1},<br>{8, 0.4, 0.5, 0.83, 1, 0.2, 0.2, 0, 1, 1, 1} |
| `vaccinebuildup` | {0, 0, 0, 0.002, 0.013, 0.038, 0.083, 0.152, 0.249, 0.379, 0.545, 0.751, 1, 1, 1, 1, 1, 1, 1, 1, 1, 1, 1, 1, 1, 1, 1, 1, 1, 1, 1} |
| `vaccineefficacybyage` | {1, 1, 1, 1, 0.6} |
| `AVEs` | 0.3 |
| `AVEi` | 0.62 |
| `AVEp` | 0.6 |
| `responsethreshhold` | 0 |
| `responsedelay` | −1 |
| `ascertainmentdelay` | 1 |

| Parameter name | Values |
|---|---|
| `ascertainmentfraction` | 0.8 |
| `essentialfraction` | 0.069 |
| `pregnantfraction` | {0,0,0.0262,0.0182,0} |
| `highriskfraction` | {0.089,0.089,0.212,0.212,0.01} |
| `seedinfected` | 100 |
| `schoolclosurepolicy` | All |

Before doing the numerical experiments, let us first test the magnitude of the noise relative to the loss function value, and the effect of each intervention method. Even though we do not know the optimal solution, from these tests we can at least have a limited idea regarding the optimal solution. We consider the intervention policies: 1) no intervention, 2) only HHAP, 3) vaccination of 50% of the people for early vaccination, 4) vaccination of 50% of the people for late vaccination, 5) vaccination of 50% of high risk people for early vaccination, 6) vaccination of 50% of high risk people for late vaccination, 7) vaccination of 50% of high risk people and essential workforce for early vaccination, 8) vaccination of 50% of high risk people and essential workforce for late vaccination, 9) one week of school closure, and 10) HHTAP and vaccination of 50% of high risk people for late vaccination. We simulate 20 replicates for each policy and get:

**Table 7.3** Noisy measurements of the loss function for the ten different intervention policies and the corresponding sample means and sample standard deviations (SD) derived from unbiased variance estimation (the values are in terms of millions of dollars).

| Repl-icate | 1<br>No interv-ention | 2<br>HHAP | 3<br>50% early vaccina-tion | 4<br>50% late vaccina-tion | 5<br>50% high risk early vaccina-tion | 6<br>50% high risk late vaccina-tion | 7<br>50% high risk and essential workfo-rce early vaccina-tion | 8<br>50% high risk and essential workfo-rce late vaccina-tion | 9<br>One week school closure | 10<br>HHTAP 50% high risk late vaccina-tion |
|---|---|---|---|---|---|---|---|---|---|---|
| 1 | 5.86 | 3.44 | 4.97 | 6.20 | 4.18 | 4.75 | 3.95 | 4.99 | 15.0 | 2.66 |
| 2 | 6.05 | 2.21 | 4.92 | 6.20 | 4.27 | 4.97 | 4.08 | 5.08 | 14.8 | 3.48 |
| 3 | 5.83 | 3.05 | 4.76 | 6.25 | 3.76 | 4.81 | 4.01 | 4.97 | 15 | 1.84 |
| 4 | 5.60 | 3.48 | 4.46 | 5.74 | 4.36 | 4.69 | 3.82 | 4.78 | 14.7 | 2.99 |
| 5 | 5.74 | 3.39 | 4.17 | 6.12 | 4.05 | 4.79 | 4.07 | 4.82 | 15.1 | 2.46 |
| 6 | 5.70 | 3.01 | 4.59 | 5.99 | 4.16 | 4.91 | 4.03 | 5.05 | 14.8 | 2.87 |
| 7 | 5.63 | 3.63 | 4.57 | 5.93 | 4.12 | 4.79 | 3.98 | 4.88 | 14.9 | 2.39 |
| 8 | 5.68 | 3.85 | 4.57 | 6.01 | 4.01 | 4.84 | 3.89 | 4.89 | 15.0 | 2.90 |
| 9 | 5.56 | 2.70 | 4.24 | 5.95 | 3.92 | 4.61 | 4.10 | 4.76 | 15.1 | 3.07 |
| 10 | 5.85 | 2.97 | 4.88 | 5.95 | 4.39 | 4.86 | 4.34 | 4.68 | 15.0 | 2.55 |
| 11 | 5.70 | 2.82 | 5.05 | 6.65 | 4.15 | 5.34 | 4.36 | 5.31 | 14.5 | 2.57 |

| Repl-icate | 1 No interv-ention | 2 HHAP | 3 50% early vaccina-tion | 4 50% late vaccina-tion | 5 50% high risk early vaccina-tion | 6 50% high risk late vaccina-tion | 7 50% high risk and essential workfo-rce early vaccina-tion | 8 50% high risk and essential workfo-rce late vaccina-tion | 9 One week school closure | 10 HHTAP 50% high risk late vaccina-tion |
|---|---|---|---|---|---|---|---|---|---|---|
| 12 | 5.65 | 2.17 | 4.41 | 6.22 | 4.18 | 4.74 | 3.88 | 4.90 | 15.0 | 2.67 |
| 13 | 5.84 | 2.73 | 5.06 | 6.34 | 4.38 | 5.03 | 4.30 | 5.04 | 15.0 | 2.81 |
| 14 | 5.76 | 3.48 | 4.75 | 6.15 | 4.31 | 5.17 | 4.19 | 5.13 | 14.8 | 1.85 |
| 15 | 5.48 | 3.34 | 5.06 | 6.16 | 4.28 | 5.05 | 4.43 | 4.89 | 14.7 | 3.16 |
| 16 | 5.87 | 3.45 | 4.98 | 6.22 | 4.00 | 4.68 | 4.36 | 4.71 | 14.8 | 2.07 |
| 17 | 5.54 | 3.10 | 5.00 | 6.43 | 4.16 | 5.08 | 4.42 | 5.11 | 15.0 | 2.98 |
| 18 | 5.59 | 3.56 | 4.80 | 6.14 | 4.20 | 5.05 | 4.15 | 5.00 | 14.9 | 2.48 |
| 19 | 5.76 | 2.04 | 5.18 | 6.53 | 4.26 | 5.14 | 4.23 | 5.26 | 14.8 | 2.86 |
| 20 | 5.71 | 3.07 | 5.10 | 6.34 | 3.98 | 4.99 | 4.01 | 4.68 | 15 | 3.28 |
| Mean | 5.72 | 3.07 | 4.78 | 6.18 | 4.16 | 4.91 | 4.13 | 4.94 | 14.9 | 2.70 |
| SD | 0.14 | 0.51 | 0.30 | 0.22 | 0.16 | 0.19 | 0.19 | 0.18 | 0.15 | 0.44 |

From the standard deviations, we see that different intervention policies may involve different levels of noise. Some noise levels are significantly large relative to the loss function values (e.g. the policy of HHTAP and vaccination of 50% of high risk people for late vaccination), and some are not (e.g. policy of one week of school closure). However, it is not enough to just discuss the noise level compared to the loss function values. It is also necessary to discuss the noise level compared to the difference between any two policies, because it affects the performance of DSPSA significantly. As we know, $\hat{\boldsymbol{g}}_k(\hat{\boldsymbol{\theta}}_k) = (y(\hat{\boldsymbol{\theta}}_k^+) - y(\hat{\boldsymbol{\theta}}_k^-))\boldsymbol{\Delta}_k^{-1}$. Thus, if the difference between two policies contains huge noise compared to the value of the difference itself, then the value of $\hat{\boldsymbol{g}}_k(\hat{\boldsymbol{\theta}}_k)$ contains significant noise, which may affect the efficiency of DSPSA. In the tests of the ten policies above, we find that for both the early vaccination case and the late vaccination case, the noise level of the policies is significantly larger than the difference between the policy of vaccinating 50% of the high risk people and the policy of vaccinating 50% of the high risk people and the essential workers (between policy 5 and policy 7, and between policy 6 and policy 8). In addition, comparing the result of all policies in Table 7.3, we see that the effect of antiviral agents is the most significant in reducing the total cost. Even though the protective effect of antiviral agents is limited, the availability of them is assumed to be infinite. Thus, antiviral agents still can help to reduce the total cost to society significantly. The effect of the vaccination is corrupted by the availability. Also, we see that if the vaccination is available 30 days earlier, the intervention effect of the vaccination is better. Moreover, we see that vaccination is more efficient if it focuses on special subgroups. Furthermore, the cost of school closure is very significant. These tests give us a brief idea of the effects of these intervention policies, which is quite similar to the sensitivity analysis discussed by other papers (e.g. Chao et al., 2011). In the following paragraphs, we start to solve the problem by using DSPSA, and after obtaining the results we will check whether the solution is consistent with the sensitivity analysis here.

First, we start to solve the discrete stochastic optimization for the case of early vaccination. The initial guess is the policy that uses no intervention, expressed as {5, 0, 0, 0, 0, 0, 0, 0, 0, 0, 0, 0, 0, 0, 0, 0}. The number of iterations is 10,000. Based on the guidelines of coefficients selection in Section 3.3, the coefficients are set to be: $a$ = 0.0000002 (this value makes the multiplication of $a_0$ and the magnitude of elements in $\hat{\boldsymbol{g}}_0(\hat{\boldsymbol{\theta}}_0)$ approximately 0.05), $A$ = 1,000, $\alpha$=0.501.

In Figure 7.4, we plot the graph in terms of the noisy loss function values (total loss to the economy) from one run of DSPSA. The plot is based on a noisy loss evaluation at each $[\hat{\boldsymbol{\theta}}_k]$; this evaluation is not needed in running the algorithm. The solution from DSPSA that we get is: {3, 0, 0, 1, 0, 0, 3, 2, 0, 3, 1, 0, 0, 0, 3, 0}, which is the expression of the policy where we give all family members antiviral agents if one member is ascertained (HHTAP), and to vaccinate 30% of the assigned groups (members of families containing infants, high risk young adults, high risk older adults, all preschoolers, and all school-age children). Among these groups, members of families containing infants and all school-age children have the highest priority. High risk older adults have the second priority. The high risk young adults and the preschoolers have the lowest priority. For the early vaccination case, vaccines are available before the peak time of H1N1. Under this case, our result is consistent with the results in Yang et al. (2009) and Chao et al. (2011), where they claim that if the vaccines are available early, it is better to vaccinate school students first. In particular, our results support Basta et al. (2009), where the authors show numerically that "School-aged children have high influenza illness attack rates and play a key role in influenza transmission."

Second, we solve the discrete stochastic optimization for the case of late vaccination (reality case). The initial guess for the reality case is the policy without using any intervention expressed as {5, 0, 0, 0, 0, 0, 0, 0, 0, 0, 0, 0, 0, 0, 0, 0}. The number of iterations is 10,000. The coefficients are set to be: $a$ = 0.00000025 (this value makes the multiplication of $a_0$ and the magnitude of elements in $\hat{\boldsymbol{g}}_0(\hat{\boldsymbol{\theta}}_0)$ approximately 0.05), $A$ = 1,000, and $\alpha$ = 0.501. In Figure 7.5, we show the total cost to the economy through iterations for late vaccination (related to "reality"). The solution that we get is: {2, 0, 0, 0, 0, 0, 1, 2, 0, 0, 3, 0, 0, 0, 3, 0}, which is the expression of the policy where we give all family members antiviral agents if one member is ascertained (HHTAP), and to vaccinate 20% of the assigned groups (high risk young adults, high risk older adults, and all school-age children). Among these groups, high risk young adults have the highest priority. High risk older adults have a lower priority. All school-age children have the lowest priority. It seems that when the vaccines are available late near the peak time of H1N1, it is better to first vaccinate the high risk adults, rather than the students. This result is also consistent with the arguments of Chao et al. (2011), where they claim that when sufficient amounts of vaccines are available near the peak of a pandemic, the recommendation of vaccination strategy would change from vaccinating children to protecting high risk persons.

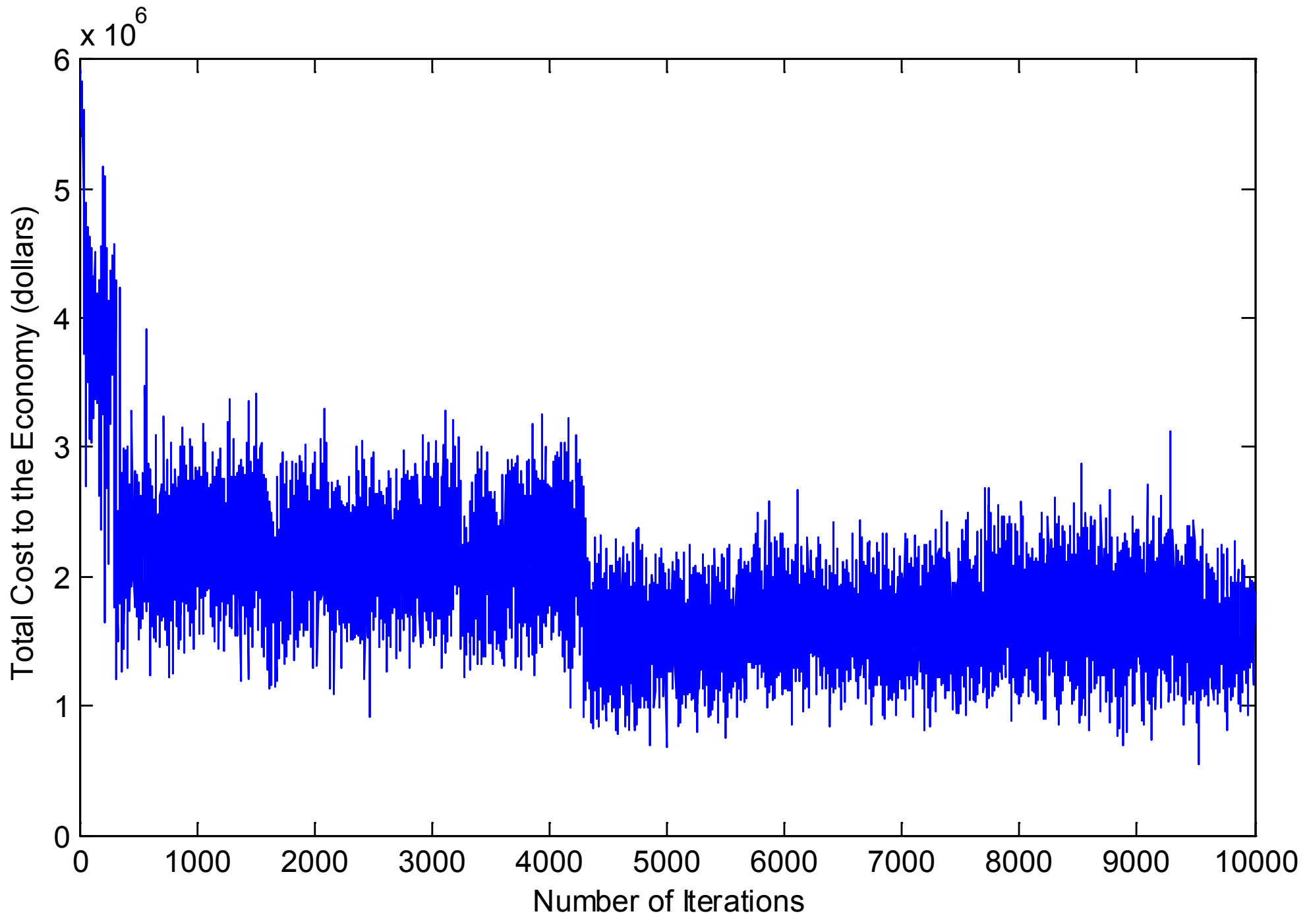

**Figure 7.4** Noisy loss values for DSPSA run for the early vaccination case. The solution indicates that for the early vaccination case, school-age children should have high priority to take the vaccination.

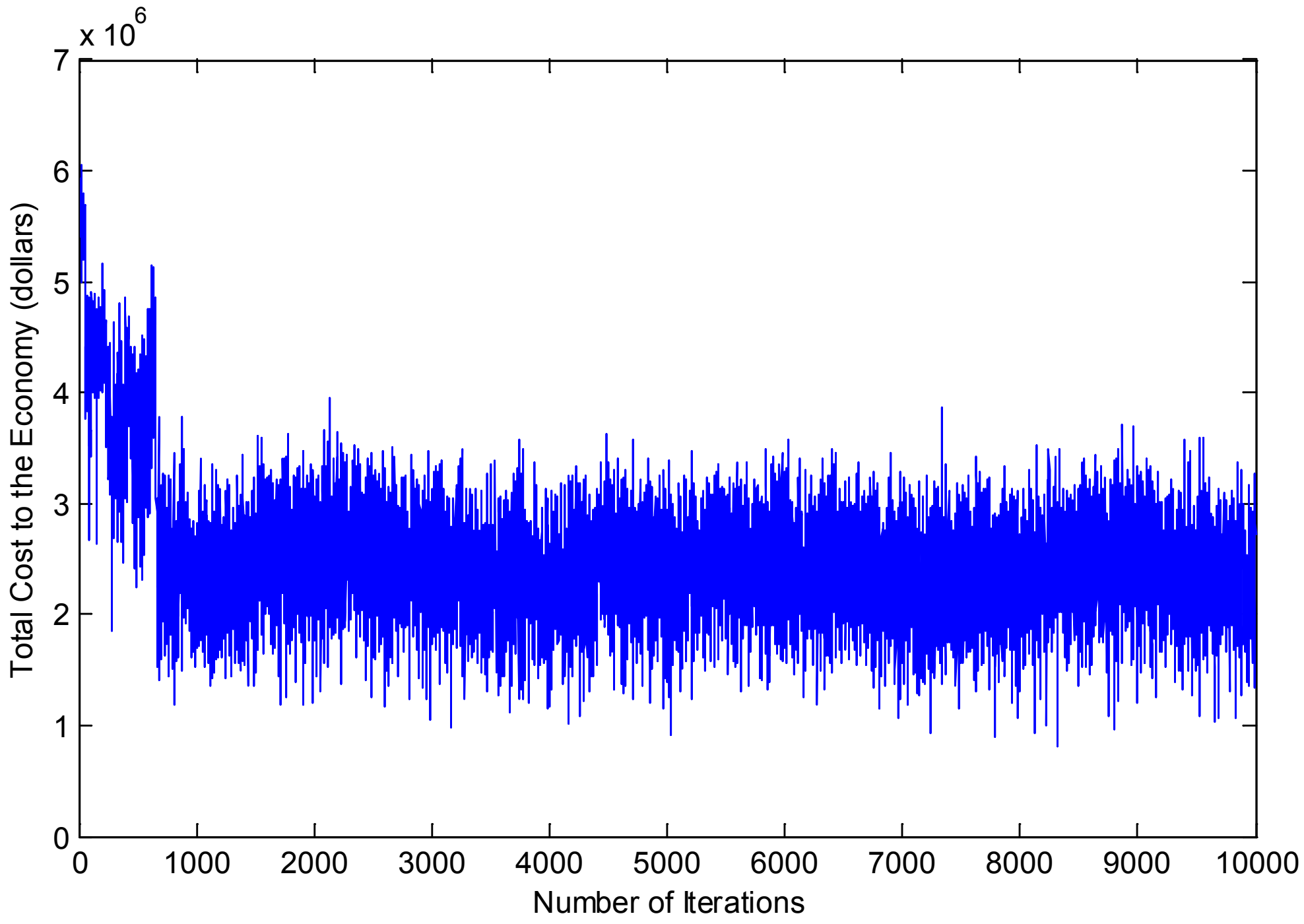

**Figure 7.5** Noisy loss values for DSPSA run for the late vaccination case (reality case). The solution indicates that for the late vaccination case high risk adults may have higher priority to take the vaccination.

In summary, because of the limited availability of the vaccines, for both early vaccination case and late vaccination case, our results show that it is preferable to use HHTAP for the antivirals policy. Even though the antiviral agents only take effect during a short period, they are still very useful when vaccination is limited. The cost of school closure is significant, so we would not use the intervention of school closure. Compared with late vaccination case, we find early vaccination can help to reduce more loss to the economy. All these results are consistent with the sensitivity analysis in Section 7.6. In regard to the priorities of the vaccinations, we would prefer to give school age children higher priority in the early vaccination case, and give high risk adults higher priority in the late vaccination case. Generally, in dealing with the epidemics intervention strategy problem, most researchers do the sensitivity analysis to achieve a relatively better solution, but they do not solve this problem by using a stochastic optimization algorithm (e.g. Halder et al., 2011, Khazeni et al., 2009). On the contrary, DSPSA is an algorithm that solves the optimal combination of interventions to achieve the least economic loss to society under the noisy measurements. Thus, DSPSA can provide a better result than just doing the sensitivity analysis.

# Chapter 8

# Conclusions

In this chapter, we summarize the results of this dissertation and discuss some related future work. In Section 8.1, we review the research problem and our contributions in this dissertation. In Section 8.2, we consider the potential opportunities for future research in DSPSA.

## 8.1 Summary

Discrete stochastic optimization problems are very important in many real-world applications, such as transmission problems in networks, facility locating problems, and resource allocation problems. Two key characteristics of discrete stochastic optimization problems are: 1) the domain is defined on the discrete points; 2) there is noise in the measurements of the loss function. Due to the discrete domain, many well developed algorithms for continuous stochastic optimization problems, such as stochastic approximation algorithms, cannot be used directly. In addition, the noise makes many algorithms designed for deterministic discrete optimization problems not applicable directly.

Nevertheless, some algorithms have been designed in order to solve the discrete stochastic optimization problems. There are three major classes of algorithms that include statistical approaches type algorithms, random search type algorithms, and stochastic approximation type algorithms. For the statistical approaches type, the algorithms can only guarantee that the solution is the optimal one with some confidence probability, which is less than 1. For the random search type, the generated sequence can converge to the optimal solution under some conditions. Multiple comparisons are needed in each step, and the number of comparisons in each iteration goes to infinity as $k \to \infty$. For the stochastic approximation type (designed for discrete problems), the basic idea is to implicitly make use of the function structure. Some stochastic approximation type algorithms focus on constructing the continuous extension of the discrete loss function, which involves many noisy measurements in each iteration to get the estimate of the gradient (the number of noisy measurements in each iteration could be proportional to the problem dimension). Then, for high-dimensional problem, the cost of noisy measurements of the loss function might be high. For the other stochastic approximation type algorithms, such as Hill et al. (2004), the number of noisy loss function measurements in each iteration is only two, but the theoretical analysis on the convergence properties for these algorithms are not well developed.

By the “No Free Lunch Theorem” (Spall 2003, Subsection 1.2.2), we know no one algorithm can beat all others across a broad range of problem types. We need to recognize and manage the trade-off between the robustness and efficiency. Here we mainly focus on the part of efficiency and introduce the algorithm of DSPSA in Section 2.1. The algorithm of DSPSA is a new discrete version of SPSA. DSPSA implicitly makes use of the function structure to provide efficient performance for the loss functions that satisfy some sufficient conditions. However, at the same time the robustness of DSPSA is sacrificed because of these conditions, which indicates that for some loss functions DSPSA may not provide sequences that converge to the optimal solution. Overall, DSPSA inherits several good properties of SPSA algorithm, but it also introduces some difficulties.

Besides introducing the new algorithm DSPSA, we also show the theoretical analysis of convergence properties of DSPSA. In Section 2.2, we show that under some conditions, the sequence generated by DSPSA converges to the optimal solution. These conditions are only sufficient but not necessary. In Section 4.5, the numerical test on the skewed quartic loss function indicates that even though the skewed quartic loss function does not satisfy the sufficient conditions, DSPSA still can provide a sequence convergent to the optimal solution.

After showing that the DSPSA algorithm converges, we further discuss in Chapter 3 the rate of convergence of DSPSA. Rate of convergence analysis is not an easy job. Different algorithms use different criteria when considering the convergence rate. Since DSPSA inherits the idea of SPSA and the elements in the sequence $\{\hat{\boldsymbol{\theta}}_k\}$ may not be multivariate integer points, it is natural to use the mean square error as the criterion to discuss the convergent rate for DSPSA. Under some conditions, we determine an upper bound for the mean square error. Based on this upper bound, we analyze the finite sample performance and asymptotical performance of DSPSA. We find the upper bound is composed of two terms: 1) a term related to the initial guess, and 2) a term related to the extra error introduced in each iteration. The first term is more significant in the early iterations and the second term starts to be significant in the later iterations, which implies that different selections of coefficients of gain sequence are preferred in different stages of DSPSA performance. We determine the practical selection and asymptotically optimal selection of values of coefficients for DSPSA based on the criterion of mean square error. These guidelines of coefficient selections can help people choose the appropriate set of coefficients in a practical problem. In the numerical tests in Chapter 4, we follow these guidelines to pick the coefficients. Furthermore, based on the upper bound, we show that the rate of convergence of DSPSA is $O(1/k^{\alpha})$, where $\alpha$ is the decaying rate of the gain sequence. This rate of convergence result (in the big-*O* sense) helps us to compare DSPSA with other discrete stochastic optimization algorithms.

Many of the discrete stochastic algorithms other than DSPSA are random search type algorithms. Moreover, for random search type algorithms, the sequence converges to the optimal solution, and the theoretical analysis is well developed. The stochastic ruler algorithm and the stochastic comparison algorithm are two basic representatives of the random search type algorithms. Therefore, we pick the stochastic ruler algorithm and the stochastic comparison algorithm to do the comparison with DSPSA. As we have mentioned, different types of algorithms use different criteria when considering the rate of convergence, so we set up a bridge between DSPSA and the other two random search type algorithms by considering the criterion $P([\hat{\boldsymbol{\theta}}_k] \neq \boldsymbol{\theta}^*)$. We calculate $P([\hat{\boldsymbol{\theta}}_k] \neq \boldsymbol{\theta}^*)$ for the three algorithms (DSPSA, stochastic ruler algorithm, and stochastic comparison algorithm) in the big-*O* sense. We compare the three algorithms theoretically in Chapter 5 and numerically in Chapter 6.

After finishing these theoretical analyses and numerical tests, we discuss the application of DSPSA towards developing optimal public health strategies for containing the spread of influenza, given limited societal resources. We use an online free simulator (FluTE) to simulate the real process of virus spreading. Due to the randomness in the synthetic social networks and transmission of the diseases, the output of FluTE involves noise. Furthermore, due to the complexity of the model in the simulator, it takes a long time to do a single simulation for a large population. We use DSPSA to solve for the optimal intervention policies in a public health decision problem related to the H1N1 virus. The objective function is defined as the total economic loss of H1N1 to society. The goal is to solve for the optimal intervention strategy, including vaccination priority, antiviral agent policy and time of school closure, to minimize the total economic cost of H1N1. Some researchers have used FluTE to do sensitivity analysis to discuss the effectiveness of some intervention policies. However, as to our knowledge, no one has used the discrete stochastic optimization algorithms to solve the problem simulated by FluTE directly. Our attempt in solving the problem by DSPSA is new for the simulator FluTE.

Overall, this dissertation introduces a new algorithm DSPSA, discusses the convergence properties, compares it with other algorithms, and uses it for an application in a public health problem. At last, let us summarize the good properties and the difficulties of DSPSA. The good properties of DSPSA include: 1) DSPSA is a simple algorithm and easily implemented in software; 2) The number of "tuning" coefficients to be picked is small; 3) DSPSA implicitly makes use of loss function structure, which leads to very efficient performance for some loss functions that satisfy some sufficient conditions; 4) The number of noisy measurements of the loss function is only two in each iteration; 5) The theoretical analysis of convergence properties are available in this dissertation. The difficulties of DSPSA include: 1) DSPSA is only a locally convergent algorithm; 2) Not all sufficient conditions are easy to check. 3) For non-integer-grid domain, we need to do reformulation before using DSPSA. These difficulties are also true for many other discrete stochastic algorithms.

In addition to the main body of the thesis, there are three appendices. In Appendix A, we consider the analysis of practical step size selection in stochastic approximation algorithms for continuous problem settings. The practical gain sequence selection is different from the theoretical optimal selection (derived from asymptotical performance). We provide a formal justification of the reasons why we choose this gain sequence in practice. In Appendix B, we consider the rate of convergence of SPSA for time-varying loss functions. One important application of time-varying loss function is in the model-free adaptive control with nonlinear stochastic systems, and model-free adaptive control is useful in many practical areas. Therefore, the results in Appendix B show the reasonable performance of SPSA in model-free control in the big-*O* sense. In Appendix C, we do the numerical experiments on the properties of the upper bound for $E\left\|\hat{\boldsymbol{\theta}}_k - \boldsymbol{\theta}^*\right\|^2$ discussed in Chapter 3. We show that the numerical results in Appendix C are consistent with the theoretical analysis in Section 3.2.

# 8.2 Future Work

In this dissertation, there are some potential issues that have not been considered or solved, and these problems can provide direction for future work.

The first issue is the consideration of mixed discrete-continuous problems. In such mixed problems, some variables are constrained to be integers while some variables can be non-integers. For these kinds of problems, we may have several research directions. First, we can divide the set of variables into two groups, with one group containing integer variables and the other group containing non-integer variables. For the group of integer variables, we can use DSPSA directly, and for the group of non-integer variables, SPSA (Spall, 1992) can be used directly. The second possible way is to discretize the continuous domain and use DSPSA directly.

The second issue is the problem that we have already discussed in Section 2.1 on the discrete mapping problem. For DSPSA, we assume that the domain is the subset of $\mathbb{Z}^p$. However, for some problems, the original domain of the loss function is discrete, but may not be the subset of $\mathbb{Z}^p$ (e.g. the original domain can be $\{\ldots, -1.5, -0.5, 0.5, 1.5, \ldots\}$ or

{…, apple, banana, lemon, … }). Therefore, we need to map the original domain to $\mathbb{Z}^p$ before using DSPSA. We have not found any paper that mainly discusses the problem of mapping an arbitrary discrete domain to $\mathbb{Z}^p$. Therefore, we think it is an interesting and important problem to be considered in the future.

The third issue that we are interested in is related to the second one, which also focuses on expanding the applicability of DSPSA to more general problems. In the second problem, we focus on mapping the discrete domain to $\mathbb{Z}^p$ to make DSPSA applicable for the discrete problem that is not defined on $\mathbb{Z}^p$. Here we can consider modifying the DSPSA algorithm description to make it applicable for more general types of domains, not only integer. The basic goal of the third problem is similar to the second one, but the two go in different research directions.

The fourth issue is the possible extension of the perturbation direction $\boldsymbol{\Delta}$. In the numerical tests of this thesis, we mainly consider the case when the components of $\boldsymbol{\Delta}$ are independent Bernoulli random variables taking the values $\pm 1$ with probability $1/2$. We do just one numerical test on the case when the components of $\boldsymbol{\Delta}$ are discrete uniformly distributed over the set $\{\pm 1, \pm 3\}$ in Section 4.6. In the future, we can analyze more non-Bernoulli distribution cases for $\boldsymbol{\Delta}$ and determine the optimal distribution for the perturbation directions. For the algorithm of SPSA on the continuous settings, Sadegh and Spall (1998) find that the asymptotically optimal distribution for the components of the simultaneous perturbation vector is a symmetric Bernoulli distribution. However, Cao (2011) shows that a non-Bernoulli distribution can achieve better finite sample performance than the Bernoulli distribution for some loss functions. Since DSPSA inherits many properties of SPSA, we believe that we can do similar analysis on DSPSA.

The fifth issue is how to use all the old information up to the current iteration. There are several possible ways that we can make use of the old information. For the first way, when the gain sequence is small, the sequence $\{\hat{\boldsymbol{\theta}}_k\}$ moves within one unit hypercube for a while, then it might be good to use old noisy measurements of the loss function in that unit hypercube to decrease the effect of noise in DSPSA. The second way could be that we can set the solution at iteration $k$ to be $\sum_{i=0}^{k} \hat{\boldsymbol{\theta}}_k \Big/ (k+1)$ (same as iterate averaging in Section 4.5.3 of Spall, 2003), which might improve the performance of DSPSA when the sequence bounces around the optimal solution. The third way could be similar to the idea of Andradottir (1999), where an extra optimization step is added to determine the current solution based on all old information (the current solution is the point with the smallest average noisy loss function based on all old information).

In summary, this dissertation introduces the DSPSA algorithm, discusses its convergence properties, compares it with other discrete stochastic optimization algorithm, and uses it in a public health problem. In the future, we want to follow the research directions discussed above to further generalize the DSPSA algorithm, and analyze more properties related to DSPSA.

# Appendix A

# Analysis of Practical Step Size Selection in Stochastic Approximation Algorithms

For many popular stochastic approximation algorithms, such as the stochastic gradient method and the simultaneous perturbation stochastic approximation method, the practical gain sequence selection is different from the optimal selection, which is theoretically derived from asymptotical performance. We provide formal justification for the reasons why we choose such gain sequence in practice.

# A.1 Introduction

Stochastic approximation (SA) algorithms are widely used in many stochastic optimization problems (e.g., Spall, 2003, Kushner and Yin, 2003, Fu, 1994). The stochastic gradient (SG) algorithm and the simultaneous perturbation stochastic approximation (SPSA) algorithm are two examples of SA methods. All SA algorithms include step-size parameters, frequently referred to as "gain sequences." The gain sequences are critical for practical implementation and asymptotic analysis. This appendix provides a rigorous basis for practical gain selection principles that appear in the literature relative to finite-sample implementation of SA algorithms. Although such principles have been presented as informal "rules of thumb," we show that they have a rigorous foundation. Hence, the finite-sample theory here complements well-known asymptotic theory for SA.

The general form of the SA algorithm is

$$\hat{\boldsymbol{\theta}}_{k+1} = \hat{\boldsymbol{\theta}}_k - a_k \hat{\boldsymbol{g}}_k(\hat{\boldsymbol{\theta}}_k),$$

where $\hat{\boldsymbol{\theta}}_k$ is the estimated optimal point from iteration *k*, $\hat{\boldsymbol{g}}_k(\hat{\boldsymbol{\theta}}_k)$ is the estimator of the gradient at $\hat{\boldsymbol{\theta}}_k$, and the gain sequence (step size) is $a_k = a/(1+A+k)^{\alpha}$. Rubinstein and Kroese (2007, p. 215) say: "The crucial question in implementations is the choice of the step sizes." We see that the choice of coefficients is a very important issue.

For the stochastic gradient algorithm, when $\alpha$ = 1, the rate of convergence of $\{\hat{\boldsymbol{\theta}}_k\}$ is maximized (Section 4.4 in Spall, 2003). But in practical problems, $\alpha$ = 1 may not be the optimal choice. Many authors use $\alpha$ = 0.501 in their numerical experiments, such as Hutchison and Spall (2009). In the Section 4.4 of Spall (2003), the author also discusses the roles of coefficients $a$ and $A$ briefly.

For SPSA, Spall (1992) shows that the theoretical optimal choice of gain sequence through asymptotical distribution is to pick $\alpha = 1$. Later Spall (1998b, 2003) discusses in detail about the rule to pick up the coefficients for real world problems with finite number of iterations, where $\alpha$ equals 0.602, $A$ equals 10% or less of the maximum number of allowed iterations, and $a$ is chosen such that $a/(1+A)^{0.602}$ $\times$[magnitude of components in gradient estimator] is equal to the smallest of the desired change magnitudes of the sequence in the early iterations. After checking many applications that use SPSA, we find that people often directly use Spall's rule of gain sequence (e.g., Bangerth et al., 2006, Maryak and Spall, 2005, and Schwartz et al., 2006).

To the best of our knowledge, for both SG and SPSA, almost no paper provides formal reasons and conditions for the practical gain sequence selection, even though many people use the practical choice directly. In this appendix, we provide a mathematical justification for why practical gain sequence selections are different from theoretical ones. In the process, we are able to provide more precise guidance with respect to particular problem characteristics. In fact, it appears that this appendix is one of very few that provide finite-sample theory for SA (another paper with finite-sample theory relative to a different aspect of SA is Cao, 2011). In Section A.2, we discuss the practical gain sequence selection for SG and SPSA, and in Section A.3, numerical results are presented. This appendix is concluded by Section A.4.

# A.2 Formal Analysis

In this section we provide mathematical justifications on the practical gain sequence selections for both SG and SPSA. For both algorithms, the loss function $L(\boldsymbol{\theta})$ that we consider is strictly convex.

## A.2.1 Stochastic Gradient

The formula for the stochastic gradient algorithm is given by

$$\hat{\boldsymbol{\theta}}_{k+1} = \hat{\boldsymbol{\theta}}_k - a_k \hat{\boldsymbol{g}}_k(\hat{\boldsymbol{\theta}}_k), \tag{A.1}$$

where $\hat{\boldsymbol{g}}_k(\hat{\boldsymbol{\theta}}_k)$ is the unbiased estimate of gradient. Suppose $\hat{\boldsymbol{g}}_k(\hat{\boldsymbol{\theta}}_k) = \boldsymbol{g}(\hat{\boldsymbol{\theta}}_k) + \boldsymbol{\varepsilon}_k(\hat{\boldsymbol{\theta}}_k)$, $\boldsymbol{g}(\hat{\boldsymbol{\theta}}_k)$ is the real gradient at point $\hat{\boldsymbol{\theta}}_k$, and $E\left(\boldsymbol{\varepsilon}_k(\hat{\boldsymbol{\theta}}_k) \middle| \hat{\boldsymbol{\theta}}_k\right) = 0$. Substituting $\hat{\boldsymbol{g}}_k(\hat{\boldsymbol{\theta}}_k)$ into eqn. (A.1), we have

$$\hat{\boldsymbol{\theta}}_{k+1} - \boldsymbol{\theta}^* = \hat{\boldsymbol{\theta}}_k - \boldsymbol{\theta}^* - a_k \boldsymbol{g}(\hat{\boldsymbol{\theta}}_k) + a_k \boldsymbol{\varepsilon}_k(\hat{\boldsymbol{\theta}}_k). \tag{A.2}$$

By calculating the norm to the squared of both sides of eqn. (A.2) and taking expectation on them, we have

$$\begin{aligned} E\left\|\hat{\boldsymbol{\theta}}_{k+1} - \boldsymbol{\theta}^*\right\|^2 &= E\left\|\hat{\boldsymbol{\theta}}_k - \boldsymbol{\theta}^*\right\|^2 - 2a_k E\left[(\hat{\boldsymbol{\theta}}_k - \boldsymbol{\theta}^*)^T \boldsymbol{g}(\hat{\boldsymbol{\theta}}_k)\right] + 2a_k E\left[(\hat{\boldsymbol{\theta}}_k - \boldsymbol{\theta}^*)^T \boldsymbol{\varepsilon}_k(\hat{\boldsymbol{\theta}}_k)\right] \\ &\quad + a_k^2 E\left\|\hat{\boldsymbol{g}}_k(\hat{\boldsymbol{\theta}}_k)\right\|^2. \end{aligned} \tag{A.3}$$

Suppose $L(\boldsymbol{\theta})$ is twice continuously differentiable, and $\Gamma(\boldsymbol{\theta})$ is the Hessian matrix of $L(\boldsymbol{\theta})$. Since $L(\boldsymbol{\theta})$ is a strictly convex function, $\Gamma(\boldsymbol{\theta})$ is a positive definite matrix. For $\boldsymbol{g}(\hat{\boldsymbol{\theta}}_k)$, there exists $\tilde{\boldsymbol{\theta}}_k$ on the line segment of $\hat{\boldsymbol{\theta}}_k$ and $\boldsymbol{\theta}^*$ such that

$$\boldsymbol{g}(\hat{\boldsymbol{\theta}}_k) = \Gamma(\tilde{\boldsymbol{\theta}}_k)(\hat{\boldsymbol{\theta}}_k - \boldsymbol{\theta}^*).$$

Then, we have

$$E\left[(\hat{\boldsymbol{\theta}}_k - \boldsymbol{\theta}^*)^T \boldsymbol{g}(\hat{\boldsymbol{\theta}}_k)\right] = E\left[(\hat{\boldsymbol{\theta}}_k - \boldsymbol{\theta}^*)^T \Gamma(\tilde{\boldsymbol{\theta}}_k)(\hat{\boldsymbol{\theta}}_k - \boldsymbol{\theta}^*)\right],$$

which is positive. By Corollary 5.4.5 of Horn and Johnson (1985), there exists $\mu_k > 0$ such that

$$E\left[(\hat{\boldsymbol{\theta}}_k - \boldsymbol{\theta}^*)^T \Gamma(\tilde{\boldsymbol{\theta}}_k)(\hat{\boldsymbol{\theta}}_k - \boldsymbol{\theta}^*)\right] = \mu_k E\left\|\hat{\boldsymbol{\theta}}_k - \boldsymbol{\theta}^*\right\|^2. \tag{A.4}$$

In addition, the third term on the right-hand side of eqn. (A.3) can be written as

$$\begin{aligned} 2a_k E\left[(\hat{\boldsymbol{\theta}}_k - \boldsymbol{\theta}^*)^T \boldsymbol{\varepsilon}_k(\hat{\boldsymbol{\theta}}_k)\right] &= 2a_k E\left[E\left((\hat{\boldsymbol{\theta}}_k - \boldsymbol{\theta}^*)^T \boldsymbol{\varepsilon}_k(\hat{\boldsymbol{\theta}}_k)\middle|\hat{\boldsymbol{\theta}}_k\right)\right] \\ &= 2a_k E\left[(\hat{\boldsymbol{\theta}}_k - \boldsymbol{\theta}^*)^T E\left(\boldsymbol{\varepsilon}_k(\hat{\boldsymbol{\theta}}_k)\middle|\hat{\boldsymbol{\theta}}_k\right)\right] \\ &= 0. \end{aligned} \tag{A.5}$$

Furthermore, the last term of the right-hand side of eqn. (A.3) is

$$\begin{aligned} a_k^2 E\left\|\hat{\boldsymbol{g}}_k(\hat{\boldsymbol{\theta}}_k)\right\|^2 &= a_k^2 E\left(\left\|\boldsymbol{g}(\hat{\boldsymbol{\theta}}_k) + \boldsymbol{\varepsilon}_k(\hat{\boldsymbol{\theta}}_k)\right\|^2\right) \\ &= a_k^2\left(E\left\|\boldsymbol{g}(\hat{\boldsymbol{\theta}}_k)\right\|^2 + 2E\left(\boldsymbol{g}(\hat{\boldsymbol{\theta}}_k)^T \boldsymbol{\varepsilon}_k(\hat{\boldsymbol{\theta}}_k)\right) + E\left\|\boldsymbol{\varepsilon}_k(\hat{\boldsymbol{\theta}}_k)\right\|^2\right). \end{aligned} \tag{A.6}$$

By the law of total expectation, we have

$$E\left(\boldsymbol{g}(\hat{\boldsymbol{\theta}}_k)^T \boldsymbol{\varepsilon}_k(\hat{\boldsymbol{\theta}}_k)\right) = E\left(E\left(\boldsymbol{g}(\hat{\boldsymbol{\theta}}_k)^T \boldsymbol{\varepsilon}_k(\hat{\boldsymbol{\theta}}_k)\middle|\hat{\boldsymbol{\theta}}_k\right)\right)$$
$$= E\left(\boldsymbol{g}(\hat{\boldsymbol{\theta}}_k)^T E\left(\boldsymbol{\varepsilon}_k(\hat{\boldsymbol{\theta}}_k)\middle|\hat{\boldsymbol{\theta}}_k\right)\right)$$
$$= 0.$$

Moreover,

$$E\left\|\boldsymbol{g}(\hat{\boldsymbol{\theta}}_k)\right\|^2 = E\left((\hat{\boldsymbol{\theta}}_k - \boldsymbol{\theta}^*)^T \Gamma^T(\tilde{\boldsymbol{\theta}}_k)\Gamma(\tilde{\boldsymbol{\theta}}_k)(\hat{\boldsymbol{\theta}}_k - \boldsymbol{\theta}^*)\right).$$

Thus, eqn. (A.6) can be written as

$$a_k^2 E\left\|\hat{\boldsymbol{g}}_k(\hat{\boldsymbol{\theta}}_k)\right\|^2 = a_k^2 E\left((\hat{\boldsymbol{\theta}}_k - \boldsymbol{\theta}^*)^T \Gamma^T(\tilde{\boldsymbol{\theta}}_k)\Gamma(\tilde{\boldsymbol{\theta}}_k)(\hat{\boldsymbol{\theta}}_k - \boldsymbol{\theta}^*)\right)$$
$$+ a_k^2 E\left\|\boldsymbol{\varepsilon}_k(\hat{\boldsymbol{\theta}}_k)\right\|^2. \quad \text{(A.7)}$$

Plugging the results of eqn. (A.4), (A.5) and (A.7) into eqn. (A.3), we have

$$E\left\|\hat{\boldsymbol{\theta}}_{k+1} - \boldsymbol{\theta}^*\right\|^2 = (1 - 2a_k\mu_k)E\left\|\hat{\boldsymbol{\theta}}_k - \boldsymbol{\theta}^*\right\|^2 + a_k^2 E\left((\hat{\boldsymbol{\theta}}_k - \boldsymbol{\theta}^*)^T \Gamma^T(\tilde{\boldsymbol{\theta}}_k)\Gamma(\tilde{\boldsymbol{\theta}}_k)(\hat{\boldsymbol{\theta}}_k - \boldsymbol{\theta}^*)\right)$$
$$+ a_k^2 E\left\|\boldsymbol{\varepsilon}_k(\hat{\boldsymbol{\theta}}_k)\right\|^2. \quad \text{(A.8)}$$

By the recursive relationship in eqn. (A.8), we get

$$E\left\|\hat{\boldsymbol{\theta}}_{k+1} - \boldsymbol{\theta}^*\right\|^2 = \prod_{j=0}^{k}(1 - 2a_j\mu_j)E\left\|\hat{\boldsymbol{\theta}}_0 - \boldsymbol{\theta}^*\right\|^2$$
$$+ \sum_{i=0}^{k}\prod_{j=i+1}^{k}(1 - 2a_j\mu_j)a_i^2 E\left((\hat{\boldsymbol{\theta}}_i - \boldsymbol{\theta}^*)^T \Gamma^T(\tilde{\boldsymbol{\theta}}_i)\Gamma(\tilde{\boldsymbol{\theta}}_i)(\hat{\boldsymbol{\theta}}_i - \boldsymbol{\theta}^*)\right)$$
$$+ \sum_{i=0}^{k}\prod_{j=i+1}^{k}(1 - 2a_j\mu_j)a_i^2 E\left\|\boldsymbol{\varepsilon}_i(\hat{\boldsymbol{\theta}}_i)\right\|^2. \quad \text{(A.9)}$$

The first term in eqn. (A.9) is corresponding to the initial guess; the second term and the third term are positive. In addition we know

$$a_i = \frac{a}{(1 + A + i)^{\alpha}},$$

and $0.5 < \alpha \le 1$. Suppose $\hat{\boldsymbol{\theta}}_0$ is far away from the optimal solution $\boldsymbol{\theta}^*$ and $\mu_k$ does not change significantly in the early iterations for different sets of coefficients, which means the change of $\mu_k$ from different sets of coefficients is insignificant compared to the value of $\mu_k$ itself. The loss functions satisfying the assumptions have smooth Hessian matrix when point is far from the optimal solution. In order to have stable performance in the early iterations, the last two terms of eqn. (A.8) should have small values. It can be achieved by small value of step size, which is equivalent to control the value of *a* not too big and the value of *A* not too small. The first term of eqn. (A.8) is significant due to the initial guess term. The rest terms of eqn. (A.8) have coefficients $a_i^2$, and it makes the last two terms less significant when the step size is small. Then, the first term on the right-hand side of eqn. (A.9) is the dominant term in the early iterations, so we prefer bigger value of $a_i$, which can be achieved by smaller value of $\alpha$ ($\alpha = 0.501$), larger value of *a* and smaller value of *A*. We see that the requirements for *a* and *A* are on the opposite sides for stability and better finite sample performance. Therefore, we need to pick *a* relative large and pick *A* relative small in the domain that satisfies the stability requirement.

Now we discuss a detailed rule for the choice of *a* and *A*. Due to the restriction of the budget, the number of iterations may be limited. In order to achieve reasonable performance for the limited number of iterations, we want *A* to be proportional to the maximum number of allowed iterations *N*: $A = \eta N$. In addition, we want that the effect of $A$ would disappear in later iterations to achieve proper decaying gain. Therefore, we prefer $A$ to be in the lower rate than $N$, which indicates $\eta$ = 0.1, 0.01, 0.001,... Moreover, as we have discussed, we want $a/(1+A+i)^{\alpha}$ to be reasonable small in the early iteration to keep the stability of the algorithm, but we also do not want the gain step too small in the later iteration to keep reasonable speed. We know the effect of $A$ disappears when $k$ is large; while the effect of $a$ is always there. It means that we prefer to pick larger $A$ instead of smaller $a$ to keep the stability of the algorithm, because we still need $a$ not to be too small for better finite sample performance. Therefore, among all the choices of $\eta$, we prefer larger value, which is $\eta = 0.1$. It indicates that $A = 0.1N$ satisfies all the requirements, including the stability and finite sample performance. After getting the rule of $A$, we can get the rule for $a$ by considering the performance in the early iterations of the algorithm. For example, we can make the multiplication of $a_0$ and approximated magnitude value of $\hat{\boldsymbol{g}}_0(\hat{\boldsymbol{\theta}}_0)$ to equal the desired change magnitude in the early iterations, and it leads the possible choice of $a$.

In contrast, for good asymptotic performance of the algorithm, we want the sequence to bounce as little as possible around the optimal solution, which is achieved by small values of the gain sequence in the very late iterations. It is equivalent to say for asymptotical performance, we prefer big $\alpha$ ($\alpha$ = 1). Thus, we see that the theoretical optimal selection of coefficients is for better asymptotical performance, while the practical selection of coefficients focuses on better finite sample performance.

## A.2.2 Simultaneous Perturbation Stochastic Approximation

The basic formula for SPSA algorithm discussed in Spall (1992) is

$$\hat{\boldsymbol{\theta}}_{k+1} = \hat{\boldsymbol{\theta}}_k - a_k \hat{\boldsymbol{g}}_k(\hat{\boldsymbol{\theta}}_k),$$

where $\hat{\boldsymbol{\theta}}_k$ is the estimate for the optimal solution $\boldsymbol{\theta}^*$ at the $k$th iteration, $\{a_k\}$ is the gain sequence and

$$\hat{\boldsymbol{g}}_k(\hat{\boldsymbol{\theta}}_k) = \frac{y(\hat{\boldsymbol{\theta}}_k + c_k \boldsymbol{\Delta}_k) - y(\hat{\boldsymbol{\theta}}_k - c_k \boldsymbol{\Delta}_k)}{2c_k} \boldsymbol{\Delta}_k^{-1},$$

where $y(\cdot)$ is the noisy measurement of the loss function $L(\cdot)$, $c_k$ is a positive scalar, $\boldsymbol{\Delta}_k = \left(\Delta_{k1},...,\Delta_{kp}\right)$ is the direction of perturbation, and $\boldsymbol{\Delta}_k^{-1} = \left(1/\Delta_{k1},...,1/\Delta_{kp}\right)$. We also assume that $y(\cdot) = L(\cdot) + \varepsilon(\cdot)$ and $E\left(\varepsilon \middle| \boldsymbol{\Delta}, \boldsymbol{\theta}\right) = 0$. Moreover, the form of $a_k$ and $c_k$ are

$$a_k = \frac{a}{(1 + A + k)^{\alpha}}, \quad c_k = \frac{c}{(1+k)^{\gamma}}.$$

The assumptions for $\alpha$ and $\gamma$ are discussed in Spall (1992): $\alpha \le 1$, $2\alpha - 2\gamma > 1$, and $3\gamma - \alpha/2 \ge 0$. Spall (1992) also shows that the asymptotical optimal coefficient choices are $\alpha = 1$ and $\gamma = 1/6$, but in practical finite sample implementations, people often use $\alpha = 0.602$ and $\gamma = 0.101$. In the following, we show the reasons for the practical choice of $\alpha$, $\gamma$, and discuss the choice of the other coefficients of $a$, $A$, $c$.

First we add and subtract $a_k \boldsymbol{g}(\hat{\boldsymbol{\theta}}_k)$ to the right-hand side of the algorithm formula and get

$$\hat{\boldsymbol{\theta}}_{k+1} - \boldsymbol{\theta}^* = \hat{\boldsymbol{\theta}}_k - \boldsymbol{\theta}^* - a_k \boldsymbol{g}(\hat{\boldsymbol{\theta}}_k) + a_k \boldsymbol{g}(\hat{\boldsymbol{\theta}}_k) - a_k \hat{\boldsymbol{g}}_k(\hat{\boldsymbol{\theta}}_k), \tag{A.10}$$

where $\boldsymbol{g}(\hat{\boldsymbol{\theta}}_k)$ is the real gradient of function $L(\cdot)$ at point $\hat{\boldsymbol{\theta}}_k$. In Spall (1992) the loss function is assumed to be strictly convex. By calculating the norm to the squared of both sides of eqn. (A.10) and taking the expectation on them, we have

$$\begin{aligned} E\left\|\hat{\boldsymbol{\theta}}_{k+1} - \boldsymbol{\theta}^*\right\|^2 &= E\left\|\hat{\boldsymbol{\theta}}_k - \boldsymbol{\theta}^*\right\|^2 - 2a_k E\left[(\hat{\boldsymbol{\theta}}_k - \boldsymbol{\theta}^*)^T \boldsymbol{g}(\hat{\boldsymbol{\theta}}_k)\right] \\ &\quad + 2a_k E\left[(\hat{\boldsymbol{\theta}}_k - \boldsymbol{\theta}^*)^T \left(\boldsymbol{g}(\hat{\boldsymbol{\theta}}_k) - \hat{\boldsymbol{g}}_k(\hat{\boldsymbol{\theta}}_k)\right)\right] + a_k^2 E\left\|\hat{\boldsymbol{g}}_k(\hat{\boldsymbol{\theta}}_k)\right\|^2. \end{aligned} \tag{A.11}$$

Suppose $\Gamma(\boldsymbol{\theta})$ is the Hessian matrix of $L(\boldsymbol{\theta})$. Since $L(\boldsymbol{\theta})$ is a strictly convex function, then $\Gamma(\boldsymbol{\theta})$ is a positive definite matrix. For $\boldsymbol{g}(\hat{\boldsymbol{\theta}}_k)$, there exists $\tilde{\boldsymbol{\theta}}_k$ on the line segment of $\hat{\boldsymbol{\theta}}_k$ and $\boldsymbol{\theta}^*$ such that

$$\boldsymbol{g}(\hat{\boldsymbol{\theta}}_k) = \Gamma(\tilde{\boldsymbol{\theta}}_k)(\hat{\boldsymbol{\theta}}_k - \boldsymbol{\theta}^*).$$

Therefore, we have

$$E\left[(\hat{\boldsymbol{\theta}}_k - \boldsymbol{\theta}^*)^T \boldsymbol{g}(\hat{\boldsymbol{\theta}}_k)\right] = E\left[(\hat{\boldsymbol{\theta}}_k - \boldsymbol{\theta}^*)^T \Gamma(\tilde{\boldsymbol{\theta}}_k)(\hat{\boldsymbol{\theta}}_k - \boldsymbol{\theta}^*)\right],$$

which is positive. By Corollary 5.4.5 of Horn and Johnson (1985), there exists $\mu_k > 0$ such that

$$E\left[(\hat{\boldsymbol{\theta}}_k - \boldsymbol{\theta}^*)^T \Gamma(\tilde{\boldsymbol{\theta}}_k)(\hat{\boldsymbol{\theta}}_k - \boldsymbol{\theta}^*)\right] = \mu_k E\left\|\hat{\boldsymbol{\theta}}_k - \boldsymbol{\theta}^*\right\|^2. \tag{A.12}$$

In addition, the third term on the right-hand side of eqn. (A.11) can be written as

$$\begin{aligned}
&2a_k E\left[(\hat{\boldsymbol{\theta}}_k - \boldsymbol{\theta}^*)^T \left(\boldsymbol{g}(\hat{\boldsymbol{\theta}}_k) - \hat{\boldsymbol{g}}_k(\hat{\boldsymbol{\theta}}_k)\right)\right] \\
&\quad = 2a_k E\left[E\left((\hat{\boldsymbol{\theta}}_k - \boldsymbol{\theta}^*)^T \left(\boldsymbol{g}(\hat{\boldsymbol{\theta}}_k) - \hat{\boldsymbol{g}}_k(\hat{\boldsymbol{\theta}}_k)\right) \middle| \hat{\boldsymbol{\theta}}_k\right)\right] \\
&\quad = 2a_k E\left[(\hat{\boldsymbol{\theta}}_k - \boldsymbol{\theta}^*)^T E\left(\boldsymbol{g}(\hat{\boldsymbol{\theta}}_k) - \hat{\boldsymbol{g}}_k(\hat{\boldsymbol{\theta}}_k) \middle| \hat{\boldsymbol{\theta}}_k\right)\right] \\
&\quad = 2a_k c_k^2 E\left[(\hat{\boldsymbol{\theta}}_k - \boldsymbol{\theta}^*)^T \boldsymbol{b}_k(\hat{\boldsymbol{\theta}}_k)\right], \qquad \text{(A.13)}
\end{aligned}$$

where

$$\boldsymbol{b}_k(\hat{\boldsymbol{\theta}}_k) = \frac{1}{c_k^2} E\left(\boldsymbol{g}(\hat{\boldsymbol{\theta}}_k) - \hat{\boldsymbol{g}}_k(\hat{\boldsymbol{\theta}}_k) \middle| \hat{\boldsymbol{\theta}}_k\right),$$

$$b_{kl} = \frac{1}{12} E\left[\sum_{i_1}\sum_{i_2}\sum_{i_3}\left(L^{(3)}_{i_1,i_2,i_3}(\tilde{\tilde{\boldsymbol{\theta}}}_k) + L^{(3)}_{i_1,i_2,i_3}(\tilde{\tilde{\tilde{\boldsymbol{\theta}}}}_k)\right)\frac{\Delta_{ki_1}\Delta_{ki_2}\Delta_{ki_3}}{\Delta_{ki_l}}\right],$$

and $\tilde{\tilde{\boldsymbol{\theta}}}_k$, $\tilde{\tilde{\tilde{\boldsymbol{\theta}}}}_k$ are on the line segment between $\hat{\boldsymbol{\theta}}_k$ and $\hat{\boldsymbol{\theta}}_k \pm c_k \boldsymbol{\Delta}_k$. By the result of Lemma 1 in Spall (1992), we have $\boldsymbol{b}_k(\hat{\boldsymbol{\theta}}_k) = O(1)$. Furthermore, the last term of the right-hand side of eqn. (A.11) is

$$a_k^2 E\left\|\hat{\boldsymbol{g}}_k(\hat{\boldsymbol{\theta}}_k)\right\|^2 = a_k^2 E\left(\left\|\left(\frac{L(\hat{\boldsymbol{\theta}}_k + c_k\boldsymbol{\Delta}_k) - L(\hat{\boldsymbol{\theta}}_k - c_k\boldsymbol{\Delta}_k) + \varepsilon_k^+ - \varepsilon_k^-}{2c_k}\right)\boldsymbol{\Delta}_k^{-1}\right\|^2\right)$$

$$= a_k^2 E\left[\frac{\left(L(\hat{\boldsymbol{\theta}}_k + c_k\boldsymbol{\Delta}_k) - L(\hat{\boldsymbol{\theta}}_k - c_k\boldsymbol{\Delta}_k)\right)^2}{4c_k^2}\boldsymbol{\Delta}_k^{-T}\boldsymbol{\Delta}_k^{-1}\right]$$

$$+ \frac{a_k^2}{2c_k^2} E\left[\left(L(\hat{\boldsymbol{\theta}}_k + c_k\boldsymbol{\Delta}_k) - L(\hat{\boldsymbol{\theta}}_k - c_k\boldsymbol{\Delta}_k)\right)\left(\varepsilon_k^+ - \varepsilon_k^-\right)\boldsymbol{\Delta}_k^{-T}\boldsymbol{\Delta}_k^{-1}\right]$$

$$+ a_k^2 E\left[\frac{\left(\varepsilon_k^+ - \varepsilon_k^-\right)^2}{4c_k^2}\boldsymbol{\Delta}_k^{-T}\boldsymbol{\Delta}_k^{-1}\right]. \tag{A.14}$$

Moreover, by the law of total expectation, we have

$$E\left[\left(L(\hat{\boldsymbol{\theta}}_k + c_k\boldsymbol{\Delta}_k) - L(\hat{\boldsymbol{\theta}}_k - c_k\boldsymbol{\Delta}_k)\right)\left(\varepsilon_k^+ - \varepsilon_k^-\right)\boldsymbol{\Delta}_k^{-T}\boldsymbol{\Delta}_k^{-1}\right]$$

$$= E\left[E\left(\left(L(\hat{\boldsymbol{\theta}}_k + c_k\boldsymbol{\Delta}_k) - L(\hat{\boldsymbol{\theta}}_k - c_k\boldsymbol{\Delta}_k)\right)\left(\varepsilon_k^+ - \varepsilon_k^-\right)\boldsymbol{\Delta}_k^{-T}\boldsymbol{\Delta}_k^{-1}\middle|\hat{\boldsymbol{\theta}}_k, \boldsymbol{\Delta}_k\right)\right]$$

$$= E\left[\left(L(\hat{\boldsymbol{\theta}}_k + c_k\boldsymbol{\Delta}_k) - L(\hat{\boldsymbol{\theta}}_k - c_k\boldsymbol{\Delta}_k)\right)\boldsymbol{\Delta}_k^{-T}\boldsymbol{\Delta}_k^{-1} E\left(\varepsilon_k^+ - \varepsilon_k^-\middle|\hat{\boldsymbol{\theta}}_k, \boldsymbol{\Delta}_k\right)\right]$$

$$= 0.$$

Therefore, eqn. (A.14) can be written as

$$a_k^2 E\left\|\hat{\boldsymbol{g}}_k(\hat{\boldsymbol{\theta}}_k)\right\|^2 = a_k^2 E\left[\frac{\left(L(\hat{\boldsymbol{\theta}}_k + c_k\boldsymbol{\Delta}_k) - L(\hat{\boldsymbol{\theta}}_k - c_k\boldsymbol{\Delta}_k)\right)^2}{4c_k^2}\boldsymbol{\Delta}_k^{-T}\boldsymbol{\Delta}_k^{-1}\right]$$

$$+ a_k^2 E\left[\frac{\left(\varepsilon_k^+ - \varepsilon_k^-\right)^2}{4c_k^2}\boldsymbol{\Delta}_k^{-T}\boldsymbol{\Delta}_k^{-1}\right]. \tag{A.15}$$

Plugging the results of eqn. (A.12), (A.13) and (A.15) into eqn. (A.11), we have

$$E\left\|\hat{\boldsymbol{\theta}}_{k+1}-\boldsymbol{\theta}^*\right\|^2=(1-2a_k\mu_k)E\left\|\hat{\boldsymbol{\theta}}_k-\boldsymbol{\theta}^*\right\|^2+2a_kc_k^2E\left[(\hat{\boldsymbol{\theta}}_k-\boldsymbol{\theta}^*)^T\boldsymbol{b}_k(\hat{\boldsymbol{\theta}}_k)\right]$$
$$+a_k^2E\left[\frac{\left(L(\hat{\boldsymbol{\theta}}_k+c_k\boldsymbol{\Delta}_k)-L(\hat{\boldsymbol{\theta}}_k-c_k\boldsymbol{\Delta}_k)\right)^2}{4c_k^2}\boldsymbol{\Delta}_k^{-T}\boldsymbol{\Delta}_k^{-1}\right]$$
$$+a_k^2E\left[\frac{\left(\varepsilon_k^+-\varepsilon_k^-\right)^2}{4c_k^2}\boldsymbol{\Delta}_k^{-T}\boldsymbol{\Delta}_k^{-1}\right]. \tag{A.16}$$

By the recursive relationship in eqn. (A.16), we get

$$E\left\|\hat{\boldsymbol{\theta}}_{k+1}-\boldsymbol{\theta}^*\right\|^2=\prod_{j=0}^{k}(1-2a_j\mu_j)E\left\|\hat{\boldsymbol{\theta}}_0-\boldsymbol{\theta}^*\right\|^2$$
$$+2\sum_{i=0}^{k}\prod_{j=i+1}^{k}(1-2a_j\mu_j)a_ic_i^2E\left[(\hat{\boldsymbol{\theta}}_i-\boldsymbol{\theta}^*)^T\boldsymbol{b}_i(\hat{\boldsymbol{\theta}}_i)\right]$$
$$+\sum_{i=0}^{k}\prod_{j=i+1}^{k}(1-2a_j\mu_j)a_i^2E\left[\frac{\left(L(\hat{\boldsymbol{\theta}}_i+c_i\boldsymbol{\Delta}_i)-L(\hat{\boldsymbol{\theta}}_i-c_i\boldsymbol{\Delta}_i)\right)^2}{4c_i^2}\boldsymbol{\Delta}_i^{-T}\boldsymbol{\Delta}_i^{-1}\right]$$
$$+\sum_{i=0}^{k}\prod_{j=i+1}^{k}(1-2a_j\mu_j)a_i^2E\left[\frac{\left(\varepsilon_i^+-\varepsilon_i^-\right)^2}{4c_i^2}\boldsymbol{\Delta}_i^{-T}\boldsymbol{\Delta}_i^{-1}\right]. \tag{A.17}$$

The first term on the right-hand side of eqn. (A.17) is corresponding to the initial guess, and the third term and the forth term are positive.

Suppose $\hat{\boldsymbol{\theta}}_0$ is far away from the optimal solution $\boldsymbol{\theta}^*$, and $\mu_k$ does not change significantly on different set of coefficients, which means the change of $\mu_k$ from different sets of coefficients is insignificant compared to the value of $\mu_k$ itself. The loss functions satisfying the assumptions have smooth Hessian matrix when point is far from the optimal solution. In order to have stable performance, the last two terms of eqn. (A.16) should have reasonable small values, which can be achieved by controlling the value of $a$ not too big and controlling the value of $A$ not too small. The first term in eqn. (A.16) is significant due to the initial guess part and the rest terms in eqn. (A.16) have corresponding coefficients $a_i c_i^2$ and $a_i^2 / c_i^2$ ($2\alpha - 2\gamma > 1$), which makes the last three terms less significant for small value of $a_i$ and $c_i$. Then, the first term on the right-hand side of eqn. (A.17) is the dominant term for finite sample case, so we prefer bigger value of $a_i$ which can be achieved by smaller value of $\alpha$ ($\alpha = 0.602$), larger value of $a$, and smaller value of $A$. We can see the stability requirements and better finite sample requirements are on the opposite sides for $a$ and $A$, so we need to pick $a$ relative large and pick $A$ relative small in the domain that satisfies the stability requirements. The detailed rule of the choice of $a$ and $A$ is the same as that for the algorithm of SG.

The feasible domain for $\alpha$ and $\gamma$ is a triangle composed of three lines $\alpha \le 1$, $2\alpha - 2\gamma > 1$, and $3\gamma - \alpha/2 \ge 0$. The triangle is equivalently defined by one achievable corner point, $(\gamma, \alpha) = (1/6, 1)$ and two unachievable corner points, (0.1, 0.6) and (1/2, 1). In a numerical implementation, we see that when $\alpha = 0.602$, we have to choose $\gamma = 0.101$ approximately. This pair of values is very close to the unachievable extreme point of the triangle domain with the smallest coefficient values, (0.1, 0.6). In addition

$$E\left[\frac{(\varepsilon_i^+ - \varepsilon_i^-)^2}{4c_i^2}\boldsymbol{\Delta}_i^{-T}\boldsymbol{\Delta}_i^{-1}\right] = E\left[\frac{(\varepsilon_i^+ - \varepsilon_i^-)^2(1+i)^{2\gamma}}{4c^2}\boldsymbol{\Delta}_i^{-T}\boldsymbol{\Delta}_i^{-1}\right].$$

Then, it is reasonable to pick $c$ at the level approximately equal to the standard derivation of the measurement noise, because it can control the effect of noise.

For good asymptotic performance of the algorithm, we want the sequence to bounce as little as possible around the optimal solution, which is achieved by small values of the gain sequence in the very late iterations. It is equivalent to say for asymptotical performance, we prefer big $\alpha$ ($\alpha = 1$). When $\alpha = 1$, we must have $\gamma = 1/6$ and it is the feasible extreme point of the triangular domain. Thus, for SPSA, we have the similar result: the theoretical optimal selection of coefficients is for better asymptotical performance; while the practical selection of coefficients focuses on better finite sample performance.

# A.3 Numerical Experiments

Obviously, in practice, only a finite number of iterations can be performed for any problem. In this section, we present two simple numerical experiments to show that the theoretical optimal gain sequence may perform worse than the practical selection for finite sample performance. We also show that when the initial condition is close to the optimum, the asymptotically optimal gains may be preferred.

We do the numerical tests on both SG and SPSA. The first example is a test with SG applied in online training to a special case of state-space model considered in Spall (2012). This special case is shown in Spall (2003, Example 13.7), and the data are distributed according to $z_i \sim N(\mu, \sigma^2 + q_i)$, where the $q_i$ are known. Here we suppose $q_i = 1$, $\mu = -1$, and $\sigma^2 = 1$ (more general forms in the above references have $q_i$ being non-identical across *i*). We use the idea of online training to estimate the parameters $\mu$ and $\sigma$. In each step we want to maximize the log likelihood function $\log\left(e^{-(z_i-\mu)^2/2(\sigma^2+q_i)}\Big/\sqrt{2\pi(\sigma^2+q_i)}\right)$ based on the observation at that time. It is equivalent to minimize the negative log likelihood function. Here $\hat{\boldsymbol{\theta}}_k=\left(\hat{\mu}_k,\hat{\sigma}_k\right)^T$, then the estimated gradient is $\hat{\boldsymbol{g}}_k(\hat{\boldsymbol{\theta}}_k)=\left((\hat{\mu}_k-z_k)\big/(\hat{\sigma}_k^2+q_i),\hat{\sigma}_k\big/(\hat{\sigma}_k^2+q_i)-(\hat{\mu}_k-z_k)^2\hat{\sigma}_k\Big/\left(\hat{\sigma}_k^2+q_i\right)^2\right)^T$ in the $k$th iteration. The dimension of the case here is 2, and $a_k=a\big/(k+1+A)^{\alpha}$. The initial guess is set to be $[5,5]^T$ in all runs. The total number of observations in each replicate is 1000 and the number of replicates is 50. By the rule of $A$, we have $A$ = 0.1×1000 = 100. For practical selection, we choose $\alpha=0.501$, and for the asymptotically optimal selection, we choose $\alpha=1$. After computing some values of $\hat{\boldsymbol{g}}_0(\hat{\boldsymbol{\theta}}_0)$, we know that the mean magnitude of the components in $\hat{\boldsymbol{g}}_0(\hat{\boldsymbol{\theta}}_0)$ is approximately 0.16. Suppose we want the elements of $\boldsymbol{\theta}$ move by a magnitude of 0.01 in the early iterations. Then, for the practical selection, $a$ = 0.65 is according to $\left(0.65/101^{0.501}\right)\times 0.16\approx 0.01$; for the asymptotically optimal selection, $a$ = 6.5 is according to $\left(6.5/101\right)\times 0.16\approx 0.01$. Thus, we compare two sets of coefficients: 1) (practical) $a$ = 0.65, $A$ = 100, $\alpha=0.501$; 2) (asymptotically optimal $\alpha$) $a$ = 6.5, $A$ = 100 $\alpha=1$. Figure A.1 shows that the practical choice of coefficients provides better performance than the asymptotically optimal one.

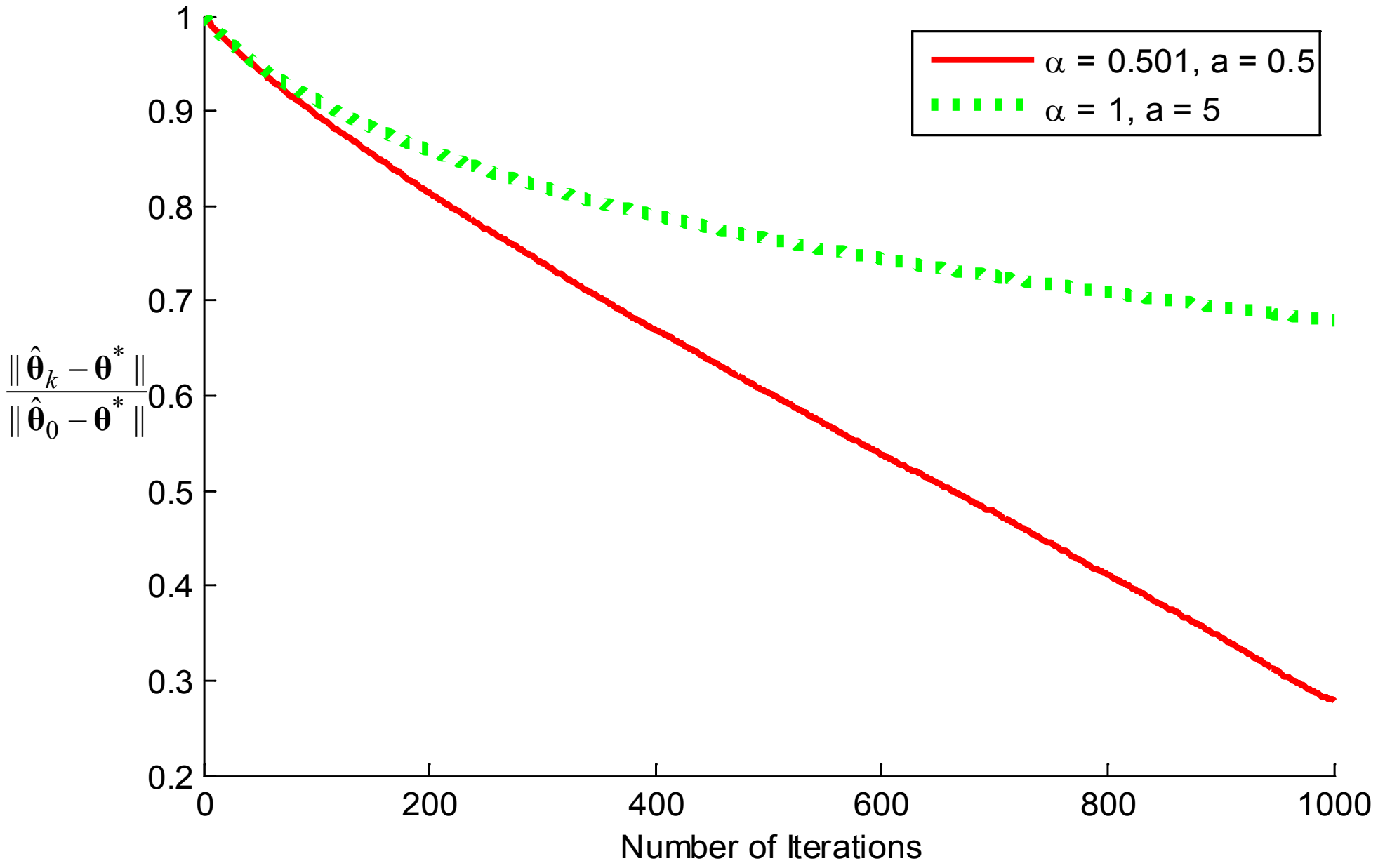

**Figure A.1** Comparison of practical selection and asymptotically optimal selection of coefficients for a simple case of state-space model. Each curve represents the sample mean of 50 independent replicates.

The second example is a test with SPSA applied to the skewed quartic loss function (Spall, 2003, Example 6.6): $L(\boldsymbol{\theta}) = \boldsymbol{\theta}^T \boldsymbol{B}^T \boldsymbol{B}\boldsymbol{\theta} + 0.1\sum_{i=1}^{p}(\boldsymbol{B}\boldsymbol{\theta})_i^3 + 0.01\sum_{i=1}^{p}(\boldsymbol{B}\boldsymbol{\theta})_i^4$, where $(\boldsymbol{B}\boldsymbol{\theta})_i$ represents the $i$th component of vector $\boldsymbol{B}\boldsymbol{\theta}$ and $p\boldsymbol{B}$ is an upper triangular matrix of 1's. We consider the case of dimension $p$ = 200, and the measurement noise ε is i.i.d $N(0,1)$. For the perturbation $\boldsymbol{\Delta}_k = [\Delta_{k1}, \Delta_{k2}, ..., \Delta_{kp}]^T$, each component $\Delta_{ki}$ is independent Bernoulli random variable taking the values $\pm 1$ with probability $1/2$ . Here $a_k = a/(k+1+A)^{\alpha}$ and $c_k = c/(k+1)^{\gamma}$ .

We consider two initial values for $\boldsymbol{\theta}$ in the search process. The first initial value is $10\times\mathbf{1}_{200}$ , where $\mathbf{1}_{200}$ is a 200-dimensional vector with each component being 1. This initial guess is far away from the optimal solution $\mathbf{0}_{200}$, where $\mathbf{0}_{200}$ is a 200-dimensional vector with each component being 0. The total number of iterations in each replicate is 1000 and the number of replicates is 50. By the rule of $A$ and $c$, we have $A = 0.1\times 1000 =$ 100 and $c$ = 1. For practical selection, we choose $\alpha = 0.602$ and $\gamma = 0.101$, and for the asymptotically optimal selection, we choose $\alpha = 1$ and $\gamma = 1/6$ . After computing some values of $\hat{g}_0(\hat{\boldsymbol{\theta}}_0)$, we know that the mean magnitude of the components in $\hat{g}_0(\hat{\boldsymbol{\theta}}_0)$ is approximately 300. Suppose we want the elements of $\boldsymbol{\theta}$ move by a magnitude of 1.5 in the early iterations. Then, for the practical selection, we have that $a$ = 0.08 is according to $\left(0.08/101^{0.602}\right)\times 300 \approx 1.5$; for the asymptotically optimal selection, $a$ = 0.5 is according to $(0.5/101)\times 300 \approx 1.5$. Thus, we compare two sets of coefficients: (1) (practical) $a$ = 0.08, $c$ = 1, $A$ = 100, $\alpha = 0.602$, $\gamma = 0.101$ ; (2) (asymptotically optimal $\alpha$ and $\gamma$) $a$ = 0.5, $c$ = 1, $A$ = 100, $\alpha = 1$, $\gamma = 1/6$ .

We also want to check the asymptotical performance for both sets of coefficients, so we also compare the performance for the initial guess $0.12\times\mathbf{1}_{200}$, which is close to the optimal solution. Therefore, plot (a) in Figure A.2 shows that the practical choice of coefficients provides better performance than the asymptotically optimal one, when the initial guess is far away from the optimal solution. We see that the asymptotically optimal selection of coefficients also provides reasonable finite sample performance for skewed quartic loss function. Plot (b) in Figure A.2 shows that when the initial guess is close to the optimal solution, the asymptotically optimal selection of coefficients provides better and more stable performance.

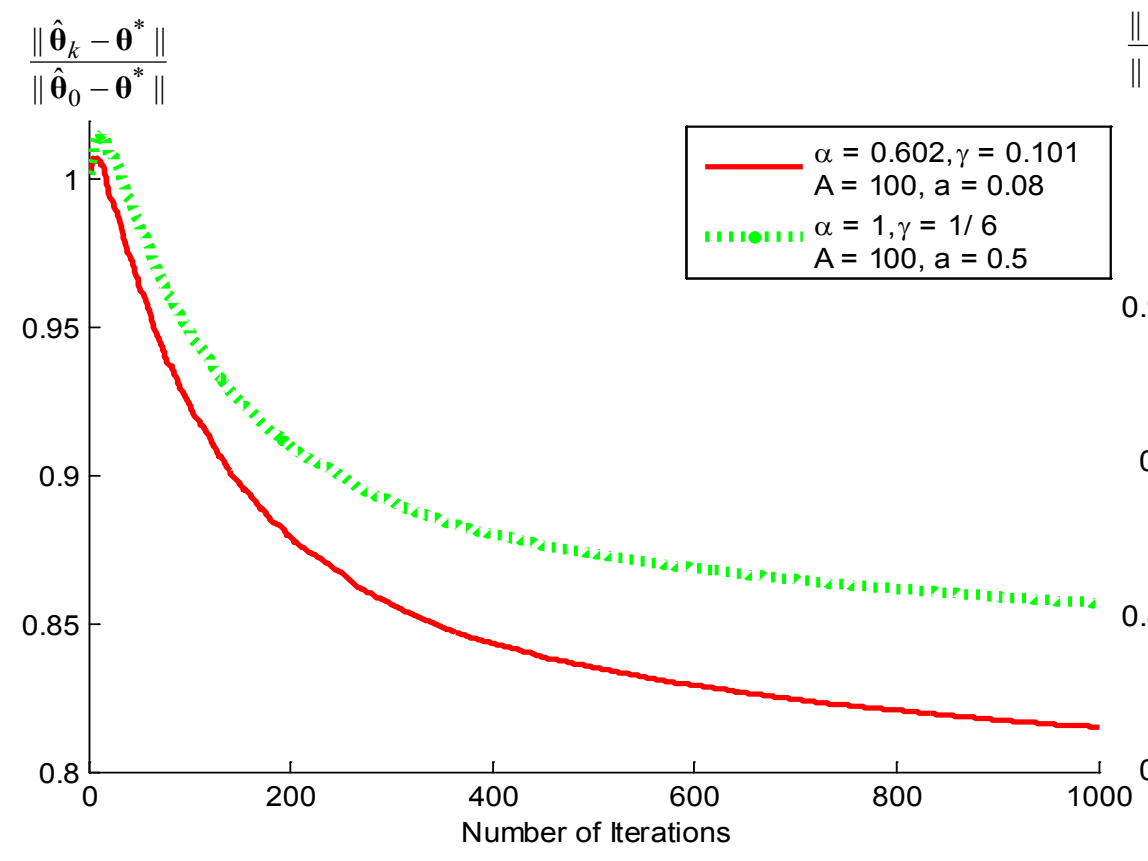

(a)

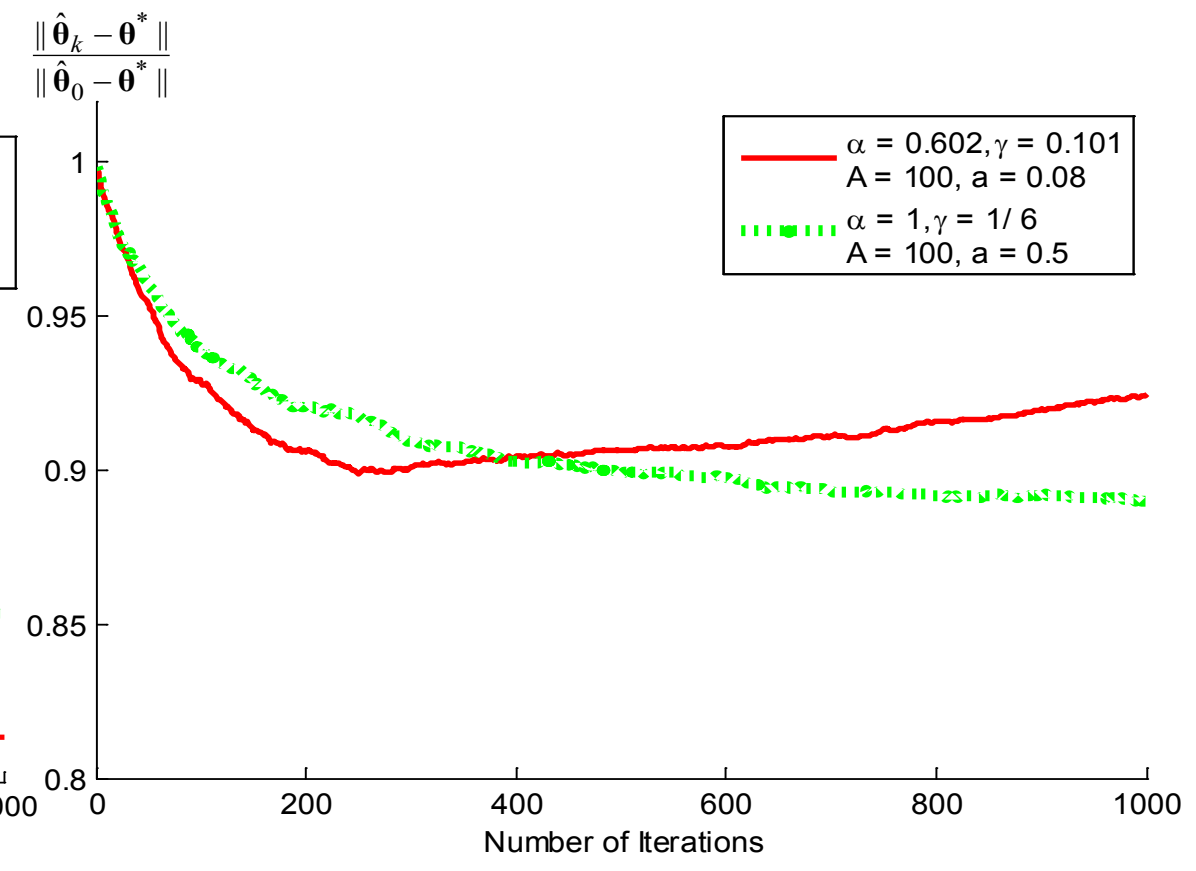

(b)

**Figure A.2** Comparisons of practical selection and asymptotically optimal selection of coefficients for the skewed quartic loss function. Plot (a) shows the comparison when the initial guess is far away from the optimal solution. Plot (b) shows the comparison when the initial guess is close to the optimal solution. Each curve represents the sample mean of 50 independent replicates.

# A.4 Conclusions

Stochastic approximation algorithms are widely used in many applications related to stochastic optimization. This appendix discusses the practical gain sequence selection for stochastic approximation algorithms, especially for SG and SPSA. For the choice of coefficients of the gain sequence (step size), people often use those popular "rules of thumb" directly, such as Spall (1998b, 2003). But to our knowledge, nobody has provided formal justifications for those practical gain sequence selection. In our appendix, we show formally that in order to have better finite sample performance the selection of gain sequence is different from the theoretical optimal gain sequence derived from the asymptotic distribution. This result provides reasons for those popular used practical rules of coefficient choices. From the proof we can see the rules for the practical coefficient selection are useful for those problems with bad initial guess and smooth Hessian matrix.

# Appendix B

# Rate of Convergence Analysis of Simultaneous Perturbation Stochastic Approximation Algorithm for Time-Varying Loss Function

A popular method for continuous stochastic optimization problem is simultaneous perturbation stochastic approximation (SPSA). Spall (1992) introduces SPSA and discusses the rate of convergence of SPSA for fixed loss functions. In this appendix, we use different criteria to discuss the rate of convergence of SPSA for time-varying loss functions. The rate of convergence result shows that SPSA is an effective algorithm for time-varying problems, such as the model-free adaptive control of nonlinear stochastic systems with unknown dynamics. The analysis here is for continuous stochastic optimization, which is different from the rate of convergence analysis in Chapter 3 on the discrete case.

# B.1 Introduction

The continuous stochastic optimization problem with time-varying loss function is widely used in the real world. Many people use SPSA directly for time-varying problems without formal justification on the efficiency of SPSA. We want to provide a formal result on the rate of convergence of SPSA on time-varying loss function.

Let us first see some examples of the applications related to time-varying problems. Neely (2006) considers the energy optimal control problem for time-varying wireless networks. The transmission rates are determined by link channel conditions from time to time, and the link conditions may change dramatically because of the environmental effects. Moore and Schneider (1996) mention the problem of simulated power-plant process on chemical process engineering, where the performance of the system depends on variations of the unsensed chemicals over time. Schwartz et al. (2006) consider simulation based optimization problem for the finical decision making and forecasting in inventory management of supply chain under the variations of supply and demand.

Now let us discuss one important application of time-varying loss function, which is in the model-free adaptive control with nonlinear stochastic systems. Model-free adaptive control is useful in many practical areas: wastewater treatment system, bio-inspired terrestrial micro-robot design, nuclear steam generator water level control, airplane design, and etc. Spall and Cristion (1998) discuss the problem of developing a model-free controller for general nonlinear stochastic systems by using SPSA. The reasons to use SPSA are that 1) the functions governing the system are not known; 2) the dimension of the problem may be high; 3) noise is involved in the system; and 4) the gradient is not generally computable. Spall and Chin (1994) consider the problem of optimal traffic light timing by using the model-free approach with SPSA. Spall and Cristion (1996) and Spall (1997) are two licensed U.S. patents related to the model-free control by using SPSA. Furthermore, many other papers have considered the model-free control problem by using SPSA or similar algorithms (Ji and Familoni, 1999, Song et al., 2008, Zhou et al., 2008, Hahn and Oldham, 2010, Dong and Chen, 2012) in many practical areas. But none of these papers offer a formal rate of convergence analysis. The formal result on the rate of convergence of SPSA algorithm can provide theoretical support on the efficiency of the algorithm, which is quite similar to the purpose of the work of Stark and Spall (2003) where the authors discuss the formal rate of convergence analysis for evolutionary computation algorithms (such as genetic algorithms).

Spall (1992) introduces the algorithm of SPSA and discusses the rate of convergence in terms of asymptotical normality for the case $\gamma \geq \alpha/6$, where $\gamma$ is the decaying rate of the finite difference approximation $c_k = c/(1+k)^{\gamma}$ ($c > 0$) and $\alpha$ is the decaying rate of the gain sequence $a_k = a/(1+k+A)^{\alpha}$ ($a > 0$, $A \geq 0$). But the asymptotical normality result requires fixed loss function and does not work for time-varying problems. We use the mean square error as the measure of rate of convergence in the analysis, and discuss the rate of convergence for time-varying loss function by using SPSA. We show that SPSA is an efficient algorithm in solving time-varying problems in terms of big-$O$ sense.

Suppose $L_k(\boldsymbol{\theta})$ is the loss function at time $k$, and $y_k(\boldsymbol{\theta})$ is the corresponding noisy measurement of the loss function with mean 0 noise. The objective function for the time-varying case here is

$$\min_{\boldsymbol{\theta}\in\Theta} L_k(\boldsymbol{\theta}) = \min_{\boldsymbol{\theta}\in\Theta} E(y_k(\boldsymbol{\theta})),$$

Suppose the optimal solution at time $k$ is $\boldsymbol{\theta}_k^*$, and there exists a $\boldsymbol{\theta}^*$ such that $\boldsymbol{\theta}_k^* \to \boldsymbol{\theta}^*$ as $k \to \infty$. Let $\{\hat{\boldsymbol{\theta}}_k\}$ be the sequence generated by SPSA.

This appendix is organized as: In Section B.2, we provide the description of SPSA. In Section B.3, we discuss the rate of convergence of SPSA under the criterion of $E\left\|\hat{\boldsymbol{\theta}}_k - \boldsymbol{\theta}^*\right\|^2$. In Section B.4, numerical results are presented. This appendix is concluded by Section B.5.

## B.2 Description of SPSA for Time-Varying Loss Function

The basic steps of SPSA for time-varying function are listed below. As we have mentioned above, $L_k(\cdot)$ is the time-varying loss function, $y_k(\cdot)$ is the noisy measurement of the loss function, and $a_k = a/(1+A+k)^{\alpha}$, $c_k = c/(1+A+k)^{\gamma}$, where $a$, $c$, $\alpha$, $\gamma$ are positive and $A$ is nonnegative.

Step 1 Pick an initial guess $\hat{\boldsymbol{\theta}}_0$.

Step 2 Generate $\boldsymbol{\Delta}_k = [\Delta_{k1}, \Delta_{k2}, ..., \Delta_{kp}]^T$, where $\boldsymbol{\Delta}_k$ has a user-specified distribution satisfying conditions discussed in Theorem B.1 below. A special case is when the $\Delta_{ki}$ are independent Bernoulli random variables taking the values $\pm 1$ with probability $1/2$.

Step 3 Evaluate noisy measurements of loss function at $\hat{\boldsymbol{\theta}}_k + c_k \boldsymbol{\Delta}_k$ and $\hat{\boldsymbol{\theta}}_k - c_k \boldsymbol{\Delta}_k$: $y_k^{(+)} = L_k(\hat{\boldsymbol{\theta}}_k + c_k \boldsymbol{\Delta}_k) + \varepsilon_k^{(+)}$ and $y_k^{(-)} = L_k(\hat{\boldsymbol{\theta}}_k - c_k \boldsymbol{\Delta}_k) + \varepsilon_k^{(-)}$. Form the estimate of $\hat{\boldsymbol{g}}_k(\hat{\boldsymbol{\theta}}_k)$

$$\hat{\boldsymbol{g}}_k(\hat{\boldsymbol{\theta}}_k) = \left( \frac{y_k^{(+)} - y_k^{(-)}}{2c_k} \right) \boldsymbol{\Delta}_k^{-1},$$

where $\varepsilon_k^{(+)}$, $\varepsilon_k^{(-)}$ are noise and $\boldsymbol{\Delta}_k^{-1} = \left[ \Delta_{k1}^{-1}, ..., \Delta_{kp}^{-1} \right]^T$.

Step 4 Update the estimator according to the recursion

$$\hat{\boldsymbol{\theta}}_{k+1} = \hat{\boldsymbol{\theta}}_k - a_k \hat{\boldsymbol{g}}_k(\hat{\boldsymbol{\theta}}_k).$$

# B.3 Rate of Convergence Analysis of SPSA in Time-Varying Loss Function

Spall and Cristion (1998) show that under some conditions, SPSA converges for time-varying loss function: $\hat{\boldsymbol{\theta}}_k \to \boldsymbol{\theta}^*$ a.s., where $\boldsymbol{\theta}^*$ is the asymptotical optimal solution. Since the asymptotical normality result may not be derived for time-varying loss function, then we consider the $E\left\|\hat{\boldsymbol{\theta}}_k - \boldsymbol{\theta}^*\right\|^2$ as the criteria of convergence rate for the algorithm. The conditions we consider here for the mean square error convergence are not consistent with the conditions in Spall and Cristion (1998), because Spall and Cristion (1998) consider the almost sure convergence.

**Theorem B.1** Assume $L_k(\boldsymbol{\theta})$ is a strictly convex function on $\mathbb{R}^p$ for all *k*, and $\boldsymbol{\theta}^*$ is the asymptotical optimal solution. Assume also (i) $a_k = a/(1+A+k)^{\alpha}$, $c_k = c/(1+A+k)^{\gamma}$, $a_k > 0$, $c_k > 0$, $a_k \to 0$, $c_k \to 0$ as $k \to \infty$, and $\sum_{k=1}^{\infty} a_k = \infty$, $\sum_{k=1}^{\infty} (a_k/c_k)^2 < \infty$; (ii) the components of $\boldsymbol{\Delta}_k$ are independently and symmetrically distributed about 0 with $|\Delta_{ki}|$ and $|\Delta_{ki}^{-1}|$ being uniformly bounded; (iii) for all *k*, $E\left[(\varepsilon_k^+ - \varepsilon_k^-) \mid \Im_k, \boldsymbol{\Delta}_k\right] = 0$ and the variance of $\varepsilon_k^{\pm}$ is uniformly bounded over *k*; (iv) $E\left(L_k(\hat{\boldsymbol{\theta}}_k + c_k\boldsymbol{\Delta}_k) - L_k(\hat{\boldsymbol{\theta}}_k - c_k\boldsymbol{\Delta}_k)\right)^2$ is uniformly bounded over *k*; (v) there exists $\mu > 0$ such that $E\left[(\hat{\boldsymbol{\theta}}_k - \boldsymbol{\theta}^*)^T \boldsymbol{g}_k(\hat{\boldsymbol{\theta}}_k)\right] \geq \mu E\left\|\hat{\boldsymbol{\theta}}_k - \boldsymbol{\theta}^*\right\|^2$ for all *k*; (vi) for all $\boldsymbol{\theta}$, $L_k^{(3)}(\boldsymbol{\theta}) \equiv$

$\partial^3 L_k \big/ \partial\boldsymbol{\theta}^T \partial\boldsymbol{\theta}^T \partial\boldsymbol{\theta}^T$ exists continuously with individual elements smaller than or equal to $\delta$, where $\delta$ is a positive constant; (vii) $2a\mu > 1 - 2\gamma$; Then

$$E\left\|\hat{\boldsymbol{\theta}}_k - \boldsymbol{\theta}^*\right\|^2 = \begin{cases} O\left(\dfrac{1}{k^{\alpha-2\gamma}}\right) & \gamma \geq \dfrac{\alpha}{6} \\ O\left(\dfrac{1}{k^{4\gamma}}\right) & \gamma < \dfrac{\alpha}{6}. \end{cases}$$

*Remark:* Since $L_k(\boldsymbol{\theta})$ is a strictly convex function for all *k*, the Hessian matrix $\Gamma_k(\boldsymbol{\theta})$ is positive definite. Thus, there exists $\tilde{\boldsymbol{\theta}}_k$ on the line segment of $\hat{\boldsymbol{\theta}}_k$ and $\boldsymbol{\theta}^*$ such that $\boldsymbol{g}_k(\hat{\boldsymbol{\theta}}_k) = \Gamma_k(\tilde{\boldsymbol{\theta}}_k)(\hat{\boldsymbol{\theta}}_k - \boldsymbol{\theta}^*)$. Then $(\hat{\boldsymbol{\theta}}_k - \boldsymbol{\theta}^*)^T \boldsymbol{g}_k(\hat{\boldsymbol{\theta}}_k) = (\hat{\boldsymbol{\theta}}_k - \boldsymbol{\theta}^*)^T \Gamma_k(\tilde{\boldsymbol{\theta}}_k)(\hat{\boldsymbol{\theta}}_k - \boldsymbol{\theta}^*)$. By Corollary 5.4.5 of Horn and Johnson (1985), we have that there exists $\psi_k > 0$ and $\tilde{\psi}_k > 0$ such that $\psi_k \left\|\hat{\boldsymbol{\theta}}_k - \boldsymbol{\theta}^*\right\|^2 \leq (\hat{\boldsymbol{\theta}}_k - \boldsymbol{\theta}^*)^T \Gamma_k(\tilde{\boldsymbol{\theta}}_k)(\hat{\boldsymbol{\theta}}_k - \boldsymbol{\theta}^*) \leq \tilde{\psi}_k \left\|\hat{\boldsymbol{\theta}}_k - \boldsymbol{\theta}^*\right\|^2$, which implies that $\psi_k \left\|\hat{\boldsymbol{\theta}}_k - \boldsymbol{\theta}^*\right\|^2 \leq (\hat{\boldsymbol{\theta}}_k - \boldsymbol{\theta}^*)^T \boldsymbol{g}_k(\hat{\boldsymbol{\theta}}_k) \leq \tilde{\psi}_k \left\|\hat{\boldsymbol{\theta}}_k - \boldsymbol{\theta}^*\right\|^2$. Therefore, condition (v) is related to the curvature of the loss function.

*Proof.* The recursive formula of SPSA is

$$\hat{\boldsymbol{\theta}}_{k+1} = \hat{\boldsymbol{\theta}}_k - a_k \hat{\boldsymbol{g}}_k(\hat{\boldsymbol{\theta}}_k). \tag{B.1}$$

Subtracting both sides of eqn. (B.1) by $\boldsymbol{\theta}^*$ and calculating the norm squared for both sides, we have

$$\left\|\hat{\boldsymbol{\theta}}_{k+1} - \boldsymbol{\theta}^*\right\|^2 = \left\|\hat{\boldsymbol{\theta}}_k - \boldsymbol{\theta}^*\right\|^2 - 2a_k(\hat{\boldsymbol{\theta}}_k - \boldsymbol{\theta}^*)^T \hat{\boldsymbol{g}}_k(\hat{\boldsymbol{\theta}}_k) + a_k^2 \left\|\hat{\boldsymbol{g}}_k(\hat{\boldsymbol{\theta}}_k)\right\|^2. \tag{B.2}$$

Adding and subtracting $2a_k(\hat{\boldsymbol{\theta}}_k - \boldsymbol{\theta}^*)^T \boldsymbol{g}_k(\hat{\boldsymbol{\theta}}_k)$ to the right-hand side of eqn. (B.2), where $\boldsymbol{g}_k(\hat{\boldsymbol{\theta}}_k)$ is the real gradient of function $L_k(\boldsymbol{\theta})$ at point $\hat{\boldsymbol{\theta}}_k$. Then, we have

$$\begin{aligned}
&\left\|\hat{\boldsymbol{\theta}}_{k+1} - \boldsymbol{\theta}^*\right\|^2 \\
&= \left\|\hat{\boldsymbol{\theta}}_k - \boldsymbol{\theta}^*\right\|^2 - 2a_k\left(\hat{\boldsymbol{\theta}}_k - \boldsymbol{\theta}^*\right)^T\left(\hat{\boldsymbol{g}}_k(\hat{\boldsymbol{\theta}}_k) + \boldsymbol{g}_k(\hat{\boldsymbol{\theta}}_k) - \boldsymbol{g}_k(\hat{\boldsymbol{\theta}}_k)\right) + a_k^2\left\|\hat{\boldsymbol{g}}_k(\hat{\boldsymbol{\theta}}_k)\right\|^2 \\
&= \left\|\hat{\boldsymbol{\theta}}_k - \boldsymbol{\theta}^*\right\|^2 - 2a_k\left(\hat{\boldsymbol{\theta}}_k - \boldsymbol{\theta}^*\right)^T \boldsymbol{g}_k(\hat{\boldsymbol{\theta}}_k) + 2a_k\left(\hat{\boldsymbol{\theta}}_k - \boldsymbol{\theta}^*\right)^T\left(\boldsymbol{g}_k(\hat{\boldsymbol{\theta}}_k) - \hat{\boldsymbol{g}}_k(\hat{\boldsymbol{\theta}}_k)\right) \\
&\quad + a_k^2\left\|\hat{\boldsymbol{g}}_k(\hat{\boldsymbol{\theta}}_k)\right\|^2. \qquad \text{(B.3)}
\end{aligned}$$

Taking expectation on both sides of eqn. (B.3), we get

$$\begin{aligned}
&E\left\|\hat{\boldsymbol{\theta}}_{k+1} - \boldsymbol{\theta}^*\right\|^2 \\
&= E\left\|\hat{\boldsymbol{\theta}}_k - \boldsymbol{\theta}^*\right\|^2 - 2a_k E\left[\left(\hat{\boldsymbol{\theta}}_k - \boldsymbol{\theta}^*\right)^T \boldsymbol{g}_k(\hat{\boldsymbol{\theta}}_k)\right] \\
&\quad + 2a_k E\left[\left(\hat{\boldsymbol{\theta}}_k - \boldsymbol{\theta}^*\right)^T\left(\boldsymbol{g}_k(\hat{\boldsymbol{\theta}}_k) - \hat{\boldsymbol{g}}_k(\hat{\boldsymbol{\theta}}_k)\right)\right] + a_k^2 E\left\|\hat{\boldsymbol{g}}_k(\hat{\boldsymbol{\theta}}_k)\right\|^2. \qquad \text{(B.4)}
\end{aligned}$$

By condition (v), we get that the recursive relationship (B.4) can be written as

$$\begin{aligned}
E\left\|\hat{\boldsymbol{\theta}}_{k+1} - \boldsymbol{\theta}^*\right\|^2 &\le (1 - 2\mu a_k)E\left\|\hat{\boldsymbol{\theta}}_k - \boldsymbol{\theta}^*\right\|^2 + 2a_k E\left[\left(\hat{\boldsymbol{\theta}}_k - \boldsymbol{\theta}^*\right)^T\left(\boldsymbol{g}_k(\hat{\boldsymbol{\theta}}_k) - \hat{\boldsymbol{g}}_k(\hat{\boldsymbol{\theta}}_k)\right)\right] \\
&\quad + a_k^2 E\left\|\hat{\boldsymbol{g}}_k(\hat{\boldsymbol{\theta}}_k)\right\|^2. \qquad \text{(B.5)}
\end{aligned}$$

Since $a_k \to 0$, then there exists large integer value $N$, such that when $k \ge N$, $1 - 2\mu a_k >$ 0. Then by the recursive formula (B.5), we have that for $k \ge N$

$$E\left\|\hat{\boldsymbol{\theta}}_{k+1}-\boldsymbol{\theta}^*\right\|^2$$

$$\leq \prod_{i=N}^{k}(1-2\mu a_i)E\left\|\hat{\boldsymbol{\theta}}_N-\boldsymbol{\theta}^*\right\|^2+\prod_{i=N+1}^{k}(1-2\mu a_i)2a_N E\left[\left(\hat{\boldsymbol{\theta}}_N-\boldsymbol{\theta}^*\right)^T\left(\boldsymbol{g}_N(\hat{\boldsymbol{\theta}}_N)-\hat{\boldsymbol{g}}_N(\hat{\boldsymbol{\theta}}_N)\right)\right]$$

$$+\ldots+(1-2\mu a_k)2a_{k-1}E\left[\left(\hat{\boldsymbol{\theta}}_{k-1}-\boldsymbol{\theta}^*\right)^T\left(\boldsymbol{g}_{k-1}(\hat{\boldsymbol{\theta}}_{k-1})-\hat{\boldsymbol{g}}_{k-1}(\hat{\boldsymbol{\theta}}_{k-1})\right)\right]$$

$$+2a_k E\left[\left(\hat{\boldsymbol{\theta}}_k-\boldsymbol{\theta}^*\right)^T\left(\boldsymbol{g}_k(\hat{\boldsymbol{\theta}}_k)-\hat{\boldsymbol{g}}_k(\hat{\boldsymbol{\theta}}_k)\right)\right]+\prod_{i=N+1}^{k}(1-2\mu a_i)a_N^2 E\left\|\hat{\boldsymbol{g}}_N(\hat{\boldsymbol{\theta}}_N)\right\|^2$$

$$+\ldots+(1-2\mu a_k)a_{k-1}^2 E\left\|\hat{\boldsymbol{g}}_{k-1}(\hat{\boldsymbol{\theta}}_{k-1})\right\|^2+a_k^2 E\left\|\hat{\boldsymbol{g}}_k(\hat{\boldsymbol{\theta}}_k)\right\|^2. \qquad \text{(B.6)}$$

Let $T_k=\prod_{i=N}^{k}(1-2\mu a_i)$, then inequality (B.6) can be written as

$$E\left\|\hat{\boldsymbol{\theta}}_{k+1}-\boldsymbol{\theta}^*\right\|^2$$

$$\leq T_k E\left\|\hat{\boldsymbol{\theta}}_N-\boldsymbol{\theta}^*\right\|^2+T_k\sum_{i=N}^{k}\frac{a_i^2 E\left\|\hat{\boldsymbol{g}}_i(\hat{\boldsymbol{\theta}}_i)\right\|^2}{T_i}+T_k\sum_{i=N}^{k}\frac{2a_i E\left[\left(\hat{\boldsymbol{\theta}}_i-\boldsymbol{\theta}^*\right)^T\left(\boldsymbol{g}_i(\hat{\boldsymbol{\theta}}_i)-\hat{\boldsymbol{g}}_i(\hat{\boldsymbol{\theta}}_i)\right)\right]}{T_i}. \qquad \text{(B.7)}$$

In the following, we consider the rate of convergence of the three terms on the right-hand side of inequality (B.7) separately. First we discuss the upper and lower bounds of $T_k$. By the definition of $T_k$, we have $\log T_k=\sum_{i=N}^{k}\log(1-2\mu a_i)$. By Taylor expansion, we have

$$\log(1-2\mu a_i)=-2\mu a_i-\frac{1}{2}(2\mu a_i)^2-\frac{1}{3}(2\mu a_i)^3-\ldots$$
$$=\sum_{j=1}^{\infty}\frac{-(2\mu a_i)^j}{j},$$

so

$$\log T_k = \sum_{i=N}^{k}\sum_{j=1}^{\infty}\frac{-(2\mu a_i)^j}{j} = \sum_{j=1}^{\infty}\sum_{i=N}^{k}\frac{-(2\mu a_i)^j}{j}. \tag{B.8}$$

For each fixed $j$, we have

$$\sum_{i=N}^{k}\frac{(2\mu a_i)^j}{j} = \sum_{i=N}^{k}\frac{(2\mu a)^j}{j(1+A+i)^{j\alpha}},$$

and it is a decreasing function on $i$. Therefore,

$$\int_{N}^{k+1}\frac{(2\mu a)^j}{j(1+A+x)^{j\alpha}}dx \le \sum_{i=N}^{k}\frac{(2\mu a)^j}{j(1+A+i)^{j\alpha}} \le \int_{N-1}^{k}\frac{(2\mu a)^j}{j(1+A+x)^{j\alpha}}dx,$$

then let us solve the lower and the upper bounds of $\sum_{i=N}^{k}(2\mu a)^j \big/ j(1+A+i)^{j\alpha}$ in the inequality above for different cases: 1) $j \ge 2$ ($j\alpha > 1$); 2) $j = 1$, $\alpha = 1$ ($j\alpha = 1$); 3) $j = 1$ $0.5 < \alpha < 1$ ($j\alpha < 1$).

Case 1: $j \ge 2$ ($j\alpha > 1$)

$$\begin{aligned}&\frac{(2\mu a)^j}{j(j\alpha-1)}\left(\frac{1}{(1+A+N)^{j\alpha-1}} - \frac{1}{(2+A+k)^{j\alpha-1}}\right)\\ &\le \sum_{i=N}^{k}\frac{(2\mu a)^j}{j(1+A+i)^{j\alpha}}\\ &\le \frac{(2\mu a)^j}{j(j\alpha-1)}\left(\frac{1}{(A+N)^{j\alpha-1}} - \frac{1}{(1+A+k)^{j\alpha-1}}\right).\end{aligned} \tag{B.9}$$

Case 2: $j = 1,\ \alpha = 1$

$$2\mu a\left(\log(2+A+k)-\log(1+A+N)\right) \le \sum_{i=N}^{k} \frac{(2\mu a)^j}{j(1+A+i)^{j\alpha}}$$
$$\le 2\mu a\left(\log(1+A+k)-\log(A+N)\right). \qquad \text{(B.10)}$$

Case 3: $j = 1,\ 0.5 < \alpha < 1$

$$\frac{2\mu a}{1-\alpha}\left((2+A+k)^{1-\alpha}-(1+A+N)^{1-\alpha}\right) \le \sum_{i=N}^{k} \frac{(2\mu a)^j}{j(1+A+i)^{j\alpha}}$$
$$\le \frac{2\mu a}{1-\alpha}\left((1+A+k)^{1-\alpha}-(A+N)^{1-\alpha}\right). \qquad \text{(B.11)}$$

Combining the results of inequalities (B.9) and (B.11), we have that for the case of $0.5 < \alpha < 1$

$$\frac{2\mu a}{1-\alpha}\left((2+A+k)^{1-\alpha}-(1+A+N)^{1-\alpha}\right)+\sum_{j=2}^{\infty}\frac{(2\mu a)^j}{j(j\alpha-1)}\left(\frac{1}{(1+A+N)^{j\alpha-1}}-\frac{1}{(2+A+k)^{j\alpha-1}}\right)$$
$$\le \sum_{j=1}^{\infty}\sum_{i=N}^{k}\frac{(2\mu a_i)^j}{j}$$
$$\le \frac{2\mu a}{1-\alpha}\left((1+A+k)^{1-\alpha}-(A+N)^{1-\alpha}\right)+\sum_{j=2}^{\infty}\frac{(2\mu a)^j}{j(j\alpha-1)}\left(\frac{1}{(A+N)^{j\alpha-1}}-\frac{1}{(1+A+k)^{j\alpha-1}}\right),$$

so by eqn. (B.8) we get the lower and upper bounds for $T_k$ as

$$\exp\left(\frac{-2\mu a}{1-\alpha}\left((2+A+k)^{1-\alpha}-(1+A+N)^{1-\alpha}\right)\right)$$

$$\times\exp\left(\sum_{j=2}^{\infty}\frac{-(2\mu a)^j}{j(j\alpha-1)}\left(\frac{1}{(1+A+N)^{j\alpha-1}}-\frac{1}{(2+A+k)^{j\alpha-1}}\right)\right)$$

$$\geq T_k$$

$$\geq\exp\left(\frac{-2\mu a}{1-\alpha}\left((1+A+k)^{1-\alpha}-(A+N)^{1-\alpha}\right)\right)$$

$$\times\exp\left(\sum_{j=2}^{\infty}\frac{-(2\mu a)^j}{j(j\alpha-1)}\left(\frac{1}{(A+N)^{j\alpha-1}}-\frac{1}{(1+A+k)^{j\alpha-1}}\right)\right). \tag{B.12}$$

Similarly combining the results of inequalities (B.9) and (B.10), we have that for the case of $\alpha = 1$

$$2\mu a\left(\log(2+A+k)-\log(1+A+N)\right)+\sum_{j=2}^{\infty}\frac{(2\mu a)^j}{j(j\alpha-1)}\left(\frac{1}{(1+A+N)^{j\alpha-1}}-\frac{1}{(2+A+k)^{j\alpha-1}}\right)$$

$$\leq\sum_{j=1}^{\infty}\sum_{i=N}^{k}\frac{(2\mu a_i)^j}{j}$$

$$\leq 2\mu a\left(\log(1+A+k)-\log(A+N)\right)+\sum_{j=2}^{\infty}\frac{(2\mu a)^j}{j(j\alpha-1)}\left(\frac{1}{(A+N)^{j\alpha-1}}-\frac{1}{(1+A+k)^{j\alpha-1}}\right),$$

so by eqn. (B.8) we get the lower and upper bounds for $T_k$ as

$$\frac{1/(2+A+k)^{2\mu a}}{1/(1+A+N)^{2\mu a}}\exp\left(\sum_{j=2}^{\infty}\frac{-(2\mu a)^j}{j(j\alpha-1)}\left(\frac{1}{(1+A+N)^{j\alpha-1}}-\frac{1}{(2+A+k)^{j\alpha-1}}\right)\right)$$

$$\geq T_k$$

$$\geq\frac{1/(1+A+k)^{2\mu a}}{1/(A+N)^{2\mu a}}\exp\left(\sum_{j=2}^{\infty}\frac{-(2\mu a)^j}{j(j\alpha-1)}\left(\frac{1}{(A+N)^{j\alpha-1}}-\frac{1}{(1+A+k)^{j\alpha-1}}\right)\right). \tag{B.13}$$

Moreover, let us further discuss the lower and upper bounds in the inequalities (B.12) and (B.13) to make them simpler. In the following, we will show that some terms in the lower and upper bounds in the inequalities (B.12) and (B.13) are uniformly bounded. First we show that

$$\sum_{j=2}^{\infty}\frac{(2\mu a)^j}{j(j\alpha-1)}\left(\frac{1}{(1+A+k)^{j\alpha-1}}\right)$$

is uniformly bounded for all $k$.

$$\begin{aligned}
&\sum_{j=2}^{\infty}\frac{(2\mu a)^j}{j(j\alpha-1)}\left(\frac{1}{(1+A+k)^{j\alpha-1}}\right)\\
&=\sum_{j=2}^{\infty}\left(\frac{\alpha}{j\alpha-1}-\frac{1}{j}\right)\left(\frac{(2\mu a)^j}{(1+A+k)^{j\alpha-1}}\right)\\
&=\sum_{j=2}^{\infty}\left(\frac{\alpha}{j\alpha-1}\right)\left(\frac{(2\mu a)^j}{(1+A+k)^{j\alpha-1}}\right)-\sum_{j=2}^{\infty}\frac{(2\mu a)^j}{j(1+A+k)^{j\alpha-1}}\\
&\le\sum_{j=2}^{\infty}\left(\frac{1}{j-1/\alpha}\right)\left(\frac{(2\mu a)^j}{(1+A+k)^{j\alpha-1}}\right)\\
&=\left(\frac{\alpha}{2\alpha-1}\right)\left(\frac{(2\mu a)^2}{(1+A+k)^{2\alpha-1}}\right)+\sum_{j=3}^{\infty}\left(\frac{1}{j-1/\alpha}\right)\left(\frac{(2\mu a)^j}{(1+A+k)^{j\alpha-1}}\right)\\
&=\left(\frac{\alpha}{2\alpha-1}\right)\left(\frac{(2\mu a)^2}{(1+A+k)^{2\alpha-1}}\right)+\sum_{j=1}^{\infty}\left(\frac{1}{j+2-1/\alpha}\right)\left(\frac{(2\mu a)^{j+2}}{(1+A+k)^{(j+2)\alpha-1}}\right)\\
&\le\left(\frac{\alpha}{2\alpha-1}\right)\left(\frac{(2\mu a)^2}{(1+A+k)^{2\alpha-1}}\right)+\sum_{j=1}^{\infty}\frac{(2\mu a)^{j+2}}{j(1+A+k)^{\alpha j}}\\
&=\frac{(2\mu a)^2\alpha}{(2\alpha-1)(1+A+k)^{2\alpha-1}}+(2\mu a)^2\sum_{j=1}^{\infty}\frac{(2\mu a/(1+A+k)^{\alpha})^j}{j}.
\end{aligned}$$

By Taylor expansion, we have

$$\log\left(1-\frac{2\mu a}{(1+A+k)^{\alpha}}\right)=\sum_{j=1}^{\infty}\frac{-(2\mu a/(1+A+k)^{\alpha})^{j}}{j}.$$

Therefore,

$$\begin{aligned}
0 &\le \sum_{j=2}^{\infty}\frac{(2\mu a)^{j}}{j(j\alpha-1)}\left(\frac{1}{(1+A+k)^{j\alpha-1}}\right) \\
&\le \frac{(2\mu a)^{2}\alpha}{(2\alpha-1)(1+A+k)^{2\alpha-1}}+(2\mu a)^{2}\sum_{j=1}^{\infty}\frac{(2\mu a/(1+A+k)^{\alpha})^{j}}{j} \\
&= \frac{(2\mu a)^{2}\alpha}{(2\alpha-1)(1+A+k)^{2\alpha-1}}-(2\mu a)^{2}\log\left(1-\frac{2\mu a}{(1+A+k)^{\alpha}}\right) \\
&\le \frac{(2\mu a)^{2}\alpha}{(2\alpha-1)(1+A+N)^{2\alpha-1}}-(2\mu a)^{2}\log\left(1-\frac{2\mu a}{(1+A+N)^{\alpha}}\right)<\infty, \qquad \text{(B.14)}
\end{aligned}$$

which indicates that

$$\sum_{j=2}^{\infty}\frac{(2\mu a)^{j}}{j(j\alpha-1)}\left(\frac{1}{(1+A+k)^{j\alpha-1}}\right)$$

is uniformly bounded on *k*. By similar calculation, we also have

$$\begin{aligned}
0 &\le \sum_{j=2}^{\infty}\frac{(2\mu a)^{j}}{j(j\alpha-1)}\left(\frac{1}{(1+A+N)^{j\alpha-1}}\right) \\
&\le \frac{(2\mu a)^{2}\alpha}{(2\alpha-1)(1+A+N)^{2\alpha-1}}-(2\mu a)^{2}\log\left(1-\frac{2\mu a}{(1+A+N)^{\alpha}}\right)<\infty, \qquad \text{(B.15)}
\end{aligned}$$

$$0 \le \sum_{j=2}^{\infty} \frac{(2\mu a)^j}{j(j\alpha-1)}\left(\frac{1}{(2+A+k)^{j\alpha-1}}\right)$$

$$\le \frac{(2\mu a)^2\alpha}{(2\alpha-1)(2+A+N)^{2\alpha-1}} - (2\mu a)^2 \log\left(1-\frac{2\mu a}{(2+A+N)^{\alpha}}\right) < \infty, \tag{B.16}$$

and

$$0 \le \sum_{j=2}^{\infty} \frac{(2\mu a)^j}{j(j\alpha-1)}\left(\frac{1}{(A+N)^{j\alpha-1}}\right)$$

$$\le \frac{(2\mu a)^2\alpha}{(2\alpha-1)(A+N)^{2\alpha-1}} - (2\mu a)^2 \log\left(1-\frac{2\mu a}{(A+N)^{\alpha}}\right) < \infty, \tag{B.17}$$

which indicates that

$$\sum_{j=2}^{\infty} \frac{(2\mu a)^j}{j(j\alpha-1)}\left(\frac{1}{(1+A+N)^{j\alpha-1}}\right) \text{ and } \sum_{j=2}^{\infty} \frac{(2\mu a)^j}{j(j\alpha-1)}\left(\frac{1}{(A+N)^{j\alpha-1}}\right)$$

are bounded, and

$$\sum_{j=2}^{\infty} \frac{(2\mu a)^j}{j(j\alpha-1)}\left(\frac{1}{(2+A+k)^{j\alpha-1}}\right)$$

is uniformly bounded on $k$. In all, for $0.5 < \alpha < 1$, by inequality (B.12), (B.14), (B.15), (B.16) and (B.17), there exist $\tilde{\upsilon}, \tilde{\upsilon}' > 0$ such that

$$\tilde{\upsilon} \exp\left(-\frac{2\mu a}{1-\alpha}(1+A+k)^{1-\alpha}\right) \le T_k \le \tilde{\upsilon}' \exp\left(-\frac{2\mu a}{1-\alpha}(1+A+k)^{1-\alpha}\right). \tag{B.18}$$

Similarly, for $\alpha = 1$, by inequality (B.13), (B.14), (B.15), (B.16) and (B.17), there exist $\upsilon, \upsilon' > 0$ such that

$$\upsilon \frac{1}{(1+A+k)^{2\mu a}} \le T_k \le \upsilon' \frac{1}{(1+A+k)^{2\mu a}}. \tag{B.19}$$

After getting the upper and lower bounds of $T_k$ from inequality (B.18) and (B.19), we know that the first term on the right-hand side of inequality (B.7) $T_k E\left\|\hat{\boldsymbol{\theta}}_N - \boldsymbol{\theta}^*\right\|^2$ goes to 0 at the rate of

$$\begin{cases} O\left(\exp\left(-\dfrac{2\mu a}{1-\alpha}(1+A+k)^{1-\alpha}\right)\right), 0.5<\alpha<1, \\ O\left(\dfrac{1}{k^{2\mu a}}\right), \alpha = 1. \end{cases} \tag{B.20}$$

Now let us start to discuss the second term on the right-hand side of inequality (B.7). By the definition of $\hat{\boldsymbol{g}}_k(\hat{\boldsymbol{\theta}}_k)$ and condition (iii), we have

$$\begin{aligned} E\left\|\hat{\boldsymbol{g}}_k(\hat{\boldsymbol{\theta}}_k)\right\|^2 &= E\left\|\frac{L_k(\hat{\boldsymbol{\theta}}_k + c_k\boldsymbol{\Delta}_k) - L_k(\hat{\boldsymbol{\theta}}_k - c_k\boldsymbol{\Delta}_k) + \varepsilon_k^+ - \varepsilon_k^-}{2c_k}\boldsymbol{\Delta}_k^{-1}\right\|^2 \\ &= E\left[\frac{\left(L_k(\hat{\boldsymbol{\theta}}_k + c_k\boldsymbol{\Delta}_k) - L_k(\hat{\boldsymbol{\theta}}_k - c_k\boldsymbol{\Delta}_k)\right)^2}{4c_k^2}\boldsymbol{\Delta}_k^{-T}\boldsymbol{\Delta}_k^{-1}\right] + E\left[\frac{\left(\varepsilon_k^+ - \varepsilon_k^-\right)^2}{4c_k^2}\boldsymbol{\Delta}_k^{-T}\boldsymbol{\Delta}_k^{-1}\right] \\ &= E\left[\frac{\left(L_k(\hat{\boldsymbol{\theta}}_k + c_k\boldsymbol{\Delta}_k) - L_k(\hat{\boldsymbol{\theta}}_k - c_k\boldsymbol{\Delta}_k)\right)^2}{4c_k^2}\sum_{i=1}^{p}\Delta_{ki}^{-2}\right] + E\left[\frac{\left(\varepsilon_k^+ - \varepsilon_k^-\right)^2}{4c_k^2}\sum_{i=1}^{p}\Delta_{ki}^{-2}\right]. \end{aligned}$$

Due to condition (ii), (iii), and (iv), there exists $\xi > 0$ such that $E\left\|\hat{\boldsymbol{g}}_k(\hat{\boldsymbol{\theta}}_k)\right\|^2 \le \xi / c_k^2$. Thus, the second term on the right-hand side of inequality (B.7) can be bounded as

$$T_k \sum_{i=N}^{k} \frac{a_i^2 E\left\|\hat{\boldsymbol{g}}_i(\hat{\boldsymbol{\theta}}_i)\right\|^2}{T_i} \le \xi T_k \sum_{i=N}^{k} \frac{a_i^2 / c_i^2}{T_i}.$$

Now let us calculate the rate of convergence of

$$T_k \sum_{i=N}^{k} \frac{a_i^2 / c_i^2}{T_i}.$$

Through inequality (B.18) and (B.19) on the lower and upper bounds of $T_k$, we have

$$\begin{cases} \tilde{\upsilon} \exp\left(-\frac{2\mu a}{1-\alpha}(1+A+k)^{1-\alpha}\right) \sum_{i=N}^{k} \frac{a^2 \Big/ \left(c^2 (1+A+i)^{(2\alpha-2\gamma)}\right)}{\tilde{\upsilon}' \exp\left(-\frac{2\mu a}{1-\alpha}(1+A+i)^{1-\alpha}\right)} \\ \frac{\upsilon}{(1+A+k)^{2\mu a}} \sum_{i=N}^{k} \frac{a^2 \Big/ \left(c^2 (1+A+i)^{(2-2\gamma)}\right)}{\upsilon' \big/ (1+A+i)^{2\mu a}} \end{cases}$$

$$\le T_k \sum_{i=N}^{k} \frac{a_i^2 / c_i^2}{T_i}$$

$$\le \begin{cases} \tilde{\upsilon}' \exp\left(-\frac{2\mu a}{1-\alpha}(1+A+k)^{1-\alpha}\right) \sum_{i=N}^{k} \frac{a^2 \Big/ \left(c^2 (1+A+i)^{(2\alpha-2\gamma)}\right)}{\tilde{\upsilon} \exp\left(-\frac{2\mu a}{1-\alpha}(1+A+i)^{1-\alpha}\right)}, \ 0.5 < \alpha < 1 & \text{(B.21a)} \\ \upsilon' \frac{1}{(1+A+k)^{2\mu a}} \sum_{i=N}^{k} \frac{a^2 \Big/ \left(c^2 (1+A+i)^{(2-2\gamma)}\right)}{\upsilon \big/ (1+A+i)^{2\mu a}}, \ \alpha = 1, & \text{(B.21b)} \end{cases}$$

and we see that the difference between the lower and upper bound is a constant scalar. In the following, we calculate the rate of convergence of the lower (upper) bound.

When $0.5 < \alpha < 1$, by ignoring the multiplier of constant scalar, the term in the inequality (B.21a) is

$$\exp\left(-\frac{2\mu a}{1-\alpha}(1+A+k)^{1-\alpha}\right)\sum_{i=N}^{k}\frac{1/(1+A+i)^{(2\alpha-2\gamma)}}{\exp\left(-\frac{2\mu a}{1-\alpha}(1+A+i)^{1-\alpha}\right)}.$$

First we will show that $\exp\left(\frac{2\mu a}{1-\alpha}(1+A+i)^{1-\alpha}\right)\Big/(1+A+i)^{(2\alpha-2\gamma)}$ is an increasing function on $i$. The derivative of $\exp\left(\frac{2\mu a}{1-\alpha}(1+A+i)^{1-\alpha}\right)\Big/(1+A+i)^{(2\alpha-2\gamma)}$ is

$$\frac{\exp\left(\frac{2\mu a}{1-\alpha}(1+A+i)^{1-\alpha}\right)\frac{2\mu a}{\left(1+A+i\right)^{\alpha+2\alpha-2\gamma}}-\exp\left(\frac{2\mu a}{1-\alpha}(1+A+i)^{1-\alpha}\right)\frac{-(2\alpha-2\gamma)}{\left(1+A+i\right)^{2\alpha-2\gamma+1}}}{(1+A+i)^{2(2\alpha-2\gamma)}}$$

$$=\frac{\exp\left(\frac{2\mu a}{1-\alpha}(1+A+i)^{1-\alpha}\right)}{(1+A+i)^{3(2\alpha-2\gamma)}}\left(\frac{2\mu a}{\left(1+A+i\right)^{\alpha}}+\frac{2\alpha-2\gamma}{1+A+i}\right).$$

Due to condition (i), we know $2\alpha-2\gamma > 1$, which implies that

$$\frac{2\mu a}{\left(1+A+i\right)^{\alpha}}+\frac{2\alpha-2\gamma}{1+A+i}>0.$$

Therefore, $\exp\left(\frac{2\mu a}{1-\alpha}(1+A+i)^{1-\alpha}\right)\Big/(1+A+i)^{(2\alpha-2\gamma)}$ is an increasing function on $i$, which follows that

$$\int_{N-1}^{k}\frac{1/(1+A+x)^{(2\alpha-2\gamma)}}{\exp\left(-\frac{2a\mu}{1-\alpha}(1+A+x)^{1-\alpha}\right)}dx\le\sum_{i=N}^{k}\frac{1/(1+A+k)^{(2\alpha-2\gamma)}}{\exp\left(-\frac{2a\mu}{1-\alpha}(1+A+k)^{1-\alpha}\right)}$$

$$\le\int_{N}^{k+1}\frac{1/(1+A+x)^{(2\alpha-2\gamma)}}{\exp\left(-\frac{2a\mu}{1-\alpha}(1+A+x)^{1-\alpha}\right)}dx,$$

then

$$\exp\left(-\frac{2\mu a}{1-\alpha}(1+A+k)^{1-\alpha}\right)\int_{N-1}^{k}\frac{1/(1+A+x)^{(2\alpha-2\gamma)}}{\exp\left(-\dfrac{2a\mu}{1-\alpha}(1+A+x)^{1-\alpha}\right)}dx$$

$$\leq \exp\left(-\frac{2\mu a}{1-\alpha}(1+A+k)^{1-\alpha}\right)\sum_{i=N}^{k}\frac{1/(1+A+k)^{(2\alpha-2\gamma)}}{\exp\left(-\dfrac{2a\mu}{1-\alpha}(1+A+k)^{1-\alpha}\right)}$$

$$\leq \exp\left(-\frac{2\mu a}{1-\alpha}(1+A+k)^{1-\alpha}\right)\int_{N}^{k+1}\frac{1/(1+A+x)^{(2\alpha-2\gamma)}}{\exp\left(-\dfrac{2a\mu}{1-\alpha}(1+A+x)^{1-\alpha}\right)}dx. \tag{B.22}$$

In the following, we will solve out the rate of convergence of the lower and upper bounds in inequality (B.22). For the lower bound, we have

$$\exp\left(-\frac{2\mu a}{1-\alpha}(1+A+k)^{1-\alpha}\right)\int_{N-1}^{k}\frac{1/(1+A+x)^{(2\alpha-2\gamma)}}{\exp\left(-\frac{2\mu a}{1-\alpha}(1+A+x)^{1-\alpha}\right)}dx$$

$$=\exp\left(-\frac{2\mu a}{1-\alpha}(1+A+k)^{1-\alpha}\right)\int_{N-1}^{k}\frac{\exp\left(\frac{2\mu a}{1-\alpha}(1+A+x)^{1-\alpha}\right)}{(1+A+x)^{\alpha+(\alpha-2\gamma)}}dx$$

$$=\exp\left(-\frac{2\mu a}{1-\alpha}(1+A+k)^{1-\alpha}\right)\int_{N-1}^{k}\frac{1}{2\mu a(1+A+x)^{(\alpha-2\gamma)}}d\exp\left(\frac{2\mu a}{1-\alpha}(1+A+x)^{1-\alpha}\right)$$

$$=\frac{\exp\left(-\frac{2\mu a}{1-\alpha}(1+A+k)^{1-\alpha}\right)}{2\mu a}\left.\frac{\exp\left(\frac{2\mu a}{1-\alpha}(1+A+x)^{1-\alpha}\right)}{(1+A+x)^{(\alpha-2\gamma)}}\right|_{N-1}^{k}$$

$$-\frac{\exp\left(-\frac{2\mu a}{1-\alpha}(1+A+k)^{1-\alpha}\right)}{2\mu a}\int_{N-1}^{k}\frac{-(\alpha-2\gamma)\exp\left(\frac{2\mu a}{1-\alpha}(1+A+x)^{1-\alpha}\right)}{(1+A+x)^{(\alpha-2\gamma)+1}}dx$$

$$=\frac{\exp\left(-\frac{2\mu a}{1-\alpha}(1+A+k)^{1-\alpha}\right)}{2\mu a}\left\{\frac{\exp\left(\frac{2\mu a}{1-\alpha}(1+A+k)^{1-\alpha}\right)}{(1+A+k)^{(\alpha-2\gamma)}}-\frac{\exp\left(\frac{2\mu a}{1-\alpha}(A+N)^{1-\alpha}\right)}{(A+N)^{(\alpha-2\gamma)}}\right\}$$

$$+\frac{\exp\left(-\frac{2\mu a}{1-\alpha}(1+A+k)^{1-\alpha}\right)}{2\mu a}\int_{N-1}^{k}\frac{(\alpha-2\gamma)\exp\left(\frac{2\mu a}{1-\alpha}(1+A+x)^{1-\alpha}\right)}{(1+A+x)^{(\alpha-2\gamma)+1}}dx$$

$$=\frac{1}{2\mu a(1+A+k)^{(\alpha-2\gamma)}}-\frac{\exp\left(\frac{2\mu a}{1-\alpha}(A+N)^{1-\alpha}-\frac{2\mu a}{1-\alpha}(1+A+k)^{1-\alpha}\right)}{2\mu a(A+N)^{(\alpha-2\gamma)}}$$

$$+\frac{1}{2\mu a}\exp\left(-\frac{2\mu a}{1-\alpha}(1+A+k)^{1-\alpha}\right)\int_{N-1}^{k}\frac{(\alpha-2\gamma)\exp\left(\frac{2\mu a}{1-\alpha}(1+A+x)^{1-\alpha}\right)}{(1+A+x)^{(\alpha-2\gamma)+1}}dx. \qquad \text{(B.23)}$$

Now let compare the rate of convergence of the third term on the right-hand side of eqn. (B.23) with the term on the left-hand side of eqn. (B.23). For the ratio of the term on the left-hand side of eqn. (B.23) over the third term on the right-hand side of eqn. (B.23), after ignoring the multipliers, by L'Hôptial's rule we have

$$\lim_{k\to\infty}\frac{\exp\left(-\frac{2\mu a}{1-\alpha}(1+A+k)^{1-\alpha}\right)\int_{N-1}^{k}\frac{1/(1+A+x)^{(2\alpha-2\gamma)}}{\exp\left(-\frac{2\mu a}{1-\alpha}(1+A+x)^{1-\alpha}\right)}dx}{\exp\left(-\frac{2\mu a}{1-\alpha}(1+A+k)^{1-\alpha}\right)\int_{N-1}^{k}\frac{\exp\left(\frac{2\mu a}{1-\alpha}(1+A+x)^{1-\alpha}\right)}{(1+A+x)^{(\alpha-2\gamma)+1}}dx}$$

$$=\lim_{k\to\infty}\frac{\int_{N-1}^{k}\frac{\exp\left(\frac{2\mu a}{1-\alpha}(1+A+x)^{1-\alpha}\right)}{(1+A+x)^{(2\alpha-2\gamma)}}dx}{\int_{N-1}^{k}\frac{\exp\left(\frac{2\mu a}{1-\alpha}(1+A+x)^{1-\alpha}\right)}{(1+A+x)^{(\alpha-2\gamma)+1}}dx}$$

$$=\lim_{k\to\infty}\frac{\frac{\exp\left(\frac{2\mu a}{1-\alpha}(1+A+k)^{1-\alpha}\right)}{(1+A+k)^{(2\alpha-2\gamma)}}}{\frac{\exp\left(\frac{2\mu a}{1-\alpha}(1+A+k)^{1-\alpha}\right)}{(1+A+k)^{(\alpha-2\gamma)+1}}}$$

$$=\lim_{k\to\infty}(1+A+k)^{1-\alpha}=\infty,$$

which indicates that the term on the left-hand side of eqn. (B.23) goes to 0 at a lower rate than the third term on the right-hand side of eqn. (B.23) that can be expressed as

$$\frac{1}{2\mu a}\exp\left(-\frac{2\mu a}{1-\alpha}(1+A+k)^{1-\alpha}\right)\int_{N-1}^{k}\frac{(\alpha-2\gamma)\exp\left(\frac{2\mu a}{1-\alpha}(1+A+x)^{1-\alpha}\right)}{(1+A+x)^{(\alpha-2\gamma)+1}}dx.$$

$$=o\left(\exp\left(-\frac{2\mu a}{1-\alpha}(1+A+k)^{1-\alpha}\right)\int_{N-1}^{k}\frac{1\big/(1+A+x)^{(2\alpha-2\gamma)}}{\exp\left(-\frac{2\mu a}{1-\alpha}(1+A+x)^{1-\alpha}\right)}dx\right).$$

Since the second term on the right-hand side of eqn.(B.23) goes to 0 at a higher rate than the term on the left-hand side of eqn. (B.23) (second term = $o$(term on the left-hand side)), then the term on the left-hand side of eqn. (B.23) and the first tem on the right-hand side of eqn. (B.23) go to 0 at the same rate, which means

$$\exp\left(-\frac{2\mu a}{1-\alpha}(1+A+k)^{1-\alpha}\right)\int_{N-1}^{k}\frac{1\big/(1+A+x)^{(2\alpha-2\gamma)}}{\exp\left(-\frac{2\mu a}{1-\alpha}(1+A+x)^{1-\alpha}\right)}dx=O\left(\frac{1}{2\mu a(1+A+k)^{(\alpha-2\gamma)}}\right).$$

Similarly, for the upper bound we get

$$\exp\left(-\frac{2\mu a}{1-\alpha}(1+A+k)^{1-\alpha}\right)\int_{N}^{k+1}\frac{1\big/(1+A+x)^{(2\alpha-2\gamma)}}{\exp\left(-\frac{2\mu a}{1-\alpha}(1+A+x)^{1-\alpha}\right)}dx=O\left(\frac{1}{2\mu a(2+A+k)^{(\alpha-2\gamma)}}\right).$$

Therefore, when $0.5<\alpha<1$, by inequality (B.22) and relationship (B.21a), we have $T_k\sum_{i=N}^{k}a_i^2\big/(c_i^2T_i)$ = $O\left(1\big/k^{(\alpha-2\gamma)}\right)$, which indicates that the second term

$$T_k\sum_{i=N}^{k}\frac{a_i^2E\left\|\hat{\mathbf{g}}_i(\hat{\boldsymbol{\theta}}_i)\right\|^2}{T_i}$$

on the right-hand side of inequality (B.7) goes to 0 at the rate of

$$O\left(\frac{1}{k^{(\alpha-2\gamma)}}\right). \tag{B.24}$$

When $\alpha$ = 1, by ignoring the multiplier of constant scalar, the term in inequality (B.21b) is

$$\frac{1}{(1+A+k)^{2\mu a}}\sum_{i=N}^{k}\frac{1/(1+A+i)^{(2-2\gamma)}}{1/(1+A+i)^{2\mu a}}=\frac{1}{(1+A+k)^{2\mu a}}\sum_{i=N}^{k}(1+A+i)^{2\mu a-(2-2\gamma)}.$$

Moreover, when $2a\mu \geq 2-2\gamma$, $(1+A+i)^{2\mu a-(2-2\gamma)}$ is an increasing function on *i*, then we have

$$\begin{aligned}\int_{N-1}^{k}(1+A+x)^{2\mu a-(2-2\gamma)}dx &\leq \sum_{i=N}^{k}(1+A+i)^{2\mu a-2-2\gamma}\\ &\leq \int_{N}^{k+1}(1+A+x)^{2\mu a-(2-2\gamma)}dx,\end{aligned}$$

which indicates

$$\begin{aligned}\frac{(1+A+k)^{2\mu a-(1-2\gamma)}-(A+N)^{2\mu a-(1-2\gamma)}}{2\mu a-(1-2\gamma)} &\leq \sum_{i=N}^{k}(1+A+k)^{2\mu a-1-2\gamma}\\ &\leq \frac{(2+A+k)^{2\mu a-(1-2\gamma)}-(1+A+N)^{2\mu a-(1-2\gamma)}}{2\mu a-(1-2\gamma)}.\end{aligned}$$

Similarly, when $2\mu a < 2-2\gamma$, $(1+A+i)^{2\mu a-(2-2\gamma)}$ is a decreasing function on *i*, then we have

$$\int_{N}^{k+1} (1+A+x)^{2\mu a-(2-2\gamma)}dx \le \sum_{i=N}^{k} (1+A+i)^{2\mu a-2-2\gamma}$$

$$\le \int_{N-1}^{k} (1+A+x)^{2\mu a-(2-2\gamma)}dx,$$

and due to condition (vii) $2\mu a > 1-2\gamma$, we have

$$\frac{(2+A+k)^{2\mu a-(1-2\gamma)}-(1+A+N)^{2\mu a-(1-2\gamma)}}{2\mu a-(1-2\gamma)} \le \sum_{i=N}^{k} (1+A+k)^{2\mu a-1-2\gamma}$$

$$\le \frac{(1+A+k)^{2\mu a-(1-2\gamma)}-(A+N)^{2\mu a-(1-2\gamma)}}{2a\mu-(1-2\gamma)}.$$

Then because of condition (vii), we have the lower and upper bounds of $T_k \sum_{i=N}^{k} a_i^2 \big/ (c_i^2 T_i)$ for both cases go 0 at the rate of

$$O\left(\frac{1}{k^{2a\mu}} \times k^{2a\mu-(1-2\gamma)}\right) = O\left(\frac{1}{k^{(1-2\gamma)}}\right)$$

Therefore, when $\alpha = 1$ we have $T_k \sum_{i=N}^{k} a_i^2 \big/ (c_i^2 T_i) = O\left(1 \big/ k^{(1-2\gamma)}\right)$, which indicates that the second term

$$T_k \sum_{i=N}^{k} \frac{a_i^2 E\left\|\hat{\mathbf{g}}_i(\hat{\boldsymbol{\theta}}_i)\right\|^2}{T_i}$$

on the right-hand side of inequality (B.7) goes to 0 at the rate of

$$O\left(\frac{1}{k^{(1-2\gamma)}}\right).$$

Overall, we have that the second term on the right-hand side of inequality (B.7) goes to 0 at the rate of

$$O\left(\frac{1}{k^{(\alpha-2\gamma)}}\right).$$

Now let us start to discuss the third term on the right-hand side of inequality (B.7). We have

$$E\left[\left(\hat{\boldsymbol{\theta}}_k-\boldsymbol{\theta}^*\right)^T\left(\boldsymbol{g}_k(\hat{\boldsymbol{\theta}}_k)-\hat{\boldsymbol{g}}_k(\hat{\boldsymbol{\theta}}_k)\right)\right]=E\left\{E\left[\left(\hat{\boldsymbol{\theta}}_k-\boldsymbol{\theta}^*\right)^T\left(\boldsymbol{g}_k(\hat{\boldsymbol{\theta}}_k)-\hat{\boldsymbol{g}}_k(\hat{\boldsymbol{\theta}}_k)\right)\middle|\hat{\boldsymbol{\theta}}_k\right]\right\}$$

$$=E\left\{\left(\hat{\boldsymbol{\theta}}_k-\boldsymbol{\theta}^*\right)^T E\left[\boldsymbol{g}_k(\hat{\boldsymbol{\theta}}_k)-\hat{\boldsymbol{g}}_k(\hat{\boldsymbol{\theta}}_k)\middle|\hat{\boldsymbol{\theta}}_k\right]\right\}.$$

By Cauchy–Schwarz inequality, we have

$$\left|\left(\hat{\boldsymbol{\theta}}_k-\boldsymbol{\theta}^*\right)^T E\left[\boldsymbol{g}_k(\hat{\boldsymbol{\theta}}_k)-\hat{\boldsymbol{g}}_k(\hat{\boldsymbol{\theta}}_k)\middle|\hat{\boldsymbol{\theta}}_k\right]\right|\le\left\|\hat{\boldsymbol{\theta}}_k-\boldsymbol{\theta}^*\right\|\left\|E\left[\boldsymbol{g}_k(\hat{\boldsymbol{\theta}}_k)-\hat{\boldsymbol{g}}_k(\hat{\boldsymbol{\theta}}_k)\middle|\hat{\boldsymbol{\theta}}_k\right]\right\|.$$

By condition (ii), (iii), (vi) and using similar arguments as in the Lemma 1 of Spall (1992), we have $\left\|E\left[\hat{\boldsymbol{g}}_k(\hat{\boldsymbol{\theta}}_k)-\boldsymbol{g}_k(\hat{\boldsymbol{\theta}}_k)\middle|\hat{\boldsymbol{\theta}}_k\right]\right\|=O(c_k^2)$. Then there exists $\zeta>0$ such that

$$\left(\hat{\boldsymbol{\theta}}_k-\boldsymbol{\theta}^*\right)^T E\left[\boldsymbol{g}_k(\hat{\boldsymbol{\theta}}_k)-\hat{\boldsymbol{g}}_k(\hat{\boldsymbol{\theta}}_k)\middle|\hat{\boldsymbol{\theta}}_k\right]\le\left\|\hat{\boldsymbol{\theta}}_k-\boldsymbol{\theta}^*\right\|\left\|E\left[\boldsymbol{g}_k(\hat{\boldsymbol{\theta}}_k)-\hat{\boldsymbol{g}}_k(\hat{\boldsymbol{\theta}}_k)\middle|\hat{\boldsymbol{\theta}}_k\right]\right\|$$

$$\le\zeta c_k^2\left\|\hat{\boldsymbol{\theta}}_k-\boldsymbol{\theta}^*\right\|,$$

which indicates that

$$E\left\{\left(\hat{\boldsymbol{\theta}}_k-\boldsymbol{\theta}^*\right)^T E\left[\boldsymbol{g}_k(\hat{\boldsymbol{\theta}}_k)-\hat{\boldsymbol{g}}(\hat{\boldsymbol{\theta}}_k)\middle|\hat{\boldsymbol{\theta}}_k\right]\right\}\le\zeta c_k^2 E\left\|\hat{\boldsymbol{\theta}}_k-\boldsymbol{\theta}^*\right\|.$$

By Jensen's inequality, we further have $E\left\|\hat{\boldsymbol{\theta}}_k - \boldsymbol{\theta}^*\right\| \le \sqrt{E\left\|\hat{\boldsymbol{\theta}}_k - \boldsymbol{\theta}^*\right\|^2}$. Then the third term on the right-hand side of inequality (B.7) is bounded by the upper bound

$$T_k \sum_{i=N}^{k} \frac{2\zeta a_i c_i^2 \sqrt{E\left\|\hat{\boldsymbol{\theta}}_k - \boldsymbol{\theta}^*\right\|^2}}{T_i},$$

which implies that inequality (B.7) can be rewritten by substituting the third term on the right-hand side by its upper bound

$$E\left\|\hat{\boldsymbol{\theta}}_{k+1} - \boldsymbol{\theta}^*\right\|^2$$

$$\le T_k E\left\|\hat{\boldsymbol{\theta}}_N - \boldsymbol{\theta}^*\right\|^2 + T_k \sum_{i=N}^{k} \frac{a_i^2 E\left\|\hat{\boldsymbol{g}}_i(\hat{\boldsymbol{\theta}}_i)\right\|^2}{T_i} + T_k \sum_{i=N}^{k} \frac{2\zeta a_i c_i^2 \sqrt{E\left\|\hat{\boldsymbol{\theta}}_k - \boldsymbol{\theta}^*\right\|^2}}{T_i}. \qquad \text{(B.25)}$$

Before showing the rate of convergence of the third term on the right-hand side of inequality (B.25), let us first show that $E\left\|\hat{\boldsymbol{\theta}}_k - \boldsymbol{\theta}^*\right\|^2$ converges to 0. By way of contradiction, suppose $E\left\|\hat{\boldsymbol{\theta}}_k - \boldsymbol{\theta}^*\right\|^2$ does not converge to 0. For the inequality (B.25), we have shown that the first term and the second term on the right-hand side converge to 0. Since the term on the left-hand side of inequality (B.25) $E\left\|\hat{\boldsymbol{\theta}}_k - \boldsymbol{\theta}^*\right\|^2$ does not converge to 0, then there exists $\eta > 0$, such that we have a relationship between the third term on the right-hand side of inequality (B.25) and the term $E\left\|\hat{\boldsymbol{\theta}}_k - \boldsymbol{\theta}^*\right\|^2$ on the left-hand side of inequality (B.25) as

$$E\left\|\hat{\boldsymbol{\theta}}_k-\boldsymbol{\theta}^*\right\|^2 \le \eta T_k \sum_{i=N}^{k} \frac{2\zeta a_i c_i^2 \sqrt{E\left\|\hat{\boldsymbol{\theta}}_k-\boldsymbol{\theta}^*\right\|^2}}{T_i} \Rightarrow 1 \le \eta T_k \sum_{i=N}^{k} \frac{2\zeta a_i c_i^2}{T_i \sqrt{E\left\|\hat{\boldsymbol{\theta}}_k-\boldsymbol{\theta}^*\right\|^2}}$$

Meanwhile, by the similar arguments for the rate of convergence of the second term on the right-hand side of inequality (B.25), we have $T_k \sum_{i=N}^{k} \frac{2\zeta a_i c_i^2}{T_i}$ goes to 0 at the rate of $O\left(1/k^{2\gamma}\right)$. But $E\left\|\hat{\boldsymbol{\theta}}_k-\boldsymbol{\theta}^*\right\|^2$ does not converge to 0, then we have

$$1 \le \eta T_k \sum_{i=N}^{k} \frac{2\zeta a_i c_i^2}{T_i \sqrt{E\left\|\hat{\boldsymbol{\theta}}_k-\boldsymbol{\theta}^*\right\|^2}} \to 0,$$

which is a contradiction. Thus, we have $E\left\|\hat{\boldsymbol{\theta}}_k-\boldsymbol{\theta}^*\right\|^2 \to 0$

Since the first, second, and the third term on the right-hand side of inequality (B.25) goes to 0 at least at the polynomial rates in the big-*O* sense, then we can assume that $E\left\|\hat{\boldsymbol{\theta}}_k-\boldsymbol{\theta}^*\right\|^2 \to 0$ in terms of $O(k^{-2t})$ for $t > 0$. Next, we solve out the value of *t*. Since $E\left\|\hat{\boldsymbol{\theta}}_k-\boldsymbol{\theta}^*\right\|^2 \to 0$ in terms of $O(k^{-2t})$, then by the similar arguments for the rate of convergence of the second term on the right-hand side of inequality (B.25), we have that the third term on the right-hand side of inequality (B.25) goes to 0 at the rate of

$$O\left(\frac{1}{k^{t+2\gamma}}\right). \tag{B.26}$$

In summary, by formula (B.20), (B.24) and (B.26), we have the convergent rate of each term on the right-hand side of inequality (B.25), which is summarized in Table B.1.

**Table B.1** The convergence rate of each term in inequality (B.25).

| | $\alpha = 1$ | $0.5 < \alpha < 1$ |
|---|---|---|
| First term | $O\left(\frac{1}{k^{2a\mu}}\right)$ | $O\left(e^{-\frac{2a\mu}{1-\alpha}(1+A+k)^{1-\alpha}}\right)$ |
| Second term | $O\left(\frac{1}{k^{(1-2\gamma)}}\right)$ | $O\left(\frac{1}{k^{(\alpha-2\gamma)}}\right)$ |
| Third term | $O\left(\frac{1}{k^{(t+2\gamma)}}\right)$ | $O\left(\frac{1}{k^{(t+2\gamma)}}\right)$ |
| $E\left\lVert\hat{\boldsymbol{\theta}}_k - \boldsymbol{\theta}^*\right\rVert^2$ | $O\left(\frac{1}{k^{2t}}\right)$ | $O\left(\frac{1}{k^{2t}}\right)$ |

From Table B.1, we know that the first term goes to 0 at a higher rate than the second and the third terms (first term = $O$(second term)). Moreover, since all terms on the right-hand side of inequality (B.25) are positive, then the multiplier of the leading term in the big-$O$ function must be positive for all terms. In Table B.1, we have that when $t + 2\gamma \geq \alpha - 2\gamma$, the third term go to 0 at a rate not slower than the rate at which the second term goes to 0 (third term = $O$(second term)), which implies the rate of convergence of $E\left\lVert\hat{\boldsymbol{\theta}}_k - \boldsymbol{\theta}^*\right\rVert^2$ =

$O\left(1/k^{(\alpha-2\gamma)}\right)$. Therefore, $2t = \alpha - 2\gamma$, which follows that $(\alpha-2\gamma)/2 + 2\gamma \ \geq \ \alpha - 2\gamma$, implying $\gamma \ \geq \ \alpha/6$. Thus, when $\gamma \ \geq \ \alpha/6$, $E\left\|\hat{\boldsymbol{\theta}}_k - \boldsymbol{\theta}^*\right\|^2 = O(1/k^{(\alpha-2\gamma)})$. On the other side, when $t + 2\gamma < \alpha - 2\gamma$, the second term go to 0 at a rate not slower than the rate at which the third term goes to 0 (second term = $O$(third term)), which implies the rate of convergence of $E\left\|\hat{\boldsymbol{\theta}}_k - \boldsymbol{\theta}^*\right\|^2 = O\left(1/k^{(t+2\gamma)}\right)$. Therefore, $2t = t + 2\gamma$, which follows that $2\gamma + 2\gamma < \alpha - 2\gamma$, implying $\gamma < \alpha/6$. Thus, when $\gamma < \alpha/6$, $E\left\|\hat{\boldsymbol{\theta}}_k - \boldsymbol{\theta}^*\right\|^2 = O(1/k^{4\gamma})$. In all, we have

$$E\left\|\hat{\boldsymbol{\theta}}_k - \boldsymbol{\theta}^*\right\|^2 = \begin{cases} O\left(\dfrac{1}{k^{\alpha-2\gamma}}\right) & \gamma \geq \dfrac{\alpha}{6} \\ O\left(\dfrac{1}{k^{4\gamma}}\right) & \gamma < \dfrac{\alpha}{6}. \end{cases} \tag{B.27}$$

Q.E.D.

Theorem B.1 shows the rate of convergence of SPSA in the big-$O$ sense for the problem of time-varying loss function. By these results, we see that when $\gamma = \alpha/6$, $\alpha - 2\gamma = 4\gamma$, which make these two forms of convergence rate in eqn. (B.27) to be consistent at the critical point. The optimal rate of convergence for $E\left\|\hat{\boldsymbol{\theta}}_k - \boldsymbol{\theta}^*\right\|^2$ can be achieved when $\alpha = 1$ and $\gamma = 1/6$, which is $O(1/k^{2/3})$. Time-varying problem is a more general case of the fixed loss function problem, and compared with the rate of convergence result in Spall (1992), we find that our result is consistent with Spall (1992) for the case of $\gamma \ \geq \ \alpha/6$. We also discuss the case of $\gamma < \alpha/6$, which is not considered in Spall (1992), and our result indicates that only using the coefficients that satisfy $\gamma \geq \alpha/6$ might be enough for practical users.

## B.4 Numerical Results

From condition (i) of Theorem B.1, we have $\alpha > 0$, $\gamma > 0$, $2\alpha - 2\gamma > 1$ and $\alpha \leq 1$. In the following, we do the numerical experiments to test the results of the rate of convergence in eqn. (B.27). We will pick two sets of $(\alpha, \gamma)$ with one set satisfying $\gamma \geq \alpha/6$ and with the other set satisfying $\gamma < \alpha/6$. We will plot figures to show that the results in eqn. (B.27) are true.

The time-varying loss function we consider is a simple periodical function $L_k(\boldsymbol{\theta})$ = $\boldsymbol{\theta}^T\boldsymbol{\theta}(1-1/(\text{mod}(k,30)+2))$ , where $\text{mod}(\cdot,\cdot)$ is a modulo function that returns the remainder of division of one number by another. The noises $\varepsilon$ are i.i.d. $N(0,1)$. We set $a$ = 0.1, $c$ = 0.1, and $A$ = 1000. The dimension of the problem is $p$ = 10. The initial guess is set as $\hat{\boldsymbol{\theta}}_0 = [1, \ldots, 1]^T$. We do 10 replicates with 100,000 iterations in each replicate. In Figure B.1, we show the result for the set $\alpha = 1$ and $\gamma = \alpha/10$, and we see that this set of $(\alpha, \gamma)$ satisfyies $\gamma < \alpha/6$. Figure B.1 plots the sample means of the values of $k^{\beta}\left\|\hat{\boldsymbol{\theta}}_k - \boldsymbol{\theta}^*\right\|^2$ for $\beta = 4\gamma = 2/5$ and $\beta = \alpha - 2\gamma = 4/5$. We see that the sample mean of $k^{4\gamma}\left\|\hat{\boldsymbol{\theta}}_k - \boldsymbol{\theta}^*\right\|^2$ goes to a constant scalar, while the sample mean of $k^{\alpha-2\gamma}\left\|\hat{\boldsymbol{\theta}}_k - \boldsymbol{\theta}^*\right\|^2$ is

diverging. These results are consistent with eqn. (B.27), which indicates under the case of $\gamma < \alpha/6$, the rate of convergence of $E\left\|\hat{\boldsymbol{\theta}}_k - \boldsymbol{\theta}^*\right\|^2 = O(1/k^{4\gamma})$ and $1/k^{\alpha-2\gamma}$ goes to 0 at a faster rate than $E\left\|\hat{\boldsymbol{\theta}}_k - \boldsymbol{\theta}^*\right\|^2$. In Figure B.2, we show the result for the set $\alpha = 1$ and $\gamma = \alpha/5$, and we see that this set of $(\alpha,\gamma)$ satisfies $\gamma \geq \alpha/6$. Figure B.2 plots the sample means of the values of $k^{\beta}\left\|\hat{\boldsymbol{\theta}}_k - \boldsymbol{\theta}^*\right\|^2$ for $\beta = 4\gamma = 4/5$ and $\beta = \alpha - 2\gamma = 3/5$. We see that the sample mean of $k^{\alpha-2\gamma}\left\|\hat{\boldsymbol{\theta}}_k - \boldsymbol{\theta}^*\right\|^2$ goes to a constant scalar, while the sample mean of $k^{4\gamma}\left\|\hat{\boldsymbol{\theta}}_k - \boldsymbol{\theta}^*\right\|^2$ goes to a really large value. These results are consistent with eqn. (B.27), which indicates under the case of $\gamma \geq \alpha/6$, the rate of convergence of $E\left\|\hat{\boldsymbol{\theta}}_k - \boldsymbol{\theta}^*\right\|^2 = O(1/k^{\alpha-2\gamma})$ and $1/k^{4\gamma}$ goes to 0 at a faster rate than $E\left\|\hat{\boldsymbol{\theta}}_k - \boldsymbol{\theta}^*\right\|^2$.

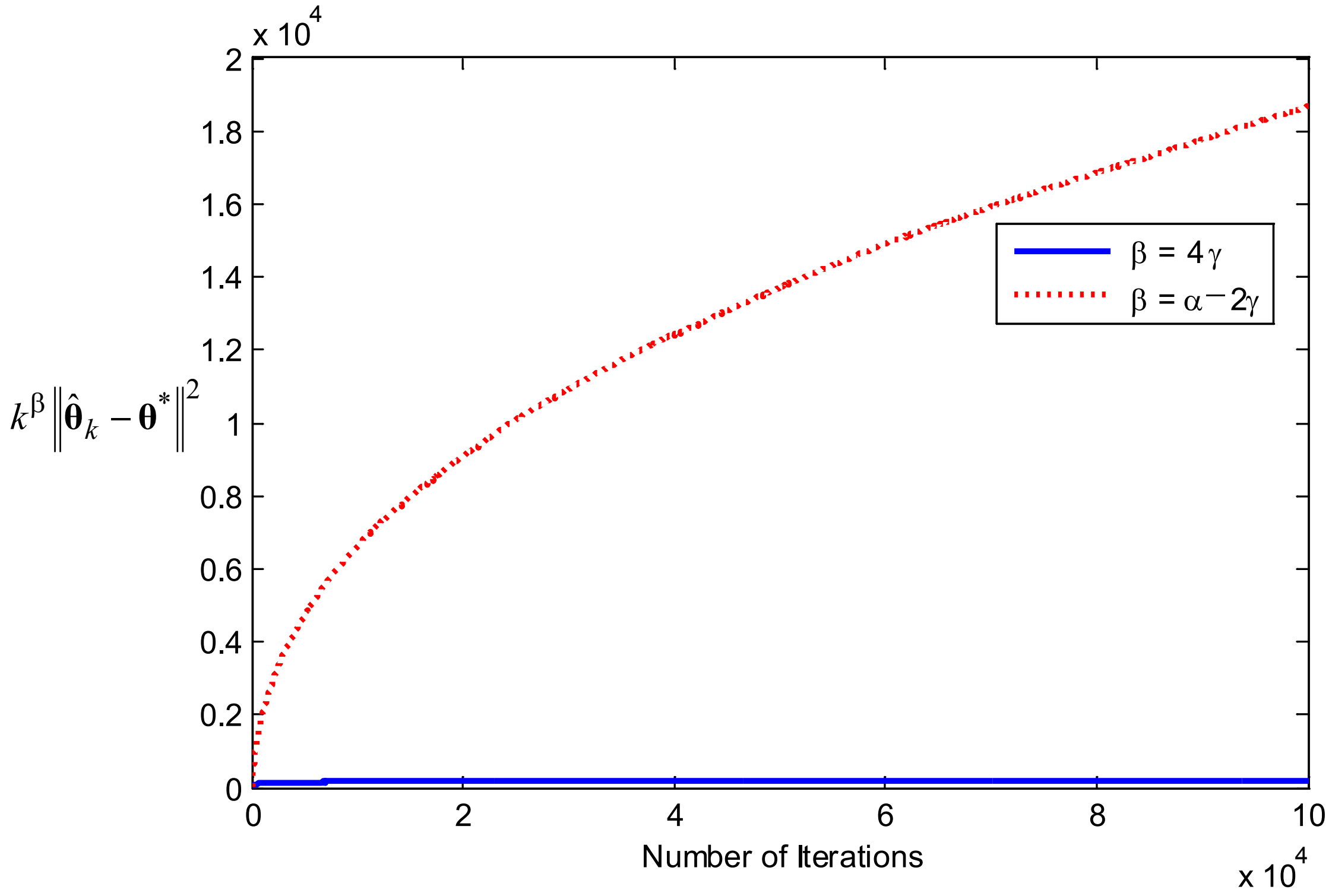

**Figure B.1** The numerical result for the case of $\alpha$ = 1, $\gamma$ = $\alpha/10$. The set $(\alpha, \gamma)$ satisfies $\gamma < \alpha/6$. The vertical axis represents the sample mean of $k^{\beta}\left\|\hat{\boldsymbol{\theta}}_k - \boldsymbol{\theta}^*\right\|^2$ for $\beta$ = $4\gamma$ = $2/5$ and $\beta$ = $\alpha - 2\gamma$ = $4/5$. We see that for the case of $\gamma < \alpha/6$, the curve for the sample mean of $k^{4\gamma}\left\|\hat{\boldsymbol{\theta}}_k - \boldsymbol{\theta}^*\right\|^2$ is flat for large *k*. However, the curve for the sample mean of $k^{\alpha-2\gamma}\left\|\hat{\boldsymbol{\theta}}_k - \boldsymbol{\theta}^*\right\|^2$ is diverging. The number of replicates is 10, and the number of iterations in each replicate is 100,000.

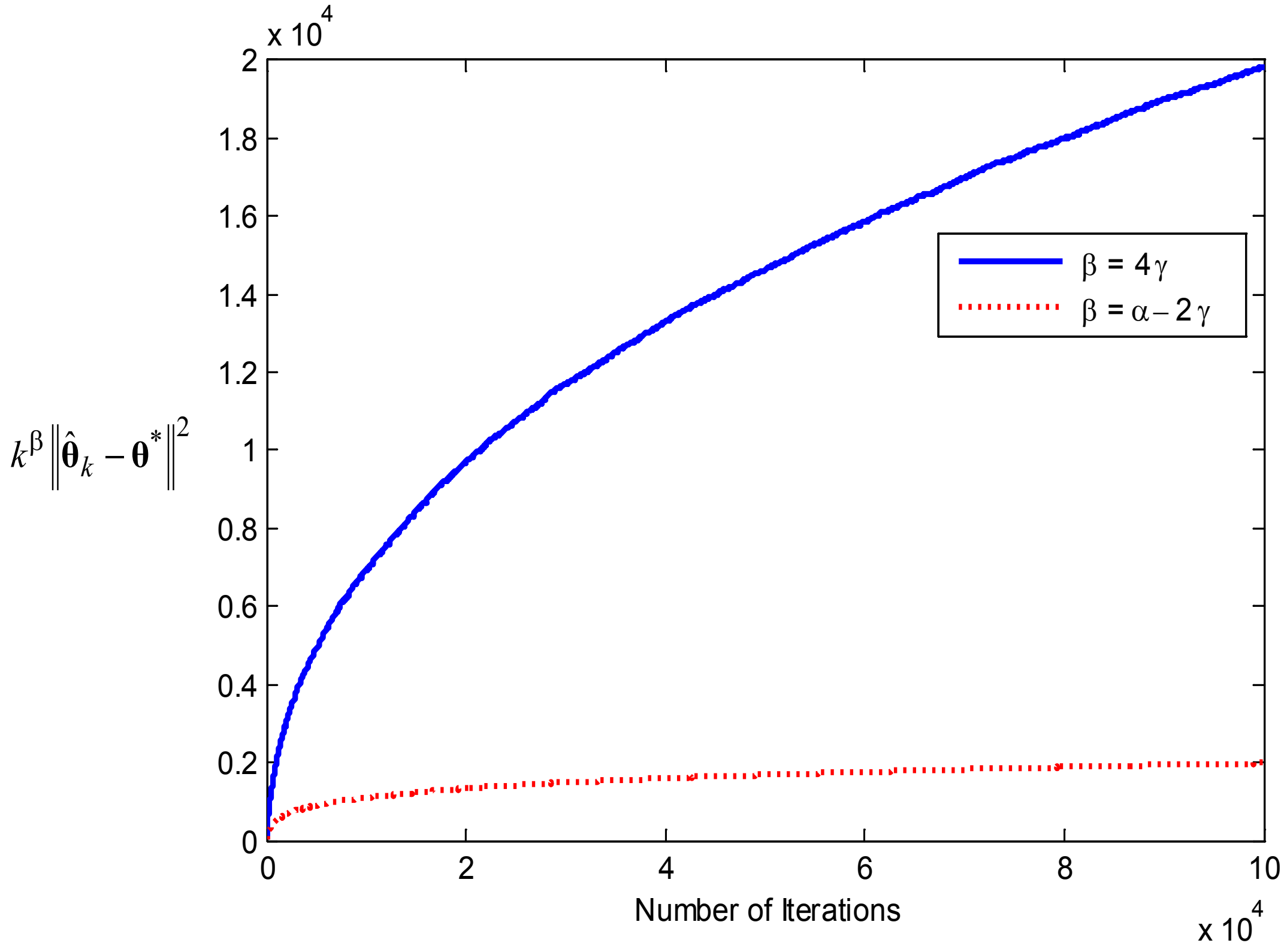

**Figure B.2.** The numerical result for the case of $\alpha = 1$, $\gamma = \alpha/5$. The set $(\alpha, \gamma)$ satisfies $\gamma \geq \alpha/6$. The vertical axis represents the sample mean of $k^{\beta}\left\|\hat{\boldsymbol{\theta}}_k - \boldsymbol{\theta}^*\right\|^2$ for $\beta = 4\gamma = 4/5$ and $\beta = \alpha - 2\gamma = 3/5$. We see that for the case of $\gamma \geq \alpha/6$, the curve for the sample mean of $k^{\alpha-2\gamma}\left\|\hat{\boldsymbol{\theta}}_k - \boldsymbol{\theta}^*\right\|^2$ is flat for large $k$. However, the curve for the sample mean of $k^{4\gamma}\left\|\hat{\boldsymbol{\theta}}_k - \boldsymbol{\theta}^*\right\|^2$ is diverging. The number of replicates is 10, and the number of iterations in each replicate is 100,000.

Overall, we see that under some conditions, the rate of convergence of $\{\hat{\boldsymbol{\theta}}_k\}$ is based on the value of $\alpha$ and $\gamma$. When the set of $(\alpha, \gamma)$ satisfies $\gamma < \alpha/6$, we have $E\left\|\hat{\boldsymbol{\theta}}_k - \boldsymbol{\theta}^*\right\|^2 = O(1/k^{4\gamma})$. When $\gamma \geq \alpha/6$, we have $E\left\|\hat{\boldsymbol{\theta}}_k - \boldsymbol{\theta}^*\right\|^2 = O(1/k^{\alpha-2\gamma})$. We see that these rates are reasonable, so the SPSA method is a quite good method for time-varying loss function.

# B.5 Conclusions

This appendix discusses the rate of convergence of SPSA for time-varying loss function rather than the fixed loss function of former rate of convergence analysis (Spall, 1992). Of course, the MSE (and closely related variance of the estimator) are popular measures of accuracy in practical problems. We show that the rate of convergence in terms of MSE $E\left\|\hat{\boldsymbol{\theta}}_k - \boldsymbol{\theta}^*\right\|^2$ is $O(1/k^{\alpha-2\gamma})$ when $\gamma \geq \alpha/6$, and $O(1/k^{4\gamma})$ when $\gamma < \alpha/6$. By the result of the convergence rate discussed here, we see that SPSA is an efficient algorithm in solving the time-varying problems in the big-*O* sense. Time-varying problem is a more general case of the fixed loss function problem considered in Spall (1992), we find that our result is consistent with Spall (1992) for the case of $\gamma \geq \alpha/6$. We also discuss the case of $\gamma < \alpha/6$, which is not discussed in Spall (1992), and our result indicates that only considering the coefficients that satisfy $\gamma \geq \alpha/6$ might be enough for practical users, who want to use SPSA to solve problems, because in our results the case of $\gamma < \alpha/6$ may not achieve higher convergence rate than the case of $\gamma \geq \alpha/6$ in the big-$O$ sense.

# Appendix C

# Numerical Experiments Results on the Properties of the Upper Bound for Mean Squared Error

After discussing the properties of the upper bound on $E\left\|\hat{\boldsymbol{\theta}}_k - \boldsymbol{\theta}^*\right\|^2$ theoretically in Chapter 3, we now present the result of some tests to see the properties of the upper bound numerically. In this appendix, we only consider the case when $\Delta_{ki}$ are independent Bernoulli random variables taking the values $\pm 1$ with probability $1/2$. Thus in the inequality (3.3a) and (3.3b), we replace $l$ with $p$.

# C.1 Effects of Coefficients in the Gain Sequence

In inequality (3.3a), we find the upper bound is a function of $\alpha$, $A$, $a$, $\mu$, $b$, $l$, and $E\left\|\hat{\boldsymbol{\theta}}_0 - \boldsymbol{\theta}^*\right\|^2$. In this section, we do the tests on 1) the importance of the first term in the upper bound in inequality (3.3a), and 2) the effects of the coefficients $\alpha$, $A$, and $a$ on the upper bound in inequality (3.3a). We use the sensitivity analysis for these tests, which means in each test we only change one coefficient and fix all the rest. In Section 3.2, we discuss the properties of the upper bound, and here we do the numerical tests to confirm these results. From Proposition 3.1, we know that the importance of the first term in the upper bound in inequality (3.3a) is strictly decreasing (monotone) with $k$. In addition, from Propositions 3.2 and 3.3, the effects of the coefficients $\alpha$, $A$, and $a$ on the upper bound in inequality (3.3a) are monotone at different stages (early iterations and late iterations). Since all these relationships are monotone, we can just pick one reasonable sets of coefficients as the base case. For example, suppose the base case of coefficients is: $\alpha = 0.65$, $A = 100$, $a = 0.1$. Other parameters of the upper bound are: $\mu = 1$, $b = 10$, $l = 100$, $E\left\|\hat{\boldsymbol{\theta}}_0 - \boldsymbol{\theta}^*\right\|^2 = 100$. The number of iterations is $m = 10000$.

Note that a potential future research direction would be to consider some or all of these SA coefficients as "nuisance parameters" (e.g., Basu, 1977; Spall, 1989; and Spall and Garner, 1990). Such a characterization might allow a formal (theoretical) quantification of the effect of the coefficients on the MSEs of the SA parameter estimates. We do not pursue this direction here, rather relying on numerical experiments.

First of all, we show in Figure C.1 about the decaying of the importance of the first term

$$\exp\left(\frac{2\mu a(1+A)^{1-\alpha}}{1-\alpha}-\frac{2\mu a(1+A+k)^{1-\alpha}}{1-\alpha}\right)E\left\|\hat{\boldsymbol{\theta}}_0-\boldsymbol{\theta}^*\right\|^2$$

over the upper bound in the inequality (3.3a). We plot the cases of $\alpha$ = 0.55 and $\alpha$ = 0.65. We see that the ratio of the first term over the upper bound (we call it the proportion of the first term) is decreasing and it decreases faster for the case when $\alpha$ is smaller. From this phenomenon, we see that in the early iterations the first term of the upper bound is more significant than the second term.

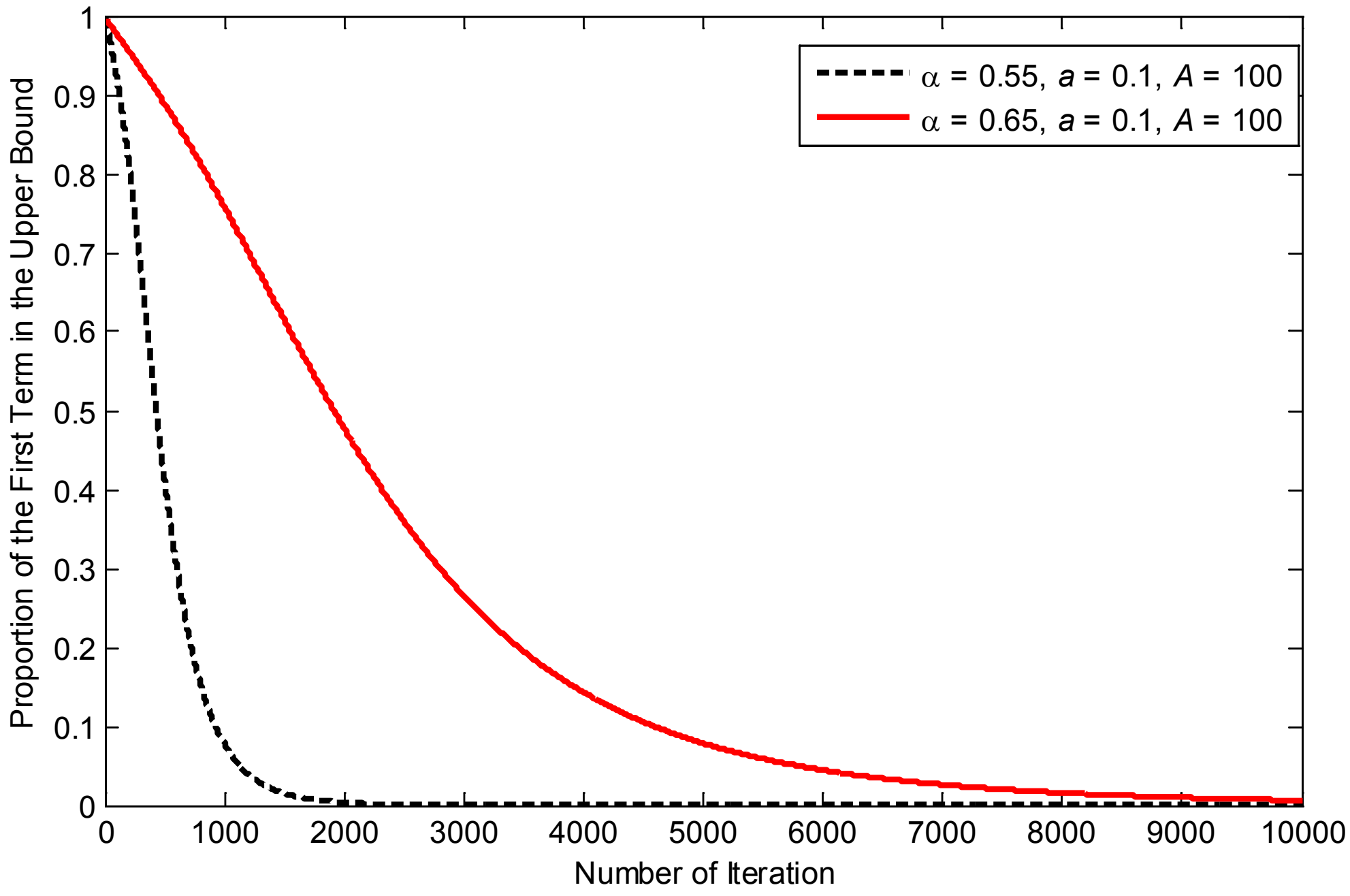

**Figure C.1** Proportion of the first term over the upper bound for $\alpha$ = 0.55 and $\alpha$ = 0.65. The proportion of the first term decreases with $k$.

The upper bound is a function on the coefficients of $\alpha$, $a$ and $A$. In the following numerical tests, we discuss the effects of these coefficients on the value of the upper bound in inequality (3.3a). These numerical tests just provide the properties of the upper bound and do not provide the exact rules of coefficients selection.

In Figure C.2, we find that a lower value of $\alpha$ provides a smaller upper bound in the early iterations; while a higher value of $\alpha$ provides a lower upper bound in the later iterations. Thus, we see that there is a crossing point for the upper bounds based on different values of $\alpha$.

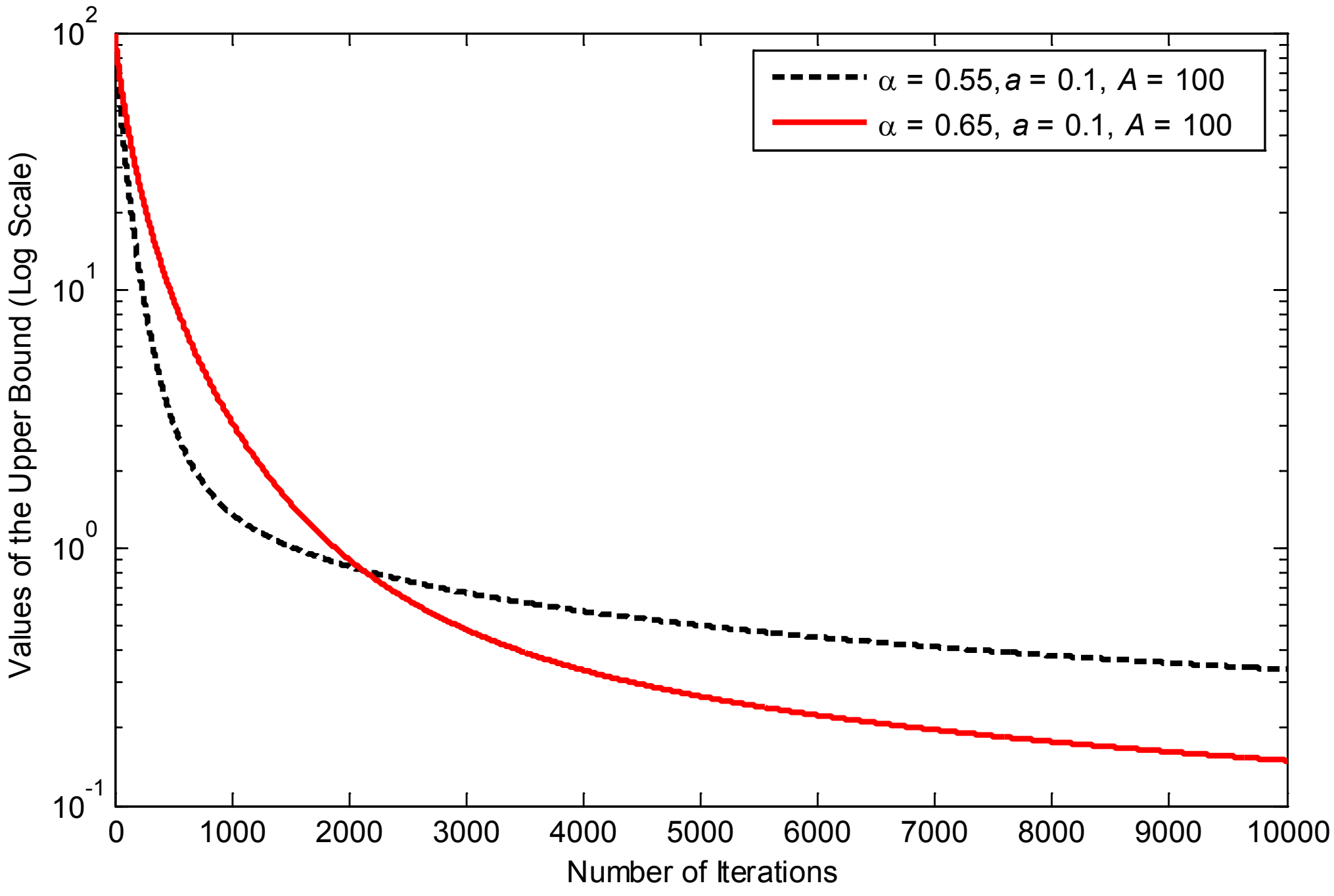

**Figure C.2** The plots of the upper bounds for the case of $\alpha = 0.55$ and $\alpha = 0.65$. Smaller $\alpha$ provides lower upper bound in the early iterations, and bigger $\alpha$ provides lower upper bound in the later iterations. The vertical axis is in the logarithm scale.

In Figure C.3, we consider the effect of $a$, and we see that larger $a$ provides lower upper bound for early iterations, but it does not lead to lower upper bound in the later iterations. Comparing the upper bound for $a = 0.1$ with the upper bound for $a = 1$, we see that the upper bound for $a = 1$ provides lower values in the early iterations, but $a = 1$ is too large for the later iterations. Thus, there is a clear crossing point between the upper bounds for $a = 0.1$ and $a = 1$, which can be seen in Figure C.3. In Figure C.4, we discuss the effect of $A$, and we see that the effect of $A$ disappears gradually with the increase of $k$.

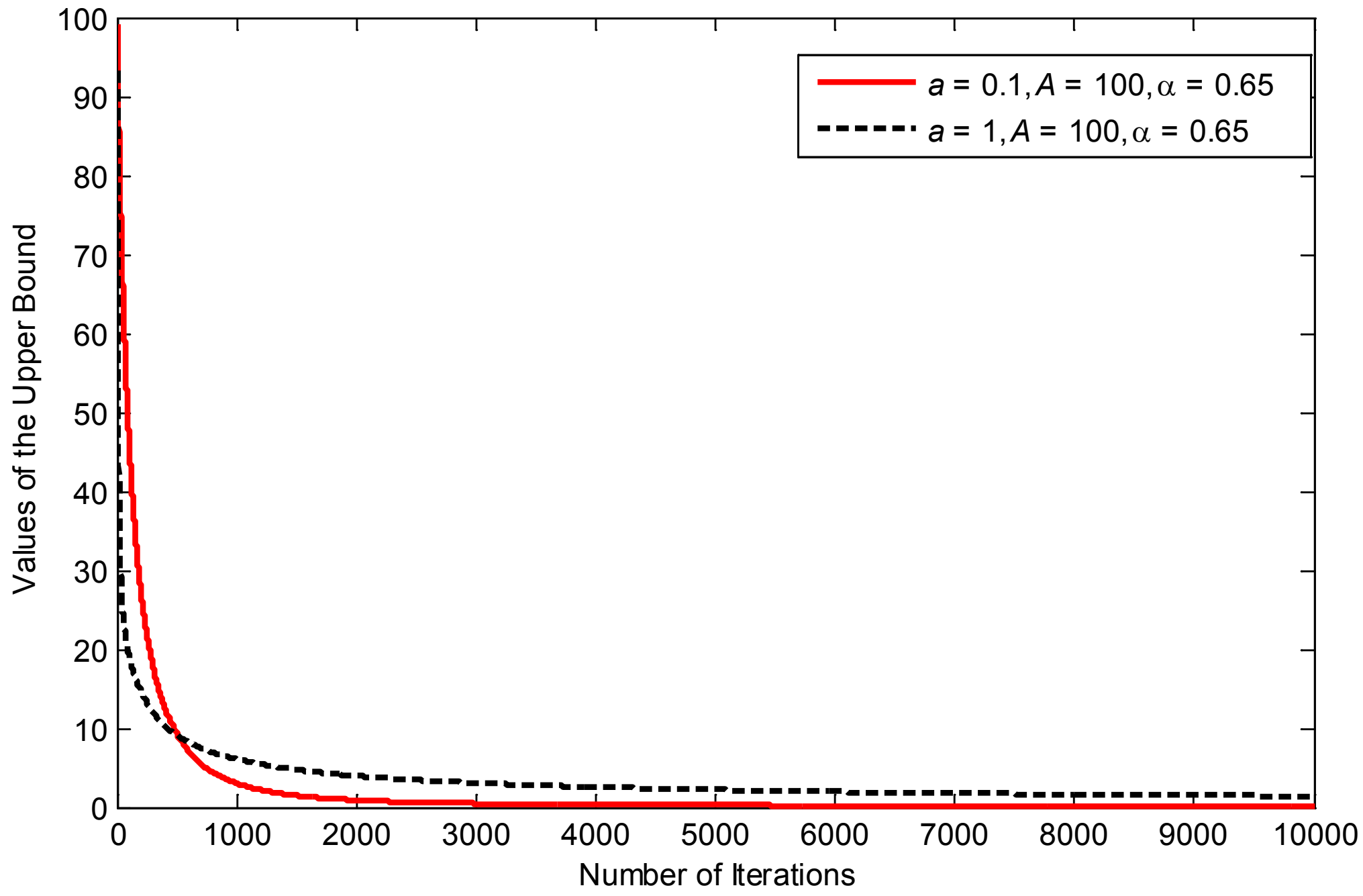

**Figure C.3** Effect of coefficient $a$ for the cases of $a = 0.1$ and $a = 1$. Bigger $a$ provides lower upper bound in the early iterations, and smaller $a$ provides lower upper bound in the later iterations.

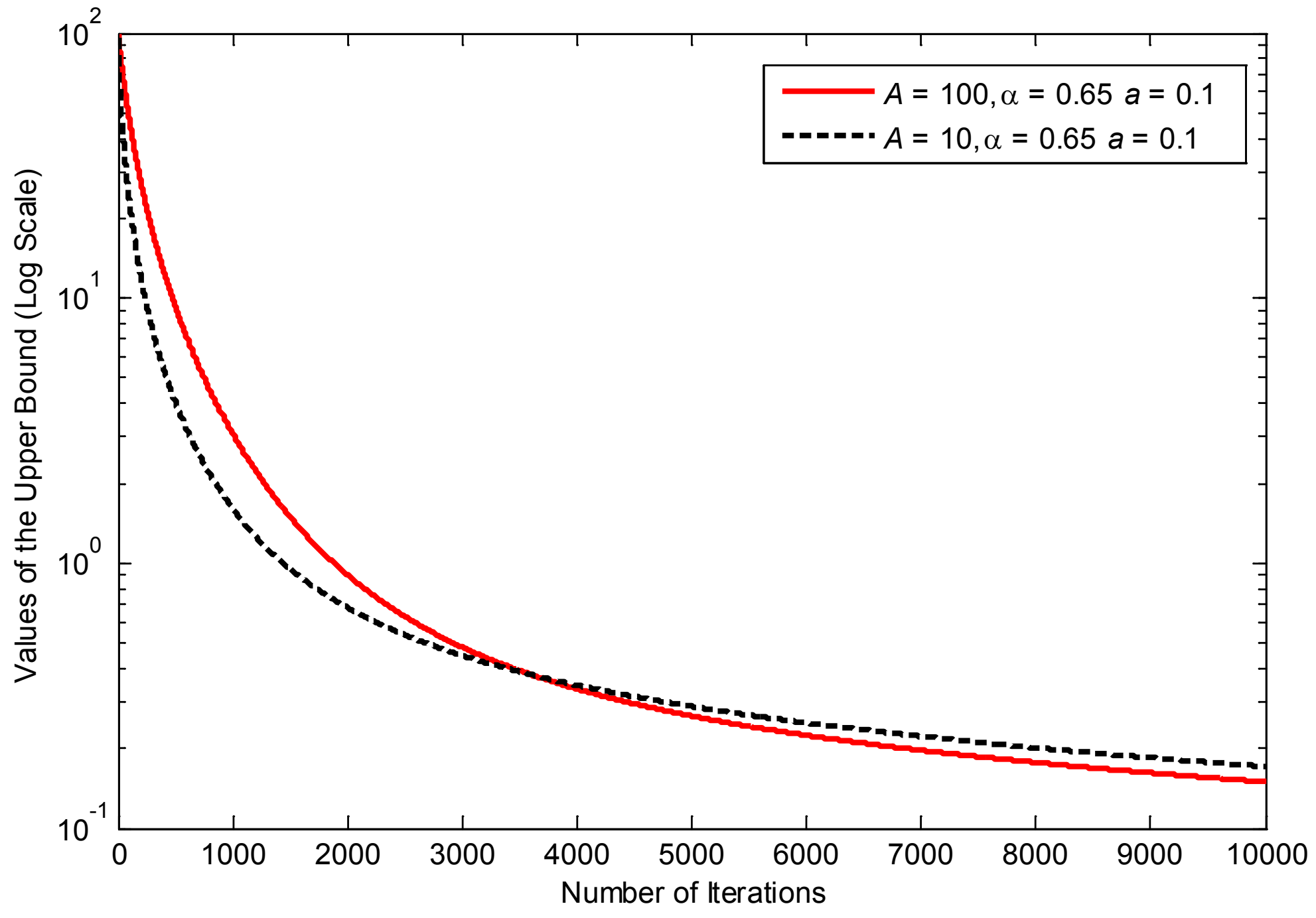

**Figure C.4** Effect of coefficient $A$ for the cases of $A$ = 10 and $A$ = 1000. Smaller $A$ provides lower upper bound in the early iterations. The vertical axis is in the logarithm scale.

In all, these numerical results are consistent with the theoretical analysis of the properties of the upper bounds in Section 3.2. Moreover, in Section 4.2 we see that the performance of DSPSA in multiple settings has the same properties as the upper bound, it indicating that the upper bound is meaningful in the sense that it captures the true properties of the performances of DSPSA.

# C.2 Discussions on the Choice of μ and *b*

After discussing the effects of $\alpha$, *A*, and *a* on the upper bound, we start to consider the other two important parameters in the upper bound: $\mu$ and *b*. In the following, we check the properties of the upper bound relative to the choice of $\mu$ and *b*. We consider the special case that $\Delta_{ki}$ are independent Bernoulli random variables taking the values $\pm 1$ with probability $1/2$. From the condition (v) of Theorem 3.1, we have the relationship that

$$0 < \mu \le E\left[(\hat{\boldsymbol{\theta}}_k - \boldsymbol{\theta}^*)^T \bar{\boldsymbol{g}}(\boldsymbol{\pi}(\hat{\boldsymbol{\theta}}_k))\right] \Big/ E\left[(\hat{\boldsymbol{\theta}}_k - \boldsymbol{\theta}^*)^T (\hat{\boldsymbol{\theta}}_k - \boldsymbol{\theta}^*)\right]$$

for all $k \ge 0$. Furthermore, *b* is a uniform upper bound for

$$E(\varepsilon_k^+ - \varepsilon_k^-)^2 + E\left[L(\hat{\boldsymbol{\theta}}_k^+) - L(\hat{\boldsymbol{\theta}}_k^-)\right]^2$$

when $k \ge 0$. The choices of $\mu$ and *b* affect the tightness of the upper bound on $E\left\|\hat{\boldsymbol{\theta}}_k - \boldsymbol{\theta}^*\right\|^2$. Generally, it is not easy to find the good values of $\mu$ and *b*.

Different from the numerical tests of the effects of $\alpha$, *A*, and *a* on the upper bound in inequality (3.3a), here we want to discuss the effects of $\mu$ and *b* on the tightness of the upper bound in inequality (3.3a). Due to the definition of $\mu$ and *b*, we need to specify the loss function in the tests below, while in the numerical tests in Section C.1 we do not need to specify the form of the loss function. Here we use the discrete separable quadratic function $\boldsymbol{\theta}^T\boldsymbol{\theta}$ as an example to discuss the choice of $\mu$ and $b$. Suppose the noises $\varepsilon_k^{\pm}$ are i.i.d. normal distributed $N(0,\sigma^2)$,

By the same arguments in Proposition 2.3, we have

$$\bar{\boldsymbol{g}}(\boldsymbol{\pi}(\hat{\boldsymbol{\theta}}_k))^T(\hat{\boldsymbol{\theta}}_k-\boldsymbol{\theta}^*)=2(\hat{\boldsymbol{\theta}}_k-\boldsymbol{\theta}^*)^T(\hat{\boldsymbol{\theta}}_k-\boldsymbol{\theta}^*)+2(\boldsymbol{\pi}(\hat{\boldsymbol{\theta}}_k)-\hat{\boldsymbol{\theta}}_k)^T(\hat{\boldsymbol{\theta}}_k-\boldsymbol{\theta}^*).$$

Then, by condition (v) in Theorem 3.1 we have for all $k$

$$\mu\le 2\left(\frac{E\left[(\hat{\boldsymbol{\theta}}_k-\boldsymbol{\theta}^*)^T(\hat{\boldsymbol{\theta}}_k-\boldsymbol{\theta}^*)\right]+E\left[(\boldsymbol{\pi}(\hat{\boldsymbol{\theta}}_k)-\hat{\boldsymbol{\theta}}_k)^T(\hat{\boldsymbol{\theta}}_k-\boldsymbol{\theta}^*)\right]}{E\left[(\hat{\boldsymbol{\theta}}_k-\boldsymbol{\theta}^*)^T(\hat{\boldsymbol{\theta}}_k-\boldsymbol{\theta}^*)\right]}\right). \tag{C.1}$$

From inequality (C.1), we know $\mu$ is the uniform lower bound for all $k\ge 0$. When the initial guess is far away from the optimal solution $\boldsymbol{\theta}^*$, then in the early iterations, the magnitude of $\boldsymbol{\pi}(\hat{\boldsymbol{\theta}}_k)-\hat{\boldsymbol{\theta}}_k$ is significantly small compared to that for $\hat{\boldsymbol{\theta}}_k-\boldsymbol{\theta}^*$. Thus, for the case of far-away initial guess, we ignore the second term of the numerator on the right-hand side of inequality (C.1), which indicates in the early iterations, we have $\mu\le 2$. In the later iterations, $\hat{\boldsymbol{\theta}}_k$ is very close to the optimal solution $\boldsymbol{\theta}^*$, and we know the magnitude of $\boldsymbol{\pi}(\hat{\boldsymbol{\theta}}_k)-\hat{\boldsymbol{\theta}}_k$ is much more significant than that of $\hat{\boldsymbol{\theta}}_k-\boldsymbol{\theta}^*$, which indicates that the upper bound for $\mu$ in inequality (C.1) is very large in the late iterations. Overall, based on the discussions above, we can choose $\mu$ as

$$\mu\approx 2. \tag{C.2}$$

Now let us consider the value of $b$, and $b$ is a uniform upper bound for

$$E(\varepsilon_k^+ - \varepsilon_k^-)^2 + E\left[L(\hat{\boldsymbol{\theta}}_k^+) - L(\hat{\boldsymbol{\theta}}_k^-)\right]^2 .$$

Since the noises $\varepsilon_k^{\pm}$ are i.i.d. normal distributed $N(0,\sigma^2)$, we have $E(\varepsilon_k^+ - \varepsilon_k^-)^2 = 2\sigma^2$. Furthermore, for the discrete separable quadratic function of interest, we have

$$\begin{aligned}
&E\left[L(\hat{\boldsymbol{\theta}}_k^+) - L(\hat{\boldsymbol{\theta}}_k^-)\right]^2 \\
&= E\left[L\left(\boldsymbol{\pi}(\hat{\boldsymbol{\theta}}_k) + \frac{1}{2}\boldsymbol{\Delta}_k\right) - L\left(\boldsymbol{\pi}(\hat{\boldsymbol{\theta}}_k) - \frac{1}{2}\boldsymbol{\Delta}_k\right)\right]^2 \\
&= E\left[E\left(\left(L\left(\boldsymbol{\pi}(\hat{\boldsymbol{\theta}}_k) + \frac{1}{2}\boldsymbol{\Delta}_k\right) - L\left(\boldsymbol{\pi}(\hat{\boldsymbol{\theta}}_k) - \frac{1}{2}\boldsymbol{\Delta}_k\right)\right)^2 \middle| \hat{\boldsymbol{\theta}}_k\right)\right] \\
&= 4E\left[\boldsymbol{\pi}(\hat{\boldsymbol{\theta}}_k)^T E\left(\boldsymbol{\Delta}_k \boldsymbol{\Delta}_k^T \middle| \hat{\boldsymbol{\theta}}_k\right) \boldsymbol{\pi}(\hat{\boldsymbol{\theta}}_k)\right]. \qquad \text{(C.3)}
\end{aligned}$$

Since the components of $\boldsymbol{\Delta}_k$ are independently Bernoulli $\pm 1$ distributed, then we have $E\left(\boldsymbol{\Delta}_k \boldsymbol{\Delta}_k^T \middle| \hat{\boldsymbol{\theta}}_k\right) = \boldsymbol{I}_p$, which implies that eqn. (C.3) can be written as

$$E\left[L(\hat{\boldsymbol{\theta}}_k^+) - L(\hat{\boldsymbol{\theta}}_k^-)\right]^2 = 4E\left(\boldsymbol{\pi}(\hat{\boldsymbol{\theta}}_k)^T \boldsymbol{\pi}(\hat{\boldsymbol{\theta}}_k)\right).$$

As we know that $b$ is the uniform upper bound for all $k \geq 0$. Therefore, when the initial guess is far away from the optimal solution, then the choice of $b$ is restricted by some large values of $E\left[L(\hat{\boldsymbol{\theta}}_k^+) - L(\hat{\boldsymbol{\theta}}_k^-)\right]^2$ in the early iterations, which implies that in order to fit the performance of DSPSA well in the early iterations, the value of $b$ might be very large, which might be too big for the later iterations. Thus, tight performance of the upper bound in (3.3a) in the early iterations might not indicate the tight performance for later iterations. Under the case of an initial guess far from $\boldsymbol{\theta}^*$ and the separable quadratic loss function, we assume the initial guess provides the largest value of $E\left[L(\hat{\boldsymbol{\theta}}_k^+)-L(\hat{\boldsymbol{\theta}}_k^-)\right]^2$, then we estimate the value of $b$ as

$$b \approx 2\sigma^2 + 4\boldsymbol{\pi}(\hat{\boldsymbol{\theta}}_0)^T \boldsymbol{\pi}(\hat{\boldsymbol{\theta}}_0). \tag{C.4}$$

By using the estimated value of $\mu$ and $b$ in formula (C.2) and (C.4), we test the effects of $\mu$ and $b$ on the upper bound for the separable quadratic function numerically. Let us set the coefficients to be $\alpha = 0.501$, $A$ = 1000, $a$ = 0.05, and based on this setting we have $2\mu a/(1+A+k)^{\alpha} < 1$ for all $k \geq 0$.We test both low-dimensional ($p$ = 5) and high-dimensional cases ($p$ = 200). The initial guess are $10\times\mathbf{1}_5$ and $10\times\mathbf{1}_{200}$, respectively. The noises of the loss function are i.i.d. $N(0,1)$. The number of replicates is 20 and in each replicate the number of iterations is 10000.

Here the initial guess is far from the optimal solution; the noise is not large; the gain step size is small. Then, the sequence $\{\hat{\boldsymbol{\theta}}_k\}$ has a low probability to go far more away than the initial guess. Thus, the unit hypercube centered by $\boldsymbol{\pi}(\hat{\boldsymbol{\theta}}_0) = 10.5\times\mathbf{1}_p$ may have the largest value among all $E\left[L(\hat{\boldsymbol{\theta}}_k+\boldsymbol{\Delta}_k/2)-L(\hat{\boldsymbol{\theta}}_k-\boldsymbol{\Delta}_k/2)\right]^2$ with $k \geq 0$. Then, we can set $b$= 2 + $4\boldsymbol{\pi}(\hat{\boldsymbol{\theta}}_0)^T \boldsymbol{\pi}(\hat{\boldsymbol{\theta}}_0)$ = 2207 for $p$ = 5, and $b$= 2 + $4\boldsymbol{\pi}(\hat{\boldsymbol{\theta}}_0)^T \boldsymbol{\pi}(\hat{\boldsymbol{\theta}}_0)$ = 88202 for $p$ = 200. In Figure C.5, we see that the upper bound is fairly accurate for the actual performance of DSPSA for the separable quadratic loss function in the low dimensional case. However, in Figure C.6, for the high dimensional case, the bound is really tight for the early iterations and not very tight for the later iterations. In both Figure C.5 and Figure C.6, we divide the mean square error and the upper bound of the mean square error by $\left\|\hat{\boldsymbol{\theta}}_0 - \boldsymbol{\theta}^*\right\|^2$ to normalize them. The vertical axis in both Figure C.5 and Figure C.6 is the sample mean of $\left\|\hat{\boldsymbol{\theta}}_k - \boldsymbol{\theta}^*\right\|^2$ over $\left\|\hat{\boldsymbol{\theta}}_0 - \boldsymbol{\theta}^*\right\|^2$, where the sample mean of $\left\|\hat{\boldsymbol{\theta}}_k - \boldsymbol{\theta}^*\right\|^2$ is the arithmetic mean of the observed values of $\left\|\hat{\boldsymbol{\theta}}_k - \boldsymbol{\theta}^*\right\|^2$ across independent replicates.

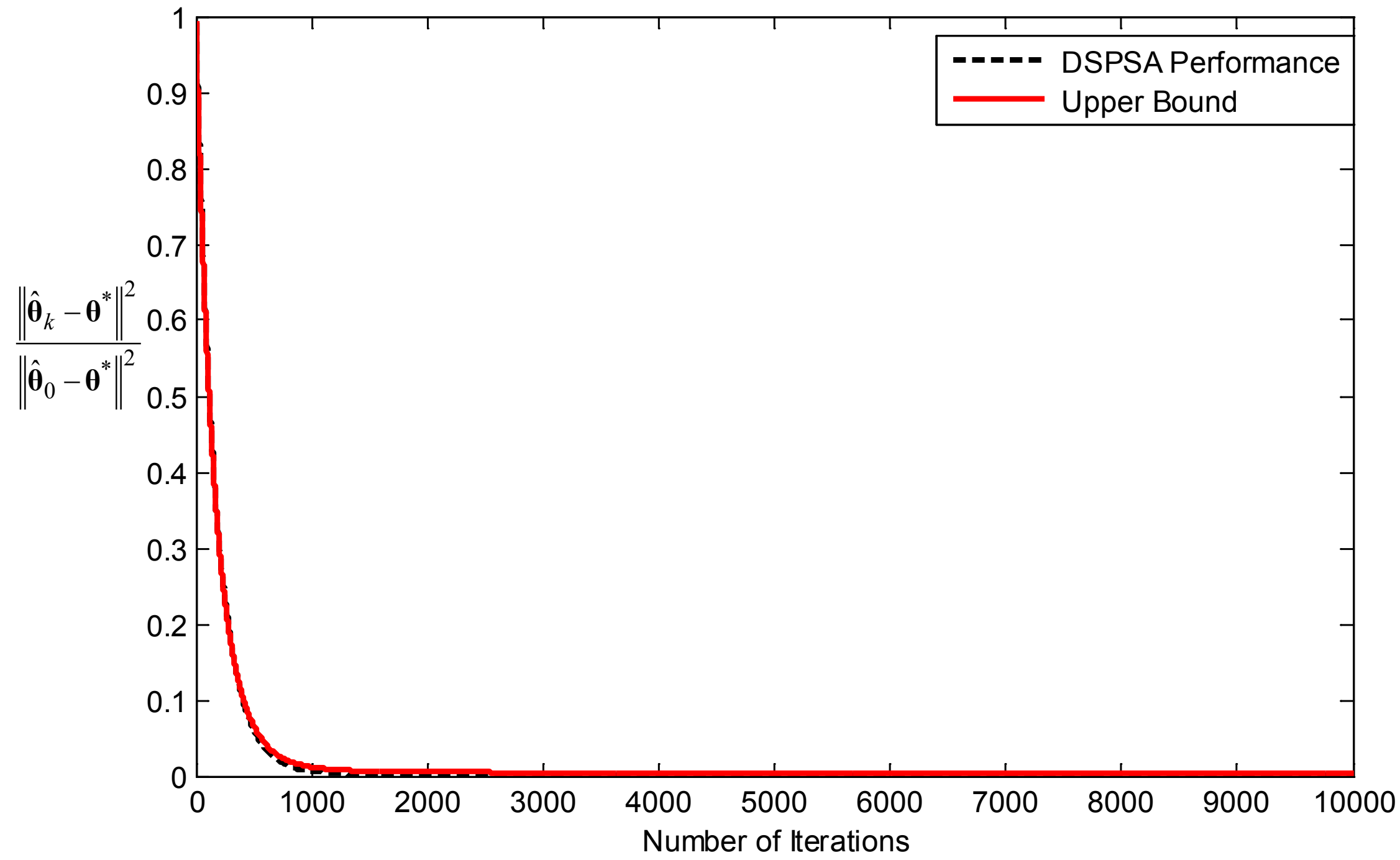

**Figure C.5** Performance of DSPSA and the upper bound for the separable loss function for the case of dimension 5 ($\mu$=2 and $b$ = 2207). The upper bound is really tight for this case. Each curve represents the sample mean of 20 independent replicates.

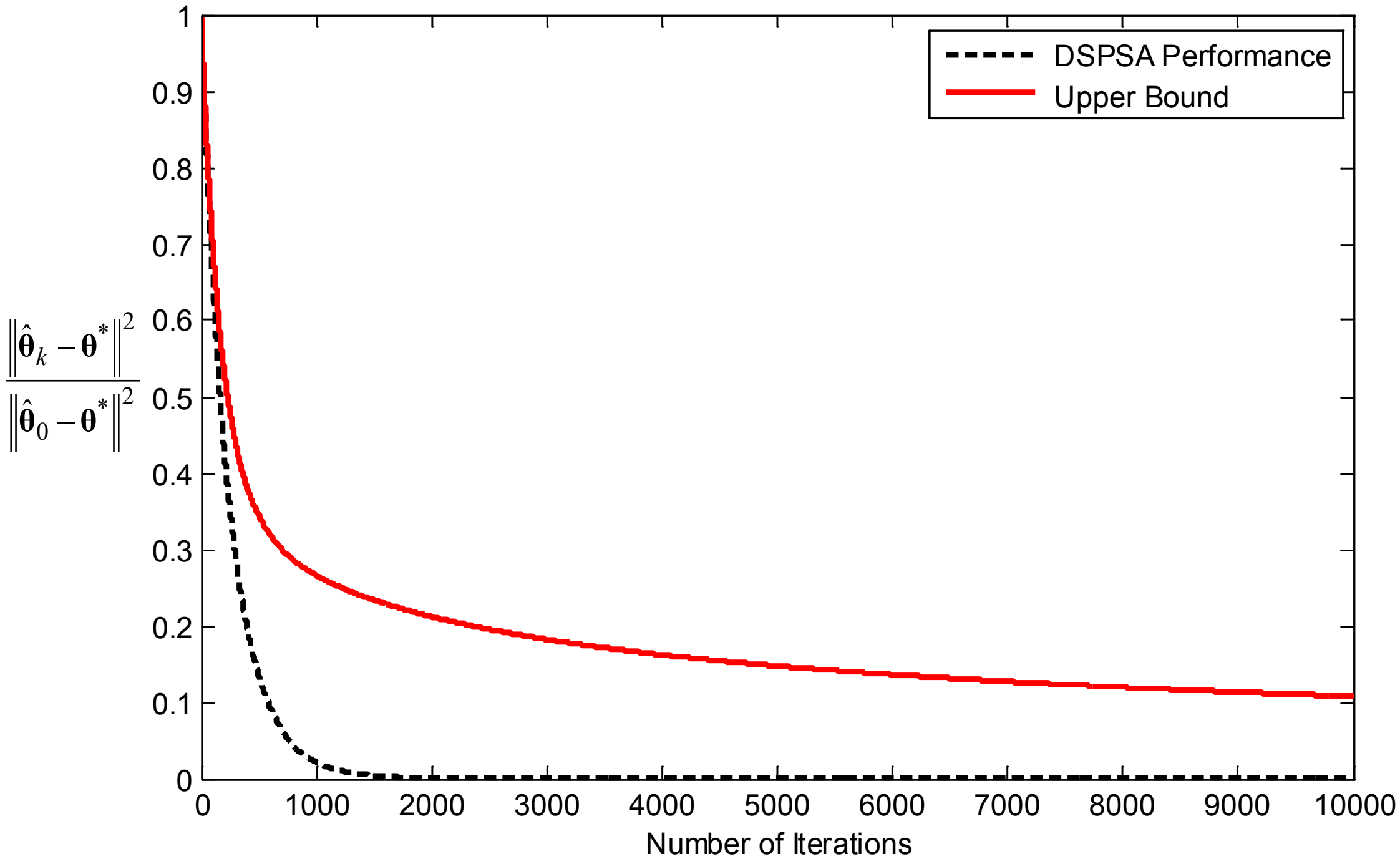

**Figure C.6** Performance of DSPSA and the upper bound for the separable quadratic loss function on dimension of 200 ($\mu$=2 and $b$ = 88202). The upper bound is tight in the early iterations, and the upper bound is not very tight in the later iterations. Each curve represents the sample mean of 20 independent replicates.

As we know the upper bound in inequality (3.3a) is composed of two terms. The first term (initial guess term) is corresponding to performance in the early iterations and the second term (integral term) is relative to performance in the later iterations. We have shown that the proportion of the first term in the upper bound is decreasing; while the proportion of the second term is increasing. From Figure C.5 and C.6, we see that for both cases (lower-dimensional and higher-dimensional) the upper bounds are really tight in the early iterations, and it means that the first term of the upper bound fits the performance of DSPSA really well in the early iterations.

From Figure C.7, we see that the proportion of the first term decreases faster for the high-dimensional case than for the low-dimensional case. For the low-dimensional case, the proportion of the first term is still higher than 50% when the error is already very close to 0. Thus, the upper bound does a good job for the low-dimensional case. However, for the high-dimensional case, the proportion of the first term decreases too fast. When the second term starts to dominate the upper bound, the value of *b* estimated by the initial guess might be too large for the later iterations.

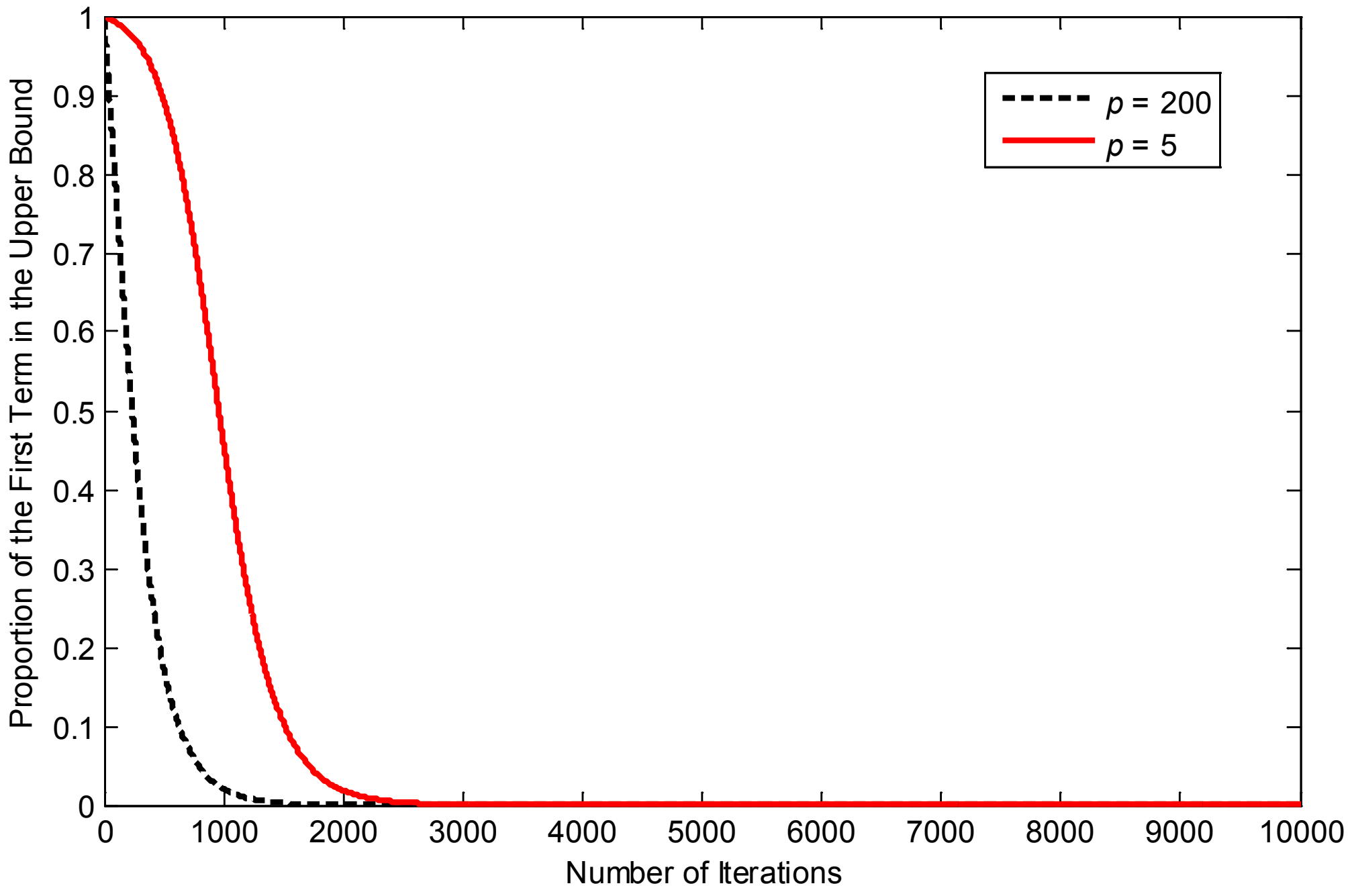

**Figure C.7** Proportion of the first term (initial guess term) in the upper bound for both low-dimensional case and high-dimensional case. The proportion of the first term decreases faster under high-dimensional case.

Here we propose one way to improve the tightness of the upper bound in the later iterations. We have picked $\mu = 2$ in the numerical test above, but by condition (v) in Theorem 3.1 we know all positive values smaller than 2 are also valid. Thus, in our proposed way, we pick a smaller value for $\mu$ and at the same time we decrease the value of *b*. We explain the reasoning in the following. The first term of the upper bound (3.3a) only contains $\mu$, so decreasing the value of $\mu$ hurts the tightness of the upper bound in the early iterations. Due to this sacrifice in the early iterations, in practice we can decrease the value of *b* to improve the tightness of the upper bound in the later iterations, even though by the definition of *b*, decreasing the value of *b* is not allowed.

In Figure C.8, we show the graph for high-dimensional case when we set $\mu$=1.2 and *b* = 5000. Similar as Figure C.5 and C.6, we divide the mean square error and the upper bound of the mean square error by $\left\|\hat{\boldsymbol{\theta}}_0 - \boldsymbol{\theta}^*\right\|^2$ to normalize them. The vertical axis in Figure C.8 is sample mean of $\left\|\hat{\boldsymbol{\theta}}_k - \boldsymbol{\theta}^*\right\|^2$ over $\left\|\hat{\boldsymbol{\theta}}_0 - \boldsymbol{\theta}^*\right\|^2$, where the sample mean of $\left\|\hat{\boldsymbol{\theta}}_k - \boldsymbol{\theta}^*\right\|^2$ is the arithmetic mean of the observed values of $\left\|\hat{\boldsymbol{\theta}}_k - \boldsymbol{\theta}^*\right\|^2$ across independent replicates.

Comparing Figure C.6 with Figure C.8, we see that when decreasing the values of $\mu$ and *b* at the same time, the upper bound in early iterations is still reasonable and the tightness in later iterations has been significantly improved.

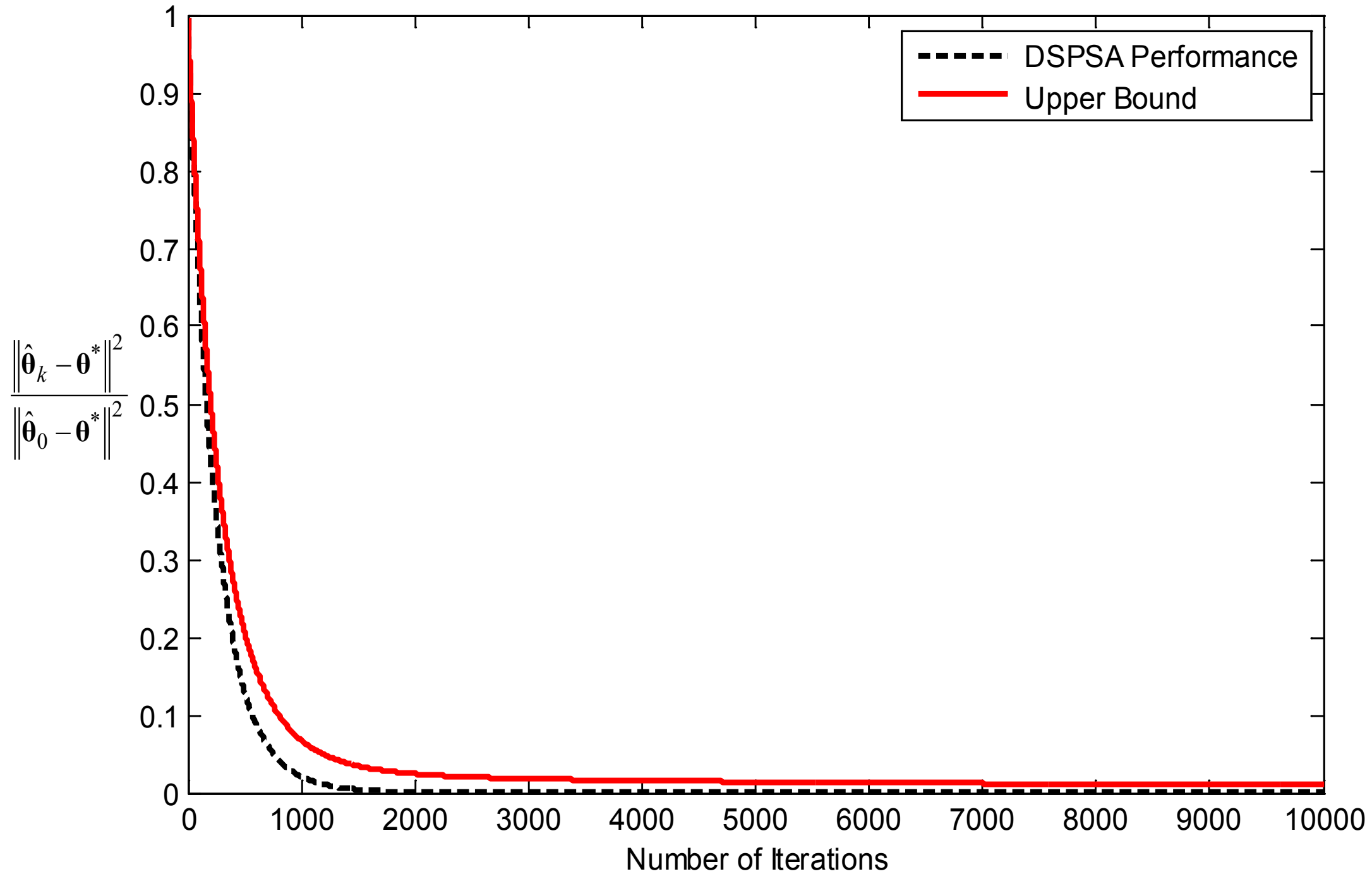

**Figure C.8** Performance of DSPSA and the upper bound for the separable loss function on dimension of 200 ($\mu$=1.2 and $b$ = 5000). The upper bound is reasonably tight in both early iterations and late iterations. Each curve represents the sample mean of 20 independent replicates.